\documentclass[11pt,oneside]{book}
\usepackage{amsthm}

\input{mystyle_arxiv.sty}

\usepackage{tabularx}
\usepackage{float}
\usepackage{mathtools}
\usepackage{microtype}
\usepackage{xurl}
\usepackage{geometry}
\usepackage{makeidx}
\usepackage{index}

\usepackage{sectsty}
\sectionfont{\bfseries}
\usepackage{tocloft}

\hypersetup{
    colorlinks=true,
    linkcolor=blue,
    filecolor=blue,
    urlcolor=blue,
    citecolor=blue,
    pdftitle={Mean Field Reinforcement Learning: Theory and Algorithms for Mean Field Control},
    pdfauthor={René Carmona and Mathieu Laurière},
}

\DeclareMathOperator{\Tr}{Tr}
\DeclareMathOperator{\Proj}{Proj}

\newcommand{\norm}[1]{\left\|#1\right\|}
\newcommand{\inner}[1]{\left\langle #1 \right\rangle}
\newcommand{\abs}[1]{\left|#1\right|}

\usepackage[noautomatic]{imakeidx}

\makeindex                 %
\makeindex[name=sub,title=Subject Index]
\makeindex[name=not,title=Notation Index]

\begin{document}

\renewcommand{\chaptermark}[1]{\markright{#1}}

\frontmatter

\title{Mean Field Reinforcement Learning\\[0.5em]
{\large Theory and Algorithms for Mean Field Control}}
\author{
Ren\'e Carmona \thanks{Department of Operations Research and Financial Engineering \& Program in Applied and Computational Mathematics, Princeton NJ 08544, USA,\\
{\tt email: rcarmona@princeton.edu}}
\and 
Mathieu Lauri\`ere \thanks{Shanghai Frontiers Science Center of Artificial Intelligence and Deep Learning; NYU-ECNU Institute of Mathematical Sciences at NYU Shanghai; NYU Shanghai, 567 West Yangsi Road, Shanghai, 200126, People’s Republic of China\\
{\tt email: mathieu.lauriere@nyu.edu}}
}

\date{}
\hypersetup{pageanchor=false}
\maketitle
\hypersetup{pageanchor=true}
\chapter*{Preface}
\addcontentsline{toc}{chapter}{Preface}

 {\small
Reinforcement learning (RL) has become a central paradigm for learning decisions directly from interaction with an uncertain environment. Its successes in games, robotics, operations, economics, and finance illustrate both the flexibility of the approach and the breadth of its possible applications. Although most of its early accomplishments concerned single-agent problems, many important environments are not single-agent systems.  Traffic
networks, energy markets, communication platforms, financial markets, and
large online communities all involve many decision makers whose choices affect one another.  A direct multi-agent description can quickly become too large to analyze or compute with, while a single-agent description misses the feedback loop between individual decisions and the collective behavior they generate.

Mean field reinforcement learning (MFRL) addresses this
difficulty by combining RL,
stochastic control, and mean field models.  Its guiding idea is to replace a
large population of weakly interacting agents by a representative agent
coupled to the distribution of the population.  This is more than a dimension reduction device. It changes the state of the decision problem: probability
measures become part of the dynamics, conditional laws may have to be tracked,
and common sources of randomness can keep the population correlated even in
the infinite-population limit. These features are what make the subject both
powerful and mathematically delicate.

There is now a substantial literature on multi-agent reinforcement learning
(MARL), mean field games (MFGs), mean field control
(MFC), and mean field Markov decision processes
(MFMDPs).  Yet these areas are often presented with
different conventions, different levels of mathematical detail, and different
algorithmic priorities.  This monograph is meant to help close that gap.  It
develops a self-contained discrete-time route from rigorous MFC
models to concrete RL algorithms, with particular care
for the distinction between finite populations and their mean field limits,
the role of common noise and common randomization, and the relation between
open-loop, closed-loop, and randomized controls.

The exposition is mathematical, but it is guided by computation.  We do not seek the greatest mathematical generality.  
Instead, the
emphasis is on formulations that are precise enough to derive dynamic programming principles, propagation of chaos arguments, and convergence
statements, while remaining close to implementable learning procedures.

The monograph is organized as follows. It first places MFRL in the broader RL and MARL context, and recalls
the probabilistic and control-theoretic tools needed for their analysis. We then study two complementary models that cover different aspects of the framework, and analyze them at a rigorous mathematical level.  
A first abstract model with common noise is used to clarify the different roles of open and closed loop policies and the relation between finite-population systems and representative-agent formulations. It also connects optimization by a central planner with the lifted MFMDP viewpoint.
Then, a linear
quadratic model provides an explicit benchmark in which stability,
optimality, and policy gradient methods can be analyzed in detail.  The final
part proposes numerical implementations of these formulations: tabular and
discretized algorithms for the abstract model, and policy gradient methods
with several simulator viewpoints for the linear quadratic model.

The intended reader can be a graduate student or an early researcher in applied mathematics coming from probability, stochastic control,
reinforcement learning, or multi-agent systems.  The guiding message is that
MFRL is not simply RL with a larger state space, nor MFC
with numerical experiments appended to it.  It is a framework in which large
population limits, information structures, and learning procedures have to be
designed together.

\vskip 6pt
This monograph is an expanded version of the lecture notes prepared by the authors for a set of tutorial lectures at the University of Chicago NSF Institute for Mathematics, Statistics and Innovation (IMSI) in March 2026, in preparation for the week-long workshop \emph{Theoretical Foundations of Multi-agent and Mean Field Reinforcement Learning} held at IMSI May 18-22, 2026.\footnote{See \url{https://www.imsi.institute/activities/theoretical-advances-in-reinforcement-learning-and-control/foundations-of-multi-agent-and-mean-field-reinforcement-learning/}}

\vskip 6pt\noindent\textbf{Acknowledgements:}

We would like to thank IMSI for a superb job at hosting the workshop. Our gratitude also goes to Zongjun Tan with whom we had numerous interactions during his PhD studies at Princeton, and with whom we co-authored several papers on RL and MFRL; some of the results of these papers are used in this monograph. R.C. was partially supported by AFOSR grant FA-9550-23-1-0324. We acknowledge the use of NYU's and NYU Shanghai's HPC for some of the numerical experiments.

\vskip 12pt
\begin{flushright}
Ren\'e Carmona\\
Princeton, N.J.\\
\vskip 2pt
Mathieu Lauri\`ere\\
Shanghai, China
\vskip 2pt
June 29, 2026
\end{flushright}
}

\pagebreak

\tableofcontents

\mainmatter

\part{Stochastic Control and Reinforcement Learning Landscape}

\chapter[MARL versus MFRL]{Multi Agent Reinforcement Learning versus Mean Field Reinforcement Learning}
\label{ch:review}

\section{\textbf{Introduction}}

The relevance and economic importance of applications of
sequential decision making  by multiple interacting agents has emerged as a major impetus for the development of the theory of Multi-Agent Reinforcement Learning (MARL). \index[not]{MARL}\index[sub]{multi-agent reinforcement learning}The survey-style papers \cite{BusoniuBabuskadeSchutter2010} and \cite{ZhangYangBasar2021} give an overview of the challenges and some of the existing algorithms involved. Noticeably, \cite{CuiTahirEkinci2022} offers an extensive bibliography of $320$  references. The applications come from domains as diverse as autonomous vehicle coordination, robotic swarms, financial markets, smart grids, and distributed resource allocation.

While the classical theory of single-agent Reinforcement Learning (RL) \index[not]{RL}\index[sub]{reinforcement learning} focuses on the approximation of agent optimal policies that depend only on the state of the agent, their own actions and their subsequent potential rewards, and possibly some random shocks, the search for optimal policies in  multi-agent systems is much more involved.
Indeed, each agent's decision influences the time evolution (classically prescribed by transition probabilities) of the states and the rewards of the other agents, creating intricate dependencies that compound the complexity of MARL models. In fact, the theoretical complexity of the model is not the only roadblock to a straightforward generalization of the theory and algorithms developed for single agent RL. The formulation of the optimization problem becomes a delicate matter as the computational complexity tends to scale with the number of agents. As modeling practical applications increasingly requires the control of tens, hundreds, or even thousands of individual agents, the need for sound models of the interactions between the agents, and for scalable and theoretically grounded algorithms becomes of critical importance.

\vskip 4pt
The goal of this chapter is to review the pros and cons of the theoretical models and their implementations as they appear in the published literature for finitely many agents RL. Highlighting some of their limitations will serve as a motivation for the introduction of Mean Field Reinforcement Learning (MFRL) models, \index[not]{MFRL}\index[sub]{mean field reinforcement learning} to which this monograph is devoted. Our contention is that these models provide tractable solutions to some of the main difficulties encountered in the analysis of large populations of agents.

\vskip 4pt
We first recall separately the main definitions and properties of MARL and MFRL. We then highlight the main differences in a comparison intended to firmly establish the basis of our mathematical treatment of MFRL in the rest of the monograph.

\vskip 4pt\noindent
\emph{Note to the reader:} while most economists and engineers set up models so that the optimization part is the maximization of some form of utility or reward, most of the models treated by mathematicians involve the minimization of energy or cost. Since these two forms of models are mathematically equivalent, we choose, because of our own backgrounds and biases, the latter approach. Throughout the monograph, we shall minimize costs instead of maximizing rewards, unless this choice is detrimental to comprehension of the application.

\section{\textbf{Multi-Agent Reinforcement Learning}}

This section is intended to set the stage for a rigorous mathematical discussion of the foundations and paradigms of MARL as found in the published literature on the subject. 

\subsection{Problem Formulation}

The premise of MARL is to extend the technology developed in classical single agent reinforcement learning to  $N\ge 2$ agents interacting within a shared environment. The natural generalization is to replace the standard mathematical setting of a Markov Decision Process (MDP) \index[not]{MDP}\index[sub]{Markov decision process} by some sort of joint MDP, and because of the multivariate nature of the optimizations, the most convenient way to formulate the model is to treat it  as a stochastic game with $N$ players. Although attractive at first, this option requires specifying the kind of information available to each agent when making their decisions. The mathematical set-up we use to describe the models includes the following components: 
\begin{itemize}
\item  A joint state space $S$ \index[sub]{state space}that captures the global system configuration. It will be assumed to be a \emph{nice} measurable space, typically a Borel or Polish space.
\item Individual action spaces \index[sub]{action space}$A^1,\cdots, A^N$, one for each agent, which we shall also assume to be Borel or Polish spaces, even if in most practical applications we assume that they are closed convex subsets of a finite dimensional Euclidean space, say $\RR^k$. We shall also use the space of action profiles, \index[sub]{action profile} that is the product space $A = A^1\times\cdots\times A^N$.
\item A state transition probability kernel $p(\,\cdot\,|s,\ba)=\PP[\cdot\,|\,s, a^1, \cdots, a^N]$ which gives the probability distribution of the next state of the system when the previous  state was $s\in S$, and the actions taken by the $N$ agents were given by the joint action profile $\ba=(a^1,\cdots,a^N)$. It is a probability measure on $S$. We denote by $\cP(S)$ \index[not]{$\cP(S)$} the space of probability measures on $S$. Notice that when the state space $S$ is finite, say with $n$ elements, $\cP(S)$ can be identified with the $n-1$ dimensional simplex \index[sub]{simplex}
$$
\Delta_{n-1}=\{(p^1,\cdots,p^n)\,|\, \sum_{i=1}^np^i=1, \; p^i\ge 0,\,i=1,\cdots,n\}.
$$
\begin{remark}
To conform with the vast majority of existing articles on the subject, we define state dynamics with a \index[sub]{transition probability}transition probability. Still, as we explain in the following subsection, the dynamics of the state can alternatively be given by a \index[sub]{system function} system function. We explain below that these two formulations are essentially equivalent and that we shall choose one over the other when convenient.
\end{remark}

\item  Individual reward functions $S\times A\ni (s,a^1,\cdots,a^N)\mapsto r^i(s, a^1,\cdots, a^N)$ that are typically assumed to be bounded measurable functions of the state and the entire action profile (i.e. the actions of all the agents);

\begin{remark}
Most of the MARL literature relates to models in which rewards are maximized, hence our choice to introduce reward functions in this section. But as explained at the end of the introduction, despite the fact that the early RL literature maximizes rewards, we shall switch to the mathematical tradition to minimize cost functions, and starting from the next chapter,  we shall use cost functions $f^i$ instead of reward functions $r^i$.
\end{remark}

\item Individual policies $\pi^i: S \to \cP(A^i)$ that each agent uses to choose which action to take in order to optimize their objective. A reader familiar with the stochastic control or MDP literature could expect that the policies are \index[sub]{strict policy}\emph{strict}, i.e. measurable functions $\pi^i: S \to A^i$, however, as it is standard in the theoretical models and the implementation algorithms of RL, one privileges \emph{mixed} policies, \index[sub]{mixed policy}\index[sub]{randomized policy} also called \emph{randomized} policies for which, given the state $s\in S$ of the system, agent $i$ will choose an action $a^i\in A^i$ drawn at random from the probability distribution $\pi^i(s)\in\cP(A^i)$. 
\end{itemize}

\begin{remark}
In most of the models considered in the MARL literature, the state $s\in S$ is, in fact, a state profile $s=(s^1,\cdots,s^N)$ of private states $s^i\in S^i$ of the individual agents, so $S=S^1\times\cdots\times S^N$. Actually, in many cases the state of the system also includes an extra factor $s^0\in S^0$ common to all the individual agents, also evolving over time.
\end{remark}

Given this set-up, each agent $i$ tried to choose a course of action to optimize their objective. However, the mathematical formulation of this optimization depends greatly on the collective behavior of the agents which encompasses cooperative settings (where agents share a common reward), competitive scenarios captured by stochastic game models as well as a wide spectrum of intermediate settings. We shall give more details and references for the differences between these models later on.
For the purpose of the present discussion, we assume that agent $i$ 
aims to maximize (resp. minimize) their expected cumulative discounted reward (resp. cost) 
\begin{equation}
    \label{fo:ith_objective}
J^i(\pi)=\EE\Bigl[\sum_{n\in\cN} \gamma^n r^i(X_n,a^1_n,\cdots,a^N_n)\Bigr]
\end{equation}
while other agents simultaneously optimize their own objectives. Here 

\begin{itemize}\itemsep=-1pt
\item $\cN$ is the time horizon over which the optimization is performed. $\cN=\{0,1,\cdots,T\}$ for a finite integer $T$ for finite horizon models, and $\cN=\{0,1,\ldots\}$ for infinite horizon models.
\item $\gamma\in(0,1]$ is a discount factor that is used to \emph{actualize} the values of the reward or the cost over time. It is often chosen to be $1$ for finite horizon models, and it is most often chosen to be strictly smaller than $1$ for infinite horizon models, helping with the convergence of the infinite summation.
\item for each $i\in\{1,\cdots,N\}$, $a^i_n\sim\pi^i(X_n)$ is the action taken by agent $i$ at time $n$ by sampling from the probability distribution $\pi^i(X_n)$, and $X_n$ represents the random realization of the state of the system at time $n$.
\end{itemize}
According to our introduction of the transition probability $p(\,\cdot\,|s,\ba)$, the time evolution of the state of the system is given by a stochastic process $\bX=(X_n)_{n\in\cN}$ whose dynamics are given by
$$
X_{n+1}\sim \int_A p(\,\cdot\,|\,X_n,\ba)\,\pi^1(da^1|X_n)\cdots\pi^N(da^N|X_n)
$$
where throughout the monograph we use the notation $Y\sim\mu$ \index[not]{$Y\sim\mu$} to state that the random variable $Y$ is a sample realization 
from the probability distribution $\mu$, implying that $\mu$ is the distribution of $Y$.

\begin{remark}
Notice that we assumed that the reward functions $r^i$ as well as the policies $\pi^i$ are independent of time, i.e. they do not change over time. This assumption is unrealistic in many applications, including those based on finite horizon models. However, it is made to avoid encumbering the notation at this early stage, and because it leads to stationary models that  are easier to solve, theoretically as well as numerically.
\end{remark} 

\subsection{Transition Probability or System Function?}

It is often convenient to define the time evolution of the state of the system using a \emph{system function}  \index[sub]{system function} $F:S\times A\times E\to S$
and a dynamical equation
\begin{equation}
\label{fo:system}
    X_{n+1}=F(X_n,\ba_n,\epsilon_{n+1})
\end{equation}\index[sub]{system function}
for an independent identically distributed ( i.i.d. for short) sequence $(\epsilon_n)_{n\ge 1}$, of innovation random variables taking values in a measurable space $E$. This is the approach used in \cite{motte2022mean,CarmonaLauriere_AAP} and \cite{CarmonaLauriereTan2019} to analyze the optimization of mean field MDPs with general \emph{open loop} and \emph{closed loop Markovian} control processes, even in the presence of common noise.
Because of their ease of implementation, and because of their ability to preserve the desirable Markov property, \emph{closed loop} policies have been studied and implemented to a larger extent than the more general \emph{open loop} policies whose implementation in practical applications is problematic. We shall limit ourselves to closed loop policies in this introductory chapter. We shall provide a comparative study of the optimizations with respect to open and closed loop policies in Chapter \ref{ch:AAP}.

\vskip 2pt
If the dynamics of the state are prescribed via a system function like in \eqref{fo:system}, it is straightforward to retrieve the transition probability. Indeed, in this case
\begin{equation}
\label{fo:transition}
p(\,\cdot\,|\,x,\ba)=\PP[F(x,\ba,\epsilon_{n+1})\in\,\cdot\,],
\end{equation}
in other words, the probability $p(\,\cdot\,|\,x,\ba)$ is merely the image (i.e. push-forward) of the common law of the random shocks $\epsilon_{n+1}$ under the mapping $E\ni\epsilon\mapsto F(x,\ba,\epsilon)$.

\begin{remark}
It is important to keep in mind  that the value of the state $X_{n+1}$ depends on the actual realizations of the past action profiles $(a^i_m)_{1\le i\le N}$ for all prior times $m\le n$, and not just on their distributions. This will come into play in some of our analyses in Chapter \ref{ch:AAP}.
\end{remark}

\vskip 4pt
On the other hand, if the dynamics of the state are given by a transition probability and we wish to study the model as if it were given by a system function, the most convenient way to introduce the appropriate system function is to use the Blackwell-Dubins theorem\index[sub]{Blackwell-Dubins theorem} \cite{BlackwellDubins}. This important result from measure theory is a generalization of the famous Skorohod representation theorem, which is particularly well suited to our needs in developing a rigorous mathematical theory of MDPs and MFRL. We will give a precise statement in Lemma~\ref{le:BlackwellDubins} for later references.

\subsection{Key Differences between MARL and Single Agent RL}

MARL differs fundamentally from single-agent RL in several critical aspects which we briefly highlight below. The main differences arise from the multitude of possible information structures for individual agents. Specifically, the main determining factor in the construction of a model and the development of optimization algorithms is the answer to the following question:\emph{what kind of information is available to an individual agent about to make a decision concerning the best action they should take}? Usually, individual agents are cognizant of their own individual states, but do they know anything about the states of the other agents, the common sources of random shocks, their histories? In order to use the information available to them, do agents collaborate? form coalitions? compete among themselves?  Also, when it comes to deciding what action to take, should the individual agent perform their own optimization? or should they implement policies determined by a central unit? These are some of the typical issues of multivariate optimization theory and multi-objective decision theory. 
Not surprisingly, they complicate the development of a uniform theory for MARL based on a straightforward generalization of the successes of the single agent theory. Here we mention a small number of obvious consequences of the plethora of differences with the single agent theory.

\begin{itemize}
\item  \emph{Growth of the Size of the Joint State-Action Space}: The effective state-action space grows exponentially with the number of agents. Even if we assume that all the agents can choose their actions from a common finite space, say that all the action spaces $A^i$ are equal and have $k$ elements, the size of the set of action profiles, i.e.  the joint action space, is $k^N$, which makes full tabular representations exponential in $N$.

\item  \emph{Non-Stationarity}: In most MARL models, the optimization problem faced by each individual agent is more complex than in single agent RL models. This extra complexity comes in no small part from the non-stationary learning environment created by the other agents concurrently updating their own policies during learning. From each individual agent's perspective, the time evolution of the environment and reward structure they are experiencing shifts as other agents learn, violating the stationarity assumption underlying classical convergence guarantees in single agent RL.

\item  \emph{Decentralization and Partial Observability}: Many MARL models impose decentralized execution, requiring individual  agents to choose their actions on the sole basis of the information they have access to, typically local observations without access to global state information or other agents' private observations. This serious restriction on the nature of the admissible individual policies is akin to the problem of partial observation in stochastic control, and it introduces additional complexity.
\end{itemize}

Beyond these structural differences, MARL also faces challenges in credit assignment (i.e., determining each agent's contribution to the aggregate reward) which becomes prohibitive in large cooperative games, as exemplified by the computation of Shapley indices. Models requiring communication between the agents add yet another layer of complexity; see, e.g., \cite{FoersterAssaelFreitasWhiteson} and \cite{SukhbaatarSzlamFergus2016}.

\subsection{Major MARL Paradigms}

As suggested earlier, MARL frameworks are similar to those of multiple player games. Indeed, a crucial component of the definition of a MARL model is the description of the interactions between the agents. Deciding which types of interaction are admissible is of crucial importance because the nature of these interactions plays a fundamental role in the theoretical analysis of these models, as well as the development of  algorithmic solutions to address the challenges raised by the multitude of agents trying to optimize simultaneously their respective objectives. The following bullet points describe in broad strokes some of the most frequently found types of interaction in the published literature.

\begin{itemize}

\item  \emph{Independent Learning}: \cite{MatignonLaurentLeFortPiat2007} is an early example of decentralized
RL in cooperative multi agent systems where a team
of learning robots tries to coordinate their
individual behavior to reach a coherent joint behavior. The authors of the paper propose an extension of Q-learning which they call \emph{hysteretic Q-learning}, algorithm that does not require any additional
communication between robots who can learn independently of each other and maximize an aggregate reward.
More generally, independent learning covers multi agent models in which
each agent treats the actions and rewards of the other agents as part of the environment, and learns independently using single-agent RL algorithms. While simple to implement, this strategy suffers from non-stationarity and falls short of providing expected results when the interactions between the agents are strong.

\item  \emph{Centralized Training with Decentralized Execution (CTDE)}: This paradigm leverages centralized information (such as global state or other agents' actions) while deploying decentralized execution policies. 
\cite{RashidSamvelyanSchroeder2018} and \cite{LoweWuTamar2017} are typical examples of this strategy. However, this approach
does not always resolve the scalability issue when the number $N$ of agents is large. Indeed, algorithms still depend on joint actions and global state information.

\item  \emph{Communication-Based Methods}: These approaches enable agents to exchange messages or learn communication protocols to facilitate coordination under partial observability. Although effective for small teams, communication overhead and the challenge of learning meaningful protocols limit scalability to large populations. See, for example,
\cite{FoersterAssaelFreitasWhiteson} for algorithms to learn communication protocols and still fit into the framework of centralized learning and decentralized execution. In a separate line of attack, \cite{SinghJainSukhbaatar}  develops an algorithm to help agents learn when to communicate based on scenarios, and profitability in cooperative and competitive MARL settings. See also \cite{chen2021communication} for a policy gradient algorithm skipping some steps to avoid communication overheads. 

\item  \emph{Network-Based Methods}: Beyond learning if and when to communicate, some algorithms treat the network as part of the model.
For example, \cite{MaLiDuDongYang} leverages local observation
and proposes a decentralized mechanism for scalable reinforcement learning in a large-scale network control.
 \cite{ZhangYangLiuZhangBasar} also considers the maximization of a globally expected return over a sparse and possibly time varying network. The authors propose two fully decentralized
actor-critic algorithms, which are applicable
to large-scale MARL problems in an online fashion.
The model studied in \cite{QuWiermanLi} also involves a network of agents whose states and actions interact in a local manner, and the objective is to find localized policies to maximize an expected 
global reward. Again, the authors propose a scalable actor-critic algorithm providing localized policies approximating a stationary point of the objective.
\end{itemize}

The underpinnings of the first two bullet points will resurface in our discussion of Mean Field RL strategies, though in a slightly different form, allowing us to use an important extra symmetry assumption.

\subsection{The Scalability Bottleneck}

Tabular MARL algorithms for general (unstructured) problems face fundamental scalability barriers when the number of agents becomes too large. Even if individual action spaces have the same finite cardinality, say $k=|A^i|$ for all $i\in\{1,\cdots,N\}$, and even if the size of the state space does not grow with $N$,  %
a full state--action value table has $|S||A| = |S|k^N$ entries (when $S$ is finite).
 Also, storing full state--action value tables scales poorly with the number of agents, and memory becomes a serious issue. 
 Neural networks can reduce storage costs, but exploring joint actions and learning coordinated behavior can still require substantial data and computation as $N$ grows. The dependence on $N$ is model- and algorithm-specific. 
 This is a major impediment for applications such as traffic control, robotic swarms, or crowd simulation, where hundreds or thousands of agents interact. %

The mean field approach offers a clear strategy to face the scalability curse. Obviously, it is not the only one, and while far from being universal, it allows us to tackle a specific class of models which we identify below.
Networked and localized MARL methods can exploit sparse interaction or communication graphs while retaining agent identities and allowing heterogeneous agents. In the mean field models studied here, interactions are described through population distributions under homogeneity or exchangeability assumptions.

\section{\textbf{Mean Field RL as a Scalable Alternative}}

The idea of \emph{mean field} approximation of a large system of agents acting similarly (at least in a statistical sense)  comes from statistical physics, and especially the study of large systems of \emph{identical} particles with \emph{symmetric} interactions. In MFRL, the way each agent is affected by, and interacts with, the other agents is captured by \emph{aggregate} quantities defining the behavior of the other agents, these aggregate quantities being quantified using the statistical distributions of the state of the system and the other agent actions. So instead of a complicated network of interactions, each individual agent interacts with a probability distribution, and the symmetry assumption is used to guarantee that this distribution does not depend upon the particular agent that we are tagging for the analysis.

The way we envision MFRL is captured by the statement \emph{centralized learning leading to decentralized execution}.
The policies of the $N$ agents are coordinated by a central controller. See \cite{gu2024meanmarldecentralized,gu2021meanQ}
In the cooperative MFC setting studied in this monograph, MFRL means learning a common population policy to minimize an expected social cost, with the population distribution evolving under that policy. More generally, MFC also covers control problems for a single decision maker in which the dynamics or objective depend explicitly on the law of the controlled state, as in mean--variance portfolio optimization~\cite{AnderssonDjehiche2011}. In a mean field game (MFG), a representative agent instead optimizes against a given population flow, and equilibrium additionally requires consistency of that flow with the distribution generated by an optimal response. Non-cooperative games between central planners, each accounting for the effect of its policy on population distributions, have been considered in the framework of mean field type games (MFTGs). For a broader treatment of MFTGs, see~\cite{BasarDjehicheTembine2026}.
Table~\ref{tab:learning-formulations} compares the decision makers and learning targets in these and related formulations.

\begin{table}[H]
\centering
\small
\begin{tabularx}{\textwidth}{|p{0.19\textwidth}|X|X|}
\hline
\textbf{Setting} & \textbf{Decision makers} & \textbf{Learning target} \\
\hline
MARL
& Finitely many agents with shared or individual objectives.
& An optimal cooperative policy or an equilibrium, depending on the model. \\
\hline
MFC
& A single controller, such as a central planner for a population.
& An optimal policy accounting for its effect on the state distribution. \\
\hline
MFG
& A representative agent optimizing against a given population flow.
& An optimal response and a population flow consistent with that response. \\
\hline
MFTG
& Strategic central planners, each accounting for the effect of its policy on population distributions.
& An equilibrium between the planners. \\
\hline
\end{tabularx}%
\caption{Decision makers and learning targets in multi-agent and mean field formulations.}
\label{tab:learning-formulations}
\end{table}

\subsection{Core Concept and Theoretical Foundation}

As argued in \cite{LaurierePerrinGirgin2022} and \cite{CuiTahirEkinci2022}, MFRL addresses the scalability challenge of MARL by replacing explicit modeling of $N$ individual agents with a representative agent that interacts with a probability distribution—the \emph{mean field}. Such a clean formulation of the mathematical set-up comes from the theory of  mean field games \cite{LasryLions2007,HuangCainesMalhame_NCE,CarmonaDelarue_book_I}.

At an intuitive level, one should argue that in large homogeneous populations in which interactions between agents are either weak or asymptotically symmetric, the effect of all other agents on a given agent can be well-approximated by sufficient statistics of the population distribution. Instead of keeping track of the state of the system and the behavior of the other $N-1$ agents, a good approximation of the system is to assume that each individual agent interacts with a \emph{mean field} $\mu_n\in\cP(S \times A)$, which represents the empirical distribution of the states and the actions of the agents at time $n$. This reformulation is justified by the \emph{propagation of chaos} which has been known for a long time in statistical physics. See for example \cite[Chapter 5]{CarmonaDelarue_book_I} for a detailed proof for system dynamics driven by stochastic differential equations, or Proposition \ref{pr:chaos} in Chapter \ref{ch:prelims} for a version tailored to the abstract RL frameworks of this monograph.

\vskip 4pt\noindent
At a minimum, the formulation of a MFRL should include:
\begin{itemize}
\item  A \emph{Representative Agent}: A single agent whose state dynamics and rewards depend on their own state-action pair $(s, a)$ and the mean field $\mu$;
\item  The \emph{Mean Field Dynamics}: The time evolution of the population distribution which we can formalize as $\mu_{n+1} = \cT(\mu_n,\pi)$, where $\cT$ is a form of transition operator and $\pi$ is a policy:
\item  The \emph{Optimization of a Mean Field MDP}: 
The central planner optimizes an MDP whose state is the population distribution and whose actions are feasible joint state--action distributions, as defined in Chapter~\ref{ch:prelims}.
\end{itemize}

When the solution of a MFRL reformulated in this way is used as an approximation of an $N$-agent MARL, 
directly evaluating all pairwise interactions costs $O(N^2)$ operations per population time step. If interactions depend only on a fixed number of empirical averages, forming these averages and updating all $N$ agents costs $O(N)$ operations per population time step, with constant cost per pairwise evaluation, contribution to an average and individual update. These are simulation costs; solving the control problem also depends on distribution representation and the desired accuracy
\cite{gu2024meanmarldecentralized,LaurierePerrinGirgin2022}.

\subsection{Key Assumptions}

In order to reinforce the highlights of the previous subsection, we emphasize once more that the limiting argument at the core of the MFRL's tractability relies on the following assumptions:

\begin{itemize}
\item  \emph{Agent Homogeneity and Symmetry}: in a statistical sense, agents are assumed to be \emph{exchangeable} or statistically permutation-invariant, implying that the identity of individual agents  is irrelevant in the sense that only the probability distribution of the  state-action couples matters. This is the basis on which we can characterize the agent population by a single distribution $\mu$ rather than by the state-action pairs of $N$ individual agents.

\item  \emph{Weak Coupling}: Implied in part by the above assumption is the fact that the influence of any single agent on the others is sufficiently small that aggregate effects dominate. This assumption justifies rejection of individual agent identities and focus on population-level statistics.

\item  \emph{Large Population Asymptotics}: A limiting theorem of the type known as \emph{propagation of chaos} must hold. While the de Finetti extension to exchangeable systems of the classical law of large numbers ensures that the empirical distribution of the state of the system and the actions of the agents converges to the mean field distribution $\mu$, and the contribution of any single agent to the empirical distribution becomes negligible, the propagation of chaos  theoretical result describes how the dynamical equation of the state of the finite agent system needs to be modified to describe the dynamics of the state interacting with the limiting mean field.
\end{itemize}
\subsection{Mean Field Approximation}

In the mathematical mean field theory of stochastic systems, whether we work on a Mean Field Game (MFG) \index[not]{MFG}\index[sub]{mean field game} model or a Mean Field Control (MFC) \index[not]{MFC}\index[sub]{mean field control} problem, the basic approximation result is a form of \emph{propagation of chaos} formalization and quantification. See, for example, \cite{HuangCainesMalhame_NCE} or \cite{CarmonaDelarue_book_I} for precise statements.
Typically, the mean field limit justifies the replacement of the original $N$-agent system with a limiting system where each agent faces an evolving mean field which can be deterministic or stochastic in the presence of a  \index[sub]{common noise} \emph{common noise}, i.e. random shocks affecting all the agents. The rate of convergence in the various forms of the propagation of chaos theorems provides a quantification of the approximation error based on the difference between characteristics of the $N$-agent system and their corresponding mean field limits. 
The rate depends on the quantity, metric, dimension and model assumptions; empirical Wasserstein rates are generally dimension-dependent. For conditionally i.i.d. particles, the empirical average of a fixed bounded measurable function of the state has root mean square error $O(1/\sqrt{N})$ relative to its conditional mean. 
See, for example, \cite{HuangCainesMalhame2007} or \cite{GomesMohrSouza2013} for specific examples.

\vskip 2pt
For example, in the finite-state, finite-action setting without common noise of \cite[Theorem~6.3]{gu2021meanQ}, the absolute difference of expected discounted policy values from the same initial empirical distribution is $O(1/\sqrt{N})$, for fixed regularity and discount constants satisfying that theorem's hypotheses.
This type of analysis requires detailed control of the various approximations used in the
algorithmic implementations of the theoretical limiting MFRL model to connect it to the original $N$-agent cooperative MARL. Typically, successful analyses do estimate the mean field:
\begin{itemize}
\item  by approximating the \emph{empirical distributions} by histograms or kernel density estimates of the agent states and actions;
\item  using \emph{subsampling} to approximate the mean field by appropriately selected subsets of agents, with explicit approximation error bounds.
\end{itemize}
See, for example, \cite{CarmonaLauriereTan2019} and \cite{yang2019provably} for instances of implementations of this strategy in specific models and the second part of this monograph devoted to numerical implementations.

\section{\textbf{MARL vs. MFRL: Key Differences and Trade-offs}}

Faced with the task of solving a practical problem, the choice between MARL and MFRL depends on the characteristics of the problem and modeling priorities:

\begin{table}[H]
\centering
\resizebox{0.95\textwidth}{!}{%
\begin{tabularx}{\textwidth}{|>{\hsize=.6\hsize\linewidth=\hsize}X|>{\hsize=1.2\hsize\linewidth=\hsize}X|>{\hsize=1.2\hsize\linewidth=\hsize}X|}
\hline
\textbf{Criterion } & \textbf{MARL} & \textbf{Mean Field RL} \\
\hline
\textbf{Applicability} & Finite and small $N$ for unstructured systems; possibly large $N$ for structured systems; strong pairwise couplings; heterogeneous agents;  strategic asymmetries & Very large $N$; homogeneous/exchangeable agents; interactions captured by aggregate statistics \\

\hline
\textbf{Scalability} & Full joint tables grow exponentially in $N$; structured methods can avoid them in some cases & Designed for large N; complexity scales with mean field representation \\

\hline
\textbf{Computational Complexity} & Direct all-pairs evaluation: $O(N^2)$ per population time step & Fixed-size empirical averages: $O(N)$ per population time step, including aggregation and all agent updates \\

\hline
\textbf{Representation} & Explicit joint state action or factored decentralized policies; agent-indexed & Representative agent + population distribution (histogram, parametric summary) \\

\hline
\textbf{Modeling Fidelity} & High fidelity for heterogeneity and few agents & Approximation error when agents are heterogeneous or $N$ is small \\

\hline
\textbf{Theoretical Analysis} & Convergence guarantees limited by non-stationarity & Strong convergence and approximation results under assumptions \\
\hline
\end{tabularx}
}
\caption{Main Differences between MARL and MFRL}
\end{table}

MARL offers higher modeling freedom for problems with a very small number of agents, especially when the population is heterogeneous and interactions are agent specific. However, unstructured joint representations can become computationally prohibitive at scale. On the other hand, MFRL may sacrifice some faithfulness to some of the population characteristics by assuming exchangeability (and weak coupling), but can support tractable learning with convergence guarantees under suitable assumptions in large-population settings.

\subsection{Application Domains}

MFRL emerged to make learning and control feasible in settings with large populations of agents by trading exact replication of multi-agent interactions for a tractable population-level description. This trade-off enables both algorithmic scalability and cleaner theoretical analysis, addressing the fundamental limitations of conventional MARL in large-population regimes.

MFRL is particularly well-suited to models for which assuming that agents are \emph{exchangeable} \index[sub]{exchangeable} (in a statistical sense) is reasonable. When agents are in large numbers, %
modeling the interactions with aggregate characteristics provides realistic results. Such models have been used for the following applications:
\begin{itemize}
\item  \emph{Robotic Swarms}: Coordinate hundreds of simple robots for coverage, formation control, or collective transport. See, for example, \cite{PerrinLaurierePerolatGeistEliePietquin} for an instance of MFRL application to flocking, possibly with multiple populations, or the more recent \cite{cui2023learningpomfc} providing numerical comparisons with applications of state-of-the-art MARL algorithms to Kuramoto and Vicsek swarming models; see also \cite{achdoulasry2019meancrowd,burger2013mean,aurell2019modeling,nourian2010synthesisflocking,grover2018meanhomogeneousflocking}.
\item  \emph{Traffic and Transportation}: Managing vehicle routing, traffic signal control, and ride-sharing with thousands of participants. See, for example, \cite{jusup2024safe} for an application to the vehicle positioning problem; related mean field models for routing and traffic include \cite{huang2019game,bauso2016densitynetwork,salhab2018meanroute,tanaka2020linearly,cabannes2021solving}.
See also \cite{BrunnbauerLemmelBabaieeNeubauerGrosu} for the use of off-line algorithms for crowd navigation purposes.
\item  \emph{Smart Grids and Resource Allocation}: Distributed energy management and load balancing across many devices. See, for example, \cite{HeLiu} for a hybrid model with an MFC and MFG framework for integrating Distributed Energy Resources (DER) aggregators into wholesale electricity markets, or \cite{SubramanianSerajMahajan} for a note containing an application to demand response in smart grids. Other mean field game/control applications to energy systems include \cite{alasseur2020extended,elie2019mean,bagagiolo2014mean,kizilkale2019integral,li2016mean,graber2018existence}.
\item  \emph{Epidemiology and Social Dynamics}: Modeling disease spread, opinion dynamics, or crowd behavior. See, for example, \cite{laguzet2015individualsir,hubert2018nash,elie2020contact,lee2021controllingepidemics,aurell2022optimalincentives,doncel2022meansir,stella2013opinion,bauso2016opinion,parise2023graphon}.
\item  \emph{Financial Markets}: High-frequency trading and market making with many competing agents. See, for example, \cite{BernasconiVittoriTrovoRestelli} for a model of dealer market, or \cite{QiuWangDongWangStrbac} for an example of energy trading; see also \cite{lachapelle2016efficiency,cardaliaguetlehalle2018mean,gomes2020mean,carmona2020applicationsmfgfinanceeconsurvey,achdou2014pde,achdou2017income}.
\item \emph{Security and communication networks}: botnet defense, wireless networks, distributed caching, and interference control. See, for example, \cite{meriaux2012mean,samarakoon2015energy,hamidouche2016mean,yang2017mean,kolokoltsov2016mean,kolokoltsov2018corruption}.
\end{itemize}

\subsection{Theoretical Guarantees}

On the theoretical front, MFRL takes full advantage of two decades of remarkable developments in stochastic analysis, optimal transport and analysis on Wasserstein spaces of probability measures, as well as the probabilistic approach to MFGs and MFC. In the following, we highlight some of these developments and their relevance to progress in MFRL.

\begin{itemize}
\item  \emph{Dynamic Programming and Optimality Characterization}: MFC formulations lead to dynamic programming principles and characterizations of value functions and optimal policies via Hamilton-Jacobi-Bellman equations and Pontryagin maximum principles in Wasserstein spaces of probability distributions. See for example \cite{Cardaliaguet2013} and \cite{CarmonaDelarue_book_I} for a textbook form expos\'e.

\item  \emph{Approximation Bounds}: Finite-population approximation results quantify the error between N-agent systems and mean field limits, with rates depending on the quantity, metric, dimension and model assumptions. See, for example, \cite{HuangCainesMalhame2007} or \cite{GomesMohrSouza2013}, and \cite[Chapter 5]{CarmonaDelarue_book_I} for propagation-of-chaos results under suitable conditions and in specified metrics.

\item  \emph{Convergence of Learning Algorithms}: Rigorous bounds on learning algorithm complexity and convergence guarantees are highly technical and notoriously difficult to prove in classical MARL settings, and these difficulties are even greater in the MFRL setting. 
\cite{CarmonaLauriereTan2019, yang2019provably, PerrinLaurierePerolat2020,gu2021meanQ,subramanian2019reinforcement,elie2020convergence,anahtarci2020qregu,wang2021globallqmfcg,fu2019actorcriticMFG,frikha2025actor,gu2024meanmarldecentralized,wang2020breaking,li2021permutationmarl,wu2024population} are recent works establishing sample complexity bounds and convergence guarantees for mean field Q-learning, policy gradient methods, and actor-critic algorithms in the mean field setting.

\item  \emph{Extensions Beyond Homogeneity}: The rare models of MFRL for which the fundamental assumption of homogeneity has been relaxed draw inspiration from the developments in the theory of MFGs. They include several sub-population (typically $2$) models, and major-minor agent formulations. See, for example, \cite{BonginiFornasierJunge2017} and \cite{ZamanLauriereKoppel2024}.

\item  \emph{Inclusion of Constraints}: Despite the presence of constraints, the error in the MFC approximation of a MARL can still be controlled in some specific models. See, for example, \cite{MondalAggarwalUkkusuri} for an example.
\end{itemize}

\part{Mathematical Models}

\chapter{Mathematical Preliminaries and Model Foundations}
\label{ch:prelims}

\begin{abstract}
\it The purpose of this chapter is to consolidate the notational conventions used in the subsequent chapters, and to outline the model features that underpin the classes of stochastic control problems considered in the sequel. Beyond notation, we collect the definitions and selected reference results used in later chapters. While no proof is given, precise references are provided in the last section of \emph{Notes and Complements} for the reader's convenience. 

The first class of problems is the object of Chapter \ref{ch:AAP}. It concerns \textbf{mean field control} (MFC) models in which the optimization is performed by a central unit (which we sometimes call central planner) for the benefit of an infinite population of agents subject to both idiosyncratic and common random shocks. 

The second class of problems is treated in Chapter \ref{ch:LQMFC}. It addresses \textbf{linear-quadratic (LQ) mean field control}, examining both finite-agent systems and their mean field limits. 

Although these models arise from different motivations and possess distinct structural properties, they share fundamental mathematical tools and conceptual frameworks that we develop here.
\end{abstract}

\section{Probability Spaces and Measure Theory}

We first introduce the notations and the definitions used to define and characterize the stochastic dynamical systems we plan to optimize.
\subsection{Basic Probability Structure}

Throughout the monograph, unless stated otherwise, all the random quantities are assumed to be defined on a \textbf{probability space} $(\Omega, \cF, \PP)$, where:
\begin{itemize}[leftmargin=*]
    \item $\Omega$ is the sample space,
    \item $\cF$ is a $\sigma$-algebra of events,
    \item $\PP$ is a probability measure on $(\Omega, \cF)$.
\end{itemize}

We denote by $\EE[X]$ the expectation with respect to $\PP$ of a random element $X$, by $\PP_X$ or $\cL(X)$ its distribution (also called its law), and by $\EE[X|\cG]$ and $\PP_{X|\cG}$ or $\cL(X|\cG)$ its conditional expectation and its conditional distribution given a sub-$\sigma$-algebra $\cG \subseteq \cF$.

\subsection{Borel and Polish Spaces}
\label{sub:Polish}
Our goal is to provide a rigorous mathematical presentation of concepts and theories that are too often discussed only at an intuitive level. For that reason, most of our arguments will be grounded in abstract measure theory. Throughout the monograph, a standard Borel space means a measurable space isomorphic to a Borel subset of a Polish space. For simplicity, whenever possible, we shall work directly with \emph{Polish} spaces, \index[sub]{Polish space} namely complete separable metrizable spaces. They offer tractable analysis tools and are widely adopted as a convenient setting for stochastic control problems. We mentioned the more general and more arcane notion of Borel spaces only because we will need in some rare instances to appeal to abstract measure theoretical results whose statements and proofs require such a level of generality.

For a Polish space $S$, we denote by:
\begin{itemize}[leftmargin=*]
    \item $\cB_S$ its Borel $\sigma$-field,\index[not]{$\cB_S$}
    \item $\cP(S)$ the space of probability measures on $(S, \cB_S)$.
\end{itemize}\index[not]{$\cP(S)$}
We implicitly assume that the space $\cP(S)$ is equipped with the topology of the weak convergence \index[sub]{weak convergence} of probability measures, according to which a sequence $(\mu_n)_{n\ge 0}$ in $\cP(S)$ converges toward $\mu\in\cP(S)$ if and only if
\[
\lim_{n\to\infty}\int_S \varphi \, d\mu_n = \int_S \varphi \, d\mu
\]
for every bounded continuous function $\varphi: S \to \RR$.
With this uniform structure, the space $\cP(S)$ is itself a Polish space. We shall denote by $d_w$ a distance for this Polish space structure, in other words, a distance for the weak convergence of probability measures.\index[not]{$d_w$}

\subsection{Probability Kernels and Disintegration}
\label{sub:kernels}
In most reinforcement learning studies, the dynamics of the states are prescribed by \emph{probability kernels}.\index[sub]{probability kernel}

\begin{definition}[Probability Kernel]
A probability kernel from a measurable space $(S, \cS)$ to a measurable space $(T, \cT)$ is a mapping $\kappa: S \times \cT \to [0,1]$ such that:
\begin{enumerate}[label=(\roman*)]
    \item For each $x \in S$, the map $\cT \ni B \mapsto \kappa(x, B)$ is a probability measure on $(T, \cT)$,
    \item For each $B \in \cT$, the map $S \ni x \mapsto \kappa(x, B)$ is $\cS$-measurable.
\end{enumerate}
\end{definition}

Given a probability measure $\mu$ on $S$ and a kernel $\kappa$ from $S$ to $T$, their \textbf{composition} $\mu \circ \kappa$ \index[not]{$\mu \circ \kappa$} is the probability measure on $T$ defined by
$$
(\mu \circ \kappa)(B) = \int_S \kappa(x, B) \, \mu(dx), \quad B \in \cT.
$$
Intuitively, the probability measure $\mu \circ \kappa$ is the probability on $T$ obtained by transitioning from a random distribution $\mu$ on $S$ to a distribution on $T$, $\kappa$ providing the rules of the transition.
Their \textbf{product} $\mu \measprod \kappa$ is the probability measure on $S\times T$ defined by \index[not]{$\mu \measprod \kappa$}
$$
\int_{S\times T}\varphi(x,y)(\mu \measprod \kappa)(dx,dy) = \int_S \mu(dx)\int_T \varphi(x,y)\kappa(x, dy)
$$
for all bounded measurable functions $\varphi$ on $S\times T$.

\vskip 4pt
Conditional expectations and conditional probabilities are best understood and manipulated when they can be handled with probability kernels. If $\cA\subset\cF$ is a sub $\sigma$-field of $\cF$ and $X$ is a random element in a measurable space $(S,\cS)$, by a \emph{regular conditional distribution} of $X$ given $\cA$ \index[sub]{regular conditional distribution} we mean a version of the function $\PP[X\in\cdot|\cA]$ on $\Omega\times S$ which is a probability kernel from $(\Omega,\cA)$ to $(S,\cS)$, hence a $\cA$-measurable random probability measure on $(S,\cS)$.

\begin{proposition}[Existence of regular conditional distributions]
\label{pr:regular_version}
For any Borel space $S$ and measurable space $T$, and for any random elements $X$ and $Y$ in $S$ and $T$ respectively, there exists a probability kernel $\kappa$ from $T$ to $S$ satisfying 
\[
\PP[X\in\cdot|Y]=\kappa(Y,\cdot) \qquad a.s
\]
and such a $\kappa$ is unique $\PP_Y$-a.s.
\end{proposition}
We use the notation $\PP[\cdot|Y]$ for the conditional probability given the $\sigma$-field $\cF_Y$ generated by $Y$, namely the subset of events depending only on $Y$. Similarly $\EE[\cdot|Y]$ stands for $\EE[\cdot|\cF_Y]$.
The following result should be interpreted as a generalization of the classical Fubini theorem of integration theory.

\begin{proposition}[Disintegration]
\label{pr:disintegration}
If $(S,\cS)$ and $(T,\cT)$ are measurable spaces, $\cA$ is a sub $\sigma$-field of $\cF$, and $X$ is a random element in $S$ such that the conditional probability $\PP[X\in\cdot|\cA]$ has a regular version $\nu$, then for any $\cA$-measurable random element $Y$ in $T$ and any bounded measurable function $f$ on $S\times T$, we have
$$
\EE[f(X,Y)|\cA]=\int_Sf(x,Y)\nu(dx)\qquad \PP-a.s.
$$ 
\end{proposition}

The following theorem is the most general form of disintegration of probability measures which we shall use.
\begin{theorem}[Universal Disintegration]
\label{th:universal_disintegration}
For any Borel spaces $S$ and $T$, there exists a probability kernel $\kappa: S\times \cP(S\times T)\times \cP(S)\mapsto T$ such that for any $\nu\in\cP(S)$ and $\rho\in\cP(S\times T)$ such that $\rho(\cdot\times T)$ is absolutely continuous with respect to $\nu$ it holds
\[
\rho=\nu\measprod\kappa(\cdot,\rho,\nu).
\]
Here $\kappa(\cdot,\rho,\nu)$ is unique $\nu$ a.e. for fixed $\rho$ and $\nu$.
\end{theorem}
We use these two forms of disintegration in the sequel.

\subsection{Skorohod and Blackwell-Dubins Tools}

The weak convergence of probability measures whose definition we recalled earlier is a fundamental tool in the analysis of the regularity and the asymptotic properties of stochastic structures. In the late fifties, Skorohod made the following discovery: if one is willing (or able) to change the probability space on which random elements are defined, weak convergence of random elements is not any different from almost sure convergence. This remarkable result turned out to be a powerful tool allowing all sorts of results from real analysis to be applied to random elements whose distributions converge weakly. To be specific:

\begin{theorem}[Skorohod Representation]
\label{th:Skorohod}
If $(\mu_n)_{n\ge 0}$ is a sequence of probability measures on a Polish space which converges weakly to $\mu$, there exist random variables $X_n$ and $X$, defined on a common probability space, such that:
\begin{enumerate}[label=(\roman*)]
    \item $\cL(X_n) = \mu_n$ for $n\ge 0$ and $\cL(X) = \mu$,
    \item $\lim_{n\to\infty}X_n = X$ almost surely.
\end{enumerate}
\end{theorem}\index[sub]{Skorohod representation}
This result is instrumental in the proof of convergence of value functions in stochastic control and reinforcement learning, and of continuity properties in the weak topology.

\vskip 6pt
Generalizing some of the ideas at the heart of the proof of Skorohod's representation theorem, Blackwell and Dubins developed a fundamental tool for constructing random variables with prescribed distributions. We shall use it repeatedly in the sequel.

\begin{lemma}[Blackwell-Dubins]
\label{le:BlackwellDubins}
For any Polish space $S$, there exists a measurable function 
\[
\rho_S : \cP(S) \times [0,1] \to S
\]
such that:
\begin{enumerate}[label=(\roman*)]
    \item For each $\nu \in \cP(S)$ and each uniform random variable $U \sim \text{Uniform}[0,1]$, the random variable $\rho_S(\nu, U)$ has distribution $\nu$,
    \item For almost every $u \in [0,1]$, the mapping $\nu \mapsto \rho_S(\nu, u)$ is continuous with respect to the weak topology on $\cP(S)$.
\end{enumerate}
\end{lemma}\index[sub]{Blackwell-Dubins lemma}

We shall use such a function $\rho_S$ and call it "the" \index[sub]{Blackwell-Dubins function} \emph{Blackwell-Dubins function} of the space $S$ even if it is not unique. \index[not]{$\rho_S$}Note that the atomless randomization space $([0,1],\mathcal B([0,1]),\mathrm{Leb})$ can be replaced by another atomless probability space if needed.

\section{Functional Analysis and Function Space Preliminaries}

\subsection{Spaces of Bounded Functions}

For a Polish space $S$, we denote by:
\begin{itemize}[leftmargin=*]
    \item $\cB(S)$ the vector space of bounded Borel measurable functions $\varphi: S \to \RR$. We shall implicitly equip this space with the supremum norm $\norm{\varphi}_\infty = \sup_{x\in S} \abs{\varphi(x)}$, even if the corresponding normed space does not have many desirable properties. For example, it is complete, but not necessarily separable.
    \item $\cC_b(S)$ the subspace of $\cB(S)$ formed by bounded continuous functions. This is a Banach space for the sup norm.
    \item $\cL_{sc}(S)$ the set of bounded lower semi-continuous functions.
\end{itemize}
\index[not]{$\cB(S)$}\index[not]{$\cC_b(S)$}\index[not]{$\cL_{sc}(S)$}
\begin{definition}[Lower Semi-Continuity]
A function $\varphi: S \to \RR$ is said to be lower semi-continuous if for every $\alpha \in \RR$, the level set $\{x \in S ; \varphi(x) \leq \alpha\}$ is closed. \index[sub]{lower semi-continuous}
\end{definition}

Lower semi-continuity is preserved under pointwise suprema, and is crucial for establishing existence of optimal policies.

\subsection{Contraction Mappings and Fixed Points}

\begin{definition}[Contraction Mapping]
A mapping $\tau: S \to S$ on a metric space $(S, d_S)$ is said to be a contraction if there exists $\kappa \in (0,1)$ such that
\[
d_S(\tau x, \tau y) \leq \kappa \, d_S(x, y) \quad \text{for all } x, y \in S.
\]
\end{definition}\index[sub]{contraction}

\begin{theorem}[Banach Fixed Point Theorem]
\label{th:banach}
Any contraction $\tau$ on a complete metric space $S$ has a unique fixed point, which can be obtained as the limit of iterates $\tau^n x_0$ starting from any initial point $x_0\in S$.
\end{theorem}\index[sub]{Banach fixed point theorem}

When the Dynamic Programming Principle (DPP) \index[not]{DPP} holds, optimal value functions are often characterized as solutions of an equation called the Bellman equation. Their existence and uniqueness are typically proven by identifying solutions as fixed 
points of a contraction map. So if one can show that a certain  \emph{Bellman operator} is in fact a  contraction on a Banach space of functions,  existence and uniqueness of an optimal value function will follow. We shall use this strategy in Chapter \ref{ch:AAP} and we describe its steps in broad brushstrokes in Section \ref{se:abstract_model} below.

\subsection{Matrix Spaces and Norms}

Linear Quadratic (LQ)\index[not]{LQ}\index[sub]{linear-quadratic model} models play an essential role in stochastic control, the theory of MDPs, and as a result in Reinforcement Learning (RL). The main reason is their simplicity and the fact that their optimization is solved via explicit (or at least semi-explicit) formulas. Their popularity spreads from theoretical circles to real-world applications and practitioners.  For this reason, one of the two main classes of models we consider in this monograph will be the class of LQ models. In such a framework, spaces of finite-dimensional matrices play an important role, so in preparation for their analyses, we introduce notations and definitions for some of the matrix spaces we use in the sequel. For convenience, we shall use several notions of norm. For a matrix $U \in \RR^{d \times d}$, we shall use the:
\begin{itemize}[leftmargin=*]
    \item \textbf{Spectral norm:} $\norm{U} = \sigma_{\max}(U)$, the largest singular value,
    \item \textbf{Frobenius norm:} $\norm{U}_F = \left(\sum_{i,j} U_{ij}^2\right)^{1/2} = \left(\Tr(U^T U)\right)^{1/2}$.
\end{itemize}
We shall also use the \textbf{Trace inner product:} $\inner{U, V}_{\text{tr}} = \Tr(UV^T)$.
\index[sub]{Spectral norm}\index[sub]{Frobenius norm}\index[sub]{Trace inner product}
\begin{definition}[Symmetric Matrix Spaces]
We denote by $\text{SM}_d$ the space of $d \times d$ real symmetric matrices. For $U \in \text{SM}_d$, we write:
\begin{itemize}[leftmargin=*]
    \item $U \geq 0$ if $U$ is positive semi-definite (all the eigenvalues are non-negative),
    \item $U > 0$ if $U$ is positive definite (all the eigenvalues are  strictly positive),
    \item $\lambda_{\min}(U)$ and $\lambda_{\max}(U)$ for the minimum and maximum eigenvalues, respectively.
\end{itemize}
\end{definition}
\index[not]{$SM_d$}\index[not]{$\lambda_{\min}(U)$}\index[not]{$\lambda_{\max}(U)$}
\begin{definition}[Kronecker Product]
The \textbf{Kronecker product} $U \otimes V$ of matrices $U \in \RR^{m \times n}$ and $V \in \RR^{p \times q}$ is the $mp \times nq$ block matrix
\[
U \otimes V = \begin{bmatrix}
U_{11}V & \cdots & U_{1n}V \\
\vdots & \ddots & \vdots \\
U_{m1}V & \cdots & U_{mn}V
\end{bmatrix}.
\]
\end{definition}
\index[sub]{Kronecker product}\index[sub]{Dilatation matrix}
\begin{definition}[Dilatation Matrix]
For a matrix $U \in \RR^{\ell\times d}$, its \textbf{dilatation matrix} $\tilde U$ is the $(\ell+d)\times(\ell+d)$ symmetric matrix
\[
\tilde U = \begin{bmatrix}
0 & U \\
U^\top & 0
\end{bmatrix}.
\]
It satisfies $\|\tilde U\|=\|U\|$ for the spectral norm.
\end{definition}

\section{Stochastic Dynamical Systems and Information Structures}
\label{sec:finite-agent-models}

We aim at optimizing systems evolving in time and for the purpose of implementation we shall assume that time is discrete. In other words, observations, measurements, actions, $\ldots$ , are performed at regularly spaced times labelled by natural integers $n=0,1,\cdots$. A decision at time $n$ can only be made on the basis of information available at that time. Typically, the information structure is given by a filtration \index[sub]{filtration} $\FF=(\cF_n)_{n\ge 0}$, namely a non-decreasing sequence of sub $\sigma$-fields of $\cF$: $\cF_n$ contains the events known to have happened (or not) at time $n$ or before. Since we investigate dynamical systems 
whose  states are random, we rely on models based on discrete time stochastic processes.

\begin{definition}[Stochastic Process]
A \textbf{stochastic process} is a sequence $(X_n)_{n\geq 0}$ of random variables taking values in a measurable space $(S, \cS)$. The process is said to be \textbf{adapted} to a filtration $\FF=(\cF_n)_{n\geq 0}$ if $X_n$ is $\cF_n$-measurable for each $n$, in other words, if the values of $X_n$ only depend upon events in $\cF_n$.
\end{definition}
\index[sub]{stochastic process}\index[sub]{adapted}
We shall say that the process is a noise process if the random variables $X_n$ are independent of each other. \index[sub]{noise}
A central feature of both models subsequently studied in this monograph is the presence of two types of noise processes.

\begin{definition}[Idiosyncratic Noise]
In a model with $N$ agents, a double sequence $(\varepsilon_n^i)_{n\geq 0, i=1,\cdots,N}$ of i.i.d.\ random variables with common distribution $\nu_\varepsilon \in \cP(E)$, where $(E,\cE)$ is a measurable space, is called an idiosyncratic noise sequence if the random shocks $(\varepsilon_n^i)_{n\geq 0}$ affect individual agent $i$ independently of each other.
\end{definition}
\index[sub]{idiosyncratic noise}\index[sub]{common noise}
\begin{definition}[Common Noise]
Irrespective of the number of agents in the model, a sequence $(\varepsilon^0_n)_{n\geq 0}$ of i.i.d.\ random variables with common distribution $\nu_{\varepsilon^0} \in \cP(E^0)$, where $(E^0,\cE^0)$ is a measurable space, is called a common noise if the random shocks it creates affect all agents simultaneously.
\end{definition}

In models including both idiosyncratic and common noise processes, we shall assume that 
the sequences $(\varepsilon^i_n)_{n\geq 0, i=1,\cdots,N}$ and $(\varepsilon^0_n)_{n\geq 0}$ are mutually independent.

\vskip 2pt
The goal of this dual noise structure is to capture both individual heterogeneity and aggregate uncertainty, which are essential features of many economic and engineering systems.

\subsection*{Filtrations and Information Structures}

The information available to the agents evolves over time as captured by the choice of filtrations for information structure. 
These filtrations vary according to the model and the type of decision maker.
\begin{itemize}[leftmargin=*]
        \item \textbf{Open loop policies:} When the model is set up with open-loop controls, the $\sigma$-field $\cF_n$ will typically include the realizations of all past random noise terms (both idiosyncratic and common, when appropriate) and the sources of randomness used to randomize past actions.
    
    \item \textbf{Agent-level filtration:} In the case of finitely many agents, the $\sigma$-field $\cF^i_n$ representing the information available to agent $i$ at time $n$, includes their own state history (or the current value of their state in the case of Markovian policies), and possibly observations of the population distribution. In the mean field limit,  the $\sigma$-field $\cF^c_n$, representing the information available to a generic agent at time $n$ includes the history up to time $n$ of a generic state and observations of the population distribution.
\end{itemize}
For each of the specific models studied in the sequel, we shall specify explicitly the various sources of randomness, including noise sequences and randomization devices, used in the definitions of the filtration we use.

\subsection{Dynamics Prescription and Objective Optimization}

If the state of the system is given by a generic element $x$ of a Polish space $S$ and the actions agents are allowed to take are elements of a Polish space $A$, the stochastic dynamics of the state of the system are given by a stochastic process $(X_n)_{n\ge 0}$ in $S$. Historically, and especially in the existing RL literature, the rules governing the time evolution of the state are given by an initial distribution $\mu_0\in\cP(S)$ and a sequence $(p_n)_{n\ge 0}$ of probability kernels from $S\times A$ into $S$. Intuitively speaking, $\mu_0=\PP_{X_0}$ is the initial distribution of the state, and for each integer $n\ge 0$, state-action couple $(x,a)$ and measurable set $B\in\cB_S$, the quantity $p_n\bigl((x,a),B\bigr)$ gives the probability that the state $X_{n+1}$ is in the set $B$ whenever (i.e. conditioned by the fact that) the state at time $n$ is $x$, i.e. $X_n=x$ and the action taken at that time was $a\in A$. 

\vskip 2pt
Mostly for mathematical convenience, some form of stationarity in time is assumed, so that the probability kernel $p_n$ is independent of $n$.

\vskip 2pt
In order to define the dynamics of the state using a system function instead of a transition probability kernel, we use the Blackwell-Dubins function $\rho_S$ of the Polish space $S$ (recall Lemma \ref{le:BlackwellDubins}) and define, independently of the time $n\ge 0$, the function $F:S\times A\times [0,1]\mapsto S$ by 
\[
F(x,a,u)=\rho_S\bigl(p_n\bigl((x,a),\,\cdot\,\bigr),u\bigr), 
\]
so that, 
\begin{itemize}
\item if $(a_n)_{n\ge 0}$ is a prescribed sequence of actions in $A$,
\item if $(U_n)_{n\ge 0}$ is an i.i.d. sequence of uniform random variables in $[0,1]$, 
\item the stochastic process $(X_n)_{n\ge 0}$ defined by
\[
X_0=\rho_S(\mu_0,U_0)\quad\text{and}\quad X_{n+1}=F(X_n,a_n,U_{n+1}),\qquad n=0,1,\cdots
\]
\end{itemize}
satisfies all the requirements prescribed above in the description of the dynamics of the state of the system.
Since the original description of the state dynamics did not specify the information structure, the above use of the Blackwell-Dubins map allows us to specify the state dynamics by means of a sequence of independent random shocks (which we shall sometimes call \emph{innovations}) \index[sub]{innovation} and a (deterministic) function mapping a triple $(x,a,u)$ of a state, an action and a random shock into another state. As explained in the previous chapter, we call such a function $F$ a \emph{system function}.\index[sub]{system function}

\vskip 4pt
When the state of the system is affected by only one controller, the goal of the latter is to choose a sequence $(a_n)_{n\ge 0}$ of actions in order to minimize an expected cost aggregated over time.
In practice, the sequence $(a_n)_{n\ge 0}$ is random since at each time $n$,  the individual action $a_n$ chosen by the controller, is picked as a function of the information available to the controller at that time. Specifying the exact nature of the available information will force the measurability of the $A$-valued random variables $\alpha_n$ from which the actions $a_n$ are sampled. So the objective function can be viewed as being of the form
\[
J = \EE\left[\sum_{n=0}^\infty \gamma^n f(X_n, \alpha_n, \epsilon_n)\right],
\]
where 
\begin{itemize}
\item $f:S\times A\times E \mapsto \RR$ is a bounded measurable function, $f(x,a,\epsilon)$ giving the one stage cost incurred when the state is $x\in S$, the action taken by the controller is $a\in A$ and the realization of the noise (random shock) is $\epsilon\in E$;
\item $\gamma\in (0,1)$ is a discount factor guaranteeing that the infinite sum converges, 
\end{itemize}
the expectation being over the random values of the states $X_n$, all the sources of randomness: initial conditions, the realizations $\epsilon_n$ of the random shocks, and the random actions $\alpha_n$ from which the actions $a_n$ are sampled, $\ldots$.

\subsection*{Multi-Agent Model Definitions}

When the system involves $N$ agents ($N$ is a finite integer $N\ge 1$), each agent $i$ chooses an action $a^i_n\in A^i$ at time $n$, the action space in which agent $i$ chooses their action being possibly different from agent to agent. Accordingly, each agent may have a different one stage cost function and agent $i$ aims at minimizing the cost
\[
J^i(\balpha) = \EE\left[\sum_{n=0}^\infty \gamma^n f^i(X_n, \balpha_n, \epsilon^i_n)\right],
\]
where 
\begin{itemize}
\item $(\epsilon^i_n)_{n\ge 0, i=1,\cdots,N}$ is an idiosyncratic noise process,
\item $\bA=A^1\times\cdots\times A^N$ is the set of action profiles,
\item $\balpha=(\balpha_n)_{n\ge 0}$ is a stochastic process in the space $\bA$ of action profiles,
\item $f^i:S\times \bA\times E \to \RR$ is a bounded measurable function, $f^i(x,\ba,\epsilon)$ giving the one stage cost incurred by agent $i$ when the state is $x\in S$, the actions taken by all the agents are  $\ba=(a^1,\cdots,a^N)\in \bA$ and the realization of the noise (random shock) is $\epsilon^i\in E$;
\item and as before, $\gamma\in (0,1)$ is a discount factor guaranteeing that the infinite sum converges.
\end{itemize}
Notice that the interactions between the agents enter the costs $J^i$ through the values $X_n$ of the state of the system and the actions $a^j_n$ of all the agents. We mentioned earlier that we were assuming some form of stationarity in time allowing us to use one stage cost functions independent of the stage index $n$. Now, and again purely for mathematical convenience,  we shall assume a form of symmetry among the agents and use the same one-stage cost function $f$ for all the agents. As emphasized earlier, this exchangeability between the agents is crucial in the models we want to study, and it will be especially useful when we consider the mean field limit $N\to\infty$.

\vskip 4pt
In many $N$-agent models, the state $X_n$ of the system at time $n$ is best described by an $(N+1)$-tuple $X_n=(X^0_n,X^1_n,\cdots,X^N_n)$ where the stochastic process $(X^0_n)_{n\ge 0}$ with values in a Polish space $S^0$ captures stochastic factors common to all the agents, and for each $i=1,\cdots,N$, the stochastic process $(X^i_n)_{n\ge 0}$ with values in a Polish space $S^i$ captures stochastic factors specific to agent $i$. When this is the case, $X^i_n$ can be interpreted as the value at time $n$ of the \emph{private state} of agent $i$, and the state space of the system can be factorized in the form $S=S^0\times S^1\times\cdots\times S^N$. Accordingly, the system function $F$ splits in a form $F=(F^0, F^1,\cdots, F^N)$, the components being interpreted as individual system functions $F^j$ for $j=0,1,\cdots,N$ 
describing the time evolutions of the common factor and the individual private states 
\[
X^j_{n+1}=F^j(X_n,a_n,\epsilon^j_n),\qquad n=0,1,\cdots,\quad j=0,1,\cdots,N.
\]
For the desirable assumption of \emph{symmetry} among the agents to hold, we shall often assume that the individual system functions $F^i$ are the same for $i=1,\cdots,N$.

\subsection*{Competition or Cooperation?}

Even though we are not yet ready to give the precise details of how the agents choose their actions, in other words
the precise definitions of the filtrations $\FF^i=(\cF^i_n)_{n\ge 0}$ with respect to which the control processes $\balpha^i=(\alpha^i_n)_{n\ge 0}$ need to be adapted for the sequences of sampled actions $\ba^i=(a^i_n)_{n\ge 0}$
to be admissible, we can still describe the major two forms of optimization.

\vskip 4pt
Whether or not the individual system functions are identical, each agent $i$ has their own objective cost $J^i$ and this form of multivariate optimization requires formalizing the behavior of the agents.

\vskip 2pt
If the agents compete among themselves and optimize their objectives independently, the most frequently searched outcome is given by the notion of \emph{Nash equilibrium}. \index[sub]{Nash equilibrium}
\begin{definition}[Competitive Equilibrium]
\label{de:competitive}
A control process profile  $\balpha=(\balpha^i)_{i=1,\cdots,N}$ of admissible control processes $\balpha^i=(\alpha^i_n)_{n\ge 0}$ is said to form a Nash equilibrium if for any $i\in\{1,\cdots,N\}$ and any alternative control process $\bbeta=(\beta_n)_{n\ge 0}$ for player $i$,
\[
J^i(\balpha^{-i},\bbeta)\ge J^i(\balpha).
\]
\end{definition}
Here (and in the following) we use the notation $(\balpha^{-i},\bbeta)$ for the control process profile where the $i$-th component of $\balpha$ is replaced by $\bbeta$. In other words, when the $N$ agents use the control processes $\balpha=(\balpha^i)_{i=1,\cdots,N}$, no single agent $i$ can be better off by deviating alone from using their control $\balpha^i=(\alpha^i_n)_{n\ge 0}$.

\vskip 4pt
An alternative is to assume that all the agents cooperate and aim at minimizing an aggregate or average cost
\begin{equation}\label{fo:social_cost_N}
    \tilde J(\balpha)=\frac1N\sum_{i=1}^NJ^i(\balpha),
\end{equation}
which is often referred to as \emph{social cost}. \index[sub]{social cost} The factor $1/N$ does not play any role when the number of agents is fixed. However, it will facilitate the analysis of the limit $N\to\infty$ when we study large populations of agents.

\begin{definition}[Cooperative Equilibrium]
\label{de:cooperative}
A control process profile  $\hat\balpha=(\hat\balpha^i)_{i=1,\cdots,N}$ of admissible control processes $\hat\balpha^i=(\hat\alpha^i_n)_{n\ge 0}$ is said to form a cooperative equilibrium if 
\[
\hat\balpha = \argmin_{\balpha} \tilde J(\balpha).
\]
\end{definition}
Because the objective function to minimize is an aggregate quantity, one should envision the actual optimization to be performed by a single entity that would identify the best strategy profile and dispatch the individual optimal strategies to the agents. We shall call this entity the central unit, in line with the terminology \emph{central planner} used in macro-economics. \index[sub]{central planner}\index[sub]{central unit}

\vskip 4pt
For all the reasons already articulated in Chapter~\ref{ch:review}, the thrust of this monograph is the study of models of large populations of cooperative agents. As a result, we shall concentrate on cooperative equilibrium, and only mention Nash equilibria when appropriate, essentially to provide insightful comparisons.

\subsection*{Including a Common Noise}

Many realistic situations cannot be modeled with idiosyncratic noise as the only source of random shocks. In order to 
add random factors and random shocks affecting simultaneously all the agents, we can use the flexibility of the system function to add an extra variable from another Polish space, say $E^0$, and use an i.i.d. sequence of random variables $\bepsilon^0=(\epsilon^0_n)_{n\ge 0}$ to include these common random effects in the dynamics of the state of the system. Accordingly, the dynamics of the private state of the $i$-th agent is now given by:
\begin{equation}
\label{fo:dynamics_N}
X^i_{n+1}=F^i(X_n,\alpha_n,\epsilon^i_n,\epsilon^0_n),\qquad n=0,1,\cdots.
\end{equation}
Note that the important Proposition~\ref{pr:chaos} below, referred to as the \emph{propagation of chaos} \index[sub]{propagation of chaos} is stated in this very setting including a common noise.
Notice that in the same way, the one-stage cost function of agent $i$ can also be made dependent upon the common noise and read:
\begin{equation}
\label{fo:cost_N}
f^i(X_n,\alpha_n,\epsilon^i_n,\epsilon^0_n).
\end{equation}

\subsection*{Randomization and Mixed Strategies}

Whether the motivation is to help the theoreticians have more tools to prove existence of optimal strategies, or to help the individual agents explore better the state space in hope of finding more favorable regions with lower costs, randomization is a time honored trick used very frequently in stochastic optimization. In the context of continuous time stochastic control, it relies on the introduction of \emph{relaxed controls} \index[sub]{relaxed controls}, altering the dynamics while keeping the same optimal values. In the discrete time setting considered in this monograph, the reinforcement learning practice is to choose actions by independently sampling from probability distributions on the action space. 
So instead of using adapted stochastic processes with values in the action spaces $A^i$, the controls processes $(\alpha^i_n)_{n\ge 0}$ are replaced by adapted stochastic processes $(\pi^i_n)_{n\ge 0}$ with values in the spaces $\cP(A^i)$ of probability measures on the action spaces $A^i$. These measure-valued stochastic processes are called \emph{mixed controls}. \index[sub]{mixed control} Consequently, the action values $\alpha^i_n$ actually used in the state update formula \eqref{fo:dynamics_N} and the cost evaluation formula \eqref{fo:cost_N} are random samples from the probability distributions $\pi^i_n$. Including sampling from probability distributions in the theoretical models naturally calls for the Blackwell-Dubins lemma, and indeed, the implementation of this idea of randomization in the theoretical analysis of the models will rely on this fundamental result from measure theory.

\paragraph{Running example: cybersecurity.}
We consider a cybersecurity model which will be revisited in detail in Chapter~\ref{ch:numeric_I}. A computer is defended or undefended (D or U), and infected or susceptible (I or S), so the individual state space is $S=\{DI,DS,UI,US\}$ and the action space is $A=\{0,1\}$. Action $1$ requests a change of protection status; action $0$ requests no change.
For $N$ computers, write $\mu_n^N(x)=N^{-1}\sum_{i=1}^N\mathbf{1}_{\{X_n^i=x\}}$. The planner chooses four probabilities $p_n^{\mathrm{sw}}(x)\in[0,1]$: conditional on the states and this shared choice, computers draw their actions independently, with $\alpha_n^i\sim\mathrm{Bernoulli}(p_n^{\mathrm{sw}}(X_n^i))$.
Let $P_{xy}^a(\mu,z)$ be a discrete-time transition kernel for the computer state, where $z$ is a common shock to attack intensity. Given the current states, actions and the fresh common shock $\varepsilon_{n+1}^0$, the next states are conditionally independent, with probabilities
\[
\PP\bigl(X_{n+1}^i=y\mid (X_n^j,\alpha_n^j)_{j=1}^N,\varepsilon_{n+1}^0\bigr)
=P_{X_n^i y}^{\alpha_n^i}(\mu_n^N,\varepsilon_{n+1}^0).
\]
With protection cost $k_D>0$ and infection penalty $k_I>0$, the average one-stage cost is
\begin{equation}
\label{eq:running-cyber-cost}
c(\mu)=k_D\bigl(\mu(DI)+\mu(DS)\bigr)+k_I\bigl(\mu(DI)+\mu(UI)\bigr),
\end{equation}
where $\mu(x)=\mu(\{x\})$. The planner minimizes $\EE[\sum_{n\ge0}\gamma^n c(\mu_n^N)]$.

\section{Mean Field Models and their Limiting Behavior}

In this section, we still consider a multi-agent model with $N$ agents, and we assume stationarity in time and symmetry among the agents as we mentioned earlier. First, we review the main properties which we require for the class of \emph{Mean Field} models.

\subsection{Symmetry Requirements}
\label{sub:symmetry}

A first step in enforcing the symmetry among the agents is to assume that all the individual state spaces are the same, say $S^i=S$, that all the action spaces are the same, say $A^i=A$, that all the 
system functions are the same, i.e. $F^i=F$, and all the one-stage cost functions are the same, i.e. $f^i=f$ for all $i\in[N]$. In fact, we require a stronger form of symmetry by assuming that the interactions between the individuals are symmetric in the following sense: for each agent $i\in[N]$, we assume that the dynamics of the state $X^i_n$ is given by equation \eqref{fo:dynamics_N} with a common system function $F$, the dependence of the right hand side upon $X_n=(X^1_n,\cdots,X^N_n)$ being given by a function of $X^i_n$ and a symmetric function of the $X^j_n$ for $j\ne i$.
Following the intuition that symmetric functions of a large number of variables are essentially functions of their empirical distribution (see for example \cite[Lemma 1.2]{CarmonaDelarue_book_I} for a mathematical justification), we can imagine that because of the symmetry in the interactions between the states of the agents, the dependence of the system function $F$ on $X_n=(X^1_n,\cdots,X^N_n)$ could be a function of $X^i_n$ and the empirical measure 
\[
\mu^X_n=\frac1N\sum_{j=1}^N\delta_{X^j_n}.
\]
 Similarly, one can imagine that the dependence of  $F$ on $\alpha_n=(\alpha^1_n,\cdots,\alpha^N_n)$ could be through a function of $\alpha^i_n$ and the empirical measure 
\[
\mu^\alpha_n=\frac1N\sum_{j=1}^N\delta_{\alpha^j_n}.
\]
In fact, when we say that the agents are in mean field interaction, we shall mean that the dynamics of the private states of the agents are given by an equation of the form
\[
X^i_{n+1}=F\Bigl(X^i_n, \alpha^i_n, \frac1N\sum_{j=1}^N\delta_{(X^j_n,\alpha^j_n)}, \varepsilon^i_{n+1}, \varepsilon^0_{n+1}\Bigr)
\] 
for a system function $F$ defined on the product space $S\times A\times \cP(S\times A)\times E\times E^0$. De facto, the interactions between the agents are modeled via the empirical distribution of the state-action couples $(X^j_n,\alpha^j_n)$. Similarly, we can assume that the common individual one-stage cost functions $f$ have the same structure.

\subsection{Limit \texorpdfstring{$N\to\infty$}{N -> infinity} and Propagation of Chaos}
\label{sub:chaos}

Although we added significant sources of coupling between the agents experiences, exchangeability among the agents persists, and in the limit $N\to\infty$, 
a conditional form of the propagation of chaos should lead to conditional independence of  individual dynamics.
This asymptotic conditional independence is quite intuitive if one is familiar with the standard theory of the propagation of chaos for It\^o stochastic differential equations. However, stating the discrete time conditional form of this result precisely with full mathematical rigor is rather difficult, and it requires new notation.

\vskip 6pt\noindent
\textbf{Special Notation:} For any sequence $(\zeta_n)_{n\ge 0}$ we shall use the notation $\underline{\zeta}_n=(\zeta_0,\zeta_1,\cdots,\zeta_n)$.

\vskip 6pt
Even if we do not call them by their names, in what follows, we introduce the elements of the definition of open-loop policies which is given later on in Definition \ref{de:open_loop}. 

We fix a sequence $\bpi=(\pi_n)_{n\ge 0}$ of deterministic measurable functions $\pi_n:\Xi_n\to \cP(A)$  where $\Xi_n = \cY\times (\Theta\times\Theta^0\times E\times E^0)^n\times \Theta^0$ where $\cY$, $\Theta$ and $\Theta^0$ are Polish spaces from which samples are drawn to randomize the initial condition, the individual action choices and the common mixed action.
Next, for each integer $i\ge 1$, which we may want to think of as the index of an agent in an infinite population, we consider a random variable $Y^i$ in $\cY$ to randomize the initial condition, and i.i.d. sequences $(\theta^i_n)_{n\ge 0}$ to randomize the action choices of the individual agents, $(\theta^0_n)_{n\ge 0}$ to randomize the action choices of the central unit, $(\varepsilon^i_n)_{n\ge 0}$ as realizations of the idiosyncratic random shocks, and $(\varepsilon^0_n)_{n\ge 0}$ as realizations of the common random shocks. The families indexed by $i$ are independent across agents and independent of the shared common sequences. 

Motivated by the classical form of the theory of propagation of chaos, we anticipate the form of the statistical behavior of the population of agents in the limit $N\to\infty$. So for each $i\ge 1$ fixed, we consider the stochastic process $(X^{\bpi,i}_n)_{n\ge 0}$ describing the time evolution of the $i$-th state given by: 

\begin{equation}
\label{fo:MFC_individual_dynamics}
X^{\bpi,i}_{n+1}=F(X^{\bpi,i}_n,\alpha^{\bpi,i}_n,\PP^0_{(X^{\bpi,i}_n,\alpha^{\bpi,i}_n)},\varepsilon^i_{n+1},\varepsilon^0_{n+1}),\quad n\ge 0,
\end{equation}
where the random action $\alpha^{\bpi,i}_n$ is chosen as
\begin{equation}
    \label{fo:alpha_n}
\alpha^{\bpi,i}_n=\pi_n\bigl(Y^i,\underline\theta^i_n,\underline\theta^0_n,\underline\varepsilon^i_n,\underline\varepsilon^0_{n+1}\bigr),\qquad n\ge 0.
\end{equation}
Note that all processes $(X^{\bpi,i}_n)_{n\ge 0}$ are conditionally independent given the shared common sequences and they follow the same dynamical equation
\begin{equation}
\label{fo:MFC_dynamics}
X_{n+1}=F(X_n,\alpha_n,\PP^0_{(X_n,\alpha_n)},\varepsilon_{n+1},\varepsilon^0_{n+1}),\quad n\ge 0,
\end{equation}
where the notation $\PP^0_{(X_n,\alpha_n)}$ stands for the conditional law of the state-action couple $(X_n,\alpha_n)$ given the common shocks $(\varepsilon^0_k)_{1\le k\le n}$ and common policy randomization $(\theta^0_k)_{0\le k\le n}$, and where the action process $(\alpha_n)_{n\ge 0}$ is obtained from $\bpi=(\pi_n)_{n\ge 0}$ and a sequence of randomization random variables and idiosyncratic and common noise sequences as in \eqref{fo:alpha_n}. 

Accordingly, each agent $i$ of this infinite population of conditionally independent agents, should solve the following optimization problem:
\begin{equation}
\label{fo:MFC_optimization}
\inf_{\balpha=(\alpha_n)_{n\ge 0}} \EE\Bigl[\sum_{n=0}^\infty \gamma^n f\bigl(X_n,\alpha_n,\PP^0_{(X_n,\alpha_n)}\bigr)\Bigr].
\end{equation}
The intuition that such set-up is the proper asymptotic regime for a symmetric population of $N$ agents when $N\to\infty$ is substantiated by the following mathematical result. By constructing the $N$-agent model on the same probability structure (i.e. with the same sequences of randomization random variables and idiosyncratic and common random shocks), it is possible to prove that this convergence is actually in the almost sure sense.

\begin{proposition}
\label{pr:chaos}
In the set-up described above, if we assume that for each fixed $(x_0,a,\nu^0, e^0)\in S\times A\times\cP(S\times A)\times E^0$, the map
$$
S\times \cP(S\times A)\ni (x,\mu)\mapsto F(x,a,\mu,e,e^0)\in S
$$
is continuous for almost every $e\in E$ with respect to $\nu$, and if for each integer $N$ we define the $N$-agent state dynamics by:
\begin{equation}
\label{fo:MFC_N_individual_dynamics}
X^{\bpi,i,N}_{n+1}=F(X^{\bpi,i,N}_n,\alpha^{\bpi,i}_n,\frac1N\sum_{j=1}^N\delta_{(X^{\bpi,j,N}_n,\alpha^{\bpi,j}_n)},\varepsilon^i_{n+1},\varepsilon^0_{n+1}),\quad n\ge 0,
\end{equation}
then for any $\bpi$, $i$ and $n$, we have
\begin{itemize}
\item  $\lim_{N\to\infty}X^{\bpi,i,N}_n=X^{\bpi,i}_n$ a.s.
\item $\lim_{N\to\infty}d_w\Bigl(\frac1N\sum_{j=1}^N\delta_{(X^{\bpi,j,N}_n,\alpha^{\bpi,j}_n)},\PP^0_{(X^{\bpi,i}_n,\alpha^{\bpi,i}_n)}\Bigr)=0$ a.s.
\end{itemize}
Moreover, if we also assume that for each fixed $a\in A$, the map
$
S\times \cP(S\times A)\ni (x,\mu)\mapsto f(x,a,\mu)\in \RR
$
is continuous and that $f$ is bounded, then $\lim_{N\to\infty}J_N^{\bpi}=J^{\bpi}$, where $J_N^{\bpi}$ is the expected discounted social cost $\tilde J$ defined in~\eqref{fo:social_cost_N} of the $N$-agent model under $\bpi$, and $J^{\bpi}$ is the expected discounted cost inside the infimum in~\eqref{fo:MFC_optimization}, evaluated under $\bpi$.

\end{proposition}
Recall that the distance $d_w$ is the metric of the weak convergence of probability measures introduced in Subsection \ref{sub:Polish}. In the above result, for each fixed policy and time, the $i$-th tagged-state and weak empirical-measure limits hold almost surely under the above coupling, also conditionally almost surely for almost every realization of the common shocks and common policy randomization.

\subsection*{More Measurability Considerations}
The form of Mean Field discrete time Control problem with common noise given by \eqref{fo:MFC_dynamics} and \eqref{fo:MFC_optimization} is now very close to the model we investigate in the next chapter.
The last question we would like to address before turning to the theoretical analysis of the model is:
\begin{quote}  
\emph{Could the conditional distributions $\PP^0_{(X_n,\alpha_n)}$ depend upon some extra sources of randomness beyond the randomness of the common noise?}
\end{quote}
Indeed, if individual agents and the central unit performing the social cost optimization are allowed to use mixed strategies, they need independent sources of randomness to randomize their actions and decisions at each time. So since the cost optimization and the choice of action strategy are performed by such a central unit,
\begin{itemize}
\item the central unit needs a source of randomness to \emph{randomize} the choice of a mixed policy to be dispatched to the individual agents, and
\item the individual agents need to sample their actions from the policy sampled for them by the central unit.
\end{itemize}
In the limit $N\to\infty$, we expect that the \emph{idiosyncratic} randomizations of the individual agents average out. On the other hand, we do not see why the \emph{central unit randomization} should also average out, so we should expect that it remains.
The motivation for the common randomization is two-fold. First, in full analogy with single-agent RL, we expect the common randomization to help exploring the state space and the action space during learning. This point will be of crucial importance when we show that the mean field Q-learning algorithm converges under the assumption that the state-action sequence covers repeatedly the whole state-action space. Moreover, this form of randomization will also prove useful for technical reasons.

As a result, the conditioning in $\PP^0_{(X_n,\alpha_n)}$ should be with respect to the past values of the common noise $\varepsilon^0_1,\cdots,\varepsilon^0_n$ \emph{as well as} all the sources of randomness used by the central unit to randomize the policy handed out to the individual agents for implementation. We shall provide rigorous mathematical definitions for this conditioning in the next chapter.

\paragraph{Running example: the mean field state.}
In the cybersecurity example, the population state becomes a probability vector $\mu_n\in\cP(S)$, conditional on the common shocks and the planner's common randomization. Its update is
\begin{equation}
\label{eq:running-cyber-update}
\mu_{n+1}(y)=\Phi_y(\mu_n,p_n^{\mathrm{sw}},\varepsilon_{n+1}^0),\qquad
\Phi_y(\mu,p^{\mathrm{sw}},z)=\sum_{x\in S}\mu(x)\bigl[(1-p^{\mathrm{sw}}(x))P_{xy}^0(\mu,z)+p^{\mathrm{sw}}(x)P_{xy}^1(\mu,z)\bigr].
\end{equation}
The independent action and transition draws have been averaged out. The common shock still changes the whole probability vector, and a randomized choice of $p_n^{\mathrm{sw}}$ is also shared by the population. The limiting social cost is $\EE[\sum_{n\ge0}\gamma^n c(\mu_n)]$, with the same function $c$ as in~\eqref{eq:running-cyber-cost}.

\section{A General Abstract MFC Model}
\label{se:abstract_model}

The next two chapters are devoted to the complete mathematical analyses of two fundamental models
whose definitions we make explicit in this section for the sake of future reference. Both models will be set in infinite horizon for the sake of definiteness, and we shall use the notation $\gamma\in(0,1)$ for the \emph{discount factor} \index[sub]{discount factor} in both cases.
We start with the model that we study in Chapter \ref{ch:AAP}.

\begin{definition}[MFC Model]
\label{de:mfc_model}
Our first model of \textbf{infinite-horizon discounted mean field control (MFC) with common noise} is specified by a tuple
\[
(S, A, E, E^0, F, f),
\]
where:
\begin{itemize}[leftmargin=*]
    \item $S$ is a Polish space (the state space for individual agents),
    \item $A$ is a Polish space (the action space for individual agents), 
    \item $E$ is a Polish space (the idiosyncratic noise space),
    \item $E^0$ is a Polish space (the common noise space),
    \item $F: S \times A \times \cP(S \times A) \times E \times E^0 \to S$ is the system function, Borel measurable, describing the state transition dynamics,
    \item $f: S \times A \times \cP(S \times A) \to \RR$ is the one-stage cost function, bounded and Borel measurable,
\end{itemize}
\end{definition}

The model aims at describing the time evolution of a \textbf{homogeneous population of agents}, the dynamics of the state of a generic agent being given by:
\[
X_{n+1} = F(X_n, a_n, \mu_n, \varepsilon_{n+1}, \varepsilon^0_{n+1}),
\]
where at each time $n$:
\begin{itemize}[leftmargin=*]
    \item $X_n \in S$ is the state of the generic agent,
    \item $a_n \in A$ is the action taken by the agent,
    \item $\mu_n \in \cP(S \times A)$ is the joint distribution of states and actions across the population,
    \item $(\epsilon_n)_{n\ge 1}$ is an independent identically distributed (i.i.d. for short) sequence of random elements in $E$ whose common distribution we denote by $\nu$, $\varepsilon_{n+1} \in E$ representing the next value of the idiosyncratic random shock,
    \item $(\epsilon^0_n)_{n\ge 1}$ is an i.i.d. sequence of random elements in $E^0$ whose common distribution we denote by $\nu^0$, $\varepsilon^0_{n+1} \in E^0$ is the next value of the common noise.
\end{itemize}

The one-stage cost incurred at time $n$ by the generic agent is $f(X_n, a_n, \mu_n)$. Notice that it may depend on the population distribution $\mu_n$, capturing \textbf{mean field interactions}.

The goal of the analysis is the minimization of the expected discounted social cost
\begin{equation}
    \label{fo:J}
J = \EE\left[\sum_{n=0}^\infty \gamma^n f(X_n, a_n, \mu_n)\right],
\end{equation}
where the expectation is over all the sources of randomness: initial conditions, noise sequences, and sources of randomness used for policy randomization. 
Notice that we did not specify yet the space of actions - controls - policies over which the optimization is performed. We shall discuss the admissible information structures later on.

\subsection*{Assumptions for MFC with Common Noise}
The results stated and proven in Chapter \ref{ch:AAP} hold under various sets of assumptions. Still, for future reference, we list here the assumptions that will be in force throughout the analysis of the model. In some instances, more restrictive assumptions will be added in order to derive specific results.

\begin{assumption}[Continuity]\label{ass:H1}
\begin{enumerate}[label=(\roman*)]
    \item For each fixed $(e, e^0) \in E \times E^0$, the function $F(\cdot, \cdot, \cdot, e, e^0)$ is continuous on $S\times A\times \cP(S\times A)$.
    \item The one-stage cost function $f: S \times A \times \cP(S \times A) \to \RR$ is continuous.
\end{enumerate}
\end{assumption}
Since $S$ and $A$ are assumed to be Polish spaces, $S\times A$ and $\cP(S\times A)$ are also Polish spaces, and the notion of continuity alluded to above is understood with respect to the corresponding weak topology on probability measures.

\begin{assumption}[Compactness]\label{ass:H2-mfc}
The state space $S$ and action space $A$ are compact metric spaces.
\end{assumption}

\subsection{Policies and Policy Spaces}
\label{sub:policies}

Policies are functions with values in the action space $A$, or possibly in the space $\cP(A)$ of probability distributions on $A$ when using randomized actions i.e. mixed strategies. Their measurability properties offer a convenient way to specify the information available to the agent making a decision. Indeed, this information is coded by the arguments of the policy function.
Because we shall use several classes of policies, we state their definitions for future reference.
\vskip 6pt\noindent
Recall the special notation introduced earlier: for any given sequence $(\zeta_n)_{n\ge 0}$ we use  $\underline{\zeta}_n=(\zeta_0,\zeta_1,\cdots,\zeta_n)$.
\vskip 6pt
We want to think of minimization of the objective \eqref{fo:J} as the job of a central unit. As suggested by applications in economics, we shall sometimes use the terminology \emph{central planner}. Still, the decisions made by the central unit will have to be implemented by the individual agents, typically the \emph{generic agent} in the mean field limit. Both the central unit and the generic agent should be able to randomize their choices, so for this purpose, we introduce spaces $\Theta^0$ and $\Theta$ and i.i.d. sequences $(\theta^0_n)_{n\ge 0}$ and $(\theta_n)_{n\ge 0}$  of random variables in $\Theta^0$ and $\Theta$ respectively which they will use to randomize their action choices.

\begin{definition}[Open-Loop Policies]
\label{de:open_loop}
The space $\bPi^{\text{OL}}$ \index[not]{$\bPi^{OL}$}of open-loop policies is the set of sequences $\bpi=(\pi_n)_{n\ge 0}$ of deterministic measurable functions $\pi_n:\Xi_n\to \cP(A)$  where $\Xi_n = \cY\times (\Theta\times\Theta^0\times E\times E^0)^n\times \Theta^0$.
A control process $\bfa=(\fa_n)_{n\ge 0}$ is said to be generated by the open loop policy $\bpi$ if $\fa_n=\pi_n(\xi_n)$ $\PP$-almost surely where $\xi_0=(y,\theta_0^0)$ and, for $n\ge1$, $\xi_n=(y,(\theta_k,\theta_k^0,\epsilon_{k+1},\epsilon^0_{k+1})_{k=0,\ldots,n-1},\theta_n^0)$. Here $y\in\cY$ is a random variable providing the randomization of the initial condition.
\end{definition}

Typically, an action taken from an open loop policy is a sample from a distribution on the action space $A$, which depends on the the initial condition, the past sources of randomization of the central unit and the generic agent, the past values of the idiosyncratic noise and the common noise observed up to time $n$ and the current common randomization.

\begin{definition}[Closed-Loop Markov Policies]
\label{de:close_loop}
The space $\bPi^{\text{CL}}$ \index[not]{$\bPi^{CL}$} of closed loop Markov policies is the set of sequences $\bpi=(\pi_n)_{n\ge 0}$ where each $\pi_n$ is a measurable function from $S\times\cP(S)\times\Theta^0$ into $\cP(A)$.  
\end{definition}
In other words, an action taken from a closed-loop Markov policy is a sample from a distribution on the action space $A$, which depends on the current state $X_n$, the current population distribution $\mu_n$, and the current common randomization $\theta^0_n$.

\vskip 2pt
We shall say that a policy $\pi = (\pi_n)_{n\geq 0}$ is \textbf{stationary} if $\pi_n = \hat{\pi}$ for all $n$, for some fixed function $\hat{\pi}$ independent of $n$.\index[sub]{stationary policy}
Also we shall say that a policy is pure if its components take values in $A$ instead of $\cP(A)$. \index[sub]{pure policy}

\subsection{The Lifted Mean Field MDP}

The analysis of McKean-Vlasov dynamics, \index[sub]{McKean-Vlasov dynamics} is typically based on stochastic differential equations whose coefficients depend upon the distribution of the solution, sometimes called nonlinear stochastic differential equations. These equations are also called mean field equations as they appear in the limit of the description of the time evolution of particle systems interacting in a mean field manner.
Most theoretical analyzes of these equations rely on a lifting procedure replacing the original problem by a (deterministic) dynamical system that describes the time evolution of the probability distribution of the original state of the system. As can be seen in the monograph \cite{BertsekasShreve}, this idea is already present in early works on the control of  Markov Decision Processes (MDPs). \index[not]{MDP} However, due to the presence of two independent sources of random shocks, a simple lifting of our model will not lead to the optimization of a deterministic dynamical system. It will lead  to a simpler MDP, which is still stochastic, though driven only by the common noise. 
We will work directly with this stochastic lifted MDP.

\begin{definition}[Mean Field MDP]\label{de:mfmdp}
The result of the first lifting is the \textbf{mean field MDP (MFMDP)} defined by the tuple
\[
(\hat{S}, \hat{A}, \hat{\Gamma}, \hat{P}, \hat{f}),
\]
where:
\begin{itemize}[leftmargin=*]
    \item the state space is $\hat{S} = \cP(S)$ the space of probability distributions over agent states,
    \item the action space is $\hat{A} = \cP(S \times A)$ the space of joint distributions over states and actions,
    \item the constraint space is $\hat{\Gamma} \subseteq \hat{S} \times \hat{A}$ defined by $\hat{\Gamma} = \{(\mu, \bar{a}) : \bar{a} \in \hat{U}(\mu)\}$ \index[not]{$\hat{\Gamma}$} where \index[not]{$\hat{U}(\mu)$} $\hat{U}(\mu) = \{\bar{a} \in \hat{A} : \text{pr}_1(\bar{a}) = \mu\}$.
Here and in the following we denote by  $\text{pr}_1: \cP(S \times A) \to \cP(S)$ \index[not]{$\text{pr}_1$} the projection onto the first marginal. For later purposes we note that the set $\hat{\Gamma}$ is closed and analytic in the product space $\hat{S} \times \hat{A}$.    
    \item the time evolution of the lifted MDP is given by the transition probability kernel $\hat{P}: \hat{\Gamma} \to \cP(\hat{S})$
    \[
	\hat{P}(\cdot|\mu, \bar{a}) = \cL\big(\hat F(\mu,\bar a,\varepsilon^0)\big),
    \]
    where $(X, a) \sim \bar{a}$, $\varepsilon \sim \nu$, $\varepsilon^0 \sim \nu^0$, all independent.  The stochastic dynamics could be specified as well via the system function $\hat F:\hat\Gamma\times E^0\mapsto\hat S$  defined by:
    \begin{equation}
        \label{fo:Fhat}
        \hat F(\mu,\bar a,e^0)=\cL_{\bar a\otimes \nu}\bigl(F(\cdot,\cdot,\bar a,\cdot,e^0)\bigr)
    \end{equation}
    by which we mean the image (push forward) of the probability measure $\bar a\otimes \nu$ under the map $S\times A\times E\ni (x,a,e)\mapsto F(x,a,\bar a,e,e^0)$. Recall that $\mu$ is the first marginal of $\bar a$.
    \item the one-stage cost function $\hat{f}: \hat{\Gamma} \to \RR$ is defined by
    \[
    \hat{f}(\mu, \bar{a}) = \int_{S \times A} f(x, a, \bar{a}) \, \bar{a}(dx, da).
    \]
\end{itemize}
\end{definition}

\begin{remark}[Connection with Open-Loop Policies]
A key insight is that the introduction of randomization of the decisions made by the central unit will enable a one-to-one correspondence between open-loop policies in the original MFC and Markov policies in the lifted MFMDP. This connection will allow us to leverage standard MDP theory for the analysis of MFC problems with common noise.
\end{remark}

The lifted MFMDP is a standard MDP in the space $\hat S$ of probability measure on the agent state space $S$, the new state of which evolves according to the recursion
\[
\mu_{n+1} \sim \hat{P}(\cdot|\mu_n, \bar{a}_n),\quad\text{or equivalently}\quad  
\mu_{n+1} = \hat{F}(\mu_n, \bar{a}_n,\varepsilon^0),\quad \bar{a}_n \in \hat{U}(\mu_n),
\]
with one-stage cost $\hat{f}(\mu_n, \bar{a}_n)$. We shall use tools from the standard theory of MDPs, relying heavily on \cite{BertsekasShreve}, to analyze this lifted MFMDP in the next chapter. Under suitable assumptions, the value function of this lifted MDP  should satisfy the Bellman equation
\[
J(\mu) = \inf_{\bar{a} \in \hat{U}(\mu)} \left\{\hat{f}(\mu, \bar{a}) + \gamma \, \EE[J(\hat{F}(\mu, \bar{a}, \varepsilon^0))]\right\},
\]
where the expectation is over the common noise $\varepsilon^0 \sim \nu^0$. Recall that $\hat F$ is the system function of the lifted MDP, namely $\hat{F}(\mu, \bar{a}, \varepsilon^0)$ is the distribution of $F(X, a, \bar{a}, \varepsilon, \varepsilon^0)$ with $(X, a) \sim \bar{a}$, $\varepsilon \sim \nu$.
In order to guarantee solvability of the Bellman equation, we shall prove:

\begin{proposition}[Contraction Property]\label{prop:contraction}
Under Assumptions~\ref{ass:H1} and a couple of mild assumptions on $\hat F$ and $\hat f$, namely, Assumptions \ref{ass:basic_assumption_MFMDP} and \ref{ass:compactness_MFMDP} introduced in the next chapter, the Bellman operator $T$ defined by
\[
[TJ](\mu) = \inf_{\bar{a} \in \hat{U}(\mu)} \left\{\hat{f}(\mu, \bar{a}) + \gamma \, \EE[J(\hat{F}(\mu, \bar{a}, \varepsilon^0))]\right\},
\]
is a strict contraction on the space of bounded lower semi-continuous functions (or bounded continuous functions) equipped with the supremum norm.
\end{proposition}
Indeed, among other things, we shall prove in Theorem \ref{th:MFMDP_DPP_BOREL} that:
\[
\norm{TJ - TJ'}_\infty \leq \gamma \, \norm{J - J'}_\infty.
\]
Consequently, by the Banach fixed point theorem (Theorem~\ref{th:banach}), $T$ has a unique fixed point $J^* = TJ^*$, which is the optimal value function. Moreover, $J^*$ can be obtained as the limit of value iterations:
\[
J^* = \lim_{n \to \infty} T^n J_0,
\]
for any initial bounded lower semi-continuous function $J_0$.

\subsection*{State-Action Value Functions and Existence of Optimal Policies}

Because of our interest in RL, we shall use the state-action value function \index[sub]{state-action value function} also known as the \textbf{Q-function} \index[sub]{Q-function} defined by
\[
Q^*(\mu, \bar{a}) = \hat{f}(\mu, \bar{a}) + \gamma \, \EE[J^*(\hat{F}(\mu, \bar{a}, \varepsilon^0))],
\]
With its help, one will characterize optimal states for the lifted MDP:

\begin{proposition}[Optimality Condition]\label{prop:optimality}
An action $\bar{a}^* \in \hat{U}(\mu)$ is optimal in state $\mu$ if and only if
\[
Q^*(\mu, \bar{a}^*) = \min_{\bar{a} \in \hat{U}(\mu)} Q^*(\mu, \bar{a}) = J^*(\mu).
\]
\end{proposition}

Finally, using the fact that $Q^*$ is the unique solution of the Bellman equation for the Q-functions, i.e.
\[
Q^*(\mu, \bar{a}) = \hat{f}(\mu, \bar{a}) + \gamma \, \EE_{\mu' \sim \hat{P}(\cdot|\mu,\bar{a})} \left[\min_{\bar{a}' \in \hat{U}(\mu')} Q^*(\mu', \bar{a}')\right].
\]
one will derive existence of pure stationary optimal policies for the lifted MDP model:

\begin{theorem}[Existence of Optimal Policies]\label{th:existence-optimal}
Under Assumptions~\ref{ass:H1}--\ref{ass:H2-mfc}, the optimal value function $J^*$ is bounded and lower semi-continuous, and there exists a \textbf{pure stationary optimal policy} $\pi^* = (\hat{\pi}^*, \hat{\pi}^*, \ldots)$ such that
\[
J^*(\mu) = \EE^{\pi^*}\left[\sum_{n=0}^\infty \gamma^n \hat{f}(\mu_n, \bar{a}_n) \;\Big|\; \mu_0 = \mu\right].
\]
\end{theorem}

\section{The Linear-Quadratic Mean Field Control Model}
\label{se:prelims_LQ}
We now switch to the second MFC model studied in this monograph.
In some sense, the LQ MFC model which we analyze in detail in Chapter \ref{ch:LQMFC} can be viewed as a particular case of the model of Chapter \ref{ch:AAP} for which the system function is linear and the one-stage cost function is quadratic. To be more specific, we assume that the state and action spaces are Euclidean spaces, say $S= \RR^d$ and  $A= \RR^\ell$
\begin{equation}
\label{eq:lq_sys_function_F}
    F(x,a,\mu,\epsilon,\epsilon^0)=Ax +Ba+\bar A\bar\mu^x+\bar B \bar\mu^a+\epsilon+\epsilon^0,
\end{equation}
where the matrices $A \in \RR^{d \times d}$ and $B \in \RR^{d \times \ell}$ drive the dynamics and the matrices $\bar{A} \in \RR^{d \times d}$ and $\bar{B} \in \RR^{d \times \ell}$ capture the mean field coupling as we denote by $\bar\mu^x$ the mean of the first marginal of $\mu$, namely $\bar\mu^x=\int_{S\times A}x\mu(dx,da)$
and by $\bar\mu^a$ the mean of the second marginal of $\mu$, namely $\bar\mu^a=\int_{S\times A}a\mu(dx,da)$.
Obviously, the marginal distributions $\nu$ and $\nu^0$ of the idiosyncratic and common noises are now probability measures on the Euclidean space $\RR^d$. They are most commonly assumed to be mean-zero Gaussian distributions, but we shall show that many results still hold for more general distributions.

\vskip 2pt
Similarly, we assume that the one-stage cost function is of the form 
\begin{equation}
\label{fo:LQ_cost}
    f(x,a,\mu) = (x-\bar{\mu}^x)^T Q (x-\bar{\mu}^x) + (a-\bar{\mu}^a)^T R (a-\bar{\mu}^a) + (\bar{\mu}^x)^T (Q+\bar{Q}) \bar{\mu}^x + (\bar{\mu}^a)^T(R+\bar{R}) (\bar{\mu}^a),
\end{equation}
where we use the exponent ${}^T$ to denote the transpose of a matrix or a vector and where:
\begin{itemize}[leftmargin=*]
    \item $Q, \bar{Q} \in \text{SM}_d$ are state cost matrices with $Q \geq 0$, $Q + \bar{Q} > 0$,
    \item $R, \bar{R} \in \text{SM}_\ell$ are control cost matrices with $R > 0$, $R + \bar{R} > 0$.
\end{itemize}
While it would be tempting to apply the results of Chapter \ref{ch:AAP} to the present model, most of the proofs rely on the compactness of some spaces or the boundedness of the system and cost functions, hypotheses which are obviously not satisfied here. Still the linear-quadratic structure of the model makes the problems amenable to explicit or semi-explicit solutions and we shall take full advantage of this possibility.

\subsection{Finite-Agent LQ Model}

In the $N$ agent version of the model,  each agent state $X^i_n \in \RR^d$ at time $n$ evolves according to the recursive equation
\begin{equation}
\label{fo:N_LQ_agent_dynamics}
X^i_{n+1} = A X^i_n + B a^i_n + \bar{A} \bar{X}_n + \bar{B} \bar{a}_n + \varepsilon^i_{n+1} + \varepsilon^0_{n+1}, \quad i = 1, \ldots, N,
\end{equation}
where:
\begin{itemize}[leftmargin=*]
    \item $\bar{X}_n = \frac{1}{N} \sum_{i=1}^N X^i_n$ and $\bar{a}_n = \frac{1}{N} \sum_{i=1}^N a^i_n$ are the empirical mean state and control,
    \item $\varepsilon^i_{n+1} \sim \nu$ are i.i.d.\ idiosyncratic noise (mean zero, covariance $\Sigma_\varepsilon$),
    \item $\varepsilon^0_{n+1} \sim \nu^0$ is common noise (mean zero, covariance $\Sigma_{\varepsilon^0}$), independent of all $(\varepsilon^i_{n+1})$.
\end{itemize}

At time $n$, the one-stage cost for agent $i$ is
\begin{equation}\label{eq:finite-agent-cost}
f^i(X_n,a_n)= (X^i_n-\bar X_n)^T Q (X^i_n-\bar X_n) + (a^i_n-\bar a_n)^T R (a^i_n -\bar a_n)+ \bar{X}_n^T (Q+\bar{Q}) \bar{X}_n + \bar{a}_n^T (R+\bar{R}) \bar{a}_n,
\end{equation}
where the cost matrices $Q, \bar{Q}, R, \bar{R}$ satisfy the same conditions as in \eqref{fo:LQ_cost}.

Since we are restricting ourselves to cooperative agents, the objective is to minimize the expected discounted social cost
\[
J^N = \EE\left[\sum_{n=0}^\infty \gamma^n \frac{1}{N} \sum_i f^i(X_n,a_n)\right].
\]

\subsection{LQ MFC Model (Mean Field Limit)}

Given what we have seen so far, it is not surprising that in the asymptotic regime of large homogeneous populations of agents, i.e. when $N \to \infty$, the finite-agent model  \emph{converges} to the \textbf{LQ MFC model}, where the state of a generic agent evolves as
\begin{equation}\label{eq:mfc-dynamics}
X_{n+1} = A X_n + B \alpha_n + \bar{A} \bar{X}_n + \bar{B} \bar{\alpha}_n + \varepsilon_{n+1} + \varepsilon^0_{n+1},
\end{equation}
with:
\begin{itemize}[leftmargin=*]
    \item $\bar{X}_n = \EE[X_n\mid\sigma\{\underline\varepsilon^0_n\}]$ the mean field state (population average),
    \item $\bar{\alpha}_n = \EE[\alpha_n\mid\sigma\{\underline\varepsilon^0_n\}]$ the mean field control.
\end{itemize}
\noindent
The one-stage cost is
\begin{equation}\label{eq:mfc-cost}
f(X_n,\alpha_n,\bar X_n,\bar\alpha_n)= (X_n-\bar X_n)^T Q (X_n-\bar X_n) + (\alpha_n-\bar\alpha_n)^T R (\alpha_n-\bar\alpha_n)+ \bar{X}_n^T (Q+\bar Q) \bar{X}_n + \bar{\alpha}_n^T (R+\bar R) \bar{\alpha}_n,
\end{equation}
and the objective is to minimize
\[
\EE\left[\sum_{n=0}^\infty \gamma^n f(X_n, \alpha_n, \bar{X}_n, \bar{\alpha}_n)\right].
\]

\subsection{Assumptions for LQ MFC}
We state the main assumptions we shall use in the proofs in Chapter \ref{ch:LQMFC} for future reference.
\begin{assumption}[Stability and Discounting]\label{ass:H2-lq}
\begin{enumerate}[label=(\roman*)]
    \item The system matrices satisfy $\;\gamma \norm{A}^2 < 1$ and $\gamma \norm{A + \bar{A}}^2 < 1$.
\end{enumerate}
\end{assumption}

\begin{assumption}[Cost Convexity]\label{ass:H3-lq}
\begin{enumerate}[label=(\roman*)]
    \item The cost matrices $Q, \bar{Q}, R, \bar{R}$ are symmetric.
    \item $Q \geq 0$, $Q + \bar{Q} > 0$, $R > 0$, $R + \bar{R} > 0$.
\end{enumerate}
\end{assumption}

We shall also impose an extra assumption requiring the optimal feedback gains to satisfy the discounted norm conditions defining the admissible class below (see Assumption~\ref{ass:optimal-feedback-in-Theta}).

\subsection{Admissible Control Parameters and Optimality Characterization}

The mainstay of our analysis of the optimization of the LQ MFC problem is to show that one can reduce the search for optimal control to the special class of controls in \textbf{linear closed-loop Markov (feedback) form}:
\begin{equation}\label{eq:linear-feedback}
\alpha_n = -K (X_n-\bar X_n) - L \bar{X}_n,
\end{equation}
for a couple $(K,L)$ of (deterministic) matrices:
\begin{itemize}[leftmargin=*]
    \item $K \in \RR^{\ell \times d}$ representing the individual feedback gain,
    \item $L \in \RR^{\ell \times d}$ representing the mean field feedback gain.
\end{itemize}

\begin{definition}[Admissible Control Parameters]\label{de:admissible-controls}
The pair $\theta = (K, L)$ of control coefficients is said to be \textbf{admissible} if it stabilizes the closed-loop system, in other words, if it belongs to the following set of admissible parameters:
\[
\Theta_{\text{ad}} = \left\{(K, L) \in \RR^{\ell \times d} \times \RR^{\ell \times d} : \gamma \norm{A - BK}^2 < 1, \; \gamma \norm{A + \bar{A} - (B + \bar{B})L}^2 < 1\right\}.
\]
\end{definition}

The rationale for these admissibility conditions is that they ensure that both the individual dynamics ($X_{n+1} \approx (A - BK)X_n + \text{noise}$) and the mean field dynamics ($\bar{X}_{n+1} \approx (A + \bar{A} - (B + \bar{B})L)\bar{X}_n + \text{noise}$) are contractive under discounting.

\vskip 4pt
Subsequently, the main result allowing us to characterize optimality in this class of admissible \emph{linear} controls will be stated at this stage in an informal form as 
\begin{theorem}[Optimal LQ Control]\label{th:optimal-lq}
Under appropriate conditions, the optimal control has the form
\[
u_n^* = -K^* (x_n-\bar x_n) - L^* \bar{x}_n,
\]
where $(K^*, L^*) \in \Theta_{\text{ad}}$ are characterized as the solution of a system of coupled matrix Riccati equations. The optimal value function at the individual level is quadratic:
\[
J(x, \bar{x}) = (x-\bar x)^T P (x-\bar x) + \bar{x}^T \bar{P} \bar{x} + \text{constant},
\]
for symmetric matrices $P, \bar{P} \geq 0$.
\end{theorem}
Here $J(x,\bar x)$ denotes the conditional expected discounted cost under the socially optimal feedback, given individual state $x$ and conditional population mean $\bar x$. For a conditional state law $\mu$ with mean $\bar x$ and finite second moment, the corresponding population value is $\int J(x,\bar x)\,\mu(dx)$.

Note that in the case of the LQMFC model, the Bellman equation reads:

\begin{theorem}[Bellman Equation for LQ MFC]\label{th:bellman-lq}
The optimal value function $J$ for the LQ MFC satisfies
\[
J(x, \bar{x}) = \inf_{a,\bar{a}} \left\{
\begin{aligned}
& (x-\bar x)^T Q( x-\bar x) + (a-\bar a)^T R (a-\bar a) \\
& {}+ \bar{x}^T (Q+\bar{Q}) \bar{x} + \bar a^T (R+\bar{R}) \bar a + \gamma \, \EE[J(x_+, \bar{x}_+)]
\end{aligned}
\right\},
\]
where $x_+ = Ax + Ba + \bar{A}\bar{x} + \bar{B}\bar{a} + \varepsilon + \varepsilon^0$ and $\bar{x}_+ = (A + \bar{A})\bar{x} + (B + \bar{B})\bar{a} + \varepsilon^0$.
\end{theorem}
For a conditional state law $\mu$ with mean $\bar x$ and a square-integrable feedback $\varphi$, admissibility requires $a=\varphi(x)$ and $\bar a=\int\varphi(y)\,\mu(dy)$. The quadratic minimization above yields $\varphi^*(x)=-K^*(x-\bar x)-L^*\bar x$ and $\bar a^*=-L^*\bar x$, which satisfy $\int\varphi^*(y)\,\mu(dy)=\bar a^*$.

\section{Reinforcement Learning Framework}

\subsection{Model-Free Learning}

In many applications, the system function $F$ and/or the cost function $f$ are not explicitly known. \textbf{Reinforcement learning (RL)} provides a framework for learning approximate optimal policies through interaction with the environment, treating it as a black box.

\begin{definition}[RL Paradigms]
Key RL paradigms include:
\begin{itemize}[leftmargin=*]
    \item \textbf{Value-based methods:} Estimate the value function $J^*$ or Q-function $Q^*$, then derive a policy by greedy action selection.
    \item \textbf{Policy gradient methods:} Parameterize the policy and optimize it directly via gradient descent on expected costs.
    \item \textbf{Actor-critic methods:} Combine value estimation (critic) with policy optimization (actor).
\end{itemize}
\end{definition}

\subsection{Q-Learning}

In the particular case of the lifted MDP, \textbf{Q-learning} proposes a model-free algorithm that estimates $Q^*$ by iteratively updating estimates based on observed transitions:
\begin{equation}\label{eq:q-learning}
Q_{k+1}(\mu, \bar{a}) = (1 - \alpha_k) Q_k(\mu, \bar{a}) + \alpha_k \left[r + \gamma \min_{\bar{a}'} Q_k(\mu', \bar{a}')\right],
\end{equation}
where:
\begin{itemize}[leftmargin=*]
    \item $\alpha_k \in (0, 1)$ is the learning rate,
    \item $r = \hat{f}(\mu, \bar{a})$ is the observed cost,
    \item $\mu' \sim \hat{P}(\cdot|\mu, \bar{a})$ is the next state.
\end{itemize}

The popularity of Q-learning is due to the fact that under appropriate conditions (diminishing step sizes, sufficient exploration), $Q_k \to Q^*$ almost surely. Such an algorithmic search for an approximation of an optimal solution works surprisingly well when the state space and the action space are finite \emph{with a very small number of elements}. As argued in Chapter \ref{ch:review}, unfortunately, Q-learning performance deteriorates very quickly when the cardinality of any of these spaces increases. So when the state or action spaces are large or continuous, $Q$ is approximated by a parametric function $Q_\theta$ (e.g., a neural network), and parameters $\theta$ are learned via updates provided by stochastic gradient descent.

\subsection{Mean Field Reinforcement Learning (MFRL)}

\begin{definition}[MFRL]
\textbf{MFRL} extends single-agent RL to the MFC setting, where:
\begin{itemize}[leftmargin=*]
    \item The agent learns to optimize in the presence of a population distribution $\mu$,
    \item The interactions and the environment lead to a \emph{MFC} problem,
    \item Algorithms must handle the infinite-dimensional lifted state space $\cP(S)$.
\end{itemize}
\end{definition}
Common approaches include:
\begin{itemize}[leftmargin=*]
    \item \textbf{Finite-dimensional approximations:} we discretize $\cP(S)$, or represent distributions via finite-dimensional sufficient statistics (e.g., moments, histograms).
    \item \textbf{Particle-based methods:} we approximate $\mu$ by empirical distributions of finite particle systems.
    \item \textbf{Neural network function approximation:} we use deep networks to represent value functions or policies over distribution spaces.
\end{itemize}
We shall provide examples of applications of all these methods.
\subsubsection*{Simplex Discretization}

For state or action spaces that are probability simplexes (e.g., discrete distributions), we approximate continuous distributions by \textbf{finite-dimensional simplexes}.

\begin{definition}[Projection Operator]
Given a finite subset $\hat{S}_{M} \subset \hat{S}$ with $|\hat{S}_{M}| = M$, we define a projection operator
\[
\Proj_{\hat{S}_{M}} : \hat{S} \to \hat{S}_{M}
\]
that maps each $\mu \in \hat{S}$ to the closest element in $\hat{S}_{M}$ (with respect to a chosen metric $d$, e.g., Wasserstein distance).
\end{definition}

\begin{remark}[Approximation Error]
To analyze the convergence of the discretized MDP as the grid is refined, we will use the assumption that
\(
h_M:=\sup_{\mu\in\hat{S}}d\bigl(\mu,\Proj_{\hat{S}_M}(\mu)\bigr)\longrightarrow0,
\)
where $d$ is the metric used above. For finite $S$ and $A$, Theorem~\ref{th:main-cv-tabular} gives Q-function error bounds in terms of state and action covering radii under suitable assumptions.
\end{remark}

\subsubsection*{Histogram Approximations}

For continuous state spaces, \textbf{histogram approximations} partition the state space into bins and represent distributions by histograms (piecewise constant densities). This reduces the infinite-dimensional problem to a finite-dimensional one, enabling standard RL algorithms.

\subsubsection*{Monte Carlo Simulation and Particle Systems}

To implement environments for RL algorithms, we simulate the mean field dynamics using \textbf{particle systems}: a large but finite collection of $N$ agents whose empirical distribution approximates the population distribution $\mu$.

\begin{enumerate}[leftmargin=*]
    \item Initialize $N$ agents with states $\{x^1_0, \ldots, x^N_0\}$ sampled from the initial distribution.
    \item At each time step, compute the empirical distribution $\mu_N = \frac{1}{N} \sum_i \delta_{x^i}$.
    \item Each agent takes action $a^i$ according to the policy, and transitions to $x^i_+$ via the dynamics $F(x^i, a^i, \mu_N, \varepsilon^i, \varepsilon^0)$.
    \item Repeat for $T$ time steps, collecting cost observations.
\end{enumerate}

\begin{remark}
In the fixed-policy setting of Proposition~\ref{pr:chaos}, the empirical state measure converges weakly almost surely, at each fixed time, to the conditional state law.
\end{remark}
The total error can be split by the triangle inequality into population error (finite-agent versus exact mean field model), discretization error (exact versus discretized mean field model), and learning error (learned versus exact solution of the discretized model). Theorem~\ref{th:main-cv-tabular} bounds the discretization and learning contributions for an exact mean field simulator; finite-population simulation requires a separate population-error analysis.

\subsubsection*{Neural Network Function Approximation}

For complex, high-dimensional problems, \textbf{deep neural networks} are used to approximate:
\begin{itemize}[leftmargin=*]
    \item the \textbf{Value functions} $J_\theta(\mu)$ or $Q_\theta(\mu, \bar{a})$,
    \item the \textbf{Policies:} $\pi_\theta(a \mid x, \mu)$.
\end{itemize}

Common architectures include:
\begin{itemize}[leftmargin=*]
    \item \textbf{Feedforward fully connected networks:} Standard multilayer perceptrons.
    \item \textbf{Convolutional networks:} For spatial or image-based state representations.
    \item \textbf{Recurrent networks:} For temporal dependencies (though less common in our discrete-time Markov settings).
\end{itemize}

\subsection{Policy Gradient Methods for LQ MFC}

In the particular case of the LQ framework, the control parameters $\theta = (K, L)$ can be optimized via \textbf{policy gradient descent}:
\begin{equation}\label{eq:policy-gradient}
\theta_{k+1} = \theta_k - \beta_k \nabla_\theta J(\theta_k),
\end{equation}
where $\beta_k > 0$ is the step size and $\nabla_\theta J(\theta)$ is the gradient of the expected social cost with respect to $\theta$.

In our subsequent analysis in Chapter \ref{ch:LQMFC}, we derive gradient estimates from finite differences by perturbing  $\theta$ and estimating gradients from cost differences. Our control of the convergence depends upon conditions such as $\gamma \norm{A - BK}^2 < 1$ which ensure that the closed-loop system is %
stable after discounting, so that the discounted effect of state perturbations decays over time and the expected discounted cost is finite.

\subsubsection*{Sub-Gaussian Concentration Inequalities}

LQ models are tailor-made for the theory of Gaussian processes. However, while convenient mathematically, assuming that random shocks are Gaussian is often unrealistic. So in order to improve the realm of application of theoretical results, a valuable endeavor is to extend their reach to sub-Gaussian and sub-exponential random variables.
Accordingly, we employ sub-Gaussian and sub-exponential norms to control tail behavior.

\begin{definition}[Sub-Gaussian Norm]
\label{de:subGaussian}
A random variable $\eta$ is \textbf{sub-Gaussian} \index[sub]{sub-Gaussian} with parameter $s$ if $\EE[\exp(\eta^2/s^2)] \leq 2$. Its \textbf{sub-Gaussian norm} is
\[
\norm{\eta}_{\psi_2} = \inf\left\{s > 0 : \EE[\exp(\eta^2/s^2)] \leq 2\right\}.
\]
\end{definition}

\begin{definition}[Sub-Exponential Norm]
A random variable $\zeta$ is \textbf{sub-exponential} \index[sub]{sub-exponential} with parameter $s$ if $\EE[\exp(|\zeta|/s)] \leq 2$. Its \textbf{sub-exponential norm} is
\[
\norm{\zeta}_{\psi_1} = \inf\left\{s > 0 : \EE[\exp(|\zeta|/s)] \leq 2\right\}.
\]
\end{definition}

We also use the standard definition and notation of $L^p$-norms for a random vector $\xi \in \RR^d$:
\[
\norm{\xi}_{L^p(\RR^d)} = \left(\EE[\norm{\xi}^p]\right)^{1/p},
\]
where $\norm{\cdot}$ denotes the Euclidean norm on $\RR^d$.
For a random vector $\xi \in \RR^d$, we use the vector sub-Gaussian norm
\[
\norm{\xi}_{\psi_2} = \sup_{\substack{u\in\RR^d\\ \norm{u}=1}} \norm{\langle u,\xi\rangle}_{\psi_2},
\]
where the norm on the right-hand side is the scalar norm of Definition~\ref{de:subGaussian}.

In the LQ framework, \textbf{sub-Gaussian concentration inequalities} \index[sub]{concentration inequality} turn out to be a powerful tool providing tail bounds for sums of independent random vectors. A standard form, see e.g.~\cite{vershynin2018high}, is the following.

\begin{theorem}[Sub-Gaussian Concentration]
\label{th:subgaussian-concentration}
Let $\xi_1, \ldots, \xi_n$ be independent, mean-zero random vectors in $\RR^d$ such that $\norm{\xi_i}_{\psi_2} \leq K$ for every $i$. Then there exist universal constants $c,C>0$ such that, for every $t>0$,
\[
\PP\Bigl[\bigl\|\sum_{i=1}^n \xi_i\bigr\| > t\Bigr]
\leq 2 \exp\left(Cd-\frac{c t^2}{nK^2}\right).
\]
\end{theorem}
This inequality controls deviations of sample averages from their expectations, with the expected dependence on the dimension.

\section{Summary and Roadmap}

This chapter has laid the mathematical groundwork for the MFC problems and learning methods studied in the subsequent chapters:

\begin{enumerate}[leftmargin=*]
    \item \textbf{Chapter \ref{ch:AAP}: Mean Field Control with Common Noise}
    \begin{itemize}
        \item Focuses on the abstract MFC model with both idiosyncratic and common noise.
        \item Analyzes the lifted MFMDP and establishes equivalence between open-loop and closed-loop formulations.
    \end{itemize}

    \item \textbf{Chapter \ref{ch:LQMFC}: Linear-Quadratic Mean Field Control}
    \begin{itemize}
        \item Specializes to the LQ setting with quadratic costs and linear dynamics.
        \item Derives explicit solutions and studies finite-agent approximations.
    \end{itemize}

    \item \textbf{Chapter \ref{ch:numeric_I}: Numerical Methods for Mean Field Reinforcement Learning}
    \begin{itemize}
        \item Develops Q-learning and deep RL algorithms for model-free learning.
        \item Presents approximation and discretization methods for computational implementation.
    \end{itemize}

    \item \textbf{Chapter \ref{ch:numeric_II}: Policy Gradient Methods for LQ Mean Field Control}
    \begin{itemize}
        \item Develops policy-gradient methods for learning feedback gains.
        \item Analyzes sample complexity, approximation errors, and numerical experiments.
    \end{itemize}
\end{enumerate}

\section{Bibliographical Notes and Complements}

The notations  and the definitions from the measure theoretical part at the beginning of the chapter are generic in stochastic analysis. For consistency, we tried to follow \cite{Kallenberg2002foundations} and \cite{Kallenberg_RM}. In particular, a proof of the disintegration result stated in Proposition~\ref{pr:disintegration} can be found in \cite[Theorem~6.4]{Kallenberg2002foundations}, and the universal disintegration result stated in Theorem~\ref{th:universal_disintegration} is from \cite[Corollary~1.26]{Kallenberg_RM}.

A proof of the original Blackwell-Dubins lemma \ref{le:BlackwellDubins} can be found in the original paper \cite{BlackwellDubins} by the authors. See also \cite[Lemma 5.29]{CarmonaDelarue_book_I} for a self-contained explicit construction. A detailed discussion of the randomization practice in reinforcement learning can be found in \cite{CarmonaLauriere_SIREV}. There, the similarities and the differences with the use of relaxed controls in continuous time control problems are identified and put in perspective.

The proof of propagation of chaos for continuous-time dynamical systems driven by stochastic differential equations can be found in many textbooks on stochastic analysis. For convenience, we refer the interested reader to \cite[Chapter~5]{CarmonaDelarue_book_I} and the references given therein, including classical propagation of chaos results in statistical physics, see \cite{Sznitman1991}. Discrete-time versions are few and far between. Even though it is only proven for compact state and action spaces, the result of Motte and Pham~\cite{motte2022mean}, quoted in Proposition~\ref{pr:chaos}, is particularly well suited to our setting because it treats mean field Markov decision processes with common noise and open-loop controls. For a recent discrete-time propagation-of-chaos result for nonlinear Markov chains, see also \cite{vuckovic2026propagation}.

\vskip 2pt
Randomization via the use of mixed controls is ubiquitous in the practice of reinforcement learning. A recent wave of interest in the connection between the theory of continuous-time stochastic control and the applications of discrete-time reinforcement learning has been the source of numerous publications and some confusion. While speculation on these connections can be enlightening in some cases, it has some limitations, especially when it comes to handling the similarities and the differences between the randomization procedures used in both fields as the use of mixed controls in discrete time reinforcement learning, and relaxed controls in continuous time stochastic controls have very different impacts on the system being controlled. The reader interested in this issue may want to consult the review article \cite{CarmonaLauriere_SIREV} for a discussion of the parallels and the differences between the two approaches in the algorithmic implementations of these two classes of models.

\chapter{A General Mathematical MFRL Model with Common Noise}
\label{ch:AAP}

\begin{abstract}\emph{
This chapter may be more theoretical than the other ones in this monograph.
The mathematical framework is of an infinite cooperative homogeneous population of agents whose state dynamics are modeled by a Mean Field MDP, hence a stochastic Mean Field Control (MFC) problem in discrete time, with an infinite horizon discounted cost in the objective function. Not only are the dynamics subject to idiosyncratic and common noise, but the actions are also subject to both idiosyncratic and {\it common randomness}. Because of the assumed cooperation of the agents, the optimization problem can be interpreted as one posed to a central unit helping a very large population of agents to minimize the aggregate social cost. To this end, the central unit can first sample a random policy for the whole population before letting each agent sample their action based on this common policy.  Another distinctive feature of this chapter is our desire to illustrate the impact of an enlargement of the information available to each agent when making their own decision. We develop in parallel the theory of the optimization over closed loop policies leading to Markovian structures, as well as the theory of optimization over open loop policies allowing the use of information accumulated over the entire past values of all the sources of randomness. The main originality of our approach is to allow both the \emph{fictitious central unit} and the individual agents to randomize their policies. Although this extra  randomization does facilitate \emph{exploration}, it also helps our mathematical analysis. Typically,  the extra source of randomness provided by the common randomization helps to connect open-loop and closed-loop policies for the MFC problem to a Mean Field MDP (MFMDP) whose state is the population distribution. This \emph{lifting} of the model from the analysis of the dynamics of the agent states and actions to the dynamics of a new MDP on the Wasserstein space of probability measures is a time-honored method to analyze mathematically mean field control problems, and we think it is important to see it in action in the RL context. Defining properly this new MFMDP and dealing rigorously with common noise and common policy randomization lead us to carry out a careful probabilistic analysis of this type of problem. }
\end{abstract}

\section{Introduction}

This chapter provides a rigorous mathematical analysis of the first model introduced in Chapter \ref{ch:prelims}.
After a quick review of the model, we study the lifted mean-field MDP (MFMDP)\index[not]{MFMDP}\index[sub]{mean-field MDP} and its connections with the original mean field control (MFC) \index[not]{MFC}\index[sub]{mean field control} problem (Section~\ref{se:relations-models}), and we introduce the state-action value function (Q-function) of the MFMDP (Section~\ref{se:Q_learning}).

The thrust of this chapter is twofold.
First, we study a setting which allows for common randomization, i.e., mixed policies at the mean-field level. To the best of our knowledge, this is a departure from the standard literature on mean field games (MFGs) and MFC problems. We believe that common randomization is of fundamental importance in the RL setup to ensure exploration in the space of the central unit's actions.
Second, we prove equality of the closed-loop and open-loop value functions for \emph{any} given policy, and not just for the optimal ones: in Section~\ref{se:back-to-MFC-2}, we prove that for any open-loop policy, there is a closed-loop policy achieving the same value, and vise versa. 

In the second part of the monograph, we describe the implementation of several RL methods based on the theoretical results of this chapter,  and we discuss the results: a tabular Q-learning algorithm relying on a discretization of the mean-field state simplex for which we prove convergence (Theorem~\ref{th:main-cv-tabular}), and a deep RL method to deal with continuous state and action spaces, which allows us to avoid simplex discretization and to deal with randomized actions.

\section{Model Description and Notations}
\label{se:model_description}

We consider the model introduced in Section \ref{se:abstract_model} of Chapter \ref{ch:prelims}. We shall repeat the notations and the definitions already stated there only when necessary for the clarity of our statements.

\vskip 4pt
Notice that because we work with a \emph{homogeneous} population of agents, we use a single action space $A$ common to all the agents, and a common cost function $f$ for all the agents. This is in contrast with our original discussion of MARL in  Chapter \ref{ch:review} where the action spaces and the cost functions were specific to the agent.

\vskip 2pt
As we argued in Chapter \ref{ch:review}, the dynamics of a MARL or MFRL model can be given equivalently by a set of transition probabilities or a system function. Since we plan to devote a significant part of this chapter to dependencies with respect to the whole past via open loop policies, we find it more convenient to work with the stochastic dynamical equations of the systems as defined by a system function and innovation processes, both idiosyncratic and common. Our system function $F$ is used to describe the evolution of the state process based on a state, an action, a mean-field interaction term\footnote{The interactions are through the joint state-action distribution, which is sometimes referred to as \emph{``extended'' MFC} or \emph{MFC of controls}.} and two noise terms. 

Even if we are willing to postpone the regularity conditions on $F$ and $f$, the definition of the model cannot be complete without explaining how the actions are taken and how these decisions affect the evolution of the system over time.  The following features of our model are important:
\begin{enumerate}[label=(\roman*)]
	\item The mean-field interactions are conditioned on all the information shared.
	\item The actions are randomized  by additional sources of randomness.
\end{enumerate} 
Thet were clearly motivated in Chapter \ref{ch:prelims} by our discussion of the propagation of chaos and how taking the limit when the number of agents goes to infinity could help identify the right mathematical features of the mean field limit.

\begin{remark}[Finite horizon setting]
In this chapter, we focus on the infinite horizon, discounted setting, and we look for stationary policies. In a finite horizon problem, it is generally not reasonable to look for stationary policies, and optimal policies can be computed by backward induction, based on a finite-horizon version of the DPP.  We refer to~\cite[Chapter 3]{BertsekasShreve} for more details. Notice, however, that the fact that the policy becomes non-stationary in itself is not crucial: we could still restrict our attention to stationary policies after extending the state space by incorporating the time component. See~\cite[Section 10.1]{BertsekasShreve} for more details. However, some issues would remain, others would arise, especially those related to the regularity of the transition and cost functions, requiring specific treatment as in~\cite[Chapter 3]{BertsekasShreve}.
\end{remark}

\begin{remark}[Boundedness of the cost function]
A key point in the definition of the problem is to ensure that the overall (aggregate expected) cost is well defined. If $f$ is bounded, this is automatically true despite the infinite horizon thanks to the discount parameter $\gamma \in (0,1)$. 
Interesting cases come naturally with an unbounded one-stage cost function $f$. This is the case for Linear Quadratic (LQ) models studied in Chapter \ref{ch:LQMFC}.
\end{remark}

As explained in Chapter \ref{ch:prelims}, the formulation of the model is motivated by the limit of an $N$-agent model when the size of the population $N$ tends to infinity. Still, because we plan to use randomized controls and policies, we need to introduce the various sources of randomization at play. Beyond the i.i.d. sequences $(\varepsilon_{n+1})_{n\ge 0}$ with distribution $\nu \in \cP(E)$  modeling the idiosyncratic random shocks
and $(\varepsilon^0_{n+1})_{n\ge 0}$ with distribution $\nu^0 \in \cP(E^0)$ modeling the common noise, we shall assume the existence of the following
 sources of randomness on the basic probability space $(\Omega,\cF,\PP)$:
\begin{enumerate}[label=(\roman*)]
	\item a random variable $\mathscr{U}$ with distribution $\PP_{\mathscr{U}} $ in a Polish space $(\Upsilon, \cB_{\Upsilon})$ providing the randomization for the initial state;
	\item  an i.i.d. sequence $(\vartheta_n)_{n\ge 0}$ of random variables in a Polish space  $(\Theta, \cB_{\Theta})$ with distribution $\PP_{\vartheta}$ providing the generic agent with a source of randomization for their action choices;
	\item  an i.i.d. sequence $(\vartheta^0_n)_{n\ge 0}$ of  random variables in a Polish space  $(\Theta^0, \cB_{\Theta^0})$ with distribution $\PP_{\vartheta^0}$ providing the central unit with a source of randomization for the choices of policies.
\end{enumerate} 
We assume that all these random sequences are independent of each other. We also assume that  $\PP_{\vartheta}$ and $\PP_{\vartheta^0}$ are both atomless. This guarantees the existence of Borel measurable functions $h:\Theta\to[0,1]$ and $h^0:\Theta^0\to[0,1]$ which are uniformly distributed when viewed as random variables on the probability spaces $(\Theta,\cB_\Theta,\PP_\vartheta)$ and $(\Theta^0,\cB_{\Theta^0},\PP_{\vartheta^0})$  respectively.

The uniform random variables constructed with the functions $h$ and $h^0$ will be used repeatedly with the Blackwell-Dubins Lemma.

 \subsection{Probabilistic Framework and Classes of Policies}
\label{se:proba-framework}
We now introduce the framework that will be used to rigorously study the MFC problem before we implement the theory so developed through RL algorithms. As explained intuitively in the introduction, we will distinguish several types of randomness. We will also distinguish between actions (elements of $A$), and controls (probability measures on $A$). They are the building blocks needed to define later the notion of policy (see \S~\ref{subsec:open-closed-policies}).

\vskip 4pt
We recall the following useful notation which was already introduced in Chapter \ref{ch:prelims}. For a given sequence $(\zeta_n)_{n\ge 0}$ we use the notation $\underline{\zeta}_n=(\zeta_0,\zeta_1,\cdots,\zeta_n)$.

\vskip 6pt
Next, we introduce the terminology which will help us characterize the information available to a generic agent and the central unit to make their choices of actions and policies. Incidentally, we start referring to everything related to agents (e.g. actions, controls, policies, costs, etc)  as \textbf{level-0}. This is in contrast with the lifted stochastic optimization model which we introduce in Section \ref{se:MFMDP} whose elements will be referred to as \textbf{level-1}. For reasons which will become clear later, the central unit controlling the agents will be called the \emph{level-1 controller}. This new terminology frees our presentation of the theoretical results from the details of the application we used to motivate the set-up.
For later use, set
\[
\cF_n^0=\sigma\{\underline{\varepsilon}_n^0,\underline{\vartheta}_{n-1}^0\},
\qquad
\cG_n^0=\sigma\{\underline{\varepsilon}_n^0,\underline{\vartheta}_n^0\},
\]
with the convention that $\underline{\vartheta}_{-1}^0$ is empty. The sigma-field $\cF_n^0$ represents the common information before sampling the current common randomization $\vartheta_n^0$, while $\cG_n^0$ includes it.

\begin{definition}
\label{de:action-control-proc}
A \defi{level-0 action} is an element of $A$. A \defi{level-0 (mixed) control} is a random probability measure on $(A, \cB_A)$, that is, a random variable with values in the Borel space $(\cP(A), \cB_{\cP(A)})$.  A \defi{level-0 action process} is a sequence of random variables $\balpha = (\alpha_n)_{n \geq 0}$ with values in $A$, $\alpha_n$ being $\sigma\{ \sU, \underline{\vartheta}_{n-1}, \underline{\vartheta}_n^0, \underline{\varepsilon}_{n}, \underline{\varepsilon}_{n}^0, \vartheta_n \}$-measurable. The set of such action processes is denoted by $\AA$.
A \defi{level-0 control process} is a sequence $\bfa = (\fa_n)_{n \geq 0}$ of level-0 controls for which  $\fa_n$ is $\sigma\{ \sU, \underline{\vartheta}_{n-1}, \underline{\vartheta}_n^0, \underline{\varepsilon}_{n}, \underline{\varepsilon}_{n}^0 \}$-measurable.
Finally, a level-0  action process $\balpha = (\alpha_n)_{n \geq 0}$ is said to be a \defi{realization} of a level-0 control process $\bfa = (\fa_n)_{n \geq 0}$ if 
$\cL \big( \alpha_n \, | \, \sigma\{ \sU, \underline{\vartheta}_{n-1}, \underline{\vartheta}_n^0, \underline{\varepsilon}_{n}, \underline{\varepsilon}_{n}^0 \} \big) = \fa_n$,
$\PP-a.s.$,  for every $n \geq 0$.	
\end{definition}
Intuitively, an action process is the realization of a control process, where the sampling is done using $(\vartheta_n)_{n \ge 0}$.
Generally speaking, we use the term \emph{mixed} for probability measures and \emph{randomized} for the corresponding random variable. Typically, a level-0 action is denoted by $a$ and a level-0 control is denoted by $\fa$. Unless specified otherwise, all the controls we consider at the level-0 are implicitly assumed to be mixed. 

For notational convenience, we use both notations $\PP_\xi$ and $\cL(\xi)$ interchangeably for the distribution of a random element $\xi$. We also use natural extensions to denote conditional distributions. 
The following simple measure theory result will prove to be useful in the sequel.

\begin{lemma}
	\label{le:identitly_level_0_action_realization}
	Given a level-0 control process $\fa$, there exists a level-0 action process $\balpha$ which is a realization of $\fa$. Moreover, if another level-0 action process $\balpha'$ is also a realization of $\fa$, then for every $n \geq 0$ and for every bounded Borel measurable function $h: A \to \RR$, if we introduce the notation $\xi_n$ for the random element $(\sU, \underline{\vartheta}_{n-1}, \underline{\vartheta}_n^0, \underline{\varepsilon}_{n}, \underline{\varepsilon}_{n}^0 )$, we have
	\begin{equation}
		\label{eq:identity_level_0_action_realization}
		\EE \left[ h( \alpha_n' ) \, | \, \sigma\{\xi_n\}\right] = \int_A h(\alpha) \fa_n(d \alpha) = 	\EE \left[ h( \alpha_n ) \, | \, \sigma\{\xi_n\}\right], \qquad \PP-a.s..
	\end{equation}
\end{lemma}

\begin{proof}
Recall that we use the notation $\sigma(\xi)$ of the $\sigma$-field generated by $\xi$, i.e. the set of events depending only on $\xi$. By definition, for each $n \geq 0$, $\fa_n$ is of the form:
	$$
		\fa_n(d \alpha)= \kappa_n^{\fa}(\xi_n)(d \alpha), \qquad\PP-a.s. \, ,
	$$
for some measurable function $\kappa_n^{\fa}$ on $\Upsilon \times \Theta^{n-1} \times (\Theta^{0})^{ n} \times E^n \times (E^{0})^{n}$ with values in $\cP(A)$. Let us denote by  $\rho_A$ the Blackwell-Dubins function (see Lemma~\ref{le:BlackwellDubins}) of the space $A$. Let us set $U_n=h(\vartheta_n)$ and:
	$$
		\alpha_n(\omega) := \rho_A( \kappa_n^\fa(\xi_n(\omega)), U_n(\omega)).
	$$
Now since the $\sigma-$field $\sigma\{\xi_n\}$ is independent of $U_n$ (with uniform distribution $\lambda_{[0,1]}$), Fubini's theorem \cite{Kallenberg2002foundations}[Theorem 6.4] yields for every bounded measurable function $\phi$ on $A$:
\begin{align*}
\EE[ \phi(\alpha_n) \, | \, \sigma\{\xi_n\}] &= \EE[ \phi( \rho_A( \kappa_n^\fa (\xi_n), U_n ) ) \, | \, \xi_n ]\\
& = \int_0^1 \lambda_{[0,1]}(du) \phi( \rho_A(\kappa_n^\fa(\xi_n), u) )\\
& = \int_A \kappa_n^\fa(\xi_n)(d\alpha) \phi( \alpha),  
\end{align*}
and so $\cL( \alpha_n | \sigma\{\xi_n\}) = \kappa_n^\fa(\xi_n) = \fa_n$, $\PP$-almost surely.
Equality~\eqref{eq:identity_level_0_action_realization} directly follows from the definition of a conditional distribution.
\qed\end{proof}

\subsubsection{Conditional distribution and state process}
We are now in a position to describe precisely the mean-field interactions in the system function, and provide a clear definition of the state process driven by a mixed control process in the mean-field model with common noise.

\begin{definition}
	\label{de:state_process_from_control_process}
For any initial distribution $\mu_0 \in \cP(S)$ and level-0 action process $\balpha = (\alpha_n)_{n \geq 0}$, we say that  a process $\bX^{\balpha,\mu_0} = (X_n^{\balpha, \mu_0})_{n \geq 0}$ is a \defi{state process associated to} $(\balpha, \mu_0)$ for the MFC model if: 
\begin{itemize}
\item[(i)] 	$X_0^{\balpha, \mu_0} $ is an $S$-valued $\sigma\{\sU\}$ - random variable with distribution $ \mu_0$,
\item[(ii)] for every $n \geq 0$, 
			\begin{equation}
			    \label{eq:system_dynamics_level_0}
				X_{n+1}^{\balpha, \mu_0} = F\big( X_n^{\balpha, \mu_0}, \alpha_n, \PP^0_{(X_n^{\balpha, \mu_0}, \alpha_n)}, \varepsilon_{n+1}, \varepsilon_{n+1}^0 \big),
			\end{equation}
			where $\PP^0_{(X_n^{\balpha, \mu_0}, \alpha_n)}$ is a regular version of $\cL\bigl( (X_n^{\balpha, \mu_0}, \alpha_n) \, | \, \cG_n^0 \bigr)$, the  joint distribution of state-action at time $n$ conditioned by the past history of the common noise and the common randomization up to the current time.
\end{itemize}
\end{definition}
Notice that for $n\ge 1$, the state $X^{\balpha, \mu_0}_n$ is measurable with respect to \[
	\sigma\{ \sU, (\vartheta_k, \vartheta_k^0, \varepsilon_{k+1}, \varepsilon_{k+1}^0 )_{k = 0, \ldots, n-1} \}
\]
which does not include $\vartheta_n$ and $\vartheta_n^0$ as $X_n^{\balpha, \mu_0}$  does not depend on $\vartheta_n$ and $\vartheta_n^0$. These random variables are used to define the control at time $n$, which will in turn influence the state at the next time step, namely, $X_{n+1}^{\balpha, \mu_0}$.

\vskip 2pt
For each level-$0$ action process $\balpha$ and each $n \geq 0$, we denote by $\PP^0_{X_n^{\balpha, \mu_0}}$ a regular version of the conditional distribution $\cL(X_n^{\balpha, \mu_0}\,|\,\cF_n^0)$. It holds:
\begin{equation}
\label{fo:P0X_n}
    \PP^0_{X_n^{\balpha, \mu_0}}
	=\cL(X_n^{\balpha, \mu_0}\,|\,\cG_n^0),  \qquad \mathbb{P}-a.s..
\end{equation}
Indeed, conditioned on $\cF_n^0$, $X_n^{\balpha,\mu_0}$ is independent of $\vartheta_n^0$.

\subsubsection{Open-loop and closed-loop policies}
\label{subsec:open-closed-policies}
The definitions of this subsection were already introduced in Subsection \ref{sub:policies} of Chapter \ref{ch:prelims}.
We repeat them here because of their technicalities and their importance for the subsequent developments.
We first consider open-loop policies. For each $n\ge 0$, let $\Xi_n = \Upsilon \times (\Theta \times \Theta^0 \times E \times E^0)^n \times \Theta^0$, with the convention that $\Xi_0 = \Upsilon \times \Theta^0$, and let us define $\bxi=(\xi_n)_{n\ge 0}$ by: 
$$
    \xi_0 = (\sU, \vartheta_0^0), 
    \qquad\text{and}\qquad \xi_n = (\sU, (\vartheta_k, \vartheta_k^0, \varepsilon_{k+1}, \varepsilon_{k+1}^0)_{k = 0, \ldots, n-1}, \vartheta_n^0), \qquad n \ge 1.
$$

\begin{definition}%
\label{de:open_loop_policy} 
A \defi{level-0 open-loop policy} is a sequence $\bpi= (\pi_n)_{n \geq 0}$ of deterministic measurable functions $\pi_n : \Xi_n \to \cP(A)$, called the \defi{open-loop strategy functions} at time $n$ of policy $\bpi$. 
The set of all open-loop policies is denoted by $\bPi^\tinyol$. 
	A level-$0$  control process $\bm{\fa} = (\fa_n)_{n \ge 0}$ is said to be \defi{generated by} the open-loop policy $\bpi$ if:  
	\begin{equation}
		\label{eq:def_admissible_open_loop_control}
		\fa_n = \pi_n( \xi_n ), \quad \PP-a.s., \qquad n \geq 0.
	\end{equation}
\end{definition}

Because of measurability restrictions, the correspondence between level-0 control processes $\bfa$ and open-loop policies $\bpi$ is understood up to $\PP$-a.s. equality of the induced control processes, after choosing measurable versions in~\eqref{eq:def_admissible_open_loop_control}. We shall use one object or the other depending on whether we want to emphasize the stochastic process or the sequence of deterministic functions defining the process. 
We shall often short-circuit the control process $\bfa$ and say that an action process $\balpha = (\alpha_n)_{n \geq 0}$ is (a realization of the control process) \defi{generated by} $\bpi$ if $\cL( \alpha_n\, | \, \sigma\{\sU,\underline{\vartheta}_{n-1},\underline{\vartheta}^0_n,\underline{\epsilon}_n,\underline{\epsilon}^0_n\} ) = \pi_n(\xi_n),$ $\PP-a.s.$, for all $n \geq 0$. We now consider closed-loop policies. 

\begin{definition}%
	\label{de:Markovian_policy}
	A \defi{closed-loop Markov strategy function} is a measurable function from $S \times \cP(S) \times \Theta^0$ into $\cP(A)$. 
	A \defi{closed-loop Markov policy} $\bpi=(\pi_n)_{n\ge 0}$  is a sequence of such functions.  
The set of all closed-loop Markov policies is denoted by $\bPi^\tinycl$. 
\end{definition} 

The choice of this form of closed-loop Markov strategy function is suggested by the mean-field nature of the dynamics~\eqref{eq:system_dynamics_level_0} and the form of the costs~\eqref{eq:J_alpha}. Taking values in $\cP(A)$ instead of $A$ indicates that the strategy functions are mixed. Typically, they suggest that at each time $n \geq 0$, the action $\alpha_n \in  A$ taken according to such a policy should be sampled from a probability measure depending directly on the values of $X_n^{\balpha, \mu_0}$ and $\PP^0_{X_n^{\balpha, \mu_0}}$, and the random variable $\vartheta_n^0$ used by the level-$1$ controller to randomize their choice. We formalize this procedure in Definition \ref{de:admissible_state_action_processes} below.

\begin{definition}%
	\label{de:admissible_state_action_processes}
	
	For a closed-loop Markov policy $\bpi\in \bPi^\tinycl$ and an initial distribution $\mu_0 \in \cP(S)$, a pair of state and action processes $(\bX, \balpha) = ( X_n, \alpha_n)_{n \geq 0}$ is said to be \defi{generated by} $(\bpi, \mu_0)$ if 
	\begin{enumerate}[label=\roman*)]
		\item $\bX$ is a state process associated to $(\balpha, \mu_0)$ in the sense of Definition~\ref{de:state_process_from_control_process}.  
		\item The action process $\balpha$ is adapted to $\GG^a$ and satisfies
		\begin{equation}
			\label{eq:definition_action_Markov_closed_loop}
			\cL \big( \alpha_n \, | \, \sigma\{\sU,\underline{\vartheta}_{n-1},\underline{\vartheta}^0_n,\underline{\epsilon}_n,\underline{\epsilon}^0_n\} \big) = \pi_n \big(  X_n, \, \PP^0_{X_n},\, \vartheta_n^0 \big), \qquad \PP-a.s., \qquad n \geq 0.
		\end{equation}
	\end{enumerate}
\end{definition}

In Definition~\ref{de:admissible_state_action_processes}, we can view the state and action processes as constructed simultaneously by alternatively invoking the system dynamics~\eqref{eq:system_dynamics_level_0} and the sampling procedure consistent with~\eqref{eq:definition_action_Markov_closed_loop}. 
A convenient way to construct an action process $\balpha$ satisfying \eqref{eq:definition_action_Markov_closed_loop} is to use the Blackwell-Dubins Lemma~ \ref{le:BlackwellDubins}. Indeed, if $\rho_A$ is the Blackwell-Dubins function of $ A$ and the uniformly distributed random variables $U_n$ is given by $U_n =h(\vartheta_n)$, we can choose $\alpha_n= \rho_A \big( \pi_n  \big(  X_n, \,  \PP^0_{X_n} ,\, \vartheta_n^0 \big), U_n \big)$, $\PP$-a.s.,  $n \geq 0$.

\begin{remark}
Even though we call a policy $\bpi\in \bPi^{\tinycl}$ a ``Markov" policy, 
	it does not imply any Markov property for the state process $\bX$ associated to such a policy. This abuse of terminology can be explained by our intention to work with level-1 Markov policies which will imply the Markov property for a lifted measure-valued state process constructed in the next section.  Also, since we only use the term ``Markov policy'' in the closed-loop setting, we shall most often drop the term Markov hereafter and only call them simply closed-loop policies.
\end{remark} 

\subsection{Optimization and value functions}
\label{subsec:optim-valuefct}

We now complete the definition of the agent level (level-0) MFC problem by stating  the optimization problems with open-loop  and closed-loop policies.
First, we introduce the value function associated to an action process. 
\begin{definition}%
	For every level-0 action process $\balpha$, the \defi{value function}  $J^{\balpha}: \cP(S) \to \RR$ is  defined for every $\mu_0 \in \cP(S)$ by:
	\begin{equation}
		\label{eq:J_alpha}
			J^\balpha(\mu_0) := \EE \left[\sum_{n \geq 0} \gamma^n f\left(X_n^{\balpha, \mu_0}, \alpha_n,\PP^0_{(X_n^{\balpha, \mu_0},\alpha_n)} \right)  \right],
	\end{equation} 
	where the state process $\bX^{\balpha, \mu_0}$ is associated to $(\balpha, \mu_0$) according to the dynamics~\eqref{eq:system_dynamics_level_0}.
\end{definition}

The value of $J^\balpha(\mu)$ is well-defined for every $\mu \in \cP(S)$ because $f$ is measurable and bounded. Furthermore, this value depends only on the sequence of time-marginal laws of the $\cP(S\times A)$-valued process $\Bigl(\PP^0_{(X_n^{\balpha, \mu_0},\alpha_n)} \Bigr)_{n\ge 0}$.

For any open-loop or closed-loop policy $\bpi$, we can show that $\mathbb{P}^0_{(X_n,\alpha_n)}$ depends on the action process only through the policy $\bpi$ provided $(\bX,\balpha)$ is generated by $\bpi$. See Lemma~\ref{le:policy_value} for a proof for open-loop policy, and a similar proof can be derived for closed-loop policy. As a consequence, we define, for any level-$0$ action process $\balpha$  generated by $\bpi$, 
$$
    J^{\bpi}(\mu) := J^{\balpha}(\mu), \qquad \mu \in \cP(S).
$$
Accordingly, we define the optimal open-loop value function and the optimal closed-loop value function as:
	\begin{equation*}
		J^{\tinyol,*}(\mu) = \inf_{ \bpi \in \bPi^\tinyol }  J^{\bpi}( \mu ), \qquad J^{\tinycl,*}(\mu) = \inf_{ \bpi \in \bPi^\tinycl }  J^{\bpi}( \mu ), \qquad  \mu \in \cP(S),
	\end{equation*}
	which are finite because we assume that the one-stage cost function $f$ is bounded and $\gamma \in (0,1)$.

\begin{lemma}
\label{le:policy_value}
For each $n\ge 0$, the law of the random measure $\PP^0_{(X_n,\alpha_n)}$ depends only upon the open-loop policy $\bpi$ as long as $\balpha$ is a realization of the control process generated by $\bpi$. 
\end{lemma}

\begin{proof}
Let $\balpha$ be an action process which is a realization of the control process generated by $\bpi$. For each integer $n\ge 0$, we compute $\EE\bigl[\Phi\bigl(\PP^0_{(X_n,\alpha_n)}\bigr)\bigr]$ for a family of bounded measurable functions $\Phi$ on $\cP(S\times A)$ which generate the Borel $\sigma$-field of $\cP(S\times A)$. For the sake of definiteness we work with functions $\Phi$ of the form:
$$
\Phi(\mu)=\prod_{j=1}^m\int_{S\times A}\varphi_j(x,\alpha)\,\mu(dx, d\alpha)
$$
for a finite set $\varphi_1,\cdots,\varphi_m$ of bounded continuous functions on $S\times A$. We have:
\begin{equation}
\label{fo:big_phi}
\begin{split}
\EE\bigl[\Phi\bigl(\PP^0_{(X_n,\alpha_n)}\bigr)\bigr] 
&=\EE\Bigl[\prod_{j=1}^m\int_{S\times A}\varphi_j(x,\alpha)\,\PP^0_{(X_n,\alpha_n)}(dx,d\alpha)\Bigr] =\EE\Bigl[\prod_{j=1}^m\EE\bigl[\varphi_j(X_n,\alpha_n)\,| \sigma\{\underline{\vartheta}^0_n,\underline{\varepsilon}^0_n\}\bigr] \Bigr].
\end{split}
\end{equation}
Now for each $j\in\{1,\cdots,m\}$ we have 
\begin{equation}
\label{fo:phi_j}
\begin{split}
&\EE\bigl[\varphi_j(X_n,\alpha_n)\,| \sigma\{\underline{\vartheta}^0_n,\underline{\varepsilon}^0_n\}\bigr]\\
&\hskip 35pt
=\int\cdots\int \varphi_j\bigl(X_n(u,\underline{\theta}_{n-1},\underline{e}_n,\underline{\theta}^0_{n},\underline{e}^0_n),
\alpha_n(u,\underline{\theta}_{n-1},\underline{e}_n,\underline{\theta}^0_{n},\underline{e}^0_n,\theta_n)\bigr)\\
&\hskip 125pt
\PP_{\cU}(du)\PP_{\underline{\vartheta}_{n-1}}(d\underline{\theta}_{n-1})\nu^n(d\underline{e}_n)\PP_{\vartheta_{n}}(d\theta_{n})\Bigr|_{\underline{\theta}^0_n=\underline{\vartheta}^0_n,\underline{e}^0_n=\underline{\varepsilon}^0_n}\\
&\hskip 35pt
=\int\cdots\int \Bigl(\int_A \varphi_j\bigl(X_n(u,\underline{\theta}_{n-1},\underline{e}_n,\underline{\theta}^0_{n},\underline{e}^0_n),
\alpha\bigr)\pi_n(d\alpha\,| u,\underline{\theta}_{n-1},\underline{e}_n,\underline{\theta}^0_{n},\underline{e}^0_n)\Bigr)\\
&\hskip 125pt
\PP_{\cU}(du)\PP_{\underline{\vartheta}_{n-1}}(d\underline{\theta}_{n-1})\nu^n(d\underline{e}_n)\Bigr|_{\underline{\theta}^0_n=\underline{\vartheta}^0_n,\underline{e}^0_n=\underline{\varepsilon}^0_n}\\\end{split}
\end{equation}
where we made explicit the dependence of $X_n$ on $\xi_n=(\cU,\underline{\vartheta}_{n-1},\underline{\varepsilon}_n,\underline{\vartheta}^0_{n},\underline{\varepsilon}^0_n)$ and $\alpha_n$ on $(\xi_n,\vartheta_n)$. This shows that the left hand side of \eqref{fo:phi_j}, and hence the left hand side of \eqref{fo:big_phi} only depend upon the action process $\balpha=(\alpha_n)_{n\ge 0}$ through the conditional distribution $\pi_n(d\alpha| u,\underline{\theta}_{n-1},\underline{e}_n,\underline{\theta}^0_{n},\underline{e}^0_n)$. We conclude that if two action processes are realizations of control processes generated by the same open-loop policy $\bpi$, the corresponding random measures $\PP^0_{(X_n,\alpha_n)}$ have the same distribution.
\qed\end{proof}

\section{Mean-Field MDP}
\label{se:MFMDP} 

In this section, we introduce the Markov Decision Process (MDP) which we use to identify the optimal closed-loop Markov policies for our original MFC model. This MDP is obtained by \emph{lifting} the agent level (level-0) dynamical system to a level-1 dynamical system on a space of probability measures.

\subsection{Mean-field MDP framework} 

The key observation is that, for a closed-loop Markov policy $\bpi \in \bPi^\tinycl$, the associated value function $ J^{\bpi}$  can be viewed as the value function of an MDP with a state process $(\PP^0_{X_n})_{n \geq 0}$ on the state space $\cP(S)$, and an action process $\bigl(\PP^0_{(X_n,\alpha_n)}\bigr)_{n \geq 0}$ with values in a new action space $\cP(S \times A)$. This new set-up requires a consistency constraint. Indeed, actions need now to be consistent with the state in the sense that the first marginal of any action which can be taken while in a given state, has to be equal to the state itself. We provide a rigorous definition of this MDP, and we show its connection to the original MFC model. 

\vskip 6pt 
The \emph{lifted} \defi{mean-field MDP (MFMDP)} consists of a six-tuple $(\bar S, \bar A, \bar \Gamma, P, \bar f, \gamma)$ as described below:
	\begin{itemize}
		\item The state space is the Polish space $\bar S := \cP(S)$; a generic element in $\bar S$ is denoted by $\mu$.
		\item The action space is the Polish space $\bar A := \cP(S \times A)$; a generic element in $\bar A$ is denoted by $\bar a$.		
		\item The control constraint is a set-valued function ${\bar U}$ from ${\bar S}$ into the set of non-empty subsets of ${\bar A}$ defined by:
		\begin{equation}
			\label{eq:Ubar}
			{\bar U}(\mu) :=\{{\bar a}\in{\bar A};\; \text{pr}_1({\bar a})=\mu\}, \qquad \forall \, \mu\in{\bar S};
		\end{equation}
		where $\pr_1: \bar A \to \bar S$ is the projection function that maps $\bar a \in \bar A$ onto its first marginal distribution on $S$. 
		We shall also use the notation
		\begin{equation}
			\label{eq:bar_Gamma}
			\bar\Gamma :=\{(\mu,{\bar a})\in{\bar S}\times{\bar A};\; {\bar a}\in{\bar U}(\mu)\}.
		\end{equation} 	
		for the graph of the constraint.	
		\item A transition probability kernel $P: \bar \Gamma \to \cP( \bar S)$, which is a  measurable function. 		
		\item A one-stage cost function $\bar f: \bar \Gamma \to \RR$, which is a bounded measurable function. 
		\item The discount coefficient $\gamma \in (0,1)$.
	\end{itemize}	

\begin{remark}
	The \emph{projection map}  $\pr_1$ is continuous, so the constraint set $\bar U(\mu)$ is closed in $\bar A$ for every $\mu \in \bar S$. The graph $Gr(pr_1) := \{ (\bar a, \mu): pr_1(\bar a) = \mu\} \subset \bar A \times \bar S$ is closed, so $\bar \Gamma$ is also closed in $\bar S \times \bar A$. Hence $\bar\Gamma$ is an analytic subset of $\bar S \times \bar A$, and a Polish space on its own. We assume that $\bar \Gamma$ is endowed with the induced topology as well as the trace $\sigma$-field inherited from $\bar S \times \bar A$.
\end{remark}

\begin{definition}%
	\label{de:MFMDP_problem}
The six-tuple $(\bar S, \bar A, \bar \Gamma, P, \bar f, \gamma)$  
is said to be the  \defi{MFMDP  lifted from the MFC model} $(S, A, E, E^0, F, f, \gamma)$ of Definition~\ref{de:mfc_model} if it satisfies:
	\begin{itemize}
		\item The transition kernel $P$ is given by
		\begin{equation}
		\label{de:MFMDP_transition_kernel_from_system_func}
	        P( \mu, \bar a)(d \mu') = \big( \nu^0 \circ \bar F( \mu, \bar a, \cdot)^{-1})(d \mu'), \qquad (\mu, \bar a) \in \bar \Gamma,
	    \end{equation}
		where $\bar F: \bar \Gamma \times E^0 \to \bar S$ is the  system function defined in terms of $F$ by: 
		\begin{equation}
			\label{eq:Fbar_pushforward}
			{\bar F}(\mu, \bar a ,e^0)  =  (\bar a \otimes \nu)\circ F(\cdot, \cdot,\bar a, \cdot,e^0)^{-1}, \qquad (\mu, \bar a, e^0 ) \in \bar \Gamma \times E^0.
		\end{equation}
		\item The one-stage cost function $\bar f: \bar\Gamma \to \RR$ of the MFMDP satisfies:
		\begin{equation}
			\label{de:bar_f}
			\bar f(\mu, \bar a )= \int_{S\times A} f(x, \alpha, \bar a) \bar a(dx,d\alpha), \qquad (\mu, \bar a) \in \bar \Gamma.
		\end{equation}
	\end{itemize}
\end{definition}

Here and in the following we denote by $\nu\circ g^{-1}$ \index[not]{$\nu\circ g^{-1}$ } the push-forward\index[sub]{push-forward} of a measure $\nu$ by a measurable function $g$.

So, for every $(\mu, \bar a , e^0) \in \bar \Gamma \times E^0$, $\bar F(\mu, \bar a , e^0)$ is defined as a probability measure on $S$ such that for every bounded and Borel measurable function $\phi : S \to \RR$, 
	\begin{equation}
		\label{eq:pushforward_of_F}
		\int_S \bar F(\mu, \bar a, e^0) (dx' ) \phi(x') = \int_{S \times A \times E} \bar a(dx, d\alpha) \nu(de) \phi \Big( F( x, \alpha, \bar a, e, e^0) \Big).
	\end{equation}
Next, we prove that $\bar F$ is Borel measurable.

\begin{lemma}
\label{le:under_A_regu-barF-barf}
Assume that $S$ and $A$ are compact metric spaces, that for every $(e, e^0) \in E \times E^0$, the function $F(\cdot, \cdot, \cdot, e, e^0): S \times A \times \cP(S \times A) \to S$ is continuous, and that $f: S \times A \times \cP(S \times A) \to \RR$ is bounded and continuous. Then, it holds:
\begin{itemize}
    \item $\bar F$ is Borel measurable and for every $e^0 \in E^0$, $\bar F(\cdot, \cdot, e^0)$ is continuous in its remaining variables.
    \item $\bar f$ is bounded and continuous.
\end{itemize}
\end{lemma}

\begin{proof}
The measurability of $\bar F$ is straightforward, see e.g. \cite[Proposition~7.29]{BertsekasShreve}, so we only argue the continuity for $e^0\in E^0$ fixed. Let $\bigl( (\mu_n,\bar a_n)\bigr)_{n\ge 0}$ be a sequence in $\bar \Gamma$ which converges weakly toward $(\mu,\bar a)$. Since $\mu_n=\text{pr}_1(\bar a_n)$ for each $n\ge 0$, we necessarily have $\mu=\text{pr}_1(\bar a)$ so that $(\mu,\bar a)\in\bar\Gamma$. We pick a continuous bounded function $\phi: S\mapsto \RR$ and we show that:
\begin{equation}
\label{fo:F_bar_n}
\lim_{n\to\infty}\int_S \phi(x')\bar F(\mu_n,\bar a_n,e^0)(dx')= \int_S \phi(x')\bar F(\mu,\bar a,e^0)(dx').
\end{equation}
Using Skorohod's characterization of weak convergence of probability measures, we have the existence of random variables $(Y_n,\beta_n)$ converging $\PP$-almost surely toward some $(Y,\beta)$ and such that $\PP_{(Y_n,\beta_n)}=\bar a_n$ for each $n\ge 0$ and $\PP_{(Y,\beta)}=\bar a$. Consequently, the integral in the left hand side of \eqref{fo:F_bar_n} can be rewritten as:
\begin{equation*}
\begin{split}
\int_S \phi(x')\bar F(\mu_n,\bar a_n,e^0)(dx')& = \int_{S\times A\times E} \phi\bigl(F(x,\alpha,\bar a_n,e,e^0)\bigr)\bar a_n(dx,d\alpha)\nu(de)\\
&=\int_E \Bigl(\EE\bigl[\phi\bigl(F(Y_n,\beta_n,\bar a_n,e,e^0)\bigr)\bigr]
\Bigr)\nu(de),
\end{split}
\end{equation*}
which converges toward $\displaystyle \int_E \Bigl(\EE\bigl[\phi\bigl(F(Y, \beta,\bar a,e,e^0)\bigr)\bigr]\Bigr)\nu(de)$
by Lebesgue's dominated convergence theorem because $F(\cdot,\cdot,\cdot,e,e^0)$ is continuous and $\phi$ is bounded continuous.

The continuity of the one-stage cost function $\bar f$ follows directly from \cite[Proposition~7.31]{BertsekasShreve}. In particular, if $f$ is only assumed to be bounded and lower semi-continuous, then $\bar f$ is bounded and lower semi-continuous.
\qed\end{proof}

\begin{lemma}
	\label{le:MFMDP_well_definedness_value_function}
	
	Let  $\bar \bpi\in\bar\bPi$. For every $\mu \in \bar S$, let $(\bmu, \bm{\bar a})$ and $(\bmu', \bm{\bar a'})$ be two pairs of state and action processes generated by $(\bar \bpi, \mu)$. Then 
	$
		\EE \left[ \sum_{n \geq 0} \gamma^n \bar f (\mu_n, \bar a_n) \right] = \EE \left[ \sum_{n \geq 0} \gamma^n \bar f (\mu_n', \bar a_n') \right].
	$
\end{lemma}

\begin{proof}
    For any fixed initial $\mu \in \bar S$, by definition of the pair of state and action processes generated by $(\bar \bpi, \mu)$, we show by induction that, for all $n \ge 0$, $\cL(\mu_n )=\cL(\mu'_n)$ and $\cL\bigl((\mu_n, \bar a_n)\bigr)=\cL\bigl( (\mu'_n, \bar a'_n)\bigr)$.
    For $n = 0$, $\mu_0 = \mu_0' = \mu \in \bar S$.  Assume that for some $n \geq 0$, we have $\cL(\mu_n )=\cL(\mu'_n)$, then for any bounded and Borel measurable function $\phi: \bar S \times \bar A \to \RR $,
	\begin{equation}
		\begin{split}
			\EE \left[ \phi ( \mu_n, \bar a_n ) \right] 
			& = \EE \left[ \EE\left[ \phi ( \mu_n, \bar a_n ) \, | \, \mu_n \right] \right]
			\\
			& = \EE \left[ \int_{\bar A} \phi( \mu_n, \bar a).  \cL(\bar a_n \, | \, \mu_n)(d \bar a) \right] 
			 \\
			 &= \EE \left[ \int_{\bar A} \phi( \mu_n, \bar a). \bar \pi_n(\mu_n) (d \bar a) \right] 
			 \\
			 &= \EE \left[ \int_{\bar A} \phi( \mu_n', \bar a).  \bar \pi_n(\mu_n') (d \bar a) \right] 
			 = \EE \left[ \EE\left[ \phi ( \mu_n', \bar a_n' ) \, | \, \mu_n' \right] \right]
			 =  \EE \left[ \phi \big( \mu_n', \bar a_n' \big) \right].
		\end{split}
	\end{equation}
	So $\cL\bigl((\mu_n, \bar a_n)\bigr)=\cL\bigl( (\mu'_n, \bar a'_n)\bigr)$. 
	Since $\varepsilon_{n+1}^0$ is independent of $(\mu_n, \bar a_n)$ and $(\mu_n', \bar a_n')$,  we have $\cL\bigl((\mu_n, \bar a_n, \varepsilon_{n+1}^0)\bigr) = \cL\bigl( (\mu_n', \bar a_n', \varepsilon_{n+1}^0)\bigr)$, which implies that the law of $
	    \mu_{n+1} = \bar F (\mu_n, \bar a_n, \varepsilon_{n+1}^0)$   is equal to the law of  $\bar F (\mu_n', \bar a_n', \varepsilon_{n+1}^0) = \mu_{n+1}'$. 
	Hence the conclusion.
\qed\end{proof}

We chose to characterize the dynamics of the lifted process by the transition kernel $P$ (and not only by the system function $\bar F$) for two reasons: 1) to conform with the standard literature on MDPs which typically uses transition kernels; 2) to be able to restrict ourselves to Markovian policies and control processes given by feedback functions of the state for the analysis of the lifted MDP.

\begin{remark}
In anticipation for what is going to come next, we want to emphasize that the
above MFMDP satisfies the standard assumptions from the theory of MDPs. See for example  \cite[Chapter 8-9]{BertsekasShreve}. We will use the results of these two classic chapters to derive a form of DPP for the MFMDP. 
\end{remark}

To highlight the tight connections between the lifted MFMDP and the original MFC, we define the notion of mixed strategy and mixed Markov policy for MFMDP.

\begin{definition}%
	\label{de:MFMDP_non_mixed_strategy_function}
We call \defi{level-1 pure strategy function} any measurable function from $\bar S$ into $\bar A$ whose graph is contained in $\bar\Gamma$.	
We call \defi{level-1 mixed strategy function} any Borel measurable function $\bar\pi$ from $\bar S$ into $\cP(\bar A)$ satisfying
	$ \bar \pi(\mu) (\bar U(\mu) )  = 1$, for every  $\mu \in \bar S$ .
We denote by $\overline{\bPi}^{p}$ (resp.  $\overline{\Pi}$) the set of pure (resp. mixed) strategy functions. 
\end{definition}
We shall identify every $\varphi \in\overline{\Pi}^{p}$ with the corresponding level-1 mixed strategy $\bar\pi$ defined by $\bar \pi(\mu) = \delta_{ \varphi(\mu) }$ for $\mu \in \bar S$.  
For the sake of brevity, we sometimes omit the term ``mixed" or ``randomized'' when $\bar \pi \in \overline \Pi$, but we always keep the term ``pure" or ``non-randomized'' when $\bar \pi \in \overline \Pi^{p}$. We now define policies as sequences of strategy functions.

\begin{definition}%
	\label{de:MFMDP_mixed_Markov_policy}	 
    A \defi{mixed Markov policy} for MFMDP is a sequence of level-1 mixed strategy functions, i.e. an element of $\bar\bPi := (\bar\Pi)^{\NN}$.
	Similarly, a \defi{pure policy} is a sequence of level-1 pure strategy functions, i.e., an element of $\bar\bPi^{p} := (\bar\Pi^{p})^{\NN}$.
    We say that a policy $\bar\bpi=(\bar\pi_n)_{n\ge 0}$ is \defi{stationary} if the strategy functions $\bar \pi_n$ are equal for all $n$. 
\end{definition}
We restrict ourselves to these policies and refrain from using history dependent policies because we are mostly interested in optimizing value functions, and we know that for MDPs like our lifted MFMDP, for each history dependent policy, there exists a Markov policy with the same value function (as defined below in Definition \ref{de:MFMDP_value_function_lsc}). See for example \cite[Proposition~9.1]{BertsekasShreve}.

\begin{definition}%
	\label{de:MFMDP_admissible_state_action_processes}	
	A pair of state and action processes $(\bmu, \bm{\bar a} ) = ( \mu_n, \bar a_n)_{n \geq 0}$ is said to be \defi{generated by} $(\bar \bpi, \mu) \in \bar\bPi \times \bar S$ if the following conditions are satisfied: $\mu_0 = \mu$, 
		\begin{equation}
			\label{eq:MFMDP_admissible_state_process}
			\mu_{n+1} = \bar F( \mu_n, \bar a_n, \varepsilon_{n+1}^0), \qquad \PP-a.s. \quad n\ge 0,
		\end{equation}
	and $ \bm{\bar a}$ is an $\bar A$-valued process adapted to $\FF^0$ satisfying:
		\begin{equation}
			\label{eq:MFMDP_admissible_action_process}
			\cL( \bar a_n | \mu_n ) = \bar \pi_n ( \mu_n ), \qquad \PP-a.s. \quad n\ge 0.
		\end{equation} 
\end{definition}

\subsection{Assumptions and optimization problem for MFMDP}
 
Some of our proofs rely on the following assumptions which we state here for the sake of later reference. 

\begin{assumption}%
	\label{ass:basic_assumption_MFMDP}
	\begin{itemize}
		\item \textbf{Mean-field system function $\bar F$:}   For every $e^0 \in E^0$, the function $\bar F(\cdot, \cdot, e^0)$ is  
		        continuous in its remaining variables on $\bar\Gamma$.
		
		\item \textbf{One-stage mean-field cost function $\bar f$:} $\bar f$ is lower semi-continuous on $\bar\Gamma$. 
	\end{itemize}
\end{assumption}

\noindent
To show existence of an optimal policy, we will make use of the following extra assumption: 
\begin{assumption}
	\label{ass:compactness_MFMDP} \textbf{Compactness}:	 The mean-field state space $\bar S$ and the mean-field action space $\bar A$ are compact metric spaces. 
\end{assumption}

\begin{remark}
\label{rem:assumptions-mfc}
It is easy to articulate assumptions formulated in terms of the original mean field model which would imply Assumptions \ref{ass:basic_assumption_MFMDP} and \ref{ass:compactness_MFMDP}. 
For example assuming that for every $(e, e^0) \in E \times E^0$, the function $F(\cdot, \cdot, \cdot, e, e^0)$ is continuous in its remaining variables, together with assuming that $f: S \times A \times \cP(S \times A) \to \RR$ is continuous, imply Assumptions \ref{ass:basic_assumption_MFMDP}. See Lemma~\ref{le:under_A_regu-barF-barf} for a proof.
Moreover, assuming that the state space $S$ and the action space $A$ are compact metric spaces implies Assumption \ref{ass:compactness_MFMDP}.  
Let us provide a concrete example. 
Consider the MFC model $(S, A, E, E^0, F, f, \gamma)$ defined by the following elements :
	\begin{itemize}
		\item State space: $(S, \cB_S) = ([-C_x, C_x], \cB_{[-C_x, C_x]})$, where $C_x$ is a positive constant. 
		\item Action space: $(A, \cB_A) = ([-C_a, C_a], \cB_{[-C_a, C_a]})$, where $C_a$ is a positive constant. 
		\item Idiosyncratic and common noise spaces: $(E, \cB_E) = (\RR, \cB_\RR)$ and $(E^0, \cB_{E^0}) = (\RR, \cB_\RR)$. 
		\item System function: $F: S \times A \times \cP(S \times A) \times E \times E^0 \to S$, 
		$$
		    F(x, a , \nu, e, e^0) 
		    = \begin{cases}
		        \tilde F(x, a , \nu, e, e^0), &\hbox{if $|\tilde F(x, a , \nu, e, e^0)| \le C_x$},
		        \\
		        \mathrm{sign}(\tilde F(x, a , \nu, e, e^0))C_x, &\hbox{ otherwise,}
		    \end{cases}
		 $$
		 where 
		 $$
		 \tilde F(x, a , \nu, e, e^0) =  A x + B a + \bar A \bar\nu_1 +\bar B \bar \nu_2 + e + e^0
		 $$
Here, $A$, $B$, $\bar A$ and $\bar B$ are constants, $\bar\nu_1=\int \xi \nu_1(d\xi) $ and $\bar\nu_2=\int \xi \nu_2(d\xi)$, where $\nu_1$ and $\nu_2$ denote the first and second marginal of $\nu$. $F$ is a continuous function hence it is measurable. Moreover it takes values in $S = [-C_x, C_x]$. 
		\item One-stage cost function: $f: S \times A \times \cP(S \times A) \to \RR$, 
		$$
		f(x,a,\nu) = Q |x - \bar\nu_1|^2+(Q+\bar Q)\bar\nu_1^2+R|a-\bar\nu_2|^2+(R+\bar R)\bar\nu_2^2
		$$
		where $Q$, $\bar Q$, $R$ and $\bar R$ are positive constants. The function $f$ is continuous on $S \times A \times \cP(S \times A)$ and bounded thanks to the boundedness of the state space and the action space.
\end{itemize}
This model is a truncated version of the linear quadratic model we investigate in the next chapter. Clearly, if $C_x$ and $C_a$ are large enough, we  expect its solution to be close to the unconstrained model.
\end{remark}

\begin{definition}
	\label{de:MFMDP_value_function_lsc}
	For every $\bar\bpi \in \bar\bPi$, we define the \defi{value function} $J^{\bar \bpi}$ by:
	\begin{equation}
		\label{eq:MFMDP_J_BOREL}
		\bar J^{\bar \bpi}(\mu) := \EE \Big[ \sum_{n \geq 0} \gamma^n \bar f (\mu_n, \bar a_n) \Big], \qquad \mu \in \bar S,
	\end{equation}
	where $(\bmu, \bm{\bar a}) = (\mu_n, \bar a_n)_{n \geq 0}$ is any pair of state and action processes generated by $(\bar \bpi, \mu)$. If $\bar \bpi \in \bar\bPi$ is stationary with strategy function $\bar \pi \in \bar\Pi$, then we let $\bar J^{\bar \pi} := \bar J^{\bar \bpi}$. 
\end{definition}

As a consequence of  Lemma~\ref{le:MFMDP_well_definedness_value_function}, the value function  $\bar J^{\bar \bpi}$ given in~\eqref{eq:MFMDP_J_BOREL} is  well defined because the expectation in \eqref{eq:MFMDP_J_BOREL} does not depend upon the particular choice of the pair of state action processes $(\bmu,\bm{\bar a})$ generated by $(\bar \bpi, \mu)$.  

 We will show below in Theorem~\ref{th:MFMDP_DPP_BOREL} that the value function $\bar J^{\bar \bpi}$ is lower semi-continuous when the policy is an optimal policy $\bar\bpi = \bar\bpi^*$. 
With the value function for MFMDP at hand, we define the \defi{optimal value function of the MFMDP} as:  
	$$
	\bar J^{*}(\mu) = \inf_{ \bar \bpi \in \bar\bPi}{\bar  J}^{\bar \bpi}(\mu), \qquad \mu \in \bar S.
	$$

\subsection{Dynamic Programming Principle for MFMDP}

We state and prove  the DPP for the optimal value function with Borel measurable mixed Markov policies.
We stress the fact that we are working in the lifted MFMDP setting by using  \emph{overlines} in the notations of all the quantities relative to the MFMDP. Still, it is worth noticing that, due to the presence of the common noise, our lifted model is still a stochastic model while in general, one lifting is enough to get to a deterministic dynamical model. This is for example the case in the monograph \cite{BertsekasShreve} to which we refer frequently for basic results.

\begin{theorem}%
	\label{th:MFMDP_DPP_BOREL}
	Suppose that Assumptions~\ref{ass:basic_assumption_MFMDP} and~\ref{ass:compactness_MFMDP} hold. 
	Then, the function $\bar J^{*}$ is bounded and lower semi-continuous, and moreover it is the unique bounded and lower semi-continuous function satisfying the following dynamic programming equation with unknown $\bar J$:  
	\begin{equation}
		\label{eq:MFMDP_DPP_fixed_pt_BOREL}
		\bar J(\mu) = \inf_{ \bar a \in \bar U(\mu) } \left\{  \bar f(\mu, \bar a) + \gamma \EE \Big[ \bar J \big( \bar F( \mu, \bar a, \varepsilon^0) \big) \Big] \right\}, \qquad \mu \in \bar S.
	\end{equation}
	Furthermore, there exists a pure stationary $\bar \bpi^{*} = (\bar \pi^{*}, \bar \pi^{*}, \dots) \in \bar\bPi^{p}$ that is optimal, i.e., 
		$
		\bar J^{\bar \pi^*} = \bar J^{*}.
	$ 
\end{theorem}

The DPP \eqref{eq:MFMDP_DPP_fixed_pt_BOREL} is known to hold for universally measurable policies and lower semi-analytic one-stage cost function. See for example \cite[Proposition~9.8]{BertsekasShreve}. The gist of the above theorem is to show that it also holds for Borel measurable policies.

\begin{proof}

\textbf{Step 1: Bellman operator and fixed point with universal measurability. } 
We define the Bellman operator $\bar T$ by:
\begin{equation}
\label{fo:bar_T} 
    [\bar T \bar J](\mu) = \inf_{ \bar a \in \bar U(\mu) } \left\{  \bar f(\mu, \bar a) + \gamma \EE \Big[ \bar J \big( \bar F( \mu, \bar a, \varepsilon^0) \big) \Big] \right\}, \qquad \mu \in \bar S.
\end{equation}
Comparison of the infima gives $\norm{\bar T\bar J-\bar T\bar J'}_\infty\leq\gamma\norm{\bar J-\bar J'}_\infty$ for bounded universally measurable $\bar J,\bar J'$. Hence there is at most one bounded universally measurable solution of the Bellman equation. By \cite[Proposition~9.8]{BertsekasShreve}, this solution is the optimal value function for universally measurable policies, namely $\bar J^{*,Univ}$ defined as:
$$
    \bar J^{*,Univ}(\mu) = \inf_{\bar\bpi \in \bar\bPi^{Univ}} \bar J^{\bar\bpi}(\mu), \qquad \mu \in \bar S
$$
where $\bar\Pi^{Univ}$ is the set of universally measurable mixed strategy functions, $\bar\bPi^{Univ}:=(\bar\Pi^{Univ})^{\NN}$ is the corresponding class of policies, and for every $\bar\bpi \in \bar\bPi^{Univ}$, $\bar J^{\bar\bpi}$ is defined as in~\eqref{eq:MFMDP_J_BOREL} in the Borel measurable case.

\textbf{Step 2: Fixed point with Borel measurability.} We will apply the Banach fixed point theorem for $\bar T$ on $\Lsc(\bar S)$, which denotes the set of real valued,  bounded and lower semi-continuous functions on $\bar S$. This set is a closed subset of the Banach space of real valued bounded functions on $\bar S$ endowed with the sup norm.

The key point is to show that $\bar T$ leaves $\Lsc(\bar S)$ invariant. Fix $\bar J\in\Lsc(\bar S)$ and set
\[
g(\mu,\bar a)
=
\bar f(\mu,\bar a)+\gamma \EE\bigl[\bar J(\bar F(\mu,\bar a,\varepsilon^0))\bigr],
\qquad (\mu,\bar a)\in\bar\Gamma ,
\]
which is the content of the curly bracket in \eqref{fo:bar_T}. 
By Assumption~\ref{ass:basic_assumption_MFMDP}, $\bar f$ is lower semi-continuous and, for each fixed $e^0$, the map $(\mu,\bar a)\mapsto \bar F(\mu,\bar a,e^0)$ is continuous. Hence $(\mu,\bar a)\mapsto \bar J(\bar F(\mu,\bar a,e^0))$ is lower semi-continuous for each fixed $e^0$, and Fatou's lemma, applied after adding a constant if necessary, implies that the expectation term is lower semi-continuous. Therefore $g$ is lower semi-continuous on $\bar\Gamma$. Since Assumption~\ref{ass:compactness_MFMDP} makes $\bar A$ compact and since the graph $\bar\Gamma$ is closed, the selection theorem for lower semi-continuous functions over compact sets, see \cite[Proposition~7.33]{BertsekasShreve}, shows that $\mu\mapsto\inf_{\bar a\in\bar U(\mu)}g(\mu,\bar a)$ is lower semi-continuous. Thus $\bar T\bar J\in\Lsc(\bar S)$. The same proposition shows that the infimum is attained and admits a Borel measurable minimizing selector, since every $\bar U(\mu)$ is nonempty and compact.

Last, $\bar T$ is a strict contraction for the sup norm. We thus conclude by the Banach fixed point theorem that $\bar T$ has a unique fixed point in $\Lsc(\bar S)$, which we denote by $\bar J^{*,bd,lsc}$. In other words, $\bar J^{*,bd,lsc}$ is the unique bounded lower semi-continuous function satisfying the dynamic programming equation~\eqref{eq:MFMDP_DPP_fixed_pt_BOREL}.

\textbf{Step 3: $\bar J^{*,bd,lsc} = \bar J^{*}$.} Being lower semi-continuous, $\bar J^{*,bd,lsc}$ is universally measurable, so it is also a fixed point in the space of universally measurable functions. Hence, by uniqueness, it coincides with $\bar J^{*,Univ}$. Furthermore $\bar\bPi \subseteq \bar\bPi^{Univ}$. So we deduce:
\begin{equation}
\label{fo:*}
    \bar J^{*,bd,lsc}(\mu)
    =\bar J^{*,Univ}(\mu)
    \le \bar J^{*}(\mu), \qquad \mu \in \bar S.
\end{equation} 

Assumption~\ref{ass:compactness_MFMDP}  implies that the lifted MFMDP satisfies the assumptions of \cite[Corollary 9.17.2]{BertsekasShreve} 
so there exists a stationary pure (at the level-1) Borel measurable policy $\bar\bpi^*=(\bar\pi^*,\bar\pi^*,\cdots)$ which is optimal in the sense that: 
$$
    \bar J^{\bar\bpi^*}(\mu) 
    =\bar J^{*, Univ}(\mu),\qquad \mu\in\bar S.
$$
Consequently:
\begin{equation*} 
    \bar J^{*}(\mu)
    \le \bar J^{\bar\bpi^*}(\mu)
    =\bar J^{*,Univ}(\mu)
    =\bar J^{*,bd,lsc} (\mu),\qquad \mu\in\bar S,
\end{equation*}
which together with \eqref{fo:*} gives the equality  
$\bar J^{*} = \bar J^{*,bd,lsc}$.

Hence $\bar J^{*}$ satisfies the DPP~\eqref{eq:MFMDP_DPP_fixed_pt_BOREL} and by Step 2, it is the unique bounded lower semi-continuous function satisfying it. 
\qed\end{proof}

\begin{remark}
\label{remark:explain_universally_strategy_function}
	When the  mixed strategy function $\bar \pi_n \in \overline \Pi^{Univ}$ is only universally measurable, for each time $n$, as in the above proof of Theorem~\ref{th:MFMDP_DPP_BOREL}, the understanding of condition~\eqref{eq:MFMDP_admissible_action_process} requires a modicum of care.
 Let $q_n = \cL( \mu_n ) \in \cP(\bar S)$ be the distribution of the random state $\mu_n$ with values in $\bar S$.
 We consider a Borel measurable kernel, $\bar \pi_{q_n} : (\bar S, \cB_{\bar S}) \to (\cP(\bar A), \cB_{\cP(\bar A)} )$, such that $\bar \pi_{q_n}(\mu) = \bar \pi_n(\mu)$ for $q_n$-almost every $\mu \in \bar S$ (see~\cite[Lemma~7.28 (c)]{BertsekasShreve} for existence). Then condition~\eqref{eq:MFMDP_admissible_action_process} says that $\bar \pi_{q_n}$ is a regular version of the conditional probability of $\bar a_n$ given $\mu_n$. 

 Furthermore, the integration of a function $\phi(\mu_n, \cdot)$ with respect to $\bar \pi_n(\mu_n)$ should be understood in the following sense:
 	$$
	     \int_{\bar A} \phi(\mu, \bar a) \bar \pi_n( \mu_n) (d \bar a)  = \int_{\bar A} \phi( \mu, \bar a) \bar \pi_{q_n}( \mu ) (d \bar a), 
	$$
	for $q_n-$almost every $\mu \in \bar S$, where $q_n = \cL(\mu_n)$.
 \end{remark}

\paragraph{Running example: lifting and the two policy levels.}
For the cybersecurity example introduced in Section~\ref{sec:finite-agent-models}, the lifted state is $\mu\in\cP(S)$ and the four probabilities $p^{\mathrm{sw}}(x)\in[0,1]$ determine the feasible joint state--action distribution
\[
\bar a(x,1)=\mu(x)p^{\mathrm{sw}}(x),\qquad \bar a(x,0)=\mu(x)(1-p^{\mathrm{sw}}(x)).
\]
Thus $\bar a(x,0)+\bar a(x,1)=\mu(x)$. With $\Phi$ from~\eqref{eq:running-cyber-update}, and under the assumptions of Theorem~\ref{th:MFMDP_DPP_BOREL}, the Bellman equation specializes to
\[
\bar J^*(\mu)=\min_{p^{\mathrm{sw}}\in[0,1]^4}\left\{c(\mu)+\gamma\int\bar J^*(\Phi(\mu,p^{\mathrm{sw}},z))\,\nu^0(dz)\right\},
\]
where $\nu^0$ is the law of the fresh common shock. Figure~\ref{fig:running-cyber-policy-levels} distinguishes randomizing the planner's rule from the computers' private action draws. The planner's common randomization is separate from the exogenous common shock entering $\Phi$.
\begin{figure}[htbp]
\centering
\begin{tikzpicture}
\node[draw,align=center,text width=0.8\textwidth] (planner) at (0,0) {\textbf{Planner: level 1}\\Given $\mu_n$, sample one rule $p_n^{\mathrm{sw}}\in[0,1]^4$ and share it with all computers.};
\node[draw,align=center,text width=0.8\textwidth] (computers) at (0,-2.2) {\textbf{Computers: level 0}\\Each computer $i$ draws $\alpha_n^i\sim\mathrm{Bernoulli}(p_n^{\mathrm{sw}}(X_n^i))$, independently conditional on the states and the shared rule.};
\draw[->] (planner) -- (computers);
\end{tikzpicture}
\caption{Two policy levels in the cybersecurity example: one shared rule and independent individual action draws.}
\label{fig:running-cyber-policy-levels}
\end{figure}

\section{Relations between the models}
\label{se:relations-models}

In the previous section, we defined the level-1 lifted MFMDP and we proved some of its main properties. We now reconnect this model with the original agent level MFC problem and derive properties of the latter from what we just found out about the lifted model

\subsection{Relations between MFC closed-loop policies and MFMDP policies}
\label{se:relation-models-closedloop-barpi}

We start by highlighting that, intuitively, a closed-loop policy for the MFC can be viewed as sampled from a policy for the MFMDP by picking the common randomness. The following definition formalizes this idea. 

\begin{definition}%
	\label{de:correspondence_between_policies}
Let $\bpi \in \bPi^\tinycl$ and $\bar \bpi \in \bar\bPi$. We say that they \defi{correspond to each other} if  for each $\mu\in\bar S$ and $n\ge 0$, $\bar\pi_n(\mu)\in\cP(\bar A)$ is equal to the push forward of $\PP_{\vartheta^0}$ by the map:
$$
    \Theta^0 \ni \theta^0 \mapsto \mu \measprod \pi_n(\cdot, \mu, \theta^0) \in \bar A,
$$
where $\measprod$ denotes the product of a measure and a kernel.
\end{definition}

Note that if $\bpi$ and $\bar\bpi$ correspond to each other and if $\bpi$ is stationary, then so is $\bar\bpi$. Conversely, if $\bpi$ and $\bar\bpi$ correspond to each other and if $\bar\bpi$ is stationary, then $\bpi$ can be chosen to be stationary. We will provide a specific construction in equation~\eqref{eq:construction_closed-loop_Markov_policy_MFC} below.
The main result of this section is the following.

\begin{theorem}%
	\label{th:identical_value_function_MFC_MFMDP}
Assume that the continuity assumption (Assumption~\ref{ass:basic_assumption_MFMDP}) holds. Then for every $\mu \in \cP(S)$, for every $\bpi \in \bPi^{\tinycl}$, there exists $\bar \bpi \in \bar\bPi$ such that $J^{\bpi}(\mu) = \bar J^{\bar \bpi}(\mu)$. Conversely, for every $\mu \in \cP(S)$, for every $\bar \bpi \in \bar\bPi$, there exists $\bpi \in \bPi^{\tinycl}$ such that $J^{\bpi}(\mu) = \bar J^{\bar \bpi}(\mu)$.
In particular, for stationary policies, we have: for every $\mu \in \cP(S)$, for every stationary $\bpi \in \bPi^{\tinycl}$, there exists a level-1 mixed strategy function $\bar \pi \in \bar \Pi$ such that $J^{\bpi}(\mu) = \bar J^{\bar \pi}(\mu)$, and conversely, for every $\mu \in \cP(S)$, for every stationary policy $\bar \bpi \in \bar\bPi$ with strategy function $\bar \pi$, there exists a stationary $\bpi \in \bPi^{\tinycl}$ such that $J^{\bpi}(\mu) = \bar J^{\bar \pi}(\mu)$.
\end{theorem}

In preparation for the proof of this result, we start with two useful technical lemmas describing the properties of conditional distributions of the level-0 state and action processes. 

\begin{lemma}
\label{le:dynamics_of_conditional_distribution}
Suppose that Assumption~\ref{ass:basic_assumption_MFMDP} holds. Let $\balpha \in \AA$, $\mu_0 \in \cP(S)$, and let $\bX$ be the associated state process. Then:
\begin{equation}
     \label{eq:formula_dynamics_in_lemma_MFC_MFMDP}
        \PP^0_{X_{n+1}} = \bar F( \PP^0_{X_n}, \PP^0_{(X_n, \alpha_n)}, \varepsilon_{n+1}^0 ), \qquad \PP-a.s. \qquad n \geq 0.
    \end{equation}
So $\cL(\PP^0_{X_{n+1}}) = P(\PP^0_{X_n}, \PP^0_{(X_n, \alpha_n)})$, $n \geq 0$, where the transition kernel $P$ was defined in~\eqref{de:MFMDP_transition_kernel_from_system_func} as transition probability of the level-1 MFMDP.
\end{lemma}

\begin{proof}  
 Let us denote by $ \zeta_{n+1}$ the right hand side of \eqref{eq:formula_dynamics_in_lemma_MFC_MFMDP}, and let  $\phi: S \to \RR$, $h_n: ( \Theta^0\times E^0 )^n \to \RR$ and $\psi_{n+1}: E^0 \to \RR$ be arbitrary bounded measurable functions. We have:
    \begin{equation*}
    \begin{split}
    &\EE \left[  \psi_{n+1}(\varepsilon_{n+1}^0)  h_n(\underline \vartheta_n^0, \underline \varepsilon_n^0  )  \int_{S} \phi(x') \zeta_{n+1}(d x') \right]
    \\
     & \hskip 5pt
     = \EE \left[  \psi_{n+1}(\varepsilon_{n+1}^0)  h_n( \underline \vartheta_n^0, \underline \varepsilon_n^0  ) \int_{S} \phi(x') \Big(  \bar F( \PP^0_{X_n}, \PP^0_{(X_n, \alpha_n)}, \varepsilon_{n+1}^0 ) \Big) (d x') \right]
    \\
     & \hskip 5pt
 = \EE \left[  \psi_{n+1}(\varepsilon_{n+1}^0) h_n( \underline \vartheta_n^0, \underline \varepsilon_n^0  ) \int_{S \times A \times E} \PP^0_{(X_n, \alpha_n)}(d x, d\alpha) \nu(d e) \phi \Big( F( x, \alpha, \PP^0_{(X_n, \alpha_n)}, e, \varepsilon_{n+1}^0) \Big) \right]
     \\
     & \hskip 5pt
 = \int_{E \times E^0}\hskip -12pt  \nu( d e)  \nu^0(d e^0) \psi_{n+1}(e^0) \EE \left[ h_n( \underline \vartheta_n^0, \underline \varepsilon_n^0  )  
 \int_{S \times A} \hskip -12pt \PP^0_{(X_n, \alpha_n)}(d x, d\alpha) \phi \Big( F( x, \alpha, \PP^0_{(X_n, \alpha_n)}, e, e^0) \Big) \right]
     \\
     & \hskip 5pt
 = \EE \left[ \psi_{n+1}( \varepsilon_{n+1}^0)  \EE \left[ h_n(\underline \vartheta_n^0, \underline \varepsilon_n^0 )
 \phi \Big( F( X_n, \alpha_n, \PP^0_{(X_n, \alpha_n)}, \varepsilon_{n+1}, \varepsilon_{n+1}^0 \big) \Big) \, \Big| \, \cF_n^0, \varepsilon_{n+1}, \varepsilon_{n+1}^0 \right]  \right]
     \\
     & \hskip 5pt
 = \EE \left[  \psi_{n+1}( \varepsilon_{n+1}^0)  h_n( \underline \vartheta_n^0, \underline \varepsilon_n^0 ) \phi (X_{n+1}) \right],
    \end{split}
    \end{equation*}
    where the first equality is by definition of $\zeta_{n+1}$, the second equality is by the definition of $\bar F$ in terms of the system function $F$ of the original MFC, the third equality is by the fact that $\varepsilon^0_{n+1}$ is independent of all the other random quantities, the fourth equality is by definition of the conditional probability $\PP^0_{(X_n, \alpha_n)}$, and  the last equality is by the tower property of conditional expectation, the fact that $X_{n+1} = F( X_n, \alpha_n, \PP^0_{(X_n, \alpha_n)}, \varepsilon_{n+1}, \varepsilon_{n+1}^0 \big)$ and the fact that $(\varepsilon_{n+1}, \varepsilon_{n+1}^0) $ is independent of $(X_n, \alpha_n) $ and $\PP^0_{(X_n, \alpha_n)}$ is measurable with respect to  $\cF_n^0 = \sigma\{ \underline \vartheta_n^0, \underline \varepsilon_n^0 \}$. This shows that  $\zeta_{n+1}=\PP^0_{X_{n+1}}$, $\PP-$almost surely.
\qed\end{proof}

\begin{lemma}
\label{le:idenfiy_conditional_joint_dist_with_kernels}
Let us assume that Assumption~\ref{ass:basic_assumption_MFMDP} holds. Let $\balpha \in \AA$, $\mu_0 \in \cP(S)$, and let $\bX$ be the associated state process. For every $n \ge 0$, let $\kappa_n: \bar S \to \cP(\bar A)$ be the Borel measurable disintegration kernel of $\cL( \PP^0_{X_n}, \PP^0_{(X_n, \alpha_n)} )$ along its first marginal. 
Then, if $(\bzeta, \bm{\bar \eta} )$ is an $(\bar S \times \bar A)$ - valued pair of stochastic processes  which are $\sigma\{\underline{\vartheta}^0_n,\underline\epsilon^0_n\}$ - adapted, and satisfy: 
    	    $\zeta_0 = \mu_0,$ $\PP-a.s.,$ 
    	    $\zeta_{n+1} = \bar F( \zeta_n, \bar \eta_n, \varepsilon_{n+1}^0)$, $\PP-a.s.$,  $n \ge 0$,
    	and if
    	    $\cL( \bar \eta_n  | \zeta_n ) = \kappa_n( \zeta_n ),$ $\PP-a.s.$ $n \ge 0$,
we have: 
    \begin{equation}
        \label{eq:identify_conditional_joint_distribution}
\cL (\zeta_n, \bar \eta_n) =\cL \big(\PP^0_{X_n}, \PP^0_{(X_n, \alpha_n)} \big), \qquad n \geq 0.
    \end{equation}
\end{lemma}

\begin{proof} 
 We prove equation~\eqref{eq:identify_conditional_joint_distribution} by induction. Equality \eqref{eq:identify_conditional_joint_distribution} holds for $n=0$ by the assumptions on $(\bzeta, \bm{\bar \eta} )$. Let us assume that \eqref{eq:identify_conditional_joint_distribution} holds for some $n\ge 0$. We first show that $\cL(\zeta_{n+1}) = \cL(\PP^0_{X_{n+1}})$. Since $\varepsilon_{n+1}^0$ is independent of $\sigma\{\underline{\vartheta}^0_n,\underline\epsilon^0_n\}$, for every bounded Borel measurable function $\psi: \bar S \to \RR$, it holds:
    \begin{equation*}
        \begin{split}
            \EE[ \psi( \zeta_{n+1} ) ] &= \EE \left[ \EE \left[ \psi \Big( \bar F( \zeta_n, \bar \eta_n, \varepsilon_{n+1}^0 ) \Big) \, \Big| \, \zeta_n, \bar \eta_n \right] \right]
            \\
            & = \EE \left[ \int_{\bar S} \psi( \mu) P( \zeta_n , \bar \eta_n)(d \mu)   \right]
            \\
            & = \EE \left[ \int_{\bar S} \psi( \mu) P\Big( \PP^0_{X_n} , \PP^0_{(X_n, \alpha_n)} \Big) (d \mu)  \right]
            \\
            & = \EE \left[ \psi \big( \PP^0_{X_{n+1}} \big) \right], 
        \end{split}
    \end{equation*}
where the second equality is by definition of $P$ in~\eqref{de:MFMDP_transition_kernel_from_system_func}, the third equality is by the induction hypothesis, and the last equality is due to Lemma~\ref{le:dynamics_of_conditional_distribution}. So $\cL(\zeta_{n+1}) = \cL(\PP^0_{X_{n+1}})$. We then consider the joint law of $(\zeta_{n+1}, \bar \eta_{n+1})$.  
Using the condition $\cL( \bar \eta_m  | \zeta_m ) = \kappa_m( \zeta_m )$ with $m=n+1$, and the definition of $\kappa_{n+1}$ as a regular conditional distribution of $\PP^0_{(X_{n+1},\alpha_{n+1})}$ given $\PP^0_{X_{n+1}}$, the equality $\cL(\zeta_{n+1})=\cL(\PP^0_{X_{n+1}})$ implies
\[
\cL(\zeta_{n+1},\bar\eta_{n+1})
=
\cL\bigl(\PP^0_{X_{n+1}},\PP^0_{(X_{n+1},\alpha_{n+1})}\bigr).
\]
We conclude that \eqref{eq:identify_conditional_joint_distribution} holds for $n+1$ instead of $n$.
 \qed\end{proof}

\vskip 6pt
Intuitively, $(\kappa_n)_{n \ge 0}$ plays the role of the conditional law of the action process. We will come back to this interpretation in the proof of Lemma~\ref{le:construct_closed_loop}.  We are now ready to prove Theorem~\ref{th:identical_value_function_MFC_MFMDP}.

\begin{proof}[Proof of Theorem~\ref{th:identical_value_function_MFC_MFMDP}]

\textbf{Step 1:} Let $\bpi \in \bPi^\tinycl$. Let $\bar\bpi$ corresponding to $\bpi$ in the sense of Definition~\ref{de:correspondence_between_policies}.  
We now check the equality of the value functions $J^{\bpi}$ and $\bar J^{\bar\bpi}$. Note that if $\bar\bpi$ is stationary, then so is $\bpi$. Let $(\bX,\balpha)$ be a pair of state and action processes generated by $(\bpi,\mu_0)$. Then
	\begin{equation*}
		\begin{split}
			J^{\bpi}(\mu_0) 
			& = \EE \Bigl[ \sum_{n \geq 0} \gamma^n f \bigl( X_n, \alpha_n,\PP^0_{(X_n,\alpha_n)} \bigr) \Bigr]
			= \EE \Bigl[ \sum_{n \geq 0} \gamma^n \bar f \bigl( \PP^0_{X_n}, \PP^0_{(X_n, \alpha_n)} \bigr) \Bigr].
		\end{split}
	\end{equation*}
	We now check that $\mu_n=\PP^0_{X_n}$ and $\bar a_n=\PP^0_{(X_n,\alpha_n)}$ form a pair of state and action processes generated by $(\bar\bpi,\mu_0)$. This is indeed the case because \eqref{eq:MFMDP_admissible_state_process} is implied by Lemma~\ref{le:dynamics_of_conditional_distribution} and \eqref{eq:MFMDP_admissible_action_process} is implied by the definition of $\bar\pi_n$ and the fact that $\PP^0_{(X_n,\alpha_n)}=\PP^0_{X_n}\measprod\pi_n(\cdot,\PP^0_{X_n},\vartheta^0_n)$. Consequently,
	\begin{equation*}
    \bar J^{\bar \bpi}(\mu_0)  
    = \EE \Bigl[ \sum_{n \geq 0} \gamma^n \bar f \bigl( \PP^0_{X_n}, \PP^0_{(X_n, \alpha_n)}\bigr) \Bigr]
    = J^{\bpi}(\mu_0).
	\end{equation*}

\vskip 6pt
\textbf{Step 2:} Conversely, let $\bar \bpi = (\bar \pi_n)_{n \geq 0} $ in $\bar\bPi$. For every $n \geq 0$, 
	 $\bar \pi_n: \bar S \to \cP( \bar A)$ is a Borel measurable map such that for every $\mu\in \bar S$ we have   $\bar\pi_n(\mu)\bigl(\bar U(\mu)\bigr)=1$. According to the
 Universal Disintegration Theorem \ref{th:universal_disintegration} recalled in Subsection \ref{sub:kernels} of Chapter \ref{ch:prelims}, there exists a Borel measurable probability kernel  $K: S \times \cP(S \times A) \times \cP(S) \to \cP(A)$ such that for every $\rho\in\cP(S\times A) $ and $\mu \in \cP(S)$ such that $\text{pr}_1(\rho)=\mu$, we have $\rho=\mu\measprod K(\cdot,\rho,\mu)$, where $\measprod$ denotes the product of a measure and a kernel as defined in Subsection \ref{sub:kernels} of Chapter \ref{ch:prelims}.
So for every integer $n\ge 0$, $x\in S$, $\mu\in \bar S$ and $\theta^0\in\Theta^0$, we define:
	\begin{equation}
		\label{eq:construction_closed-loop_Markov_policy_MFC}
		\pi_n( x, \mu, \theta^0) := K\Bigl( x, \rho_{\bar A}\bigl(\bar \pi_n (\mu), h^0(\theta^0)\bigr),\mu\Bigr), 	
	\end{equation}
where $\rho_{\bar A}$ is the Blackwell-Dubins function of $\bar A$. Note that if $\bar\bpi$ is stationary, then so is $\bpi$.
Because the functions $K$, $h^0$ and $\rho_{\bar A}$ are measurable, so is the strategy function $\pi_n$ for every $n \geq 0$. Hence $\bpi = (\pi_n)_{n\ge 0} \in \bPi^{\tinycl}$. Recall that the function $h^0$ was introduced in Subsection~\ref{se:proba-framework}, and that $ h^0(\vartheta^0)$ is uniformly distributed on $[0,1]$ by construction. Notice that for every $\mu \in \bar S$ and for almost every $\theta^0 \in \Theta^0$, the definition of the universal disintegration kernel $K$ implies that:
	\begin{equation}
	\label{eq:pi-K-barphi}
	    \Big(\rho_{\bar A}\bigl(\bar \pi_n (\mu), h^0(\theta^0) \bigr) \Big) (d x, d\alpha)
	    =
		\mu(dx) K \Big( x,\rho_{\bar A}\bigl(\bar \pi_n (\mu), h^0(\theta^0) \bigr) , \mu \Big)(d \alpha),
	\end{equation}
and as a result, we have:
	\begin{equation}
	\label{eq:mu-pin-rhoAbar}
		\Big(\rho_{\bar A}\bigl(\bar \pi_n (\mu),h^0(\theta^0) \bigr) \Big) (d x, d\alpha)
		=\mu(dx)\pi_n( x, \mu, \theta^0)(d\alpha).
	\end{equation}
When $\theta^0$ is replaced by $\vartheta_n^0$, by Blackwell-Dubins lemma (see first point in Lemma~\ref{le:BlackwellDubins}), the left hand side of \eqref{eq:mu-pin-rhoAbar} is a random variable with values in $\bar A=\cP(S\times A)$ with distribution $\bar \pi_n (\mu)$. 

\vskip 2pt
Next, we show that $J^{\bpi} = \bar J^{\bar\bpi}$. Let $\mu_0 \in \bar S$. Let $(\bzeta, \bm{\bar{\eta}} )$ be state and action processes generated by $(\bar \bpi, \mu_0)$ (see Definition~\ref{de:MFMDP_admissible_state_action_processes}). Let $(\bX, \balpha)$ be a pair of  state and action processes generated by $\bpi$ and $\mu_0$.
Using the fact that $\PP^0_{(X_n,\alpha_n)}=\PP^0_{X_n}\measprod\pi_n(\cdot,\PP^0_{X_n},\vartheta^0_n)$, and the fact that $\vartheta^0_n$ is independent of $\PP^0_{X_n}$ by equation~\eqref{fo:P0X_n}, we have:
	\begin{equation*}
		\begin{split}
            J^{\bpi}(\mu_0)
            &= \sum_{n\ge 0}\gamma^n\EE\Bigl[\int_{S\times A}f(x,\alpha,\PP^0_{(X_n,\alpha_n)})\PP^0_{(X_n,\alpha_n)}(dx,d\alpha) \Bigr]  \\
            &= \sum_{n\ge 0}\gamma^n\EE\Bigl[\int_{S\times A}f\bigl(x,\alpha,\PP^0_{X_n}\measprod\pi_n(\cdot,\PP^0_{X_n},\vartheta^0_n)\bigr) \PP^0_{X_n}(dx)\pi_n(x,\PP^0_{X_n},\vartheta^0_n)(d\alpha)\Bigr]
            \\
            &= \sum_{n\ge 0}\gamma^n\EE\Bigl[\int_{S\times A}f(x,\alpha,\bar a)\bar a(dx,d\alpha)
\bar\pi_n\bigl(\PP^0_{X_n}\bigr)(d\bar a)\Bigr].
		\end{split}
	\end{equation*}
The last equality holds by the fact that both sides of~\eqref{eq:mu-pin-rhoAbar} with $\theta^0=\vartheta_n^0$ are random variables with values in $\bar A=\cP(S\times A)$ with distribution $\bar \pi_n (\mu)$. On the other hand, with the processes $(\bzeta, \bm{\bar{\eta}})$ generated by $(\bar \bpi , \mu_0)$ following equations~\eqref{eq:MFMDP_admissible_state_process} and~\eqref{eq:MFMDP_admissible_action_process}, we have:
	\begin{equation*}
		\begin{split}
			\bar J^{\bar \bpi}(\mu_0) 
			& = \EE \left[ \sum_{n \geq 0} \gamma^n \bar f ( \zeta_n, \bar \eta_n ) \right]
			\\
			& =  \sum_{n \geq 0} \gamma^n  \EE \left[ \int_{\bar A} \bar f ( \zeta_n, \bar a ) \cL( \bar \eta_n \, | \, \zeta_n)(d \bar a) \right]
			\\
			& = \sum_{n \geq 0} \gamma^n  \EE \left[ \int_{\bar A} \bar f ( \PP^0_{X_n}, \bar a ) \bar \pi_n ( \PP^0_{X_n} ) (d \bar a) \right],
		\end{split}
	\end{equation*}
where the last equality holds by~\eqref{eq:MFMDP_admissible_action_process} and by Lemma~\ref{le:idenfiy_conditional_joint_dist_with_kernels}, which implies $\cL(\zeta_n)=\cL(\PP^0_{X_n})$. This completes the proof.
\qed\end{proof}

\vskip 2pt
At this stage, we need to emphasize the crucial role played by the common randomization provided by the sequence $(\vartheta^0_n)_{n\ge 0}$. Its presence is what allowed us to prove that the value of a policy for the lifted MDP can always be achieved by a closed loop policy of the original MFC.

\vskip 4pt
Notice that the above result proves equality of the value functions, \textbf{policy by policy}, which is much stronger that equality of the value functions after taking the optimum values.

\begin{corollary}
\label{co:early}
Under Assumption \ref{ass:basic_assumption_MFMDP}, for every pure policy $\bar \bpi = (\bar \pi_n)_{n \geq 0} $ in $\bar\bPi^{p}$, there exists a level-0 closed-loop policy $\bpi=(\pi_n)_{n\ge 0}$ not depending upon the common randomization and such that $\bar J^{\bar\bpi}(\mu)=J^\bpi(\mu)$, for all $\mu \in \bar S$.
\end{corollary}
Saying that the policy does not depend upon the common randomization amounts to say that (by an abuse of notation) for every $n\ge 0$, the strategy function takes the form $\pi_n: S\times \cP(S)\to \cP(A)$.

\begin{proof}
$\bar \bpi$ being  pure implies that for each $n \ge 0$, $\bar \pi_n : \bar S \to \bar A$ satisfies $pr_1(\bar \pi_n(\mu) ) = \mu$ for all $\mu \in \bar S$, by Definition \ref{de:MFMDP_non_mixed_strategy_function}. 
We claim that the function $(x,\mu,\theta^0)\mapsto \pi_n(x,\mu,\theta^0)=K\bigl(x,\rho_{\bar A}\bigl(\delta_{\bar\pi_n(\mu)},h^0(\theta^0)\bigr),\mu\bigr)$ is independent of $\theta^0$ defining the function $(x,\mu)\mapsto \varphi_n(x,\mu)$ whose existence is claimed in the statement of the corollary. Indeed, since by definition of the Blackwell-Dubins map $\rho_{\bar A}$ the map $u\mapsto \rho_{\bar A}\bigl(\delta_{\bar \pi_n(\mu)},u\bigr)$ has law $\delta_{\bar \pi_n(\mu)}$ when viewed as a random variable on the unit interval, one concludes that for Lebesgue almost every $u\in[0,1]$, we have $\rho_{\bar A}\bigl(\delta_{\bar \pi_n(\mu)},u\bigr)=\bar \pi_n(\mu)$.
As a result
$$
\pi_n(x,\mu,\theta^0)=K(x,\bar \pi_n(\mu),\mu) =: \varphi_n(x, \mu), \quad (x, \mu) \in S \times \cP(S),
$$
almost surely in $\theta^0$ and consequently, $\bar \pi_n(\mu)=\mu\hat\otimes K(\cdot,\bar \pi_n(\mu),\mu)$, and we conclude the proof as in the proof of Theorem~\ref{th:identical_value_function_MFC_MFMDP}.
\qed\end{proof}

\subsection{Relations between MFC closed-loop and open-loop policies}
\label{se:back-to-MFC-2}

We first prove existence of optimal closed-loop Markov policies and then, equality of the open loop and closed loop value functions.

\begin{proposition}%
\label{pr:existence_opt_policy_mean_field}

Suppose Assumptions~\ref{ass:basic_assumption_MFMDP} and~\ref{ass:compactness_MFMDP} hold.
	There exists a pure stationary closed loop Markov policy for the original MFC that is optimal, i.e., $\bpi^{*} = (\pi^*,\pi^*,\dots) \in \bPi^{\tinycl,p}$ such that:
	$
		J^{\pi^{*}} = J^{\tinycl, *}.
	$
\end{proposition}
We recall that by ``pure'' we mean that the policy does not depend upon the \emph{common} randomization. However, the representative agent actions can still be randomized.

\begin{proof}
Let $\bar \bpi^{*}$ be an optimal pure stationary Markov policy for MFMDP whose existence is given in Theorem~\ref{th:MFMDP_DPP_BOREL}, and let $\bpi^{*}\in\bPi^{CL}$ be a closed-loop Markov policy whose value function is the same and whose existence is given in Theorem~\ref{th:identical_value_function_MFC_MFMDP} (for the case of stationary policies). Moreover, by Corollary~\ref{co:early}, we can assume that $\bpi^{*}$ is pure too. We have:
$$
    J^{\bpi^{*}}(\mu) = \bar J^{\bar \bpi^{*}}(\mu) = \inf_{\bar \bpi \in \bar\bPi} \bar J^{\bar \bpi}(\mu), \qquad \mu \in \cP(S),
$$
Using Theorem~\ref{th:identical_value_function_MFC_MFMDP} again, for every $\bpi\in \bPi^\tinycl$ for MFC, there exists $\bar \bpi \in \bar\bPi$ for MFMDP such that $    \bar J^{\bar \bpi}= J^{\bpi}$.
So, for every $\bpi \in \bPi^\tinycl$, 
$$
     J^{\bpi}(\mu) \geq J^{\bpi^{*}}(\mu) = \inf_{\bar \bpi \in \bar\bPi} \bar J^{\bar \bpi}(\mu),
$$
which concludes the proof.
\qed\end{proof}

We now show equality of the open and closed loop optimal value functions of the original agent level MFC optimization problem.

\begin{theorem}%
		\label{th:equality_value_function_open_closed_Markov}
Suppose Assumption~\ref{ass:basic_assumption_MFMDP} holds.
Then for every $\bpi^{\tinycl} \in \bPi^{\tinycl}$, there exists $\bpi^{\tinyol} \in \bPi^{\tinyol}$ such that $J^{\pi^\tinyol} = J^{\pi^\tinycl}$, and conversely, for every $\bpi^{\tinyol} \in \bPi^{\tinyol}$, there exists $\bpi^{\tinycl} \in \bPi^{\tinycl}$ such that $J^{\pi^\tinyol} = J^{\pi^\tinycl}$. As a consequence, $J^{\tinycl,*} = J^{\tinyol,*}$. If, in addition, Assumption~\ref{ass:compactness_MFMDP} holds, then there exist an optimal closed-loop policy and an optimal open-loop policy achieving this value.
\end{theorem}

Note that when both Assumption~\ref{ass:basic_assumption_MFMDP} and Assumption~\ref{ass:compactness_MFMDP} hold, then combining Theorem~\ref{th:equality_value_function_open_closed_Markov} with Proposition~\ref{pr:existence_opt_policy_mean_field} yields the existence of a \emph{pure} stationary closed loop Markov policy for the original MFC that achieves the value $J^{\tinyol,*}$. 

The first part of Theorem~\ref{th:equality_value_function_open_closed_Markov} is given by Lemma~\ref{le:construct_open_loop} below. For the converse direction, the result is obtained by combining Lemma~\ref{le:construct_closed_loop} below and Theorem~\ref{th:identical_value_function_MFC_MFMDP}. Under the additional compactness assumption, the existence of an optimal closed-loop policy stems from Theorem~\ref{th:identical_value_function_MFC_MFMDP}, which entails the existence of an optimal open-loop policy by Lemma~\ref{le:construct_open_loop}.

Again, while this type of equality between the optimal value functions could be expected to hold under different assumptions and without the central randomization, we prove it here by leveraging the equalities proven in Lemma~\ref{le:construct_open_loop} and Lemma~\ref{le:construct_closed_loop} policy by policy, before computing optima over sets of policies.

\vskip 6pt
First, it is expected that every closed-loop policy for the MFC can be viewed as an open-loop policy for the MFC, which leads to  the following result in terms of value functions. 

\begin{lemma}
\label{le:construct_open_loop}
Suppose Assumption~\ref{ass:basic_assumption_MFMDP} holds.
For every $\tilde{\bpi} \in \bPi^{\tinycl}$, there exists $\bpi \in \bPi^{\tinyol}$ such that: 
    $
        J^{ \bpi} = J^{\tilde{\bpi}}.
    $
\end{lemma}
\begin{proof}
We prove this statement by showing that there exist $\bpi \in \bPi^{\tinyol}$ and   an  open-loop action process $\balpha$ generated by $\bpi$ such that:
    $
        J^{ \bpi}= J^{ \balpha } = J^{\tilde{\bpi}}.
    $
	Let $\mu_0 \in \cP(S)$ and let $(\bX, \balpha)$ be a pair of state and action processes generated by $(\tilde{\bpi}, \mu_0)$. Let $\bfa$ be the $\cP(A)$-valued process given by:
	$$
	    \fa_n= \tilde\pi_n( X_n, \PP^0_{X_n}, \vartheta_n^0) , \qquad n \geq 0.
	$$
We recall that $\Xi_n, \xi_n$ and $\sU$ are defined in \S~\ref{subsec:open-closed-policies}. Since $\fa$ is adapted to $\sigma\{\xi_n\}$, it is 
a level-0 control process, and for every $n \geq 0$, there exists a Borel measurable function $\pi_n: \Xi_n \to \cP(A)$ satisfying
$$
     \fa_n     = \pi_n(\xi_n), \qquad \PP-a.s. \, .
$$
Moreover, since $\balpha$ is generated by $\tilde\bpi$, it is adapted to $\sigma\{\xi_n,\vartheta_n\}$ and satisfies
$$  
    \cL( \alpha_n \, | \, \sigma\{\xi_n\}) = \tilde\pi_n( X_n, \PP^0_{X_n}, \vartheta_n^0) = \pi_n( \xi_n), \qquad  \PP-a.s. \qquad n \geq 0.
$$
So $\balpha$ can be viewed as an open-loop action process generated by $\bpi$. 
Meanwhile, the state process $\bX$ constructed by equation~\eqref{eq:system_dynamics_level_0} 
is also a state process associated with $(\balpha, \mu_0)$ (see Definition~\ref{de:state_process_from_control_process}). Therefore, by definition of the value function associated to an open-loop policy $\bpi$, we have:
$$
   J^{\tilde\bpi}(\mu_0)
   = J^{\balpha}(\mu_0) 
   = \EE \left[ \sum_{n\geq 0} \gamma^n f \big(X_n, \alpha_n, \PP^0_{(X_n, \alpha_n)} \big) \right].
$$
\qed\end{proof}

\vskip 6pt
Next, every open-loop policy for the MFC corresponds to a policy for the MFMDP, as we show in the following result. 

\begin{lemma}
	\label{le:construct_closed_loop}
Suppose Assumption~\ref{ass:basic_assumption_MFMDP} holds. For every $\bpi \in \bPi^{\tinyol}$, there exists  $\bar \bpi \in \bar\bPi$ such that: 
	$
	    \bar J^{\bar \bpi} = J^{\bpi}.
	$
\end{lemma}

\begin{proof} 
Let $\bpi \in \bPi^{\tinyol}$, 
let us fix an initial distribution $\mu_0 \in \cP(S)$, let 
$ \balpha$ be an action process generated by $\bpi$, and let $\bX$ be the state process associated with $(\balpha, \mu_0)$ (recall Definition~\ref{de:state_process_from_control_process}).
For each $n\ge 0$, we consider the probability kernel $\kappa_n: \bar S \to \cP(\bar A)$ defined in the statement of Lemma~\ref{le:idenfiy_conditional_joint_dist_with_kernels}.   
We construct by induction an $(\bar S \times \bar A)$-valued pair of processes $(\bzeta, \bm{\bar \eta})$ in the following way.
For $n=0$ we set $\zeta_0 = \mu_0$ and $\bar\eta_0=\rho_{\bar A}\bigl(\kappa_0(\zeta_0),h^0(\vartheta^0_0)\bigr)$ where $\rho_{\bar A}$ is the Blackwell-Dubins function of the space $\bar A$ introduced in Lemma~\ref{le:BlackwellDubins}. Then for any $n\ge 0$ we define: 
$$
\zeta_{n+1} = \bar F( \zeta_n, \bar \eta_n, \varepsilon_{n+1}^0 ),
\quad\text{and}\quad
\bar\eta_{n+1}=\rho_{\bar A}\bigl(\kappa_{n+1}(\zeta_{n+1}),h^0(\vartheta^0_{n+1})\bigr).
$$
The process $\bm{\bar \eta}$ is adapted to $\FF^0$ and by Lemma~\ref{le:BlackwellDubins} it satisfies
$$
    \cL( \bar \eta_n \, | \, \zeta_n ) = \kappa_n( \zeta_n), \qquad \PP-a.s.\,, \quad n \geq 0, 
$$
because $\vartheta^0_{n}$ is independent of $\zeta_n$. Thus, by Lemma~\ref{le:idenfiy_conditional_joint_dist_with_kernels}:
\begin{equation}
    \label{eq:prop-zeta-bareta-law}
    \cL(\zeta_n, \bar \eta_n) =\cL \Big(\PP^0_{X_n},  \PP^0_{( X_n,  \alpha_n ) }  \Big), \qquad n \geq 0.
\end{equation}
Let $\bar \bpi = (\kappa_n)_{n \geq 0}$. Since $\kappa_n(\mu)(\bar U(\mu) ) = 1$ for every $n \geq 0$ and $\mu \in \bar S$, we see that $\bar \bpi \in \bar\bPi$. We conclude by noting that:
$$
    \bar J^{\bar \bpi}(\mu_0) = \sum_{n=0}^\infty \gamma^n \EE\left[ \bar f( \zeta_n, \bar \eta_n) \right] = \sum_{n=0}^\infty \gamma^n \EE\left[ \bar f \Big( \PP^0_{X_n},  \PP^0_{( X_n, \alpha_n ) } \Big) \right]  = J^{\bpi}(\mu_0), \qquad \mu_0 \in \bar S,
$$
where the second equality holds by~\eqref{eq:prop-zeta-bareta-law}. 
\qed\end{proof}

\begin{remark}
Obviously, each time we involve $\vartheta^0_n$, we rely on the central unit randomization. Still, the conclusion of this lemma should not be considered as obvious because it proves that we can pack all the dependence on the past carried by the open loop controls at level $0$ into $\PP^0_{X_n},  \PP^0_{( X_n, \alpha_n)}$ and  a probability measure on the space of actions at level $1$.
\end{remark}

\section{Mean-Field Q-Function}
\label{se:Q_learning}

In preparation for our discussion of \emph{learning} in the third part of the monograph, and in line with our analysis of the first theoretical model of MFC, we introduce the \index[sub]{state-action value function} state-action value function, for the lifted MFMDP. Using standard arguments, we derive its first theoretical properties. This will pave the way for the algorithmic search for approximations of the solution of the MFC problem in a model-free setting, i.e. assuming the model is unknown. 

\subsection{State-action value function}

In order to take full advantage of the strongest results proven so far, from now on we assume that both assumptions \ref{ass:basic_assumption_MFMDP} and~\ref{ass:compactness_MFMDP} hold.
Under these assumptions, recall that the DPP given by Theorem~\ref{th:MFMDP_DPP_BOREL} holds. In this section, we restrict ourselves to pure stationary policies. Recall that when $\bar \bpi = (\bar \pi, \bar \pi, \ldots) \in \bar \bPi$ is stationary, we use the notation $\bar J^{\bar \pi} = \bar J^{\bar \bpi}$. 
Also, recall that we use the notation $\bar \Pi^p$ for the set of pure stationary strategy functions.\index[not]{$\bar \Pi^p$}

\vskip 2pt
For each $\bar \pi \in \bar \Pi^p$, the mapping
$\bar S \ni \mu \mapsto \delta_{\bar \pi(\mu)} \in \cP(\bar A)$ that assigns to each $\mu \in \bar S$ the Dirac point mass at the point $\bar \pi(\mu) \in \bar A$ is measurable by definition of the Borel $\sigma$-field of $\cP(\bar A)$.

Now, for each $\bar\pi \in \bar \Pi^p$, we introduce the state-action value function $\bar Q^{\bar \pi} : \bar\Gamma \to \RR$ defined by: \index[sub]{state-action value function}\index[not]{$\bar Q^{\bar \pi}$}
\begin{equation}
\label{eq:Q_function_stationary_control}
    \bar Q^{\bar \pi} (\mu, \bar a) := \bar f(\mu, \bar a) + \sum_{n \geq 1} \gamma^{n} \EE [\bar f(\mu_n, \bar\pi(\mu_n))], \qquad (\mu, \bar a) \in \bar\Gamma,
\end{equation}\index[not]{$\bar Q^{\bar \pi} (\mu, \bar a)$}where the process $(\mu_n)_{n\geq 0}$ starting at $\mu_0 = \mu$ satisfies $\mu_1 = \bar F(\mu, \bar a, \varepsilon^0)$ with $\cL(\varepsilon^0) = \nu^0$, and for every $n \geq 1$:
$
	\mu_{n+1} = \bar F \bigl( \mu_n, \bar \pi(\mu_{n}), \varepsilon_{n+1}^0 \bigr).
$
Next we define the optimal state-action value function by:
\begin{equation}
\label{fo:optimal_Q}
    \bar Q^* (\mu, \bar a) := \inf_{\bar\pi\in\bar \Pi^p}\bar Q^{\bar \pi} (\mu, \bar a), 
\qquad
(\mu,\bar a)\in \bar\Gamma.
\end{equation}\index[not]{$\bar Q^* (\mu, \bar a)$}
\index[sub]{optimal state-action value function}
The goal of this section is to prove the following dynamic programming principle for $\bar Q^*$.

\begin{theorem} 
 	\label{prop:opt_Bellman_eq_Q}
Assume that assumptions \ref{ass:basic_assumption_MFMDP} and~\ref{ass:compactness_MFMDP} hold. The optimal state-action value function $\bar Q^*$ satisfies the \defi{Bellman equation for state-action value function}:
\begin{equation}
\label{eq:Bellman_equation_Q}
	\bar Q^*(\mu, \bar a) = \bar f(\mu, \bar a) + \gamma \EE \Bigl[ \inf_{ \bar a' \in \bar U \big(\bar F(\mu, \bar a, \varepsilon^0) \big) }  \bar Q^*\big( \bar F(\mu, \bar a, \varepsilon^0) , \bar a'\big) \Bigr], \qquad (\mu, \bar a) \in \bar\Gamma.
\end{equation}\index[sub]{Bellman equation}
\end{theorem}
\index[not]{$\bar\Gamma$}\index[not]{$\cBD_u(\bar\Gamma)$}
\index[not]{$\Lsc(\bar\Gamma)$}
We will prove this result by showing that $\bar Q^*$ is the unique fixed point of the state-action Bellman operator $T$ defined on the set $\Lsc(\bar\Gamma)$ of bounded lower semi-continuous functions on $\bar\Gamma$, by:
\begin{equation}
\label{fo:state-action_Belman}
    [T \bar Q](\mu,\bar a) := \bar f(\mu,\bar a) 
    +\gamma\EE \Bigg[ \inf_{\bar a'\in \bar U \big( \bar F(\mu, \bar a, \varepsilon^0) \big) } \bar Q\bigl(\bar F(\mu,\bar a,\varepsilon^0),\bar a' \bigr) \Bigg],\qquad
(\mu,\bar a)\in \bar\Gamma.
\end{equation} 
We first justify in Lemma~\ref{le:strict_contraction} below the fact that the operator $T$ is well-defined.  
Since $\bar\Gamma$ is a closed subset of the Polish space $\bar S \times \bar A$, it is a Borel space on its own, and the space $\cBD_u(\bar\Gamma)$ of bounded universally measurable real-valued functions on $\bar\Gamma$ endowed with the sup norm $\| f \|_{\infty} = \sup_{(\mu, \bar a) \in \bar\Gamma} | f(\mu, \bar a) |$ is a Banach space. While the set $\Lsc(\bar\Gamma)$ is not a vector space, it is a closed subset of $\cBD_u(\bar\Gamma)$, 
hence a complete metric space for the metric $d_\infty(f, f') = \| f - f' \|_{\infty}$.

\begin{lemma}
\label{le:strict_contraction}
Let us assume that assumptions \ref{ass:basic_assumption_MFMDP} and~\ref{ass:compactness_MFMDP} hold. The set $\Lsc(\bar\Gamma)$ is invariant under the state-action Bellman operator $T$, which is a strict contraction on this metric space.
\end{lemma}

\begin{proof}
We first prove that the set  $\Lsc(\bar\Gamma)$ is invariant under $T$. We need to show that $T\bar Q$ is lower semi-continuous whenever $\bar Q$ is. To wit, by the projection property for infima of lower semi-continuous functions (see for example \cite[Proposition~7.33]{BertsekasShreve}), the function $\bar S  \ni \mu' \mapsto \inf_{\bar a' \in \bar U(\mu') }\bar Q(\mu',\bar a')$ is lower semi-continuous.
Since $\bar\Gamma \ni (\mu,\bar a)\mapsto \mu'=\bar F(\mu,\bar a,e^0)\in \bar S$ is continuous for $e^0\in E^0$ fixed, the infimum in formula \eqref{fo:state-action_Belman} is then a lower semi-continuous function of $(\mu,\bar a)$ for fixed $e^0\in E^0$. Finally, Fatou's lemma implies that the expectation in \eqref{fo:state-action_Belman} is a lower semi-continuous function of $(\mu,\bar a) \in \bar\Gamma$. Since $\bar f$ is lower semi-continuous by Assumption~\ref{ass:basic_assumption_MFMDP}, $T\bar Q$ is lower semi-continuous.
Now if $\bar Q_1$ and $\bar Q_2$ are elements of $\Lsc(\bar\Gamma)$ we have:
\begin{align*}
	&\| T \bar Q_1 - T \bar Q_2 \|_{\infty} \\
	&\hskip 5pt
	 \le   \gamma  \EE \Bigg[ \sup_{(\mu, \bar a) \in \bar \Gamma} \ \Bigl|\inf_{\bar a'\in \bar U(\bar F(\mu,\bar a,\varepsilon^0) )} \bar Q_1\bigl(\bar F(\mu,\bar a,\varepsilon^0),\bar a'\bigr) - \hskip -10pt\inf_{\bar a'\in\bar U(\bar F(\mu,\bar a,\varepsilon^0) )} \bar Q_2\bigl(\bar F(\mu,\bar a,\varepsilon^0),\bar a'\bigr) \Bigr| \ \Bigg]
	\\
	&\hskip 5pt 
	\le   \gamma  \EE \Bigg[ \sup_{(\mu, \bar a) \in \bar \Gamma}
	\sup_{\bar a'\in \bar U(\bar F(\mu,\bar a,\varepsilon^0) )}\ \Bigl| \bar Q_1\bigl(\bar F(\mu,\bar a,\varepsilon^0),\bar a'\bigr) - \bar Q_2\bigl(\bar F(\mu,\bar a,\varepsilon^0),\bar a'\bigr) \Bigr|\ \Bigg]\\
	&\hskip 5pt
	 \le \gamma \| \bar Q_1 - \bar Q_2 \|_{\infty}.
\end{align*}
Since $\gamma <1$, this proves that $T$ is a strict contraction on $\Lsc(\bar\Gamma)$. We conclude the proof of the result using the Banach fixed point theorem.
\qed\end{proof}

Using the Markov property, we can rewrite the state-action value function $\bar Q^{\bar \pi}$
in terms of the state value function $\bar J^{\bar \pi}$: 
\begin{equation}
    \label{eq:connection_J_Q}
    \bar J^{\bar \pi}(\mu) = \bar f(\mu, \bar \pi(\mu) ) + \gamma \EE [ \bar J^{\bar \pi}( \bar F(\mu, \bar \pi(\mu), \varepsilon^0) ) ] = \bar Q^{\bar \pi}\bigl( \mu, \bar \pi(\mu)\bigr), \qquad  \bar \pi \in \bar \Pi^p, \mu \in \bar S.
\end{equation}
Now, for the optimal value functions, we have the following.

\begin{lemma}
	\label{le:relationship_Q*_J*}
Suppose that assumptions~\ref{ass:basic_assumption_MFMDP} and~\ref{ass:compactness_MFMDP} hold. Then for all $(\mu, \bar a) \in \bar\Gamma$ 
$$
\bar Q^*(\mu, \bar a) = \bar f(\mu, \bar a) + \gamma \EE[  \bar J^{*}(\bar F(\mu, \bar a, \varepsilon^0)) ].
$$
\end{lemma}

\begin{proof}
We show the inequalities in both directions. By the definition of $\bar Q^*$ and of $\bar Q^{\bar\pi}$,
\[
\bar Q^*(\mu,\bar a)
=
\bar f(\mu,\bar a)
+
\gamma \inf_{\bar\pi\in\bar\Pi^p}
\EE\left[ \bar J^{\bar\pi}(\mu_1) \right].
\]
For every $\bar\pi\in\bar\Pi^p$, we have $\bar J^*(\mu')\leq \bar J^{\bar\pi}(\mu')$ for every $\mu'\in\bar S$. Hence
\[
\bar f(\mu,\bar a)+\gamma\EE\left[\bar J^*(\mu_1)\right]
\leq
\bar Q^*(\mu,\bar a).
\]
Conversely, by Theorem~\ref{th:MFMDP_DPP_BOREL}, there exists a stationary pure optimizer $\bar\pi^*\in\bar\Pi^p$ such that $\bar J^{\bar\pi^*}=\bar J^*$. Therefore
\[
\bar Q^*(\mu,\bar a)
\leq
\bar f(\mu,\bar a)+\gamma\EE\left[\bar J^{\bar\pi^*}(\mu_1)\right]
=
\bar f(\mu,\bar a)+\gamma\EE\left[\bar J^*(\mu_1)\right].
\]
The two inequalities prove the result.
\qed\end{proof}

\begin{lemma}
\label{le:Qstar_is_lsc}
Let us assume that assumptions \ref{ass:basic_assumption_MFMDP} and~\ref{ass:compactness_MFMDP} hold. Then, $\bar Q^*$ is lower semi-continuous and, as a result, there exists $\tilde \pi\in\bar \Pi^p$ such that for every $\mu \in \bar S$,
$
\tilde\pi(\mu) \in {\arg\inf}_{ \bar a \in \bar U(\mu) } \bar Q^*(\mu,\bar a).
$
\end{lemma}

\begin{proof}
For each fixed $e^0 \in E^0$, $(\mu,\bar a)\mapsto \bar J^*\bigl(\bar F(\mu,\bar a,e^0)\bigr)$ is lower semi-continuous by lower semi-continuity of $\bar J^{*}$ (see Theorem~\ref{th:MFMDP_DPP_BOREL}) and continuity of $\bar F(\cdot,\cdot,e^0)$. 
As in the proof of Lemma~\ref{le:strict_contraction}, Fatou's lemma implies that the function $(\mu,\bar a)\mapsto \EE\bigl[\bar J^*\bigl(\bar F(\mu,\bar a,\varepsilon^0)\bigr)\bigr]$ is also lower semi-continuous.
Since $\bar f$ is lower semi-continuous by Assumption~\ref{ass:basic_assumption_MFMDP}, Lemma~\ref{le:relationship_Q*_J*} implies that $\bar Q^*$ is a lower semi-continuous function on $\bar\Gamma$.

Since $\bar\Gamma$ is a closed subset of $\bar S \times \bar A$ and $\bar A$ is compact, by applying a selection theorem for lower semi-continuous function~\cite[Proposition~7.33]{BertsekasShreve} on $\bar Q^*: \bar\Gamma \to \RR$, we obtain the existence of a Borel measurable function $\tilde \pi \in \bar \Pi^p$ whose graph is contained in $\bar\Gamma$ and 
    $\bar Q^* \big(\mu, \tilde \pi(\mu) \big) = \inf_{\bar a \in \bar U(\mu) } \bar Q^*(\mu, \bar a),$
    for all $\mu \in \bar S.$
\qed\end{proof}

\begin{lemma}
\label{le:infimum_action_opt_barQ_and_opt_barJ}
Let us assume that assumptions \ref{ass:basic_assumption_MFMDP} and~\ref{ass:compactness_MFMDP} hold. Then, for all $\mu \in \bar S$, 
$$
\inf_{\bar a\in \bar U(\mu)} \bar Q^*(\mu, \bar a) = \bar J^{*}(\mu).
$$
\end{lemma}

\begin{proof}
    We first show the inequality $ \inf_{\bar a\in \bar U(\mu)} \bar Q^*(\mu, \bar a) \geq \bar J^*(\mu)$. 
    Let us denote by ${\tilde\pi} \in \bar \Pi^p$ the strategy function obtained in Lemma~\ref{le:Qstar_is_lsc}. By definition of $\bar Q^*$, 
    $$
    	\inf_{\bar a \in \bar U(\mu) } \bar Q^*(\mu, \bar a) = \bar Q^*(\mu, {\tilde\pi}(\mu) ) = \inf_{\bar \pi \in \bar \Pi^p } \bar Q^{\bar \pi} (\mu, {\tilde\pi}(\mu) ), \qquad \mu \in \bar S.
	$$ 
If for each $\bar \pi \in \bar \Pi^p$, we denote 
$
         \bar \bpi^{{\tilde\pi}} = ( {\tilde\pi}, \bar \pi, \bar \pi, \bar \pi, \ldots ) \in \bar\bPi
$, then  
    $$
        \bar Q^{\bar \pi}(\mu, {\tilde\pi}(\mu) ) 
        = \bar f(\mu, {\tilde\pi}(\mu) ) + \EE\left[ \sum_{n=1}^\infty \gamma^n \bar f( \mu_n, \bar \pi(\mu_n)  ) \right] 
        = \bar J^{ \bar \bpi^{{\tilde\pi}} } (\mu) \geq \bar J^*(\mu), \qquad \  \mu \in \bar S,
    $$
    which provides the first inequality. 
	To prove the converse inequality, let $\bar \pi^* \in \bar \Pi^p$ be an optimal non-randomized stationary Markov policy whose existence is given by Theorem~\ref{th:MFMDP_DPP_BOREL}, and  notice that for every $\mu \in \bar S$, 
    $$
        \bar J^{*}(\mu) = \bar J^{\bar \pi^*}(\mu) = \bar Q^{\bar \pi^*}(\mu, \bar \pi^*(\mu) ) \geq \bar Q^*(\mu, \bar \pi^*(\mu) ) \geq \inf_{\bar a \in \bar U(\mu)} \bar Q^*(\mu, \bar a),
    $$
    by equation~\eqref{eq:connection_J_Q} and the result of Lemma~\ref{le:relationship_Q*_J*}.  
This concludes the proof.
\qed\end{proof}

We can now complete the proof of Theorem~\ref{prop:opt_Bellman_eq_Q}. 
\begin{proof}[Proof of Theorem~\ref{prop:opt_Bellman_eq_Q}]
Bellman equation~\eqref{eq:Bellman_equation_Q} is a direct consequence of Lemma~\ref{le:relationship_Q*_J*} and 
 Lemma~\ref{le:infimum_action_opt_barQ_and_opt_barJ}. Since $T$ is a strict contraction on $\Lsc(\bar\Gamma)$ by Lemma~\ref{le:strict_contraction}, and since $\Lsc(\bar\Gamma)$ is closed in the Banach space $\cBD_u(\bar\Gamma)$, by the Banach fixed point theorem we conclude that $\bar Q^*$ is the unique fixed point of $T$ on $\Lsc(\bar\Gamma)$.
\qed\end{proof}

\section{Notes and Complements}

MFMDPs have been studied at least since the work of Gast et al.~\cite{gast2011mean,gast2012mean} who studied the convergence of a finite-agent MDP towards a MFMDP. More recently, Gu et al. studied in~\cite{gu2024meanmarldecentralized,gu2021meanQ} dynamic programming for the state-action value function arising in MFC problems without common noise, and they proposed RL methods for this setting. Following the posting on {\tt arxiv}  of the original version of our paper \cite{CarmonaLauriere_AAP},  Motte and Pham tackled in~\cite{motte2022mean} the open loop case in the presence of common noise, in a model very similar to the model considered in this chapter. They proved a propagation of chaos result in the spirit of Proposition~\ref{pr:chaos} of Chapter \ref{ch:prelims}, a DPP for the MFMDP, as well as the equality of optimal value functions for open-loop and closed-loop policies in their model. In this chapter, we proposed a finer model by paying close attention to all the sources of randomization, and by relying on fine measure theoretical arguments. Comparing to~\cite{motte2022mean}, while they prove equality of the optimal value functions without the central randomization, the latter allows us to prove the identity of the value functions, policy by policy, even before taking the optimum values, allowing for a deeper understanding of the agent level MFC, as done for example in Corollary \ref{co:early}.

\vskip 2pt
Note that when the state space $S$ is discrete, assuming continuity of $F(\cdot, \cdot, \cdot, e, e^0)$ implies that the dependence on the mean-field term is basically constant when the other inputs are fixed. This is very much in line with the fact that many MFG and MFC models appearing in the literature do not have mean-field interactions in the dynamics of the states, the interactions appearing explicitly only in the cost functions. See for example the original models in~\cite{LasryLions2007} as well as the applications to Bertrand and Cournot competitions~\cite{GueantLasryLions,ChanSircar}, and the models of crowd motion~\cite{achdoulasry2019meancrowd} or flocking~\cite{NourianCainesMalhame,BardiCardaliaguet}. To cover models with interactions in the dynamics of the states in a finite state space setting, one should use, instead of the assumptions mentioned in Remark~\ref{rem:assumptions-mfc}, other assumptions that are sufficient to guarantee Assumption~\ref{ass:basic_assumption_MFMDP}. In particular, it is expected that the idiosyncratic noise has a regularizing effect that helps ensuring continuity of $\bar F$ even when $F$ is discontinuous. For instance~\cite[Proposition 2.1]{motte2022mean} does not require continuity of $F$ but only some kind of almost everywhere continuity.

\chapter{A Linear Quadratic Mean Field Model with Common Noise}
\label{ch:LQMFC}

\begin{abstract}
\emph{
In this chapter, we provide a complete theoretical analysis of the Linear Quadratic (LQ) model of Mean Field Control (MFC) with common noise which was introduced in Section \ref{se:prelims_LQ} of Chapter \ref{ch:prelims}. We use the notation LQMFC \index[not]{LQMFC} to denote such a model.
We characterize optimal policies and prove their existence, preparing for the numerical analysis and RL implementations of Chapter \ref{ch:numeric_II}.
}
\end{abstract}

\section{Introduction}
\label{se:LQ_intro}
The first goal of the chapter is to properly introduce the notations and the definitions of the Mean Field (MF) Linear Quadratic  (LQ) stochastic optimization problem which was introduced in broad strokes in Chapter \ref{ch:prelims}. After arguing that the model appears naturally as the limit of a large multi agent model with a common noise, we 
prove that the optimal admissible controls exist, are unique, and  linear in the state and its conditional mean processes. We formulate and prove the necessary and the sufficient conditions for a discrete-time form of the Pontryagin maximum principle for our LQMFC model. This result implies existence and uniqueness of the optimal control process. Incidentally, this provides a crucial link between the MFC problem and learning by a finite population of agents, as it gives an approximately optimal control for the problem with a finite number of learners.

\section{Models}
\label{se:def-model}

We work with the model introduced at the beginning of Section \ref{se:prelims_LQ} of Chapter \ref{ch:prelims} and repeat the meaning of the notations only when necessary.
Because the original probability space $(\Omega, \cF,\PP)$ is assumed to be rich enough to provide for all the randomizations needed, and hence to facilitate the analysis of our LQ models, we assume the existence of two initial random variables, denoted by $\varepsilon_0^0$ and $\varepsilon_0$, called \emph{initial perturbations}, that are independent of each other and independent of the noise sequences $\bvarepsilon^0=(\varepsilon_n)_{n\ge 1}$ and $\bvarepsilon=(\varepsilon^0_n)_{n\ge 1}$. The random variable $\varepsilon_0^0$ represents the common shock at initial time from the environment. 
We distinguish the two initial perturbations $\varepsilon_0$ and $\varepsilon_0^0$ from the aforementioned noise terms. They are not assumed to follow the mean-zero distributions of the noise sequences. Also, they can be constant, i.e. random variables with a  Dirac measure as distribution. This is the case when we want the LQ model to start from a deterministic initial condition. 

\vskip 2pt
In this chapter, we assume that the initial perturbations $\varepsilon_0$ and $\varepsilon_0^0$, and all noise terms $(\varepsilon_{n+1})_{n \geq 0}$ and $(\varepsilon_{n+1}^0)_{n \geq 0}$ are \textbf{sub-Gaussian random vectors}. 
\index[sub]{sub-Gaussian}
Recall Definition \ref{de:subGaussian} at the end of Chapter \ref{ch:prelims}.

\vskip 2pt
We shall  extensively use the filtrations $(\sigma\{\underline\epsilon_n,\underline\epsilon^0_n\})_{n \geq 0}$ and $(\sigma\{\underline\epsilon^0_n\})_{ n \geq 0}$ that did not play much of a role in the analysis of Chapter \ref{ch:AAP}. Recall the meaning of the underline, namely $\underline\varepsilon_n=(\varepsilon_0,\varepsilon_1,\cdots,\varepsilon_n)$ and similarly for $\underline\varepsilon_n^0$.

\subsection{Model with a finite number of agents}
\label{sub:OC-N-Agents}

As stated above, we justify the mean field model by first considering the corresponding finite agent model and taking the limit as the number of agents goes to $\infty$.

\vskip 2pt
We assume that the state of the system is given by the $N$-tuple of individual agent states, and that the impact of the common environment is captured by a common noise included in the system functions by a transformation \emph{\`a la Blackwell-Dubins}. 
Our desire for linear dynamics and  \emph{mean field} interactions dramatically restricts the form of the system function.
The individual state dynamics in the stochastic optimization problem for $N$ agents with linear dynamics and symmetric interactions are given by the system of equations: 
\begin{align}
\label{fo:N-multi_state}
    X_{n+1}^i = {\mathrm A}  X_n^i + \bar{\mathrm A}  \bar{X}_n^{N}  + {\mathrm B}  a^i_n + \bar{{\mathrm B} } \bar{a}_n^N + \comnoise{n+1} + \idynoise{n+1}^{i},\qquad i=1,\cdots,N,\quad n\ge 0
\end{align}
with initial conditions $X_0^{i} = \varepsilon_0^0 + \varepsilon_0^{i}$, for $i = 1, \ldots, N$, where 
$$
\bar{X}^N_n = \frac1N \sum_{i=1}^N X_n^{i}
\qquad\text{and}\qquad 
\bar{a}^N_n = \frac1N\sum_{i=1}^N a_n^{i}
$$
denote the sample averages of the individual states and actions, and $( \varepsilon_0^{i}, (\idynoise{n+1}^{i})_{n \geq 0})_{i=1,\ldots, N}$ are $N$ identical and independent copies of the idiosyncratic noise sequence $( \varepsilon_0, (\varepsilon_{n+1})_{n \geq 0} )$.

\vskip 2pt
Recall that when considering a LQ model, the letter $A$ does not represent the set of possible actions, but now a deterministic matrix entering the state dynamics.

Averaging the $N$ dynamical equations above gives the time evolution of the sample average of the state as:
\begin{equation}
\label{fo:N-mean_state_dynamics}
     \bar{X}^N_{n+1} = ({\mathrm A}  + \bar{\mathrm A} ) \bar{X}^N_n  +  ({\mathrm B}  + \bar {\mathrm B} ) \bar{a}^N_n + \varepsilon^0_{n+1} + \frac{1}{N}\sum_{i=1}^N \varepsilon^{i}_{n+1},
\end{equation}
with initial value $\bar{X}^N_0 = \comnoise{0} + \frac{1}{N} \sum_{i=1}^N \idynoise{0}^{i}$.
Note that for $n$ fixed, the classical law of large numbers says that the last term in~\eqref{fo:N-mean_state_dynamics} should converge when $N\to\infty$ toward the common expectation of the $\varepsilon^{i}_{n+1}$ which is $0$. On the other hand, the random shock $\varepsilon^0_{n+1}$ affects all the agents and its impact should not disappear in the limit $N\to\infty$.

Letting $\underline{X}_n= Vect\big(X_n^{1}, \cdots, X_n^{N} \big)$ and $\underline{a}_n= Vect \big( a_n^{1}, \cdots, a_n^{N} \big)$ for $n\ge 0$, we can rewrite the system \eqref{fo:N-multi_state} of $N$ dynamical equations in a single dynamical equation: 
\begin{equation}
\label{fo:N-vector_state}
    \underline{X}_{n+1} = {\mathrm A} ^N \underline{X}_n  + {\mathrm B} ^N \underline{a}_n + \underline{\varepsilon}^0_{n+1} + \underline{\varepsilon}_{n+1}
\end{equation}
with 
$\underline{\varepsilon}_{n+1} = Vect( \idynoise{n+1}^{1}, \ldots ,\idynoise{n+1}^{N} )$, $\underline{\varepsilon}^0_{n+1} = Vect( \comnoise{n+1},\ldots,\comnoise{n+1} )$, $ {\mathrm A} ^{N} = I_N \otimes {\mathrm A}  + \frac{1}{N} \bone_N \otimes \bar{{\mathrm A} }$ and ${\mathrm B} ^{N} = I_N \otimes {\mathrm B}  + \frac{1}{N} \bone_N \otimes \bar{{\mathrm B} } $. In other words, the matrix ${\mathrm A} ^{N}$ is the sum of the block-diagonal matrix with $N$ blocks identically equal to ${\mathrm A} $, and of the matrix with $N \times N$ blocks $\frac{1}{N}\bar{{\mathrm A} }$, and similarly for the matrix ${\mathrm B} ^{N}$. 

Accordingly, since we assume that the agents cooperate with each other (instead of competing against each other), the goal is to minimize the \textit{social cost} of the population, defined as
\begin{equation}    \label{fo:social_cost_of_population}
    J^N( \underline{\ba} ) = \EE \Big[ \sum_{n \geq 0} \gamma^{n} \bar f^N(\underline{X}_n, \underline{a}_n) \Big]
\end{equation}
over the set of admissible control processes $\underline{\ba} = (\underline{a}_n)_{n \geq 0}$ consisting of $\RR^{\ell N}$-valued processes adapted to the filtration $(\sigma\{\underline\epsilon^0_n,\underline\epsilon_n\})_{n \geq 0}$ for which $J^N(\underline{\ba}) < \infty$.
Here, the one-stage \textit{social cost} function $\bar f^N : \RR^{dN} \times \RR^{\ell N} \to \RR$ is
\begin{equation}
    \label{fo:N-cost}
    \bar f^{N}( \underline{X}_n, \underline{a}_n) := \frac{1}{N} \sum_{i=1}^N f \big( X_n^i, \bar{X}_n^N, a_n^i, \bar{a}_n^N \big) = \underline{X}_n^\top Q^N \underline{X}_n + \underline{a}_n^\top R^N \underline{a}_n,
\end{equation}
where $Q^N = \frac{1}{N} I_N \otimes Q + \frac{1}{N^2} \boldsymbol{1}_{N} \otimes \bar{Q} \in \RR^{dN \times dN}$ and $R^N =  \frac{1}{N} I_N \otimes R + \frac{1}{N^2} \boldsymbol{1}_{N} \otimes \bar{R} \in \RR^{\ell N \times \ell N}$. 

\vskip 4pt
This is now a classical infinite-horizon discounted-cost linear-quadratic control problem with $(d N)-$dimensional state process $\underline{\bX} = (\underline{X}_n)_{n \geq 0}$ following dynamics~\eqref{fo:N-vector_state}. 
Under Assumptions~\ref{ass:finit_cost} and~\ref{ass:positivity-qr} stated below, an optimal control $\underline{\ba}^{*,N}$ exists and is necessarily of the form
$\underline{a}^{*,N}_n = \Phi^{*,N} \underline{X}_n$ 
for a deterministic constant matrix $\Phi^{*,N} \in \RR^{(\ell N) \times (d N)}$.

\subsection{The mean-field control model}
\label{subsec:OC-MKV}

Motivated by the analysis of the large $N$ limit of the above model, we envision a generic agent the state of which, say $X_n\in\RR^d$, evolves over time according to the equation:
\begin{align}
    X_{n+1}
    & = {\mathrm A}  X_n + \bar{{\mathrm A} }\bar{X}_n  + {\mathrm B}  a_n + \bar{{\mathrm B} }\bar{a}_n + \varepsilon^0_{n+1} + \varepsilon_{n+1},
\label{fo:MKV-state}
\end{align}
for $n \geq 0$, with an initial state $X_0 =\varepsilon^0_{0}+\varepsilon_{0}$.
Let $\tilde \mu^0_0 = \cL(\varepsilon_0^0)$ and $\tilde \mu_0 = \cL(\varepsilon_0)$ be their respective distributions.
The terms $\bar{X}_n = \EE[X_n|\sigma\{\underline\epsilon^0_n\}]$ and $\bar{a}_n = \EE[a_n|\sigma\{\underline\epsilon^0_n\}]$ appearing in the dynamics~\eqref{fo:MKV-state} of the state are the conditional means of $X_n$ and $a_n$ given $\sigma\{\underline\epsilon^0_n\}$, the past of the common noise up to and including time $n$. 
The Mean Field cost function is
\begin{equation}
    \label{fo:MKV-discounted_cost}
    J(\ba)=\EE \Big[ \sum_{n \geq 0} \gamma^n f(X_n, \bar{X}_n, a_n, \bar{a}_n) \Big].
\end{equation}
We choose a one-stage cost function $f$ which is quadratic in the state and the action. However, in order to make the derivations in the proofs of this chapter slightly easier to follow, we depart slightly from the notation of the previous section on $N$ agents and we choose:
\begin{equation}
    \label{fo:quadratic_cost}
f(x, \bar{x}, a, \bar{a})
= (x-\bar x)^\top Q (x-\bar x) + \bar{x}^\top (Q+\bar Q) \bar{x}
+ (a-\bar a)^\top R (a-\bar a) + \bar{a}^\top (R+\bar R) \bar{a}
\end{equation}
where we use the superscript ${}^\top$ to denote the transpose of a vector or a matrix.
Notice that the expected one stage cost at time $n$ reads
\begin{equation}
\label{fo:expression_expected_instateneous_cost}
\begin{split}
\EE[ f(X_n, \bar{X}_n, a_n, \bar{a}_n)] &= \EE[ (X_n-\bar X_n)^\top Q (X_n-\bar X_n) + \bar{X}_n^\top (Q+\bar Q) \bar{X}_n\\
&\hskip 45pt
+ (a_n-\bar a_n)^\top R (a_n-\bar a_n) + \bar{a}_n^\top (R+\bar R) \bar{a}_n ].
\end{split}
\end{equation}
The set of \textit{admissible control processes} is denoted by $\cA_{ad}$ \index[not]{$\cA_{ad}$} and defined as the collection of sequences of random variables on $(\Omega,\cF,\PP)$ adapted to the filtration $(\sigma\{\underline\epsilon_n,\underline\epsilon^0_n\})_{n \geq 0}$ and square integrable (a.k.a. $L^2-$discounted integrable):

\begin{equation}
\label{def:admissible_control}
    \cA_{ad} := \left\{ (a_n)_{n \geq 0} \left|\, a_n \in \RR^{\ell} \text{ is } \sigma\{\underline\epsilon_n,\underline\epsilon_n^0\}-\text{measurable}, \quad \EE \Big[ \sum_{n \geq 0} \gamma^n \Vert a_n \Vert^2 \Big] <\infty  \right. \right\}.
\end{equation}\index[not]{$\cA_{ad}$}
We also consider the set $\cX$ of $L^2-$discounted-integrable processes in $\RR^{d}$ defined by:
\begin{equation}
\label{def:L2_discoutned_integratble_process}
    \cX := \left\{ (X_n)_{n \geq 0} \left|\, X_n \in \RR^{d} \text{ is } \sigma\{\underline\epsilon_n,\underline\epsilon_n^0\}-\text{measurable}, \quad \EE \Big[ \sum_{n \geq 0} \gamma^n \Vert X_n \Vert^2 \Big] <\infty  \right. \right\}.
\end{equation}\index[not]{$\cX$}
Note that both $\cA_{ad}$ and $\cX$ are linear spaces.  
The MFC problem is then:
\begin{equation}
    \label{pb:MFC_L2_admissible_control}
    \inf_{\bu \in \cA_{ad}} J(\bu).
\end{equation}

The above stochastic control problem has some unique characteristics: 1) it is set in infinite horizon; 2) it includes a \emph{common noise} $\bvarepsilon^0$, 3) The optimization is set over open-loop controls, and 4) the interaction is not only through the conditional mean of the state, but also through the conditional mean of the control. 

\subsection{Assumptions and the well-definedness of the problem}

Throughout the whole chapter, we impose the following assumptions which we already stated in the introduction of the model in Chapter \ref{ch:prelims}.
\begin{assumption}
\label{ass:finit_cost}
    $\gamma \in (0,1)$ and the matrices ${\mathrm A} $ and $\bar{\mathrm A} $ satisfy $\gamma \| {\mathrm A}  \|^2 < 1$ and $\gamma \| {\mathrm A}  + \bar {\mathrm A}  \|^2 < 1$.
\end{assumption}

\begin{assumption}
\label{ass:positivity-qr}
    The matrices $Q, \bar Q, R, \bar R$ are symmetric and they satisfy $Q \succeq 0, Q + \bar Q \succeq 0, R \succ 0, R + \bar R \succ 0$.
\end{assumption}
Assumption~\ref{ass:finit_cost} is stronger than standard stabilizability / detectability conditions as it may exclude some unstable but stabilizable systems. As we shall see in the following sections, %
it is sufficient to ensure discounted square integrability of the state process and finiteness of the mean field cost for every admissible control.
Taking advantage of the strong convexity in the control variables provided by $R\succ0$ and $R+\bar R\succ0$, we show in Section~\ref{se:analysis} the existence and uniqueness of the optimal control process for the MFC problem.

We will sometimes use the notation $\bX^\ba$ if necessary to emphasize the fact that the state process following dynamics~\eqref{fo:MKV-state} is controlled by the admissible control process $\ba \in \cA_{ad}$.
Lemma~\ref{le:L2_discounted_integrable_of_state_process} and Proposition~\ref{pr:admissible_of_X} below justify the definitions of the set $\cA_{ad}$ of admissible control processes and $\cX$ of controlled state processes, and the fact that the discounted mean field cost $J(\ba)$ for admissible control process $\ba \in \cA_{ad}$ is finite.  

\vskip 4pt
The following result provides a simple procedure to identify $L^2-$discounted-integrable state processes.

\begin{lemma}
\label{le:L2_discounted_integrable_of_state_process}
    If we assume that $\gamma \| {\mathrm A}  \|^2 < 1$ and that the process $\bY= (Y_n)_{n \geq 0}$ satisfies 
    \begin{equation}
    \label{fo:sys_dyn_Y}
        Y_{n+1} = {\mathrm A}  Y_n +  \beta_{n+1}, \qquad n\ge 0,
    \end{equation}
for a $\sigma\{\epsilon^0_0\}$-measurable and square-integrable initial random variable $Y_0$, and a process $ (\beta_n)_{n \geq 0} \in \cX$,
then,  $\bY \in \cX$.
\end{lemma}

\begin{proof}
The fact that $Y_n$ is $\sigma\{\underline\epsilon_n,\underline\epsilon^0_n\}$-measurable can be read directly from the dynamics of $\bY = (Y_n)_{n \geq 0}$. Also, for every $n \geq 1$,
\begin{align*}
    \EE\left[ \gamma^n \| Y_n \|^2 \right] 
    &= \EE \bigg[ \gamma^n \Big\| {\mathrm A} ^n Y_0 + \sum_{m=0}^{n-1} {\mathrm A} ^{n-1-m} \beta_{m+1} \Big\|^2 \bigg]  
    \\
    & =  \EE \bigg[ \Big\| (\gamma^{1/2} {\mathrm A} )^n Y_0 + \sum_{m=1}^{n} (\gamma^{1/2} {\mathrm A} )^{n-m} \big(\gamma^{m/2} \beta_{m} \big)  \Big\|^2 \bigg] 
\end{align*}
Since $\gamma \| {\mathrm A}  \|^2 < 1$, there exists a real number $1 > \rho > \gamma$ such that $\xi := \rho \| A \|^2 < 1$ (such a $\rho$ exists because $\gamma \|A\|^2 < 1$). Let $\eta := \gamma \|A \|^2 / \xi $, then $\eta = \gamma / \rho < 1$. 
Using Fubini's theorem, we switch the two summations and obtain:
$$
    \sum_{n=1}^\infty \Big( \sum_{m=1}^{n} \xi^{n-m} \EE\big[ \gamma^m  \| \beta_m \|^2 \big] \Big) =  \sum_{m=1}^{\infty} \Big( \sum_{n=0}^\infty  \xi^n \Big) \EE\big[ \gamma^{m} \| \beta_m \|^2  \big] = \frac{1}{1 - \xi} \sum_{m \geq 0} \gamma^{m+1} \EE[ \| \beta_{m+1} \|^2 ] < \infty.
$$
Consequently, by summing up $\EE[ \gamma^n \| Y_n\|^2 ]$ from $n=0$ to $\infty$, we obtain:
\begin{align*}
    \sum_{n \geq 0}  \EE\left[ \gamma^n \| Y_n \|^2 \right]  
    & \leq \EE\big[ \| Y_0\|^2 \big] + \sum_{n \geq 1} \EE \bigg[ \Big( \| \gamma^{1/2} {\mathrm A}  \|^n  \| Y_0 \|  + \sum_{m=1}^{n} \| \gamma^{1/2} {\mathrm A}  \|^{n-m} \gamma^{m/2} \| \beta_m \| \Big)^2 \bigg] 
    \\
    & =  \EE\big[ \| Y_0\|^2 \big] + \sum_{n\geq 1}  \EE \bigg[ \Big( \eta^{n/2} (\xi^{n/2} \| Y_0 \| ) + \sum_{m=1}^{n} \eta^{(n-m)/2} \xi^{(n-m)/2} \gamma^{m/2} \| \beta_m \| \Big)^2 \bigg]
    \\
    & \leq  \EE\big[ \| Y_0\|^2 \big] +  \sum_{n\geq 1}  \EE \bigg[ \Big(\eta^n + \sum_{m=1}^{n} \eta^{n-m} \Big). \Big( \xi^n \|Y_0 \|^2 + \sum_{m=1}^{n} \xi^{n-m} \gamma^{m} \| \beta_m \|^2 \Big)   \bigg]
    \\
    & = \EE \big[ \|Y_0 \|^2 \big] + \frac{1}{1- \eta} \bigg( \sum_{n = 1}^\infty \xi^n \EE \big[ \|Y_0 \|^2 \big] +  \sum_{n = 1}^\infty \sum_{m=1}^{n} \xi^{n-m}  \EE \big[ \gamma^{m} \| \beta_m \|^2 \big] \bigg) 
    \\
    & \leq \frac{1}{1 - \eta} \frac{1}{1 - \xi} \bigg( \EE\big[ \| Y_0\|^2 \big] + \sum_{m \geq 0} \gamma^{m+1} \EE \big[ \| \beta_{m+1} \|^2 \big] \bigg)
    \\
    & < \infty.
\end{align*}
The first inequality is due to the matrix norm inequality $\| {\mathrm A} ^n \| \leq \| {\mathrm A}  \|^n$ for any $n \geq 1$, and the third inequality is justified by Cauchy-Schwarz inequality.
This proves that $\bY \in \cX$.
\qed\end{proof}

\begin{proposition}
\label{pr:admissible_of_X}
Under Assumption~\ref{ass:finit_cost}, for any admissible control $\ba \in \cA_{ad}$, the corresponding state process $\bX^\ba$ is in $\cX$. Moreover, the mean field cost $J(\ba)$ is finite.
\end{proposition}

\begin{proof}
    We use two auxiliary processes $(X_n^\ba - \bar X^{\ba}_n)_{n \geq 0}$ and $(\bar X_n^\ba)_{n \geq 0}$. Their dynamics satisfy
    \begin{align}
        X_{n+1}^\ba - \bar X_{n+1}^\ba &= {\mathrm A}  \big(X_n^\ba - \bar X^{\ba}_n\big) + {\mathrm B}  (a_n - \bar a_n) + \varepsilon_{n+1},
        \label{fo:Xt_minus_barXt_u}
        \\
       \bar X_{n+1}^\ba &= ({\mathrm A}  + \bar {\mathrm A} )  \bar X^{\ba}_n + ({\mathrm B}  + \bar {\mathrm B} ) \bar a_n + \varepsilon_{n+1}^0,
       \label{fo:barXt_u}
    \end{align}
    for $n \geq 0$, with the initial conditions $X_0^\ba - \bar X^{\ba}_0 = \varepsilon_0 - \EE[ \varepsilon_0]$ and $\bar X^{\ba}_0 = \varepsilon_0^0 + \EE[ \varepsilon_0]$. Let $\beta_{n+1}^y = {\mathrm B}  (a_n - \bar a_n) + \varepsilon_{n+1}$ for every $n \geq 0$ and $\beta_0^y = 0$. Because $\ba \in \cA_{ad}$ and $(\varepsilon_{n+1})_{n \geq 0}$ are i.i.d., by Cauchy-Schwarz inequality and Jensen's inequality for conditional expectations we get:
    $$
        \EE \big[ \| \beta_{n+1}^y \|^2 \big] \leq \big( \| B \|^2 + \|B \|^2 + 1 \big). \big( \EE \big[ \| a_n \|^2 \big] + \EE\big[ \EE[ \| a_n \|^2 \, | \, \sigma\{\underline\epsilon^0_n\} ] \big] + \EE\big[ \| \varepsilon_{n+1} \|^2 \big] \big) 
    $$
    for every $n \geq 0$, so that
    $$
        \sum_{n \geq 0} \gamma^n \EE\big[ \| \beta_n^y \|^2 \big] \leq  ( 2 \| {\mathrm B}  \|^2 + 1)  \Big( 2 \sum_{n \geq 0} \gamma^n \EE[ \| a_n \|^2 ]  + \sum_{n \geq 0} \gamma^n \EE[ \| \varepsilon_{n+1} \|^2 ] \Big) < \infty.
    $$
    Thus, $(\beta_n^y)_{n\geq 0} \in \cX$ because $\beta_n^y$ is $\sigma\{\underline\epsilon_n,\underline\epsilon^0_n\}$-measurable. We then apply Lemma~\ref{le:L2_discounted_integrable_of_state_process} with $(Y_n, \beta_n) = \big(X_n^\ba - \bar X^{\ba}_n, \beta_n^y \big)$ for $n \geq 0$. Under Assumption~\ref{ass:finit_cost}, we have $\gamma \| {\mathrm A}  \|^2 < 1$, and hence $(X_n^\ba - \bar{X}_n^\ba)_{n \geq 0} \in \cX$.  
    
    Similarly, let $\beta_{n+1}^z = ({\mathrm B}  + \bar {\mathrm B} ) \bar a_{n} + \varepsilon_{n+1}^0$ and $\beta_0^z = 0$. We have $(\beta_n^z)_{n \geq 0} \in \cX$. We consider $(Y_n, \beta_n) = \big( \bar X_n^\ba, \beta_n^z \big)$ for every $n \geq 0$ and replace ${\mathrm A} $ with ${\mathrm A}  + \bar {\mathrm A} $ in Lemma~\ref{le:L2_discounted_integrable_of_state_process}. Assumption~\ref{ass:finit_cost} implies that $(\bar X_n^{\ba})_{n \geq 0} \in \cX$. Because $\| X_n^\ba \|^2 \leq 2 \| X_n^\ba - \bar{X}_n^\ba \|^2  + 2 \| \bar{X}_n^\ba \|^2 $ for $n \geq 0$, we have $\bX^\ba \in \cX$.
    Therefore, from the definition of the mean field cost $J(\ba)$ in equation~\eqref{fo:MKV-discounted_cost} and the expectation of one-stage cost~\eqref{fo:expression_expected_instateneous_cost}, we conclude using Jensen's inequality that 
    $$
        J(\ba) \leq (\| Q \| +  \| \bar Q \|)   \sum_{n \geq 0} \gamma^n \EE \big[ \| X_n^\ba \|^2 \big] + ( \| R \| +  \| \bar R \|) \sum_{n \geq 0} \gamma^n \EE \big[ \| a_n \|^2 \big]  < \infty.
    $$
\qed\end{proof}

\begin{corollary}
\label{cor:uniqueness_of_controlled_state_process}
    If Assumption~\ref{ass:finit_cost} holds, for any admissible control process $\ba \in \cA_{ad}$,  there exists a unique process $\bX \in \cX$ controlled by $\ba$
    starting from $X_0 = \varepsilon_0 + \varepsilon_0^0$ that follows dynamics~\eqref{fo:MKV-state} with noise processes $(\bvarepsilon, \bvarepsilon^0)$.
\end{corollary}

\begin{proof}
    The existence of controlled state process $\bX \in \cX$ controlled by $\ba \in \cA_{ad}$ is given by Proposition~\ref{pr:admissible_of_X}. For the uniqueness, suppose that another process $\bX' \in \cX$ starts from the same initial state $X_0' = X_0$ and follows dynamics~\eqref{fo:MKV-state} with the same noise processes $\bvarepsilon$ and $\bvarepsilon^0$. Then, the difference between $X_n$ and $X'_n$ satisfies a deterministic dynamics
    $$
        X_{n+1} - X'_{n+1} = {\mathrm A}  (X_n - X'_n) + \bar{{\mathrm A} } \EE\big[ X_n - X'_n \, | \sigma\{\underline\epsilon^0_n\} \big]
    $$
    for every $n \geq 0$ and $X_0 - X_0' = 0$ almost surely. Taking the conditional expectation $\EE[ \cdot | \sigma\{\underline\epsilon^0_n\}]$ in the above dynamics  for every $n \geq 0$ and in the initial condition, we deduce that $\EE[ X_n - X'_n \, | \, \sigma\{\underline\epsilon^0_n\} ] = 0$ almost surely for every $n \geq 0$. Plugging back into the deterministic dynamics, we obtain that $X_n = X'_n$ almost surely for every $n \geq 0$ proving uniqueness. 
\qed\end{proof}

\begin{remark}
\label{re:equivalence_of_welldefineness_and_admissibility}
Under Assumption~\ref{ass:positivity-qr}, the matrices $R$ and $R + \bar R$ are (strictly) positive definite, and thus $\EE[ (a_n - \bar {a}_n)^\top R (a_n - \bar a_n) ] \geq  \lambda_{min}(R) \EE[ \| a_n - \bar a_n \|^2 ] $ and $\EE[ (\bar{a}_n)^\top (R + \bar R) \bar{a}_n] \geq \lambda_{min}(R + \bar R) \EE[ \| \bar a_n \|^2 ]$ for every $n \geq 0$. 
So if the mean field cost $J(\ba)$ is finite for a $(\sigma\{\underline\epsilon_n,\underline\epsilon^0_n\})_{n \geq 0}$-adapted control process $\ba = (a_n)_{n \geq 0}$, if follows that $(a_n - \bar a_n)_{n \geq 0}$ and $(\bar a_n)_{n \geq 0}$ are $L^2-$discounted integrable. Consequently, the control process $\ba$ is admissible in the sense that $\ba \in \cA_{ad}$.  Together with Proposition~\ref{pr:admissible_of_X}, we see that under Assumption~\ref{ass:finit_cost} and Assumption~\ref{ass:positivity-qr}, the $L^2-$discounted integrability condition in $\cA_{ad}$ is equivalent to the finiteness of the mean field cost $J(\ba)$. 
\end{remark}

In Section~\ref{se:analysis}, we will show existence and uniqueness of an optimal admissible control process for the MFC problem~\eqref{pb:MFC_L2_admissible_control} using a discrete-time version of Pontryagin's maximal principle. Moreover, the optimal control process $\ba$ will be shown to be in a feedback form and to be linear in $(\bX^\ba, \bar \bX^\ba)$, with coefficients independent of time. To rigorously establish these results with the aforementioned admissible controls, we define the following set of admissible control parameters. 

\begin{definition}
    The set of \textbf{admissible control parameters} is defined as
    \begin{equation}
        \label{fo:admissible_parameters_LQMF}
        \Theta := \left\{ \theta = (K,L) \in \mathbb{R}^{\ell \times d} \times \mathbb{R}^{\ell \times d} \  ; \  \gamma \| {\mathrm A}  - {\mathrm B}  K\|^2 < 1, \, \gamma \|  {\mathrm A}  + \bar {\mathrm A}   - ({\mathrm B}  + \bar {\mathrm B} ) L\|^2 < 1 \right\}.
    \end{equation} 
\end{definition}\index[not]{$\Theta$}
The set $\Theta$ is a convex subset of $\RR^{\ell \times d} \times \RR^{\ell \times d}$ because the spectral norm $\| \cdot \|$ is convex.
It is worth noting that in discrete-time control problems, contrary to their continuous-time analogs, there is no issue concerning the existence of the state process controlled by a control process given in feedback form. The following lemma states that when the feedback form is defined through a policy function $\pi_\theta : (x, \bar x) \mapsto \pi_\theta(x, \bar x)$ defined below, and the control parameter $\theta$ involved within $\pi_\theta$ is admissible, the control process $\bu$ is also admissible.

\begin{lemma}
\label{le:admissible_of_control_with_feedback_form_on_theta}
    For any admissible control parameter $\theta = (K, L) \in \Theta$ and any process $\ba = (a_n)_{n \geq 0}$ in $\RR^{\ell}$ in feedback form satisfying 
    \begin{equation}
    \label{fo:feedback_linear_form_of_control}
        a_n = - K (X_n^\ba - \bar X_n^\ba) - L \bar X_n^{\ba}
    \end{equation}
    for every $n \geq 0$, the control process $\ba$ is admissible in the sense that $\ba \in \cA_{ad}$, and the controlled state process $\bX^{\ba}$ is in $\cX$.
\end{lemma}

\begin{proof}
The idea of the proof is to replace $a_n - \bar a_n$ and $\bar a_n$ in the dynamics~\eqref{fo:Xt_minus_barXt_u}--\eqref{fo:barXt_u} with expressions $a_n - \bar a_n = - K (X_n^\ba - \bar X_n^\ba)$ and $\bar a_n = - L \bar X_n^\ba$, so that 
\begin{align}
    X_{n+1}^\ba - \bar X_{n+1}^\ba 
    & = \big( A - B K \big) (X_n^\ba - \bar{X}_n^\ba) + \varepsilon_{n+1}
    \\
    \bar X_{n+1}^\ba 
    & = \big( {\mathrm A}  + \bar {\mathrm A}  - ({\mathrm B}  + \bar {\mathrm B} ) L \big) \bar X_n^\ba + \varepsilon_{n+1}^0
\end{align}
for every $n \geq 0$. The initial conditions are $X_0^\ba - \bar X_0^\ba = \varepsilon_0 - \EE[ \varepsilon_0 ]$, $\bar X_0^\ba = \varepsilon_0^0 + \EE[ \varepsilon_0 ] $.  Now, we choose $(Y_n, \beta_n) = (X_n^\ba - \bar X_n^\ba, \varepsilon_n)$ for every $n \geq 1$, $\beta_0 = 0$, and $Y_0 = \varepsilon_0 - \EE[ \varepsilon_0]$ in Lemma~\ref{le:L2_discounted_integrable_of_state_process}. Because the process $(\varepsilon_{n+1})_{n \geq 0}$ is i.i.d. with square-integrable terms, so the process $(\beta_n)_{n \geq 0} \in \cX$ and the process $(Y_n)_{n \geq 0} = (X_n^\ba - \bar X_n^\ba)_{n \geq 0}$ is adapted to the filtration $(\sigma\{\underline\epsilon_n,\underline\epsilon^0_n\})_{n \geq 0}$. 
The admissibility condition of the control parameter $K$ giving $\gamma \| {\mathrm A}  - {\mathrm B}  K \|^2 < 1$, Lemma~\ref{le:L2_discounted_integrable_of_state_process} implies that $(X_n^\ba - \bar X_n^{\ba})_{n \geq 0} \in \cX$. Similar arguments yield that $(\bar X_n^\ba) \in \cX$. Hence, 
$$
    \sum_{n \geq 0} \gamma^n \EE[ \| a_n \|^2 ] \leq 2 \| K \|^2 \sum_{n \geq 0} \gamma^n \EE[ \| X_n^\ba - \bar X_n^\ba \|^2 ] + 2 \| L \|^2 \sum_{n \geq 0} \gamma^n \EE[ \| \bar X_n^\ba \|^2 ] < \infty. 
$$
We conclude that $\ba = (a_n)_{n \geq 0} \in \cA_{ad}$, and that the corresponding controlled state process $\bX^\ba = (X_n^\ba)_{n \geq 0} \in \cX$.
\qed\end{proof}

\begin{definition}
For an admissible control parameter $\theta = (K,L) \in \Theta$, we define its corresponding policy function by:
$
\pi_{\theta}: (x,\bar{x}) \mapsto -K (x - \bar{x}) - L \bar{x}.
$
We say that an admissible control process $\ba = (a_n)_{n \geq 0} \in \cA_{ad}$ is \textbf{parameterized} by $\theta$ if, for every $n \geq 0$, $a_n = \pi_{\theta}(X_n^\ba, \bar{X}_n^\ba)$.
\end{definition}
We shall from time to time use the notation $\ba^\theta$ and $\bX^\theta$ to stress the dependence on $\theta$.  
The collection of \textbf{(admissible) parameterized control processes} is denoted by $\cA^\Theta_{ad} = \{ \ba^\theta \, | \, \theta \in \Theta \} \subset \cA_{ad}$. \index[not]{$\cA^\Theta_{ad}$}\index[sub]{parameterized control process}

\begin{remark}
It is worth pointing out that our definition of the set of admissible controls $\cA_{ad}$ is slightly different from the one considered in~\cite{fazel2018global} for an LQ model without mean-field interactions. There, the authors consider control processes with a finite cumulative cost under some unspecified assumption on the model coefficients $({\mathrm A}, {\mathrm B} )$, and the assumption that $Q, R \succ 0$. Indeed, Remark~\ref{re:equivalence_of_welldefineness_and_admissibility} highlights that under our current assumptions, we have $\ba \in \cA_{ad} \Leftrightarrow J(\ba) < \infty$. In addition, these authors focus on control coefficient $K \in \RR^{\ell \times d}$ such that the \emph{spectral radius} of ${\mathrm A}  - {\mathrm B}  K$ satisfies $\rho({\mathrm A}  - {\mathrm B}  K) < 1$. The set $\{ K \,| \, \rho({\mathrm A}  - {\mathrm B}  K) < 1\}$ is not convex, and this notion of admissibility is different from the one based on the set $\Theta$ defined here. Given that $\rho({\mathrm A}  - {\mathrm B}  K) \leq \| {\mathrm A}  - {\mathrm B}  K \|$ for any matrix $K \in \RR^{\ell \times d}$, our definition of admissibility is more restrictive than the one considered in~\cite{fazel2018global}. As a consequence, a proper analysis of the admissibility of the optimal control is necessary. We will provide more details in Section~\ref{se:admissibility_of_opt_control}.
\end{remark}

\section{Analysis of the MFC Optimization Problem }
\label{se:analysis}

In this section, we provide an analysis of the MFC problem. To simplify the exposition we use the following notations: 
$$
\tilde{\mathrm A}  = {\mathrm A}  + \bar {\mathrm A},\;\;
\tilde {\mathrm B}  = {\mathrm B}  + \bar {\mathrm B}, \;\;
\tilde Q = Q + \bar Q, 
\;\;\text{and}\;\;
\tilde R = R + \bar R.
$$

To start, we construct $(P, \bar P)$ based on the solutions $(P^{*, y}, P^{*, z})$ of another pair of discrete-time algebraic Riccati equations (DAREs):\index[not]{DARE}\index[sub]{discrete-time algebraic Riccati equations}
\begin{equation}
    \label{fo:DARE_closeloop_feedback_excat_PG}
        \left\{
        \begin{array}{rcl}
            P^{*, y} &=& Q  + \gamma {\mathrm A} ^\top P^{*, y} {\mathrm A}  - \gamma^2 {\mathrm A} ^\top P^{*, y} {\mathrm B}  ( R + \gamma {\mathrm B} ^\top P^{*, y} {\mathrm B}  )^{-1} {\mathrm B} ^\top P^{*, y} {\mathrm A} ,
            \\
            P^{*, z} &=& \tilde{Q}  + \gamma \tilde{{\mathrm A} }^\top P^{*, z} \tilde{{\mathrm A} } - \gamma^2 \tilde{{\mathrm A} }^\top P^{*, z} \tilde{{\mathrm B} } ( \tilde{R} + \gamma \tilde{{\mathrm B} }^\top P^{*, z} \tilde{{\mathrm B} } )^{-1} \tilde{{\mathrm B} }^\top P^{*, z} \tilde{{\mathrm A} }.
        \end{array}
        \right.
\end{equation}
    For symmetric positive semidefinite solutions $P^{*,y},P^{*,z}$ of these DAREs, define the matrices $P, \bar P \in \RR^{d \times d}$ by
    \begin{equation}
    \label{fo:expression_of_solution_openloop_DARE}
        \left\{
            \begin{array}{ll}
            P &:= 2 \gamma P^{*, y} \big( I_d + \gamma {\mathrm B}  R^{-1} {\mathrm B} ^\top P^{*, y} \big)^{-1} {\mathrm A}   
            \\
            \bar P &:= 2 \gamma  P^{*, z} \big( I_d + \gamma \tilde{{\mathrm B} } \tilde{R}^{-1}  \tilde{{\mathrm B} }^\top  P^{*, z} \big)^{-1} \tilde{{\mathrm A} },
            \end{array}
        \right.
    \end{equation}

Define
$K^* =-\Gamma P$ and $L^* = -\Lambda \bar P$ where
\begin{equation}
\label{fo:Gamma_Lambda}
    \Gamma = - \frac{1}{2} R^{-1} {\mathrm B} ^\top, \qquad \Lambda = - \frac{1}{2} \tilde{R}^{-1} \tilde{\mathrm B}^\top,
\end{equation}
and $P$ and $\bar P$ satisfy the two matrix Riccati equations
\begin{equation}
\label{fo:Riccati}
\left\{
    \begin{array}{rcl}
        P & = & \gamma \big( {\mathrm A} ^\top P + 2Q \big) \big( {\mathrm A}  -  {\mathrm B}  R^{-1} {\mathrm B} ^\top P / 2 \big) 
        \\
        \bar P & = & \gamma \big( \tilde{{\mathrm A} }^\top \bar P + 2 \tilde Q \big) \big( \tilde{\mathrm A}  -  \tilde{\mathrm B}  \tilde{R}^{-1} \tilde{{\mathrm B} }^\top \bar P / 2 \big).
    \end{array}
\right.
\end{equation} %
\begin{assumption}[Optimal feedback admissibility in $\Theta$]
\label{ass:optimal-feedback-in-Theta}
The feedback matrices defined above satisfy $\gamma\norm{{\mathrm A}-{\mathrm B}K^*}^2<1$ and $\gamma\norm{\tilde{\mathrm A}-\tilde{\mathrm B}L^*}^2<1$, equivalently $\theta^*=(K^*,L^*)\in\Theta$. 
\end{assumption}
These norm conditions are stronger than Schur stability after discounting.
We shall argue later on the existence of $P$ and $\bar P$.
The main result of this section is the following.

\begin{theorem}
\label{th:existence_linear_control}
Under Assumptions \ref{ass:finit_cost},~\ref{ass:positivity-qr}, and~\ref{ass:optimal-feedback-in-Theta}, there exists a unique optimal control for the MFC problem~\eqref{pb:MFC_L2_admissible_control}. The optimal control $\ba = (a_n)_{n \geq 0}$ is in feedback form, and is linear in $(\bX^\ba, \bar{\bX}^\ba)$, specifically:
\begin{equation}
\label{fo:optimal_control_linear_in_x_barx}
    a_n = - K^* (X_n^\ba - \bar X_n^\ba) - L^* \bar X_n^\ba, \qquad n\ge 0.
\end{equation}
\end{theorem}

\begin{remark}
Within the proof of this result, we will show that the optimal control is in fact of the form
\begin{equation}
\label{fo:control_ut_expression_2}
    a_n = \Gamma ( p_n - \bar{p}_n) + \Lambda \bar{p}_n, \qquad \bar a_n = \Lambda \bar{p}_n,\qquad n \geq 0,
\end{equation}
where $(p_n)_{n\ge 0}$ is the adjoint process corresponding to $(X_n,a_n)_{n \geq 0}$ (see Definition~\ref{def:adjoint_process} below).
\end{remark}

\subsection{Admissibility of the optimal control}
\label{se:admissibility_of_opt_control}

Before proving optimality, we record the algebraic identities relating the Riccati construction to the feedback matrices in Theorem~\ref{th:existence_linear_control}. The membership $\theta^*\in\Theta$ is imposed in Assumption~\ref{ass:optimal-feedback-in-Theta}.

\begin{lemma}
\label{le:solution_to_ARE_openloop_from_closefeedback_DAREs}
    We assume that Assumption~\ref{ass:finit_cost} and Assumption~\ref{ass:positivity-qr} hold. Then, there exist symmetric and positive semi-definite matrices $P^{*, y}$ and $P^{*, z}$ in $\RR^{d \times d}$ satisfying \eqref{fo:DARE_closeloop_feedback_excat_PG}.
    For any such solutions, the matrices $P$ and $\bar P$ defined by~\eqref{fo:expression_of_solution_openloop_DARE} satisfy the matrix Riccati equations~\eqref{fo:Riccati}.
\end{lemma}

\begin{proof}
    Under Assumption~\ref{ass:finit_cost}, we know that the two pairs of matrices $(\sqrt{\gamma} {\mathrm A} , \sqrt{\gamma} {\mathrm B} )$ and $(\sqrt{\gamma} \tilde {\mathrm A} , \sqrt{\gamma} \tilde {\mathrm B} )$ are \emph{stabilizable} in the sense that there exist $K, L \in \RR^{\ell \times d}$ such that the matrices $\sqrt{\gamma} ({\mathrm A}  - {\mathrm B}  K)$ and $\sqrt{\gamma} (\tilde {\mathrm A}  - \tilde {\mathrm B}  L)$ are Schur stable \index[sub]{Schur stable} in the sense that all their eigenvalues fall strictly inside the unit disk. See for example \cite[Section 6.2]{Kailath}. 
    This is the stabilizability condition needed for the DARE existence result below.
    Under Assumption~\ref{ass:positivity-qr}, we have $Q, \tilde Q \succeq 0$ and $R, \tilde R \succ 0$, so by~\cite[Theorem 3.2 (iii)]{ran1988existence} we know that the (DAREs)~\eqref{fo:DARE_closeloop_feedback_excat_PG} admit solutions, and the solutions $P^{*, y}$ and $P^{*, z}$ are also symmetric and positive semi-definite.

    Now, to show that $P, \bar P$ given by~\eqref{fo:expression_of_solution_openloop_DARE} are solutions to~\eqref{fo:Riccati}, we first recall the following matrix identity: for any matrix $W \in \RR^{d \times d}$ such that $I_d + W$ is invertible,
    $$
        (I_d + W)^{-1} = I_d - W (I_d + W)^{-1} = I_d - W + W (I_d + W)^{-1} W.
    $$
    Let $W = \gamma {\mathrm B} R^{-1} {\mathrm B}^\top  P^{*, y}$, $U = {\mathrm B}R^{-1/2}$, and $V = (R^{-1/2})^\top {\mathrm B} ^\top P^{*, y}$. Because $P^{*, y} \succeq 0$ and $R \succ 0$, we have the minimum eigenvalue $\lambda_{min}(W) = \gamma \lambda_{min}(UV) = \gamma \lambda_{min}(VU) \geq 0$. So $I_d + W$ is invertible. Because $R + \gamma {\mathrm B} ^\top P^{*, y} {\mathrm B}  \succ 0$ and 
    \begin{align*}
       I_\ell & = I_\ell + \gamma {\mathrm B} ^\top P^{*, y} \big( I_d - (I_d + W)^{-1} - W ( I_d + W )^{-1} \big) {\mathrm B}  R^{-1}
        \\
        & =  (R + \gamma {\mathrm B} ^\top P^{*, y} {\mathrm B} ). \big( R^{-1} - \gamma R^{-1} {\mathrm B} ^\top P^{*,y} \big( I_d + \gamma {\mathrm B}  R^{-1} {\mathrm B} ^\top P^{*, y} \big)^{-1} {\mathrm B}  R^{-1} \big),
    \end{align*}
    we have
    \begin{align*}
        & I_d - \gamma {\mathrm B}  \big( R + \gamma {\mathrm B} ^\top P^{*, y} {\mathrm B}  \big)^{-1} {\mathrm B} ^\top P^{*, y}
        \\
        &= 
        I_d - \gamma {\mathrm B} . \big( R^{-1} - \gamma R^{-1} {\mathrm B} ^\top P^{*,y} \big( I_d + \gamma {\mathrm B}  R^{-1} {\mathrm B} ^\top P^{*, y} \big)^{-1} {\mathrm B}  R^{-1} \big). {\mathrm B} ^\top P^{*,y}
        \\
        &=
        I_d - \gamma {\mathrm B}  R^{-1} {\mathrm B} ^\top P^{*, y} + ( \gamma {\mathrm B}  R^{-1} {\mathrm B} ^\top P^{*, y}). \big( I_d +  \gamma {\mathrm B}  R^{-1} {\mathrm B} ^\top P^{*, y} \big)^{-1}. ( \gamma {\mathrm B}  R^{-1} {\mathrm B} ^\top P^{*, y})
        \\
        &=
        I_d - W + W ( I_d + W)^{-1} W
        \\
        &=
        (I_d + W)^{-1}
        \\
        &= \big( I_d +  \gamma {\mathrm B}  R^{-1} {\mathrm B} ^\top P^{*, y} \big)^{-1}.
    \end{align*}
    Consequently, from the first equation in (DAREs)~\eqref{fo:DARE_closeloop_feedback_excat_PG}, we obtain
    \begin{align*}
         P^{*, y}  &=  Q + \gamma {\mathrm A} ^\top P^{*, y}  \Big(  I_d - \gamma {\mathrm B}  \big( R + \gamma {\mathrm B} ^\top P^{*, y} {\mathrm B}  \big)^{-1} {\mathrm B} ^\top P^{*, y} \Big) {\mathrm A}  
         \\
         &= Q + \gamma {\mathrm A} ^\top P^{*,y} \big( I_d + \gamma {\mathrm B}  R^{-1} {\mathrm B} ^\top P^{*,y} \big)^{-1} {\mathrm A} 
         \\
         &= Q + {\mathrm A} ^\top P / 2.
    \end{align*}
    Moreover, we also have
    \begin{align*}
        {\mathrm A}  - {\mathrm B}  R^{-1} {\mathrm B} ^\top P / 2 
        & = \big[ I_d - \gamma {\mathrm B}  R^{-1} {\mathrm B} ^\top P^{*, y} \big( I_d + \gamma {\mathrm B}  R^{-1} {\mathrm B} ^{\top} P^{*, y} \big)^{-1}  \big] {\mathrm A}  
        \\
        & = \big(I_d + \gamma {\mathrm B}  R^{-1} {\mathrm B} ^\top P^{*, y} \big)^{-1} {\mathrm A} .
    \end{align*}
    Thus, we derive the matrix Riccati equation for $P$ by
    \begin{align*}
        P = \gamma (2 P^{*,y}). \big(I_d + \gamma {\mathrm B}  R^{-1} {\mathrm B} ^\top P^{*, y} \big)^{-1} {\mathrm A}  =  \gamma ( 2 Q + {\mathrm A} ^\top P)( {\mathrm A}  - {\mathrm B}  R^{-1} {\mathrm B}^\top  P / 2).
    \end{align*}
    We apply similar arguments with $(\tilde {\mathrm A} , \tilde {\mathrm B} , \tilde R, P^{*, z}, \bar P)$ to derive the matrix Riccati equation for $\bar P = \gamma ( 2 \tilde{Q} + \tilde{{\mathrm A} }^\top \bar P)( \tilde{{\mathrm A} } - \tilde{{\mathrm B} } \tilde{R}^{-1} \tilde{{\mathrm B} }^\top \bar P / 2)$.
\qed\end{proof}

\begin{corollary}
\label{cor:admissibility_of_theta^*}
Under Assumptions \ref{ass:finit_cost},~\ref{ass:positivity-qr}, and~\ref{ass:optimal-feedback-in-Theta},
if we define the matrices $K^{*}$ and $L^{*}$ in $\RR^{\ell \times d}$ by:
$$
    K^{*} = R^{-1} {\mathrm B} ^\top P/ 2, \qquad L^{*} = \tilde R^{-1} \tilde {\mathrm B} ^\top \bar P / 2,
$$
as in Theorem~\ref{th:existence_linear_control} where $(P, \bar P)$ are matrices defined by~\eqref{fo:expression_of_solution_openloop_DARE} based on the solutions matrices $(P^{*, y}, P^{*, z})$ to (DAREs)~\eqref{fo:DARE_closeloop_feedback_excat_PG}, %
then the following identities hold:
\begin{equation*}
\left\{
\begin{array}{rl}
    \sqrt{\gamma} ( {\mathrm A}  - {\mathrm B}  K^* ) &= \sqrt{\gamma} \big(I_d + \gamma {\mathrm B}  R^{-1} {\mathrm B} ^\top P^{*, y} \big)^{-1} {\mathrm A} ,
    \\
    \sqrt{\gamma} ( \tilde {\mathrm A}  - \tilde {\mathrm B}  L^* ) &=  \sqrt{\gamma} \big(I_d + \gamma \tilde{{\mathrm B} } \tilde{R}^{-1} \tilde{{\mathrm B} }^\top P^{*, z} \big)^{-1} \tilde{{\mathrm A} }.
\end{array}
\right.
\end{equation*}
\end{corollary}

\begin{proof}
    The definition of $(P, \bar P)$ in~\eqref{fo:expression_of_solution_openloop_DARE} implies
   
    \begin{align*}
        \sqrt{\gamma} ( {\mathrm A}  - {\mathrm B}  K^{* } ) & = \sqrt{\gamma} \big( {\mathrm A}  - {\mathrm B}  R^{-1} {\mathrm B} ^\top P / 2) 
        \\
        & = \sqrt{\gamma} \big[ I_d - \gamma {\mathrm B}  R^{-1} {\mathrm B} ^\top P^{*, y} \big( I_d + \gamma {\mathrm B}  R^{-1} {\mathrm B} ^{\top} P^{*, y} \big)^{-1}  \big] {\mathrm A} 
        \\
        & = \sqrt{\gamma} \big(I_d + \gamma {\mathrm B}  R^{-1} {\mathrm B} ^\top P^{*, y} \big)^{-1} {\mathrm A} .
    \end{align*}
	This proves the first identity; the second follows by the same computation with $(\tilde{\mathrm A},\tilde{\mathrm B},\tilde R,P^{*,z},\bar P,L^*)$.
\qed\end{proof}

\begin{proposition}
\label{pr:admissibility_of_u_with_feedback_control}
    Under Assumptions~\ref{ass:finit_cost},~\ref{ass:positivity-qr}, and~\ref{ass:optimal-feedback-in-Theta}, let $\theta^* = (K^*, L^*) = (- \Gamma P, -\Lambda \bar P)$ be the control parameter with coefficients $(K^*, L^*, \Gamma, \Lambda, P, \bar P)$ defined in Theorem~\ref{th:existence_linear_control}.
    The corresponding parameterized control process $\ba^{\theta^*}$ in a feedback form of its controlled state process and conditional mean process $\big( \bX^{\ba^{\theta^*}}, \bar{\bX}^{\ba^{\theta^*}} \big)$ with policy function $\pi_{\theta^*} : (x, \bar x)  \mapsto - K^* (x - \bar x) - L^* \bar x$ is admissible in the sense that $\ba^{\theta^*} \in \cA_{ad}$.
\end{proposition}

\begin{proof}
    Assumption~\ref{ass:optimal-feedback-in-Theta} states that the control parameter $\theta^* = (K^*, L^*)$ belongs to $\Theta$.
    Because the process $\ba^{\theta^*}$ is defined in feedback form with policy $\pi_{\theta^*}$ in Theorem~\ref{th:existence_linear_control}, Lemma~\ref{le:admissible_of_control_with_feedback_form_on_theta} implies immediately that $\ba^{\theta^*} \in \cA_{ad}$.
\qed\end{proof}

\subsection{Optimality of the control process}
In this subsection, we show existence and uniqueness of an optimal control using a discrete-time version of Pontryagin's maximal principle. 
For convenience, we use the notation $\check {\mathrm A}  = {\mathrm A}  - I_d$ and $\zeta = (x, \bar x, a, \bar a) \in \RR^{d} \times \RR^{d} \times \RR^{\ell} \times \RR^{\ell}$. We define the \emph{drift function} $b$ by
\begin{equation}
\label{fo:drift_function_b}
    b(\zeta) = b(x, \bar x, a, \bar a) =  \check{{\mathrm A} } x + \bar{{\mathrm A} } \bar x + {\mathrm B}  a + \bar{{\mathrm B} } \bar{a}
\end{equation}
and the Hamiltonian function $h$ by:
\begin{equation}
\label{fo:Hamilton_function_h}
    h(\zeta, p) = h(x, \bar x, a, \bar a, p) = b(\zeta) \cdot p + f(\zeta) - \delta\; x \cdot p
\end{equation}
where $p \in \RR^{d}$, $\delta = (1 - \gamma) / \gamma$, and $f(\zeta)=f(x, \bar x, a, \bar a)$ is the one-stage cost function defined by~\eqref{fo:quadratic_cost}. Then, for an admissible control process $(a_n)_{n \geq 0}$, the state dynamics~\eqref{fo:MKV-state} of the controlled state process $(X_n)_{n \geq 0}$ can be rewritten as
\begin{equation}
\label{fo:delta_X_t_with_b}
    X_{n+1} - X_n = b(X_n, \bar X_n, a_n, \bar a_n) +  \varepsilon_{n+1}^0 + \varepsilon_{n+1} = b(\zeta_n) + \varepsilon_{n+1}^0 + \varepsilon_{n+1}.
\end{equation}
To compute the partial derivatives of $h$ with respect to its arguments $(x, \bar x, a, \bar a, p)$, we treat $\zeta$ as a $(2d + 2\ell)$ dimension vector, $Vect(\zeta) = (x^\top, \bar x^\top, a^\top, \bar a^\top)^\top \in \RR^{(2d + 2 \ell)}$, and by a slight abuse of notation we use the same notation $\zeta$. Now, for every fixed $p \in \RR^{d}$, we have 
\begin{equation}
\label{fo:gradient_h}
\nabla_{\zeta} h (\zeta, p) = 
	\begin{pmatrix}
	\partial_x h(\zeta,p)\\
	\partial_{\bar{x}} h(\zeta,p)\\
	\partial_{a} h(\zeta,p)\\
	\partial_{\bar{a}} h(\zeta,p)
	\end{pmatrix}
	=
	\begin{pmatrix}
	(\check{{\mathrm A} }-\delta I_d)^\top p + 2Q (x-\bar{x}) \\
	\bar{{\mathrm A} }^\top p - 2 Q (x - \bar x) + 2 \tilde Q \bar x
	\\
	{\mathrm B} ^\top p + 2 R (a - \bar{a})\\
	\bar{{\mathrm B} }^\top p - 2 R (a - \bar a) + 2\tilde R \bar a
	\end{pmatrix},
\end{equation}
and the Hessian of $h$ with respect to $(a, \bar a)$ is given by 
\begin{equation}
\label{fo:hession_h_u}
   \nabla^2_{(a,\bar{a}),(a,\bar{a})} h(\zeta,p)=
	\begin{pmatrix}
	2R & -2R \\
	-2R & 2(R+\tilde{R})
	\end{pmatrix},
\end{equation}
which is positive definite because 
\begin{equation}
   \begin{pmatrix} a & \bar{a} \end{pmatrix}
    \begin{pmatrix}
	R & -R \\
	-R & R+\tilde{R}
	\end{pmatrix}
   \begin{pmatrix} a \\ \bar{a} \end{pmatrix}
=(a - \bar a) R(a - \bar a) + \bar a \tilde R \bar a \ge 0
\end{equation}
and equal to $0$ only when $a=\bar a=0$. Indeed, under Assumption~\ref{ass:positivity-qr}, $R$ and $\tilde R$ are assumed to be (strictly) positive definite. 

We now introduce the notion of adjoint process associated to a state process generated by an admissible control process.

\begin{definition}
\label{def:adjoint_process}
Given an admissible control process and the corresponding controlled state process $(\ba, \bX) = (a_n, X_n)_{n \geq 0}$, we say that an $\RR^d$-valued $(\sigma\{\underline\epsilon_n,\underline\epsilon^0_n\})_{n \geq 0}$-adapted process $\bp = (p_n)_{n \geq 0}$ is a corresponding adjoint process if it satisfies the backward equation:
\begin{equation}
\label{fo:adjoint}
    p_n = \EE \Big[ p_{n+1} +  \gamma \big[ (\check {\mathrm A}  - \delta I_d)^\top p_{n+1} + 2 Q X_{n+1} + \bar {\mathrm A} ^\top \bar p_{n+1} + 2 \bar Q \bar X_{n+1} \big] \big\vert \sigma\{\underline\epsilon_n,\underline\epsilon_n^0\} \Big],
\end{equation}
and the transversality condition:\index[sub]{transversality condition}
\begin{equation}
    \label{fo:transversality}
    \EE \Bigl[ \sum_{n \geq 0} \gamma^n \| p_n \|^2 \Bigr] < \infty.
\end{equation}
\end{definition}
The transversality condition should be viewed as a \emph{terminal condition} controlling the possible growth of $p_n$ when $n\to\infty$. Notice that equation~\eqref{fo:adjoint} can equivalently be written as:
\begin{equation}
   \label{fo:simpler_adjoint}
    p_n = \gamma \EE [ {\mathrm A} ^\top p_{n+1} + 2 Q X_{n+1} + \bar {\mathrm A} ^\top \bar p_{n+1} + 2 \bar Q \bar X_{n+1} \ | \  \sigma\{\underline\epsilon_n,\underline\epsilon_n^0\} ]. 
\end{equation}

\begin{proposition}
\label{pr:existence_and_uniqueness_of_adjoint_process}
    Under Assumption~\ref{ass:finit_cost}, for every state process generated by an admissible control process, there exists a unique adjoint process. 
\end{proposition}

\begin{proof}
Let $\ba \in \cA_{ad}$ be an admissible control process and $\bX = (X_n)_{n\geq 0}$ the corresponding controlled state process.

\textbf{(Uniqueness)} Consider two adjoint processes $(p_n)_{n \geq 0}$ and $(p'_n)_{n \geq 0}$ corresponding to $(\ba, \bX)$. Taking conditional expectation $\EE[\cdot | \sigma\{\underline\epsilon^0_n\} ]$ on both sides of equation~\eqref{fo:simpler_adjoint} and subtracting \eqref{fo:simpler_adjoint} for $\bar{\bp}$ to the corresponding equation for $\bar{\boldsymbol{p'}}$, we get 
$
    \EE[ \| \bar p_n - \bar p_n' \| ] \leq \gamma^{1/2} \EE[ \| \bar p_{n+1} - \bar p_{n+1}' \| ]
$ 
because of $\sqrt{\gamma} \| {\mathrm A}  + \bar {\mathrm A}  \| < 1$. This implies that for every $0 \leq n \leq m$,
$
    \gamma^{n/2} \EE [ \| \bar p_n - \bar p_n' \| ] \leq \gamma^{m/2} \EE [ \|\bar p_m -\bar p_m' \| ].
$
From the transversality condition and Jensen's inequality for conditional expectations, we have
$$
\lim_{m \to \infty} \gamma^{m/2} \EE \big[ \| \bar p_m - \bar p_m' \| \big]  \leq \lim_{m \to \infty} \left[(\EE[ \gamma^m \| p_m \|^2 ] )^{1/2} + (\EE[ \gamma^m \| p_m' \|^2 ] )^{1/2} \right]
 = 0.
$$
Hence $\bar p_n = \bar p_n'$, $\PP$-a.s., for every $n \geq 0$. Using similar arguments and the fact that $\gamma^{1/2} \| {\mathrm A}  \| < 1$,  we obtain $p_n - \bar p_n = p'_n - \bar p'_n$, $\PP$-a.s., for every $n \geq 0$ proving uniqueness.

\textbf{(Existence)} We construct directly an adjoint process $\bp = (p_n)_{n \geq 0}$ with the following two auxiliary processes $\by = (y_n)_{n \geq 0}$ and $\bz=(z_n)_{n \geq 0}$ defined for every $n \geq 0$, by: 
\begin{equation}
\left\{
\begin{array}{rcl}
    y_n & := & \sum_{k \geq 0}  \bigl( \gamma {\mathrm A} ^\top \bigr)^k \bigl( 2 \gamma Q \bigr)  \EE \bigl[ X_{n+1+k} - \bar{X}_{n+1+k}  \, | \, \sigma\{\underline\epsilon_n,\underline\epsilon_n^0\} \bigr]
    \\
    z_n & := & \sum_{k \geq 0} \bigl( \gamma \tilde {\mathrm A} ^\top \bigr)^{k} \bigl( 2 \gamma \tilde Q \bigr) \EE[ \bar{X}_{n+1 + k} \, | \, \sigma\{\underline\epsilon_n,\underline\epsilon_n^0\}]
    \\
    p_n & := & y_n + z_n.
\end{array}
\right.
\end{equation}
By definition, the process $\bp$ is $(\sigma\{\underline\epsilon_n,\underline\epsilon^0_n\})_{n \geq 0}$-adapted and satisfies the backward adjoint equation~\eqref{fo:adjoint}.
Next, we show that the process $\bz$ (resp. $\by$) is $L^2-$discounted integrable. Under Assumption~\ref{ass:finit_cost}, we consider a factor $\rho \in (0, 1)$ such that $\eta := \gamma / \rho < 1$ and $\xi := \rho \| \tilde {\mathrm A}  \|^2 < 1$, so $\eta \xi = \gamma \| \tilde {\mathrm A}  \|^2$. By Fatou's lemma and Cauchy-Schwarz inequality, for every $n \geq 0$, we bound $z_n$ by 
$$
\gamma^n \EE[ \| z_n \|^2 ] \leq \frac{4\gamma \| \tilde Q \|^2}{1 - \eta} \EE\big[ \sum_{k \geq 0} \xi^{k} (\gamma^{n+1+k} \| \bar{X}_{n+1+k} \|^2) \big].
$$
Summing over all times $n \geq 0$ and exchanging the two infinite summations involved, we get  
$$ 
\EE \big[ \sum \nolimits_{n \geq 0} \gamma^n \| z_n \|^2 \big] \leq \frac{ 4\gamma \| \tilde Q \|^2}{1 - \eta} \frac{1}{1 - \xi} \EE \big[ \sum \nolimits_{n \geq 1} \gamma^n \| \bar X_n \|^2 \big] < \infty.
$$
We then conclude that the process $\bp$ satisfies the transversality condition~\eqref{fo:transversality} and is an adjoint process corresponding to $(\ba, \bX)$.
\qed\end{proof}

Using the adjoint process, we compute the Gateaux derivative of $J$, and derive a necessary condition for optimality in the spirit of the classical Pontryagin's maximum principle.

\begin{proposition}
\label{pr:necessary_condition_of_PMP}
    Under Assumption~\ref{ass:finit_cost}, the Gateaux derivative of $J$ at $\ba$ in the direction $\bbeta \in \cA_{ad}$ exists and is given by
    \begin{equation}
    \label{fo:Gateau_derivative_DJ_u_beta}
        DJ(\ba)(\bbeta) = \EE \Bigg[ \sum_{n \geq 0} \gamma^n \big(  p_n^\top {\mathrm B}  + 2 a_n^\top R + \bar{p}_n^{\top} \bar {\mathrm B}  + 2 \bar{a}_n^\top \bar R \big) \beta_n \Bigg],
    \end{equation}
    where $(p_n)_{n \geq 0}$ is the adjoint process corresponding to the controlled state process $\bX^\ba = (X_n^{\ba})_{n \geq 0}$.
    Moreover, if the control process $\ba$ is optimal, namely $J(\ba) \leq J(\ba')$ for all $\ba' \in \cA_{ad}$, then we have for all $n \geq 0$,
\begin{equation}
    \label{fo:necessary_PMP}
        {\mathrm B} ^\top p_n + 2 R a_n + \bar {\mathrm B} ^\top \bar{p}_n + 2 \bar R \bar{a}_n = 0, \qquad \PP-a.s. \, .
\end{equation}
\end{proposition}

\begin{proof}
    For an admissible control process $\ba \in \cA_{ad}$ and a perturbation $\bbeta \in \cA_{ad}$, the control $\ba + \lambda \bbeta = (a_n + \lambda \beta_n)_{n \geq 0}$ is still admissible for any $\lambda \in (0,1)$ by the convexity of $\cA_{ad}$, and we denote by $\bX^{\ba + \lambda \bbeta} = (X_n^{\ba + \lambda \bbeta})_{n \geq 0}$ the corresponding controlled state starting from the same initial state $X_0^{\ba} = X_0^{\ba + \lambda \bbeta}$, and subject to the same noise sequences $(\varepsilon^0_{n+1})_{n \geq 0}$ and $(\varepsilon_{n+1})_{n \geq 0}$. Let us denote by $\bM^{\lambda, \bbeta} = (M_n^{\lambda, \bbeta})_{n \geq 0}$ the variation process defined as
    \begin{equation}
    \label{fo:Mt_lambda}
        M_n^{\lambda, \bbeta} := \frac{1}{\lambda} (X_n^{\ba + \lambda \bbeta} -  X_n^{\ba}). 
    \end{equation}
    The dynamics of $\bM^{\lambda, \bbeta}$ are deterministic and satisfy the following equation:
    \begin{equation}
    \label{fo:dyn_Mt_lambda_beta}
        M_{n+1}^{\lambda, \bbeta} = {\mathrm A}  M_n^{\lambda, \bbeta} + \bar {\mathrm A}  \bar M_n^{\lambda, \bbeta} + {\mathrm B}  \beta_n + \bar {\mathrm B}  \bar \beta_n,
    \end{equation}
    where $M_0^{\lambda, \bbeta} =0$, and $\bar M_n^{\lambda, \bbeta} = \EE[ M_{n}^{\lambda, \bbeta}  | \sigma\{\underline\epsilon^0_n\} ] $ and $\bar{\beta}_n = \EE[ \beta_n | \sigma\{\underline\epsilon^0_n\} ] $ are conditional expectations with respect to $\sigma\{\underline\epsilon^0_n\}$, the past of the common noise. Since $\bbeta$ is $(\sigma\{\underline\epsilon_n,\underline\epsilon_n^0\} )_{n \geq 0}$-adapted, equation~\eqref{fo:dyn_Mt_lambda_beta} implies directly that $M_{n+1}^{\lambda, \bbeta}$ is $\sigma\{\underline\epsilon_n,\underline\epsilon_n^0\}$-measurable for every $n \geq 0$. We also notice that the process $\bM^{\lambda, \bbeta}$ is independent of $\lambda$, so we shall denote it by $\bM^{\bbeta} = (M_n^{\bbeta})_{n \geq 0}$ from now on.
    
   Let $(p_n)_{n \geq 0}$ be  the adjoint process corresponding to state process $\bX^{\ba}$. From the admissibility of $\ba$ and $\ba + \lambda \bbeta$, we infer that both state processes $\bX^{\ba}$ and $\bX^{\ba + \lambda \bbeta}$ are $L^2$-integrable, and together with the transversality condition of the adjoint process $(p_n)_{n \geq 0}$, we deduce that
    $$
        \sum_{n \geq 0} \EE\Big[ \gamma^n \big( X_n^{\ba + \lambda \bbeta} - X_n^{\ba} \big)  \cdot p_n \Big] < \infty,
    $$ 
    and
    $$
        \sum_{n \geq 0} \EE \Big[ \gamma^n \Big(  b(X_n^{\ba+ \lambda \bbeta}, \bar{X}_n^{\ba+ \lambda \bbeta}, a_n + \lambda \beta_n, \bar a_n + \lambda \bar{\beta}_n) - b(X_n, \bar X_n, a_n, \bar a_n) \Big) \cdot p_n \Big] < \infty.
    $$
    Next, we compute the Gateaux derivative of $J$. Let 
    $
        \zeta^{\lambda}_n = (X_n^{\ba+ \lambda \bbeta}, \bar{X}_n^{\ba+ \lambda \bbeta}, a_n + \lambda \beta_n, \bar a_n + \lambda \bar{\beta}_n)
    $ 
    and 
    $
        \zeta_n = (X_n^{\ba}, \bar X_n^{\ba}, a_n, \bar a_n)
    $, together with~\eqref{fo:Hamilton_function_h}, \eqref{fo:delta_X_t_with_b}, and $X_0^{\ba + \lambda \bbeta} = X_0^{\ba}$, we have
    
    \begin{equation}
    \label{fo:difference_between_cost_function_in_hamilton_func}
    \begin{split}
        & \EE \Big[ \sum_{n \geq 0} \gamma^n \big( f(\zeta_n^\lambda) - f(\zeta_n) \big) \Big] \\
        &\hskip 25pt
        =\sum_{n \geq 0} \EE \Big[ \gamma^n \big( h(\zeta_n^\lambda, p_n) - h(\zeta_n, p_n) \big) \Big] - \sum_{n \geq 0} \EE \Big[ \gamma^n [ ( X_{n+1}^{\ba + \lambda \bbeta} - X_n^{\ba + \lambda \bbeta}  ) - ( X_{n+1}^{\ba} - X_n^{\ba} ) ] \cdot p_n \\
        &\hskip 195pt + \delta \gamma^n ( X_n^{\ba + \lambda \bbeta} - X_{n}^{\ba} ) \cdot p_{n} \Big] \\
        &\hskip 25pt
        =\sum_{n \geq 0} \EE \Big[ \gamma^n \big( h(\zeta_n^\lambda, p_n) - h(\zeta_n, p_n) \big) \Big] + \sum_{n \geq 0} \EE \Big[ \gamma^n \big( X_{n+1}^{\ba + \lambda \bbeta} - X_{n+1}^{\ba}) \cdot (p_{n+1} - p_n)  \big) \Big].
    \end{split}
    \end{equation}
    Using the formulas for the partial derivatives of the Hamilton function $\nabla_\zeta h(\zeta, p)$, and the process $(M_n^{\bbeta})_{n \geq 0}$ with $M_n^{\beta} = \lim_{\lambda \to 0} (X_n^{\ba + \lambda \bbeta} - X_n^\ba) / \lambda$, we then have 
    \begin{align}
        DJ(\ba)(\bbeta) 
        &= \lim_{\lambda \to 0} \frac{1}{\lambda}\big[ J(\ba + \lambda \bbeta) - J(\ba) \big] \nonumber \\
        &= \sum_{n \geq 0} \gamma^n \EE \Bigl[ \Big( \partial_x h(\zeta_n, p_n) \cdot M_n^{\bbeta} + \partial_{\bar x} h(\zeta_n, p_n) \cdot \bar{M}_n^{\bbeta} + \partial_a h(\zeta_n, p_n) \cdot \beta_n + \partial_{\bar a}h(\zeta_n, p_n) \cdot \bar{\beta}_n  \Big)
        \nonumber \\
        & \hskip 40pt + \Big( M_{n+1}^{\bbeta} \cdot ( p_{n+1} - p_n) \Big) \Bigr].  
        \label{fo:Gateaux}
    \end{align}
Properties of the conditional expectations give: 
\begin{equation*}
\begin{split}
    \EE\Bigl[ \partial_{\bar x} h(\zeta_n, p_n) \cdot \bar{M}_n^{\bbeta} \Bigr] 
    & = \EE \Bigl[ \big( \bar {\mathrm A} ^\top p_n - 2 Q (X_n^{\ba} - \bar X_n^{\ba}) + 2 \tilde{Q} \bar X_n^{\ba} \big) \cdot \EE[M_n^{\bbeta}|\,\sigma\{\underline\epsilon^0_n\} ] \Bigr] \\
    & = \EE\Bigl[ \EE\bigl[ \bar {\mathrm A} ^\top p_n - 2 Q (X_n^{\ba} - \bar X_n^{\ba}) + 2 \tilde{Q} \bar X_n^{\ba}|\,\sigma\{\underline\epsilon^0_n\} \bigr] \cdot \EE[M_n^{\bbeta}|\,\sigma\{\underline\epsilon^0_n\} ] \Bigr]
        \\
        & = \EE \Bigl[ \big( \bar {\mathrm A} ^\top \bar{p}_n + 2 \tilde{Q} \bar{X}_n^{\ba} \big) \cdot \EE[M_n^{\bbeta}|\,\sigma\{\underline\epsilon^0_n\} ] \Bigr] 
        \\
        & = \EE \Bigl[  \big( \bar {\mathrm A} ^\top \bar{p}_n + 2 \tilde{Q} \bar{X}_n^{\ba} \big) \cdot M_n^{\bbeta}  \Bigr].
\end{split}
\end{equation*}
Since $M_{n+1}^{\bbeta}$ is $\sigma\{\underline\epsilon_n,\underline\epsilon_n^0\}$-measurable, the backward equation~\eqref{fo:simpler_adjoint} for $(p_n)_{n\geq 0}$ implies that
\begin{equation*}
\begin{split}
        & \sum_{n \geq 0} \gamma^n\EE \Bigl[ M_{n+1}^{\bbeta}  \cdot ( p_{n+1} - p_{n} ) \Bigr]
        \\
         & \hskip 25pt
        =\sum_{n \geq 0} \gamma^n\EE \Bigl[ (M_{n+1}^{\bbeta})^\top (-\gamma) \EE \Big[  (\check {\mathrm A}  - \delta I_d)^\top p_{n+1} + 2 Q X_{n+1}^{\ba} + \bar {\mathrm A} ^\top \bar{p}_{n+1} + 2 \bar Q \bar{X}_{n+1}^{\ba} \Big| \sigma\{\underline\epsilon_n,\underline\epsilon_n^0\} \Big] \Bigr]
        \\
        &\hskip 25pt
        =-\sum_{n \geq 0} \gamma^{n+1} \EE \Bigl[  \big(  (\check {\mathrm A}  - \delta I_d)^\top p_{n+1} + 2 Q X_{n+1}^{\ba} + \bar {\mathrm A} ^\top \bar{p}_{n+1} + 2 \bar Q \bar{X}_{n+1}^{\ba} \big) \cdot  M_{n+1}^{\bbeta} \Bigr]
        \\
        &\hskip 25pt
        =-\sum_{n \geq 0} \gamma^{n} \EE \Bigl[  \big(  (\check {\mathrm A}  - \delta I_d)^\top p_{n} + 2 Q X_{n}^{\ba} + \bar {\mathrm A} ^\top \bar{p}_{n} + 2 \bar Q \bar{X}_{n}^{\ba} \big) \cdot  M_{n}^{\bbeta} \Bigr]
\end{split}
\end{equation*}
by shifting the time index and using the fact that $M_0^{\bbeta} = \bar{M}_0^{\bbeta} = 0$. Next, using the formula for $\partial_x h$ we get
    \begin{align*}
        & \sum_{n \geq 0} \gamma^{n} \EE \Bigl[   \partial_x h(\zeta_n, p_n) \cdot M_n^{\bbeta} + \partial_{\bar x} h(\zeta_n, p_n) \cdot \bar{M}_n^{\bbeta}  +   M_{n+1}^{\bbeta}  \cdot ( p_{n+1} - p_{n} )\Bigr] 
        \\
         &\hskip 25pt
         = \sum_{n \geq 0} \gamma^n \EE \Bigl[  \big( (\check {\mathrm A}  - \delta I_d )^\top p_{n} + 2 Q (X_n^\ba - \bar X_n^\ba \big)\cdot M^{\bbeta}_n + \big(\bar {\mathrm A} ^\top \bar p_{n} + 2 (Q + \bar Q) \bar X_n^\ba \big) \cdot M_n^{\bbeta}  
        \\
        & \hskip 60pt - \gamma \big(  (\check {\mathrm A}  - \delta I_d)^\top p_{n} + 2 Q X_{n}^{\ba} + \bar {\mathrm A} ^\top \bar{p}_{n} + 2 \bar Q \bar{X}_{n}^{\ba} \big) \cdot  M_{n}^{\bbeta} \Bigr]
        \\
     & \hskip 25pt
     = 0
    \end{align*}
As a result, \eqref{fo:Gateaux} reduces to
    \begin{equation*}
    \begin{split}
        DJ(\ba)(\bbeta) & = \sum_{n \geq 0} \gamma^n \EE \Bigl[ \partial_a h(\zeta_n, p_n) \cdot \beta_n + \partial_{\bar a}h(\zeta_n, p_n) \cdot \bar{\beta}_n  \Bigr] 
        \\
        & = \sum_{n \geq 0} \gamma^n \EE \Bigl[ \big(p_n^\top {\mathrm B}  + 2 a_n^\top R + \bar{p}_n^\top \bar {\mathrm B}  + 2 \bar a_n^\top \bar R \big) \beta_n \Bigr].
    \end{split}
    \end{equation*}
Since $\cA_{ad}$ is a linear space, if $\ba$ is optimal, then $DJ(\ba)({\bbeta})\geq 0$ and $DJ(\ba)(-\bbeta)\geq 0$ for every $\bbeta \in \cA_{ad}$, hence $DJ(\ba)({\bbeta})=0$, implying that, for every $n \geq 0$
\begin{equation}
\label{fo:necessary}
         {\mathrm B} ^\top p_n + 2 R a_n + \bar {\mathrm B} ^\top \bar{p}_n + 2 \bar R \bar{a}_n = 0, \qquad \PP-a.s.
\end{equation}
which is the desired result.
\qed\end{proof}

\subsection{Identification of the optimal control}

Let us now assume that $\ba$ is an admissible optimal control and show that $a_n$ can be written as a linear combination of the components of the process $(p_n, \bar p_n)_{n\ge 0}$ where $p_n$ is the adjoint process of $\bX^{\ba}$, and as usual we use a bar on the top of a random variable to denote its conditional expectation with respect to $\sigma\{\underline\epsilon^0_n\}$. By taking conditional expectation $\EE[\,\cdot\,|\,\sigma\{\underline\epsilon^0_n\}]$ in equation~\eqref{fo:necessary}, we get
$$
    ({\mathrm B}  + \bar {\mathrm B} )^\top \bar p_n + 2 (R + \bar R) \bar a_n = 0
$$
from which we derive:
$$
\bar a_n=-\frac12 (R +\bar R)^{-1}({\mathrm B} +\bar {\mathrm B} )^\top\bar p_n,
$$
and going back to \eqref{fo:necessary} we finally get:
\begin{equation}
    a_n = -\frac{1}{2} R^{-1} {\mathrm B} ^\top p_n  - \frac{1}{2} R^{-1} \big( \bar {\mathrm B}^\top  - \bar R (R + \bar R)^{-1} ({\mathrm B}  + \bar {\mathrm B} )^\top \big)  \bar p_n
\label{fo:control_ut_expression_1}
\end{equation}
Using the notations
\eqref{fo:Gamma_Lambda} we get formula~\eqref{fo:control_ut_expression_2} for the process $\ba = (a_n)_{n \geq 0}$.
Given the well-known fact that at the optimum of standard LQ models, the adjoint processes are affine functions of the state, one may wonder if optimal controls can also be shown to be linear in the controlled state process $\bX^\ba$ and the conditional mean process $\bar{\bX}^\ba$.  
In order to prove that it is indeed the case, we introduce the following system of forward-backward stochastic equations for two $L^2-$discounted integrable processes $(\bX, \bp) \in \cX \times \cX$:
\begin{equation}
    \label{fo:FBDSE_state_adjoint}
    \left\{
    \begin{array}{rcl}
        X_{n+1} - X_n &=& \check {\mathrm A}  X_n + \bar {\mathrm A}  \bar{X}_n + {\mathrm B} \Gamma (p_n - \bar{p}_n) + ({\mathrm B}  + \bar {\mathrm B} ) \Lambda \bar{p}_n + \varepsilon_{n+1}^0 + \varepsilon_{n+1}
        \\
        p_{n+1} - p_{n} &=& - \gamma \big[ (\check {\mathrm A}  - \delta I_d)^\top p_{n+1} + 2Q X_{n+1} + \bar {\mathrm A} ^\top \bar{p}_{n+1} + 2 \bar Q \bar{X}_{n+1} \big] 
        \\
            && \hskip 180pt + Z_{n+1}^0 \varepsilon_{n+1}^0  + Z_{n+1} \varepsilon_{n+1},
    \end{array}
    \right.
\end{equation}
where the adapted processes $(Z_n)_{n \geq 0}$ and $(Z_n^0)_{n \geq 0}$ are part of the solution, and where $X_n$ satisfies the initial condition at time $n=0$, and $p_n$ satisfies the transversality condition which we view as a terminal condition at $\infty$. This justifies our terminology of \emph{forward-backward} system, abbreviated as FBDTSE (Forward-Backward Discrete-Time Stochastic Equations).\index[not]{FBDTSE}\index[sub]{forward-backward discrete-time stochastic equations} Notice that the first equation in \eqref{fo:FBDSE_state_adjoint} 
is nothing but the state equation for a linear control of the form to which we want to restrict ourselves. Also, by taking the conditional expectations $\EE[\cdot | \sigma\{\underline\epsilon_n,\underline\epsilon_n^0\}]$ of both sides of the second equation in \eqref{fo:FBDSE_state_adjoint}, we recover the fact that $\bp$ is necessarily the adjoint process of $\bX$.

\noindent
The following proposition provides a solution to the forward-backward system~\eqref{fo:FBDSE_state_adjoint}.

\begin{proposition}
\label{pr:solution_to_FBSDE_with_optimal_control}
    Under Assumption~\ref{ass:finit_cost} and~\ref{ass:positivity-qr}, we consider a control process $\ba = (a_n)_{n \geq 0}$ in a feedback form of its controlled state process $\bX^\ba$ and its conditional mean process $\bar \bX^\ba$ with control parameter $\theta = (- \Gamma P, -\Lambda \bar P)$ and policy $\pi_\theta : (x, \bar x) \mapsto  \Gamma P (x - \bar x) + \Lambda \bar P \bar x$  such that
    \begin{equation}
        \label{fo:closed-loop-control-on-x}
        a_n = \pi_\theta(X_n^\ba, \bar X_n^\ba) =  \Gamma P (X_n^\ba - \bar{X}_n^\ba) + \Lambda \bar P \bar{X}_n^\ba,
    \end{equation}
    for every $n \geq 0$, where $(\Gamma,\Lambda)$ are defined in~\eqref{fo:Gamma_Lambda} and $(P,\bar P)$ in~\eqref{fo:expression_of_solution_openloop_DARE}. Let $\bp^\ba = (p_n^\ba)_{n \geq 0}$ be the process defined by:
    \begin{equation}
        \label{fo:relation_ajoint_and_state}
        p^{\ba}_n = P(X_n^\ba - \bar{X}^{\ba}_n) + \bar P \bar{X}^{\ba}_n,
    \end{equation}
    for every $n \geq 0$. Then, the pair $(\bX^\ba, \bp^\ba)$ is a solution of the forward-backward system~\eqref{fo:FBDSE_state_adjoint} with
    \begin{equation}
        \label{fo:def_Zt}
        Z_n = \gamma \big( {\mathrm A} ^\top P + 2 Q \big), 
        \qquad \text{and}\qquad
        Z_n^0 = \gamma \big[ ({\mathrm A}  + \bar{{\mathrm A} })^\top \bar P + 2 (Q + \bar Q) \big]
    \end{equation}
    (note that $Z_n$ and $Z_n^0$ do not depend on $n$).
    Moreover, the process $\bp^\ba$ is the unique adjoint process corresponding to the pair $(\ba, \bX^\ba)$.
\end{proposition}

\begin{proof} 
    Proposition~\ref{pr:admissible_of_X} and Proposition~\ref{pr:admissibility_of_u_with_feedback_control} state that the control parameter $\theta = (-\Gamma P, -\Lambda \bar P) \in \Theta$, and the control process $\ba$ are admissible as well as its controlled state process $\bX^\ba$ in the sense that $\ba \in \cA_{ad}$ and $\bX^\ba \in \cX$. The definition~\eqref{fo:relation_ajoint_and_state} of $\bp^\ba$ implies that it is adapted to $(\sigma\{\underline\epsilon_n,\underline\epsilon^0_n\})_{n \geq 0}$ and it is $L^2-$discounted integrable, so $\bp^\ba \in \cX$. 
    
    For any time $n \geq 0$,  we take the conditional expectation with respect to $\sigma\{\underline\epsilon^0_n\}$ and obtain
    $a_n = \Gamma \big( p_n^\ba - \bar{p}_n^\ba \big) + \Lambda \bar{p}_n^\ba$ and $\bar{a}_n = \Lambda \bar{p}_n^\ba$. 
    From the state dynamics~\eqref{fo:MKV-state} of $\bX^{\ba}$, we deduce the forward equation in \eqref{fo:FBDSE_state_adjoint}:
    \begin{equation*}
    \begin{split}
        X_{n+1}^\ba - X_n^\ba &= \check {\mathrm A}  X_n^\ba + \bar {\mathrm A}   \bar{X}^{\ba}_n + {\mathrm B}  \big( \Gamma P (X_n^\ba - \bar{X}^\ba_n) + \Lambda \bar P \bar{X}^{\ba}_n \big) + \bar {\mathrm B}  \Lambda \bar P \bar{X}^{\ba}_n  + \varepsilon_{n+1}^0 + \varepsilon_{n+1} 
        \\
        & = \check {\mathrm A}  X_n^\ba + \bar {\mathrm A}   \bar{X}^{\ba}_n + {\mathrm B}  \Gamma ( p^\ba_n - \bar{p}^{\ba}_n)  +  ({\mathrm B}  + \bar {\mathrm B} ) \Lambda \bar p^{\ba}_n  + \varepsilon_{n+1}^0 + \varepsilon_{n+1} .
    \end{split}
    \end{equation*}
    To show that the process $\bp^\ba$ satisfies the backward equation  in \eqref{fo:FBDSE_state_adjoint}, we use the following trick based on the fact that $\bp^\ba - \bar{\bp}^\ba = (p_n^\ba - \bar{p}_n^\ba)_{n \geq 0}$ and $\bar{\bp}^\ba = (\bar p_n^\ba)_{n \geq 0}$:
    $$
        p^\ba_{n+1} - p^\ba_n = \big( (p^\ba_{n+1} - \bar p^\ba_{n+1}) - (p^\ba_n - \bar{p}^\ba_n) \big) + ( \bar p^\ba_{n+1} - \bar{p}^\ba_n).
    $$
    Now, by taking conditional expectation of the state dynamics~\eqref{fo:MKV-state} of $\bX^{\ba}$, we obtain a dynamics for $\bar{X}^\ba$ that for every $n \geq 0$,
    $$
        \bar{X}^{\ba}_{n+1} = [ ( {\mathrm A}  + \bar {\mathrm A} ) + ({\mathrm B}  + \bar {\mathrm B} ) \Lambda \bar P ] \bar{X}^\ba_n + \varepsilon_{n+1}^0.
    $$
    Then, multiplying both sides by $Z_{n+1}^0 = \gamma \big[  ({\mathrm A}  + \bar{{\mathrm A} } )^\top \bar P + 2 (Q + \bar Q)\big]$ and using the Riccati equation satisfied by $\bar P$ on the right-hand side of the resulting equality, we obtain 
    $$
        \gamma \big[  ({\mathrm A}  + \bar{{\mathrm A} } )^\top \bar P + 2 (Q + \bar Q) \big] \bar{X}^{\ba}_{n+1} = \bar P \bar{X}^{\ba}_{n} + Z_{n+1}^0 \varepsilon_{n+1}^0.
    $$
    By re-arranging the terms in the above equation and using $\bar P \bar X_{n+1}^\ba = \bar p_{n+1}^\ba$, we deduce
    \begin{equation}
        \label{fo:inside_proof_diff_bar_p}
        \bar{p}^{\ba}_{n+1} - \bar{p}^{\ba}_n = - \gamma \Big[ ( \check {\mathrm A}  + \bar {\mathrm A}  - \delta I_d)^\top \bar{p}^{\ba}_{n+1} + 2 (Q + \bar Q) \bar{X}^\ba_{n+1} \Big] + Z_{n+1}^0 \varepsilon_{n+1}^0.
    \end{equation}
    Similarly, from the dynamics of $\bX^\ba$ and $\bar{\bX}^\ba$, we have 
    $$
        X_{n+1}^\ba - \bar{X}^{\ba}_{n+1} = ({\mathrm A}  + {\mathrm B}  \Gamma P) \big(X_n^\ba - \bar{X}^\ba_n \big) + \varepsilon_{n+1}.
    $$
    Then multiplying both sides by $Z_{n+1}$ and using the Riccati equation satisfied by $P$ on the right-hand side of the resulting equality, we obtain
    $$
        \gamma ({\mathrm A} ^\top P + 2Q) (X_{n+1}^\ba - \bar{X}^{\ba}_{n+1}) = P \big(X_n^\ba - \bar{X}^\ba_n \big) + Z_{n+1} \varepsilon_{n+1},
    $$
    We rearrange the terms and use the fact that $p^\ba_n - \bar p^\ba_n = P (X_n^\ba - \bar{X}^\ba_n)$. We get:
    \begin{equation}
        \label{fo:inside_proof_diff_p_minus_bar_p}
        \begin{split}
         & (p^\ba_{n+1} - \bar p^\ba_{n+1}) - (p^\ba_n - \bar{p}^\ba_n) 
         \\
         &\hskip 35pt
         = - \gamma \big[ ( \check {\mathrm A}  - \delta I_d )^\top (p^\ba_{n+1} - \bar p^\ba_{n+1}) + 2 Q (X_{n+1}^\ba - \bar{X}^{\ba}_{n+1}) \big] + Z_{n+1} \varepsilon_{n+1}.
         \end{split}
    \end{equation}
    Adding equations~\eqref{fo:inside_proof_diff_bar_p} and~\eqref{fo:inside_proof_diff_p_minus_bar_p}, we immediately find that $\bp^\ba$ satisfies the backward equation~\eqref{fo:FBDSE_state_adjoint}. Thus, the pair of state-adjoint processes $(\bX^\ba, \bp^\ba) \in \cX \times \cX$ is a solution of the forward-backward system \eqref{fo:FBDSE_state_adjoint}.
\qed\end{proof}

The following theorem shows that we can construct an optimal control process $\ba$ for the MFC problem~\eqref{pb:MFC_L2_admissible_control} based on a solution of the FBDTSE system~\eqref{fo:FBDSE_state_adjoint}.

\begin{proposition}
\label{pr:existence_of_optimal_control_from_FBDTSE}
      We assume that Assumption~\ref{ass:finit_cost} and~\ref{ass:positivity-qr} hold. If $(\bX, \bp) \in \cX \times \cX$ is a solution to the forward-backward system~\eqref{fo:FBDSE_state_adjoint} with $X_0 = \varepsilon_0^0 + \varepsilon_0$, then the control process $\ba = (a_n)_{n \geq 0}$ given by 
    \begin{equation}
    \label{fo:linearity_control_in_adjoint_process}
         a_n = \Gamma (p_n - \bar p_n) + \Lambda \bar{p}_n
    \end{equation}
    is an optimal admissible control process for the MFC problem~\eqref{pb:MFC_L2_admissible_control}. 
    Besides, the process $\bX$ is the controlled state process of $\ba$, and the process $\bp$ is the unique adjoint process corresponding to $(\ba, \bX)$.
\end{proposition}

\begin{proof}
From the definition of the control process $\ba=(a_n)_{n\geq 0}$ in a feedback form given by equation~\eqref{fo:linearity_control_in_adjoint_process}, we can rewrite $a_n = \Gamma p_n + (\Lambda - \Gamma) \bar p_n$ and $\bar a_n = \Lambda \bar p_n$ for every $n \geq 0$.
Because the process $\bp \in \cX$ is adapted to $(\sigma\{\underline\epsilon_n,\epsilon_n^0\})_{n \geq 0}$, $\ba$ is also adapted to $(\sigma\{\underline\epsilon_n,\epsilon_n^0\})_{n \geq 0}$.
Using Jensen's inequality for conditional expectations and Cauchy-Schwarz inequality, we have 
$$
\sum_{n \geq 0} \gamma^n \EE[ \| a_n \|^2 ]  \leq 2( \| \Gamma \|^2 + \| \Lambda - \Gamma \|^2) \sum_{n \geq 0} \gamma^n \EE[ \| p_n \|^2 ] <  \infty,
$$
so that $\ba \in \cA_{ad}$. In addition, using $\Gamma (p_n - \bar p_n) = a_n - \bar a_n$ and $\Lambda \bar p_n = \bar a_n$, we see from the forward equation in the forward-backward system~\eqref{fo:FBDSE_state_adjoint} that the process $\bX$ satisfies the evolution dynamics~\eqref{fo:MKV-state}: 
$$
    X_{n+1} = {\mathrm A}  X_n + \bar {\mathrm A}  \bar{X}_n + {\mathrm B}  a_n + \bar {\mathrm B}  \bar a_n + \varepsilon_{n+1}^0 + \varepsilon_{n+1},
$$
for every $n \geq 0$ and starts with the initial state $X_0 = \varepsilon_0^0 + \varepsilon_0$. Corollary~\ref{cor:uniqueness_of_controlled_state_process} states that the process $\bX$ coincides with the controlled state process $\bX^\ba$ of $\ba$ in the sense that $X_n = X_n^\ba$ almost surely for every $n \geq 0$.
Meanwhile, taking the conditional expectation with respect to $\sigma\{\underline\epsilon_n,\underline\epsilon^0_n\}$ for each time $n \geq 0$ in the backward equation of $\bp$ in~\eqref{fo:FBDSE_state_adjoint}, and applying arguments similar to those used in the last paragraph in the proof of Proposition~\ref{pr:solution_to_FBSDE_with_optimal_control}, we show that the solution process $\bp \in \cX$ is the unique adjoint process corresponding to $(\ba, \bX)$. 

Returning to the expression of $a_n$ and the definitions of the matrices $\Gamma$ and $\Lambda$, we easily see that the control process $\ba$ satisfies equation~\eqref{fo:control_ut_expression_1}, that is
$$
a_n = -\frac{1}{2} R^{-1} {\mathrm B} ^\top p_n  - \frac{1}{2} R^{-1} \big( \bar {\mathrm B}^\top  - \bar R \tilde{R}^{-1} \tilde B^\top \big) \bar p_n
$$
for all $n \geq 0$, from which we can conclude that $\ba$ satisfies the necessary condition~\eqref{fo:necessary_PMP} of Pontryagin's maximum principle with adjoint process $\bp$, that is 
$$
{\mathrm B} ^\top p_n + \bar {\mathrm B} ^\top \bar p_n + 2 R a_n + 2 \bar R \bar a_n = 0
$$ 
for every $n \geq 0$.

To show that the process $\ba$ is an optimal control process for the MFC problem~\eqref{pb:MFC_L2_admissible_control} defined for $(\bX, \ba)$, we consider the Hamilton function $h(\zeta, p) = b(\zeta) \cdot p + f(\zeta) - \delta x \cdot p$ introduced in~\eqref{fo:Hamilton_function_h} and a perturbation direction $\bbeta \in \cA_{ad}$ for the control process. The difference $J(\ba + \lambda \bbeta) - J(\ba)$ can then be expressed in the following way using equation~\eqref{fo:difference_between_cost_function_in_hamilton_func}:
\begin{equation}
\label{fo:difference_between_J_u_beta}
\begin{split}
    J(\ba + \lambda \bbeta) - J(\ba)
    &= \sum_{n \geq 0} \EE \Big[ \gamma^n \big( h(\zeta_n^\lambda, p_n) - h(\zeta_n, p_n) \big) \Big] \\
    &\qquad +  \sum_{n \geq 0} \EE \Big[ \gamma^n \big( X_{n+1}^{\ba + \lambda \bbeta} - X_{n+1}^{\ba}) \cdot (p_{n+1} - p_n)  \big) \Big]
    \\
    &   
    =\sum_{n \geq 0} \gamma^n \EE \big[ \nabla_{\zeta} h(\zeta_n, p_n) \cdot ( \zeta_n^\lambda - \zeta_n) \big] \\
    &\qquad + \sum_{n \geq 0} \gamma^n \EE \big[ \frac{1}{2} \nabla^2_{\zeta \zeta} h(\eta_n^{\lambda}, p_n) (\zeta_n^\lambda - \zeta_n) \cdot ( \zeta_n^\lambda - \zeta_n)  \big]  \\
     & \hskip 75pt + \lambda \sum_{n \geq 0} \gamma^n \EE \big[ M_{n+1}^{\lambda, \bbeta} \cdot ( p_{n+1} - p_n ) \big] 
     \\
    &
    =(i) + (ii) + (iii)
\end{split}
\end{equation}
where $\zeta^{\lambda}_n = (X_n^{\ba+ \lambda \bbeta}, \bar{X}_n^{\ba+ \lambda \bbeta}, a_n + \lambda \beta_n, \bar a_n + \lambda \bar{\beta}_n)$, $\zeta_n = (X_n^{\ba}, \bar X_n^{\ba}, a_n, \bar a_n)$, and $M_n^{\lambda, \bbeta} = ( X_n^{\ba + \lambda \bbeta} - X_n^{\ba}) / \lambda$. In the second equality in the above equation, we use the Taylor expansion of the Hamilton function where the terms $\nabla_\zeta h$ and $\nabla_{\zeta, \zeta}^2 h$ are the Jacobian and Hessian matrices with respect to $\zeta$, and $\eta_n^\lambda = \zeta_n + \rho_n (\zeta_n^\lambda - \zeta_n)$ for some $\rho_n \in [0, 1]$.

From the proof of Proposition~\ref{pr:necessary_condition_of_PMP} (especially the computation of the Gateaux derivative of $DJ(\ba)(\bbeta)$ in equation~\eqref{fo:Gateaux}), and the fact that the control process $\ba$ satisfies equation~\eqref{fo:necessary_PMP}, we have
$$
    (i) + (iii) = \lambda \sum_{n \geq 0} \gamma^n \EE \big[ ({\mathrm B} ^\top p_n + \bar {\mathrm B} ^\top \bar p_n + 2 R a_n + 2 \bar R \bar a_n) \cdot \beta_n \big] = 0.
$$
We shall use the fact that the Hessian $\nabla_{\zeta, \zeta}^2 h$ of the Hamilton function takes the form 
$$
    \nabla_{\zeta, \zeta}^2 h(\zeta, p) = \begin{pmatrix}
	2Q & -2Q & 0 & 0 \\
	-2Q & 2 Q + 2 \tilde Q & 0 & 0 \\
        0 & 0 & 2R & -2R \\
        0 & 0 & -2R & 2 R + 2 \tilde R
	\end{pmatrix}.
$$
We also notice that $\zeta^{\lambda}_n - \zeta_n = \lambda ( M_n^{\lambda, \bbeta}, \bar{M}_n^{\lambda, \bbeta}, \beta_n, \bar{\beta}_n )$. Thus,
\begin{equation*}
\begin{split}
    (ii) = & \lambda^2 \sum_{n \geq 0} \gamma^n \EE \Big[ ( M_n^{\lambda, \bbeta} - \bar{M}_n^{\lambda, \bbeta} )^\top Q  ( M_n^{\lambda, \bbeta} - \bar{M}_n^{\lambda, \bbeta} ) + (\bar{M}_n^{\lambda, \bbeta})^\top \tilde{Q} \bar{M}_n^{\lambda, \bbeta} 
    \\
    & \hskip 50pt + ( \beta_n - \bar{\beta}_n )^\top  R (\beta_n - \bar{\beta}_n ) + \bar{\beta}_n^\top \tilde{R} \bar{\beta}_n ) \Big].
\end{split}
\end{equation*}
Because $Q, \tilde Q \succeq 0$ and $R, \tilde R \succ 0$, we obtain that for any perturbation process $\bbeta \in \cA_{ad}$,
$$
    J( \ba + \lambda \bbeta) - J(\ba) \geq 0.
$$
Therefore, the process $\ba$ given by equation~\eqref{fo:linearity_control_in_adjoint_process} is an optimal control process for the MFC problem~\eqref{pb:MFC_L2_admissible_control}.
\qed\end{proof}

\subsection{Proof of Theorem~\ref{th:existence_linear_control}}

We now provide the proof of the main result of this section, namely
Theorem \ref{th:existence_linear_control}.

\begin{proof}
\textbf{(Existence.)}
Let $\ba = (a_n)_{n \geq 0}$ be the control process specified by equation~\eqref{fo:optimal_control_linear_in_x_barx} with parameter $\theta = (- \Gamma P, - \Lambda \bar P) = (K^*, L^*)$ where $(\Gamma, \Lambda, P, \bar P, K^*, L^*)$ are defined in the statement of the theorem, and let $\bX$ be the state process in $\cX$ controlled by $\ba$. Proposition~\ref{pr:admissibility_of_u_with_feedback_control} shows that under assumptions~\ref{ass:finit_cost},~\ref{ass:positivity-qr}, and~\ref{ass:optimal-feedback-in-Theta}, $\theta$ is an admissible control parameter and $\ba \in \cA_{ad}$.

We also consider the process $\bp$ defined by equation~\eqref{fo:relation_ajoint_and_state}, that is, $p_n = P ( X_n - \bar X_n) + \bar P \bar X_n$ for every $n \geq 0$. Then, under assumptions \ref{ass:finit_cost} and \ref{ass:positivity-qr}, Proposition~\ref{pr:solution_to_FBSDE_with_optimal_control} implies that the pair of processes $(\bX, \bp)$ is a solution to the forward-backward system~\eqref{fo:FBDSE_state_adjoint}, and the process $\bp$ is the unique adjoint process associated with $(\ba, \bX)$. 

Thus, from the solution process $\bp$ and the coefficients $(\Gamma, \Lambda)$, we construct another control process $\ba^* = (a_n^*)_{n \geq 0}$ defined by~\eqref{fo:linearity_control_in_adjoint_process}:
$$
    a_n^* = \Gamma ( p_n - \bar p_n) + \Lambda \bar p_n,
$$ 
for every $n \geq 0$. Under assumptions \ref{ass:finit_cost} and \ref{ass:positivity-qr}, Proposition~\ref{pr:existence_of_optimal_control_from_FBDTSE} states that the process $\ba^*$ is admissible and optimal for the MFC problem~\eqref{pb:MFC_L2_admissible_control} with $(\bX, \ba^*)$. Meanwhile, the process $\bX$ is also the state process controlled by $\ba^*$, and the process $\bp$ is the unique adjoint process corresponding to $(\ba^*, \bX)$.  From the definition of $\bp$ in~\eqref{fo:relation_ajoint_and_state}, we have $p_n - \bar p_n = P (X_n - \bar X_n)$ and $\bar p_n = \bar P \bar X_n$, so we obtain
$$
      a_n^* = \Gamma P ( X_n - \bar X_n) + \Lambda \bar P \bar X_n = a_n,
$$
for every $n \geq 0$. Hence, the control process $\ba$ considered here is an optimal control to the MFC problem~\eqref{pb:MFC_L2_admissible_control}.

\noindent \textbf{(Uniqueness.)} Now suppose that there exists another optimal control $\ba' \in \cA_{ad}$ for the MFC problem~\eqref{pb:MFC_L2_admissible_control}. We denote $\beta_n = a'_n - a_n$ the difference between the two control processes at time $n$. Because both $\ba$ and $\ba'$ are $L^2-$discounted-integrable, $\bbeta = (\beta_n)_{n \geq 0}$ is also $L^2-$discounted integrable by Assumption~\ref{ass:finit_cost}, thus $\bbeta \in \cA_{ad}$. Because the process $\bp$ is the adjoint process corresponding to $(\ba, \bX)$, from equation~\eqref{fo:difference_between_J_u_beta}, we get: 
\begin{equation*}
\begin{split}
    0 &= J(\ba') - J(\ba) 
    \\ 
    & = \sum_{n \geq 0} \gamma^n \EE \big[ ({\mathrm B} ^\top p_n + \bar{{\mathrm B} }^\top \bar{p}_n + 2 R a_n + 2\bar{R} \bar{a}_n ) \cdot \beta_n \big] + \sum_{n \geq 0} \gamma^n \EE \big[ f( M_n, \bar M_n, \beta_n, \bar{\beta}_n ) \big]
    \\
    & = (i) + (ii)
\end{split}
\end{equation*}
where $M_n = X^{\ba'}_n - X^{\ba}_n$ satisfies the equation 
$$
M_{n+1} = {\mathrm A}  M_n + \bar {\mathrm A}  \bar M_n + {\mathrm B}  \beta_n + \bar {\mathrm B}  \bar{\beta}_n
$$ 
for every $n \geq 0$ with $M_0 = 0$. From the optimality of the control process $\ba$, the necessary condition of the Pontryagin's maximum principle in Proposition~\ref{pr:necessary_condition_of_PMP} yields $(i) = 0$. In addition, Assumption~\ref{ass:positivity-qr} implies $\EE \big[ f( M_n, \bar M_n, \beta_n, \bar{\beta}_n ) \big] \geq 0$ for every $n \geq 0$. Thus $(ii)=0$. Since all summands in $(ii)$ are nonnegative, for every $n\geq 0$,
\begin{equation*}
\begin{split}
    0 & = \EE\big[ f(M_n, \bar M_n, \beta_n, \bar{\beta}_n ) \big] = \EE[ (M_n - \bar M_n)^\top Q (M_n - \bar M_n) ] + \EE[ \bar{M}_n^\top \tilde{Q} \bar M_n] 
    \\
    & \hskip 120pt + \EE[ ( \beta_n - \bar{\beta}_n)^\top R (\beta_n - \bar{\beta}_n) ] + \EE[ (\bar{\beta}_n)^\top \tilde{R} \bar{\beta}_n ]
\end{split}
\end{equation*}
Because $Q, \tilde Q \succeq 0$ and $R, \tilde R \succ 0$, we must have $\beta_n-\bar{\beta}_n=0$ and $\bar{\beta}_n=0$ almost surely, hence $\beta_n=0$ for every $n \geq 0$. This implies the uniqueness of the optimal control process.
\qed\end{proof}

\section{Notes and Complements}

Discrete-time LQ mean field control problems without common noise are treated in \cite{ElliottLiNi2013,NiElliottLi2015}, where the analysis leads to decoupled Riccati equations. Continuous-time LQ mean field control problems and related McKean--Vlasov control formulations are studied, e.g.,  in \cite{Yong2013SICON,HuangLiYong2015,Graber2016}; see also \cite[Chapter~6]{CarmonaDelarue_book_I} for a systematic treatment of the continuous-time theory. 

Classical LQ problems without mean-field interactions have also been studied from the reinforcement learning perspective, in particular for learning linear feedback gains; see \cite{recht2018tour,fazel2018global,yang2019provably}. The policy-gradient analysis for the LQ mean field model considered here is developed in Chapter~\ref{ch:numeric_II}.

\part{Numerical Implementations}

\chapter{Implementation of the First Abstract Model}
\label{ch:numeric_I}

\begin{abstract}
	\emph{
		This chapter presents numerical methods for solving the Mean Field Markov Decision Processes (MFMDPs) introduced in Chapter~\ref{ch:AAP}.
		We first describe a tabular Q-learning algorithm that uses simplex discretization to handle the measure-valued state space.
		A proof of convergence is provided for this approach.
		To improve scalability, we then discuss an adaptation of the Deep Deterministic Policy Gradient (DDPG) algorithm.
		This method leverages neural network approximations to handle continuous state and action spaces, as motivated by the theoretical developments in the first part of this monograph.
	}
\end{abstract}

\section{Introduction}
\label{se:ch5_intro}

This chapter is the first of two dedicated to the numerical resolution of Mean Field Control (MFC) through reinforcement learning.
It bridges the gap between the rigorous mathematical models developed in Chapter~\ref{ch:AAP} and their practical implementation.
As discussed in Chapter~\ref{ch:AAP}, finding optimal policies for mean field Markov decision processes (MFMDPs) requires overcoming significant computational hurdles, primarily due to the potential infinite-dimensional nature of the lifted state space consisting of probability measures.
We revisit the concept of the ``lifted'' (level-1) problem, where the central planner controls the flow of distributions, and contrast it with the level-0 dynamics of individual agents.
To address this challenge in the general setting, we present two distinct approaches.

First, we discuss a method based on space discretization, leading to a tabular version of the Q-learning algorithm.
By projecting the mean-field distributions onto a finite grid (a discretized simplex), we can rigorously establish the convergence of the algorithm (Theorem~\ref{th:main-cv-tabular}).
However, as we shall see, this approach suffers from the curse of dimensionality, making it impractical for problems where the level-0 state and action spaces are large.

Second, in order to circumvent the limitations of discretizations, we introduce function approximation techniques, specifically adapting the Deep Deterministic Policy Gradient (DDPG) \index[not]{DDPG} algorithm.
This approach operates directly on continuous state and action spaces by relying on artificial neural networks to approximate the Q-function (the critic) and the policy (the actor).\index[sub]{actor}\index[sub]{critic}
While theoretical convergence guarantees for deep RL in the mean field setting remain an active area of research, the numerical experiments reported in Section~\ref{se:numres} demonstrate the practical efficacy of this approach for models ranging from cybersecurity to distribution planning.

By exploring both a theoretically grounded discrete method and a scalable continuous method, this chapter equips the reader with the tools necessary to deploy Mean Field Reinforcement Learning algorithms in diverse and complex environments.
These methods prepare the way for Chapter~\ref{ch:numeric_II}, where we specifically focus on Linear-Quadratic MFC problems for which we develop a more refined policy gradient analysis and provide stronger convergence guarantees.

\section{Background on Reinforcement Learning Methods}
\label{sec:rl_background}

Before delving into the specific algorithms tailored to MFC, we briefly review the core concepts of Reinforcement Learning (RL) that underpin the methodologies presented in this and the subsequent chapter. 

\subsection{The Reinforcement Learning paradigm}
The core objective of RL is to solve an optimal control problem formulated as a \emph{Markov Decision Process (MDP)}. An MDP is defined by a state space $S$, an action space $A$, a transition probability $P(x' | x, a)$, a one stage cost function $f(x, a)$, and a discount factor $\gamma \in (0,1)$. Given a policy $\pi$ that determines the choice of actions $a_n\in A$ based on the current value of the state $x_n\in S$, the goal is to minimize the expected discounted cumulative cost:
\begin{equation}
    \label{fo:single_agent_objective}
    J(\pi) = \mathbb{E}^\pi \left[ \sum_{n=0}^\infty \gamma^nf(X_n, A_n) \Big| X_0=x_0 \right],
\end{equation}
where we use upper case characters to denote random variables and lower cases for specific values or realizations of these random variables. Thus typically $x_{n}$ will denote a realization or a sample of the random variable $X_n$ describing the state as a random variable in $S$. So the description of the MDP stipulates that $X_{n+1}$ is a random variable in $S$ with distribution $P(\cdot | x_n, a_n)$ when $X_n=x_n$ and $A_n=a_n$, and its realizations $x_{n+1}$ are sampled from $P(\cdot | x_n, a_n)$. The objective $J(\pi)$ is analogous to the multi-agent reward maximization described in Chapter~\ref{ch:review} (see for instance equation~\eqref{fo:ith_objective}), but adapted to our cost-minimization framework.

A distinguishing feature of RL is its \emph{model-free} approach. In many practical scenarios, the agent does not have access to the internal dynamics of the environment (the transition kernel $P$). Instead, they learn by interacting with the environment, observing states, taking actions, and incurring costs without a prior model of the environment's transition or cost structure. This paradigm shift requires the agent to balance the exploration of possibly new parts of the state and action spaces, with the exploitation of known low-cost regions, a challenge that lies at the heart of the algorithms discussed in this chapter.

\subsection{Value-based methods and Q-learning}
Value-based methods rely on the evaluation of the state-action value function, or Q-function, the value $Q(x,a)$ representing the expected cumulative cost incurred by starting in a state $x$, taking an initial action $a$, and then following a policy $\pi$ thereafter:
$$
    Q^\pi(x, a) = \mathbb{E}^\pi \left[ \sum_{n=0}^\infty \gamma^n f(X_n, A_n) \Big| X_0 = x, A_0 = a \right].
$$
The optimal Q-function, denoted $Q^*$, satisfies the Bellman optimality equation:
$$
    Q^*(x,a) = f(x, a) + \gamma \mathbb{E}_{X' \sim P(\cdot|x,a)} \left[ \inf_{a' \in A} Q^*(X', a') \right].
$$
A fundamental value-based algorithm is Q-learning~\cite{watkins1992q}, which iteratively updates an estimate of $Q^*$ using temporal difference learning. During training, upon observing a transition $(x_n, a_n, c_n, x_{n+1})$, the Q-value is updated via:
$$
    Q(x_n, a_n) \leftarrow (1-\eta_n) Q(x_n, a_n) + \eta_n \left( c_n + \gamma \inf_{a'} Q(x_{n+1}, a') \right),
$$
where $\eta_n$ is the learning rate. \index[sub]{learning rate} For environments with finite state and action spaces, these Q-values can be stored in a table, justifying the terminology \emph{tabular Q-learning}. \index[sub]{tabular Q-learning} A critical aspect of such algorithms is the exploration-exploitation trade-off, often managed via an $\epsilon$-greedy\index[not]{$\epsilon$-greedy} strategy, where the agent selects the action prescribed by the current policy with probability $1-\epsilon$, and with probability $\epsilon$, explores the action space by choosing an action uniformly at random. This tabular approach will be the foundation for the discretized simplex Mean Field Q-learning algorithm in Section~\ref{se:tabularQ-discrete}.

\subsection{Policy gradient and zeroth-order optimization}
When dealing with continuous action spaces or highly parameterized policies (such as neural networks), value-based methods become computationally expensive due to the computation of the infimum $\inf_{a'} Q(x', a')$ for each state $x'$. Policy gradient methods address this by directly parameterizing the policy $\pi_\theta(x)$ with weights $\theta$ and performing stochastic gradient descent on the expected total cost $J(\theta)$.

However, in many RL applications,  calculating exact gradients is not possible because the environment acts as a black box. \emph{Zeroth-order optimization} \index[sub]{zeroth-order optimization}(also known as derivative-free optimization) \index[sub]{derivative-free optimization}overcomes this obstacle by estimating gradients using only cost evaluations. By applying small random perturbations $u$ to the parameters and measuring the resulting cost $J(\theta + u)$, one can construct an unbiased estimate of the gradient of a smoothed objective function. This approach will be exploited in the case of LQ MFC models in Chapter~\ref{ch:numeric_II}.

\subsection{Actor-Critic and Deep RL}
To combine the benefits of value-based and policy-based methods, \emph{Actor-Critic} architectures maintain two separate structures: the \emph{actor}, representing the policy parameterized by $\theta^\pi$, and the \emph{critic}, estimating the value function parameterized by $\theta^Q$. The critic evaluates the current policy's performance, providing a low-variance critique that the actor uses to guide its parameter updates.

When state and action spaces are large or continuous, tabular methods are infeasible due to the curse of dimensionality. \emph{Deep RL} solves this quandary by using artificial neural networks as function approximators for the actor and the critic. A prominent algorithm in this class is DDPG~\cite{lillicrap2015continuous}, \index[not]{DDPG} an actor-critic method specifically designed for continuous action spaces. DDPG utilizes experience replay buffers and target networks to stabilize the neural network training. This methodology underpins the Mean Field DDPG algorithm presented in Section~\ref{sec:DDPG-algo}, enabling scalable numerical solutions applied to complex problems in Section~\ref{se:numres}.

\subsection{Mean Field MDP: notations and environment}
\label{subsec:mf_env_notations}

In the context of MFRL, we adapt the standard MDP paradigm to the lifted problem described in Chapter~\ref{ch:AAP}. We recall that the population distribution $\mu_n \in \mathcal{P}(S)$ (where $S$ is the individual state space) serves as the \emph{state} of the MFMDP, often referred to as a level-1 state. The \emph{action}, or level-1 action, is a joint distribution $\bar a_n \in \mathcal{P}(S \times A)$ consistent with $\mu_n$, meaning its first marginal is $\mu_n = \text{pr}_1(\bar a_n)$. The environment in this setting is a \emph{mean field simulator} that implements the lifted dynamics $\bar F$: given the current distribution $\mu_n$ and the lifted action $\bar a_n$, it returns the next distribution $\mu_{n+1}$ and the aggregate social cost $\bar f(\mu_n, \bar a_n)$.

Ideally, this environment provides a perfect simulation of the infinite-population dynamics. However, in practice, simulating a continuous distribution $\mu_n$ is often intractable, especially with high-dimensional state spaces. Consequently, numerical implementations frequently adopt a \emph{particle system simulator}. In this approach, the environment simulates a system of $N$ interacting agents $\underline{X}_n = (X^1_n, \dots, X^N_n)$ and provides the empirical distribution $\hat{\mu}^N_n = \frac{1}{N} \sum_{i=1}^N \delta_{X^i_n}$ as an observation. Under suitable propagation-of-chaos assumptions, and conditionally on the common information when common noise is present, such empirical distributions approximate the corresponding mean-field distribution as $N \to \infty$.

In the numerical examples of Section~\ref{se:numres} in this chapter, we primarily focus on simulating the mean field dynamics directly. In contrast, Chapter~\ref{ch:numeric_II} investigates more granular configurations. Specifically, Section~\ref{subsection:modelfree_MKV_simulator} discusses an idealized McKean-Vlasov simulator, while Section~\ref{subsec:PG-popsimu} introduces a finite-population simulator. The latter allows for a theoretical study and numerical illustration of the additional error induced by the size of the finite population.

\section{Mean Field Q-Learning}

\subsection{Controls for finite state and action spaces}\,

In the rest of this section, we assume that $S$ and $A$ are finite, we denote their cardinalities by $|S|$ and $|A|$ respectively, and we denote by $x^{(1)},\dots,x^{(|S|)}$ and  $\alpha^{(1)},\dots,\alpha^{(|A|)}$  their elements.
We first revisit the description of the action space and then propose two RL methods for this setting.
Subsequently, we explain in Section~\ref{sec:DDPG-algo} how to adapt these RL techniques to continuous spaces.

Before introducing the mean-field Q-learning algorithm, we first provide a representation of the set $\bar \Gamma \subseteq \bar S \times \bar A = \cP(S) \times \cP(S \times A)$ on which the $\bar Q^*$ function is defined.
Recall that the lifted state-action space $\bar \Gamma$, defined in Equation~\eqref{eq:bar_Gamma} of Chapter~\ref{ch:AAP}, captures the requirement that the marginal of a lifted action must be consistent with the current mean field state.

Since we assume that $S$ is finite, its lifted space $\cP(S)$ can be identified with a simplex $\mathfrak{S}$ in $\RR^{| S |}$.
In other words, we treat a distribution $\mu \in \cP(S)$ as an $|S|$-dimensional vector $(\mu^{(i)})_{i=1,\dots,|S|}$ whose non-negative coordinates sum up to one.
Similarly, since $A$ is finite, we identify $\cP(A)$ with a simplex $\mathfrak{A}$ in $\RR^{|A|}$.
However, representing admissible actions $\bar a \in \bar U(\mu) \subseteq \cP(S \times A) $ of the lifted MDP requires careful handling due to the underlying constraints.
A first approach is to identify $\cP(S \times A)$ with a simplex in $\RR^{|S| \times |A|}$ and to view a lifted action $\bar a$ as a $|S|\times |A|$ matrix $ \big( \bar a( x^{(i)}, \alpha^{(j)}) \big)_{1 \leq i \leq |S|, 1 \leq j \leq |A|}$ of non-negative numbers summing up to $1$.
Then a pair $(\mu, \bar a) \in \bar S \times \bar A$ is in $\bar \Gamma$ if and only if the following linear constraint is satisfied: $$
\sum_{j=1}^{|A|} \bar a( x^{(i)}, \alpha^{(j)} ) = \mu^{(i)}, \qquad i = 1, \cdots, |S|.
$$
While this transformation is straightforward, it is insufficient for our purposes.
It only provides a representation of the actions and controls from the central planner's perspective. It fails to capture the strategy functions of non-randomized stationary mixed Markovian closed-loop policies for an individual agent in our original optimization problem.

For any pair $(\mu, \bar a) \in \bar \Gamma$, we can define a kernel $k_{\mu}: S \to \cP(A)$ by setting, for $i=1,\dots,|S|$ and $j=1,\dots,|A|$,
\begin{equation*}
	k_\mu(x^{(i)})(\alpha^{(j)}) = \frac{\bar a( x^{(i)}, \alpha^{(j)} )}{\mu(x^{(i)})}, \qquad \hbox{ if } \mu(x^{(i)})>0,
\end{equation*} 
and by choosing an arbitrary probability measure on $A$ when $\mu(x^{(i)})=0$.

Note that common randomization is absent here.
As proved above (see Theorem~\ref{th:MFMDP_DPP_BOREL}), there exists a non-randomized stationary policy for the lifted MDP.
So the central planner can look for strategy functions within the set:
\begin{equation}
	\label{eq:def-calA-S-PA}
	\tilde{\cA} := \{ \tilde a: S \to \cP(A) \  | \   \tilde a \text{ Borel measurable} \}.
\end{equation}
So, because of the constraint, we can identify a probability measure $\bar a$ on the product space $S\times A$ whose first marginal is given and equal to $\mu$, with the probability kernel $\tilde a$ appearing in its disintegration with respect to its first projection. In this way, we can view the Q-function $\bar Q^*$ as a function of $\mu$ and $\tilde a$ instead of $\mu$ and $\bar a$, avoiding the redundancy imposed by the constraint.
So if we introduce the function $\tilde Q^*: \cP(S) \times \tilde{\cA} \to \RR$ defined by:
\begin{equation}
	\tilde Q^*(\mu, \tilde a) := \bar Q^*(\mu, \mu \measprod \tilde a),
\end{equation}
the Bellman equation~\eqref{eq:Bellman_equation_Q} becomes:
\begin{equation}
	\label{eq:optimal_Q_kernel_version}
	\tilde Q^*(\mu, \tilde a) = \int_{S \times A} f(x, \alpha, \mu \measprod \tilde a )  \tilde a(x, d\alpha) \mu( dx)+ \gamma \EE \left[ \inf_{\tilde a'\in \tilde{\cA}} \tilde Q^*( \mu_1, \tilde a')\right],
\end{equation}
for $(\mu, \tilde a) \in \cP(S) \times \tilde{\cA}$, where $\mu_1 = \bar F( \mu, \mu \measprod \tilde a, \varepsilon^0)$.
Here $\bar F$ denotes the lifted system function introduced in equation~\eqref{eq:Fbar_pushforward} of Chapter~\ref{ch:AAP}, representing how the population distribution evolves given the current state and the chosen level-1 action.
An important observation from Chapter~\ref{ch:AAP} is that the discrete-time evolution remains stochastic in the presence of common noise $\varepsilon^0$, as shown in Equation~\eqref{eq:MFMDP_admissible_state_process}.
Although $S$ and $A$ are finite, equation \eqref{eq:optimal_Q_kernel_version} is still characterized as a fixed point in the space of bounded lower semi-continuous functions on a closed subset of a finite-dimensional Euclidean space, consistent with the measurability considerations addressed in deriving equation~\eqref{eq:Bellman_equation_Q}.
We also introduce the function $\tilde f: \cP(S) \times \tilde{\cA} \to \RR$ defined by:
$$
	\tilde f(\mu, \tilde a) := \bar f(\mu, \mu \measprod \tilde a) = \int_{S \times A} f(x, \alpha, \mu \measprod \tilde a)  \tilde a(x, d\alpha) \mu( dx), \qquad (\mu, \tilde a) \in \cP(S) \times \tilde{\cA}.
$$
The function $\bar f$ relates back to the lifted cost function defined in Equation~\eqref{de:bar_f} of Chapter~\ref{ch:AAP}, which aggregates individual agent costs into a single scalar value for the central planner.
In the remainder of this section, we propose two model-free algorithms based on the optimal state-action value function $\bar Q^* : \bar \Gamma \to \RR$ or, equivalently, $\tilde Q^*: \cP(S) \times \tilde{\cA} \to \RR$.

\subsection{Simplex discretization and tabular MFQ-learning}\,
\label{se:tabularQ-discrete}

We consider two settings, depending on whether the controls at level-0 are mixed or pure.
In both cases, we discretize the simplexes, which leads to a finite-state finite-action MDP for which the tabular Q-learning algorithm can be used.
We establish its convergence.
When using pure controls, we can prove the convergence of the value function as well as that of the associated optimal control.

\subsubsection{Q-learning with level-0 mixed controls}
\label{se:tabularQ-mixed}

We first consider the situation where the controls are mixed at the level-0, meaning that the generic agent can use randomized actions.
Since the simplexes $\mathfrak{S}$ and $\mathfrak{A}$ are infinite, a tabular version of the Q-learning algorithm cannot be directly applied to approximate $\tilde Q^*$.
To address this, we first approximate these simplexes using finite subsets $\check{\mathfrak{S}} \subset \mathfrak{S}$ and $\check{\mathfrak{A}} \subset \mathfrak{A}$.
Let $\check{\cA} = \{ \check{a} : S \to \check{\mathfrak{A}} \}$.
In particular, $|\check{\cA}| = | \check{\mathfrak{A}} |^{|S|}$ because we identify functions in $\check{\cA}$ with $|S|$-dimensional vectors  whose entries take values in the finite set $\check{\mathfrak{A}}$.
To ensure that the mean-field term takes values in the finite set $\check{\mathfrak{S}}$, we use a projection:
At time $n$, given $\mu_n \in \check{\mathfrak{S}}$, we compute $\mu_{n+1} = \bar F( \mu_n, \mu_n \measprod \check{a}, \varepsilon_{n+1}^0)$, and then we project $\mu_{n+1}$ back on $\check{\mathfrak{S}}$ using a projection operator $\proj_{\check{\mathfrak{S}}}: \cP(S) \to \check{\mathfrak{S}}$.
Precise definitions of the discretization and the projection are provided below, after introducing a discrete version of the original MFC problem.

More precisely, we consider the \defi{projected MFC problem}\index[sub]{projected MFC problem}:
$$
	\inf_{ \check \pi \in \check \Pi} \check J^{\check \pi}(\mu_0),\qquad \mu_0 \in \check{\mathfrak{S}},
$$
where $\check \Pi = \{ \check \pi: S \times \check{\mathfrak{S}} \to \check{\mathfrak{A}} \}$ is a (finite) set of policies, and for every strategy function $\check \pi : S \times \check{\mathfrak{S}} \to \check{\mathfrak{A}}$,
$\check J^{\check \pi}: \check{\mathfrak{S}} \to \RR$ is defined by:
\begin{equation}
	\label{eq:generic-MFC-fctmeasure-cost-proj}
	\check J^{\check \pi}(\mu_0) =  \EE \left[ \sum_{n \geq 0} \gamma^n  \tilde f \Big( \mu^{\mu_0, \check \pi}_n,  \check \pi(\cdot, \mu^{\mu_0, \check \pi}_n)  \Big) \right]
\end{equation}
where
\begin{equation}
	\label{eq:generic-MFC-fctmeasure-dyn-proj}
	\mu^{\mu_0, \check \pi}_{n+1} = \proj_{\check{\mathfrak{S}}} \circ  \bar F \Big( \mu_n^{\mu_0, \check \pi}, \mu_n^{\mu_0, \check \pi} \measprod \check \pi( \cdot, \mu_n^{\mu_0, \check \pi}), \varepsilon_{n+1}^0 \Big) =: \check \Phi^{\check \pi, \varepsilon_{n+1}^0}(\mu_n^{\mu_0, \check \pi}).
\end{equation}
We will denote by $\check J^*$ and $\check Q^*$ respectively the optimal state and state-action value functions of this projected MFC problem.
Here $ \check Q: \check{\mathfrak{S}} \times \check{\cA} \to \RR$ can be represented by a matrix (also called a table) in $\RR^{ |\check{\mathfrak{S}}| \times | \check{\cA}| }$ and is viewed as an approximation of $\tilde{Q}^* : \cP(S) \times \tilde{\cA} \to \RR$ of the original MFC problem.

This problem can be viewed as an MDP with finite state and action spaces.
In this case, one can adapt in a straightforward way the classical tabular Q-learning algorithm.
The pseudo-code is provided in Algorithm~\ref{algo:Qtable-projection}.
At each iteration, the algorithm uses one sampled transition to update one table entry.
The state-action pairs can be prescribed if the learner has access to a generative simulator, which can sample a successor from the projected transition law at any queried pair $(\check\mu_n,\check a_n)$, using fresh transition noise.
Alternatively, they can be obtained along an online trajectory, in which case the next state is the sampled successor of the current state.
In either case, the pair at an update is selected before the fresh transition noise is drawn.
The quantitative bounds below require the covering condition in Assumption~\ref{hyp:covering-time}; trajectory-based sampling is also illustrated in Section~\ref{se:numres}.

\begin{remark}
	In this chapter, we adopt the perspective of a central planner who learns to control the distribution to minimize a social cost.
	For simplicity in the RL application, we assume the planner directly observes the mean-field population distribution.
	This allows us to analyze the convergence of an RL algorithm rigorously.
	However, in practice it is more likely that the planner could not observe this ideal distribution but rather a distorted version, such as the empirical distribution over a finite number of agents.
	In Section~\ref{ref:PGconv}, we will present and prove the convergence of a model-free policy gradient algorithm for the LQ setting in the case where an empirical distribution is observed instead of the mean-field one.
	In the present situation, it is reasonable to expect that the propagation of chaos could be used to design and analyze a Q-learning algorithm based on the mere observation of a finite population, but this is beyond the scope of the present discussion.
    \end{remark}

    \begin{remark}
    In the finite state setting, the set of empirical distributions with $N$ agents can be interpreted as the set of histograms where each bin can only take values that are multiples of $1/N$: $\{(\mu_x)_{x \in S} \,:\, \forall x \in S, \mu_x \in\{ k/N, k = 0,\dots,N\}, \sum_x \mu_x = 1\}$.
	This set can be used as a discretized simplex $\check{\mathfrak{S}}$, which is precisely the structure that is used in the tabular Q-learning algorithm.
\end{remark}

\begin{algorithm}[htbp]
	\DontPrintSemicolon
	\KwData{A number of updates $N_{\mathrm{upd}}$; a sequence of learning rates $(\eta_n)_{n=0,\dots,N_{\mathrm{upd}}-1}$; a sequence of state-action pairs $(\check\mu_n,\check{a}_n)_{n \ge 0} \in \check{\mathfrak{S}} \times \check{\cA}$.}
	\KwResult{$\check Q_{N_{\mathrm{upd}}}$, an approximation of $\tilde{Q}^*$ on $\check{\mathfrak{S}} \times \check{\cA}$.}
	\Begin{
	Initialize table $\check Q_0=0 \in \RR^{|\check{\mathfrak{S}}| \times |\check{\cA}|}$, $\mu_{0} \in \mathfrak{S}$ and $\check{a}_{0} \in \check{\cA}$\;
	\For{ $n =0, 1, \dots, N_{\mathrm{upd}}-1$}{
	\vskip 2pt
	Execute action $\check{a}_n$ in state $\check\mu_n$, observe $\check \mu'_{n+1} = \proj_{\check{\mathfrak{S}}} \circ \bar F(\check \mu_n, \check \mu_n \measprod \check{a}_{n}, \varepsilon_{n+1}^0)$ and cost $\tilde f(\check \mu_{n}, \check{a}_{n})$ \;
	Initialize $\check Q_{n+1} = \check Q_n$ on $\check{\mathfrak{S}} \times \check{\cA}$\;
	Set $\check Q_{n+1}(\check \mu_n , \check  a_n) = (1- \eta_n) \check Q_n( \check \mu_{n} , \check  a_{n}) + \eta_n  \left( \tilde f( \check \mu_{n}, \check{a}_{n}) + \gamma  \min_{\check{a}' \in \check{\cA}} \check  Q_n( \check \mu'_{n+1}, \check{a}' ) \right)$ \;
	}
	\KwRet{$\check Q_{N_{\mathrm{upd}}}$}
	}
	\caption{Mean-Field Q-learning (MFQ) with simplex discretization}
	\label{algo:Qtable-projection}
\end{algorithm}

Algorithm~\ref{algo:Qtable-projection} returns the table $\check Q_{N_{\mathrm{upd}}}$ after $N_{\mathrm{upd}}$ transitions and Q-updates.
We prove below that this table converges to the optimal Q-function $\tilde Q^*$ in a suitable sense.

To keep the chapter at a reasonable length, we will make the following simplifying assumptions.

We endow the simplexes $\mathfrak{S}$ and $\mathfrak{A}$ respectively with the Euclidean distances $d_{\mathfrak{S}}$ and $d_{\mathfrak{A}}$ of the spaces $\RR^{|S|}$ and $\RR^{|A|}$.
Because $S$ is finite, we can identify $\tilde{\cA}$ defined in~\eqref{eq:def-calA-S-PA} with $\cP(A)^{|S|}$  with the distance $d_{\tilde{\cA}}(\tilde a, \tilde a') = \sup_{x \in S} d_{\mathfrak{A}}( \tilde a(x), \tilde a'(x)) $ for $\tilde a, \tilde a' \in \tilde{\cA}$.
Furthermore, we consider the following discretizations of the simplexes.
Let $\varepsilon_{\mathfrak{S}}>0$ satisfying: for all $\mu \in \mathfrak{S},$ there exists $\check \mu \in \check{\mathfrak{S}}$ s.t. $d_{\mathfrak{S}}(\mu, \check \mu) \le \varepsilon_{\mathfrak{S}}$.
Similarly, let $\varepsilon_{\mathfrak{A}} > 0$ satisfying: for all $\nu \in \mathfrak{A}$, there exists $\check \nu \in \check{\mathfrak{A}}$ such that $d_{\mathfrak{A}}( \nu, \check \nu ) \leq \varepsilon_{\mathfrak{A}}$.
Because $S$ is finite and the definition of the distance $d_{\tilde{\cA}}$, we have for every $\tilde a \in \tilde{\cA}$, there exists $\check{a} \in \check{\cA}$, such that $d_{\tilde{\cA}}( \tilde a, \check{a}) \leq \varepsilon_{\mathfrak{A}}$.

\begin{remark}
	In practice, the subsets $\check{\mathfrak{S}}$ and $\check{\mathfrak{A}}$ should be chosen based on two criteria. First,  $\varepsilon_{\mathfrak{S}}$ and $\varepsilon_{\mathfrak{A}}$ should be as small as possible.
	Second, the number of points in these subsets, namely $|\check{\mathfrak{S}}|$ and $|\check{\mathfrak{A}}|$, should be as small as possible.
	These quantities appear in the bound given by Theorem~\ref{th:main-cv-tabular} below.
	Of course, to decrease $\varepsilon_{\mathfrak{S}}$, we expect that more points are needed and hence $|\check{\mathfrak{S}}|$ should increase (and likewise for $\varepsilon_{\mathfrak{A}}$ and $|\check{\mathfrak{A}}|$).
	This trade-off has to be managed carefully.
	A standard choice is to cover the probability simplex with balls of a given radius.
	If the dimension of the simplex is $d$ and the radius of these balls is $\delta$, then %
	the minimal number of centers needed is of order $\Theta(\delta^{-d})$ as $\delta\downarrow0$, with constants that may depend on $d$.
	In our setting, the number of points would be of order $\Theta(\varepsilon_{\mathfrak{S}}^{1-|S|})$ and $\Theta(\varepsilon_{\mathfrak{A}}^{1-|A|})$ respectively.
\end{remark}

\begin{assumption}\label{hyp:bdd-smooth-data}
	\textbf{Regularity of the data:} $\tilde f$ is nonnegative, bounded, and Lipschitz continuous with respect to $(\mu, \tilde a)$ with constant $L_{\tilde f}$, namely for every $(\mu, \tilde a), (\mu', \tilde a') \in \mathfrak{S} \times \tilde{\cA}$, we have
	$$
		| \tilde f( \mu, \tilde a) - \tilde f( \mu' , \tilde a')  | \leq L_{\tilde f} \left( \| \mu - \mu' \|_{d_{\mathfrak{S}}}  +  d_{\tilde{\cA}}( \tilde a , \tilde a')  \right)
		\qquad
		and
		\qquad
		0\leq\tilde f(\mu, \tilde a)\leq C_{\tilde f}, \qquad 0<C_{\tilde f}<\infty.
	$$
	Also, $\bar F$ is Lipschitz continuous with respect to $\mu$ and $\tilde a$ with constant $L_{\bar F}$ in expectation over the randomness of the common noise, namely:  for every $(\mu , \tilde a), (\mu', \tilde a')  \in \mathfrak{S} \times \tilde{\cA}$,
	\begin{align*}
		 & \EE_{\varepsilon^0} \left[ \| \bar F(\mu, \mu \measprod \tilde a, \varepsilon^0) -  \bar F(\mu', \mu' \measprod \tilde a', \varepsilon^0)  \|_{d_{\mathfrak{S}}} \right]
		\le L_{\bar F}  \left( \| \mu - \mu' \|_{d_{\mathfrak{S}}} + d_{\tilde{\cA}}( \tilde a , \tilde a') \right)
	\end{align*}
\end{assumption}
	Notice that any bounded one-stage cost can be made nonnegative by adding a constant $b$; this shifts all policy values and optimal Q-values by $b/(1-\gamma)$ and leaves optimal policies unchanged.
\begin{assumption}\label{hyp:smooth-V}
	\textbf{Regularity of the value function:} $\bar J^*$ is Lipschitz continuous w.r.t. $\mu$ with constant $L_{\bar J^*}$.
\end{assumption}

\begin{assumption}\label{hyp:covering-time}
\textbf{Covering time:} There exists an integer $T_{cov}\geq1$ such that, almost surely, for every integer $n_0\geq0$,
\[
\check{\mathfrak S}\times\check{\cA}
\subseteq\{(\check\mu_n,\check a_n):n_0\leq n<n_0+T_{cov}\},
\]
where $(\check\mu_n,\check a_n)_n$ is the state-action sequence used in Algorithm~\ref{algo:Qtable-projection}.
\end{assumption}

Under Assumption~\ref{hyp:bdd-smooth-data}, the optimal projected state and state-action values belong to $[0,C_{\tilde f}/(1-\gamma)]$.
For the finite-time estimates below, with zero initialization, we fix
\[
\check J_{bound}:=\max\left\{1,\frac{C_{\tilde f}}{1-\gamma}\right\}.
\]
The cost-minimization updates are equivalent to reward-maximization updates for $r=C_{\tilde f}-\tilde f\in[0,C_{\tilde f}]$, using $Q_n^r=C_{\tilde f}/(1-\gamma)-\check Q_n$.
\begin{lemma}\label{lem:lip-value-tabular}
Suppose Assumption~\ref{hyp:bdd-smooth-data} holds and $\gamma L_{\bar F}<1$.
Then $\bar J^*$ is Lipschitz continuous on $\mathfrak S$, and Assumption~\ref{hyp:smooth-V} holds with
\(
L_{\bar J^*}=\frac{L_{\tilde f}}{1-\gamma L_{\bar F}}.
\)
\end{lemma}
\begin{proof}
Set $V_0\equiv0$ on $\mathfrak S$ and, for $m\geq0$, define
\[
V_{m+1}(\mu)=\inf_{\tilde a\in\tilde{\cA}}
\left\{\tilde f(\mu,\tilde a)
+\gamma\EE_{\varepsilon^0}\left[
V_m\left(\bar F(\mu,\mu\measprod\tilde a,\varepsilon^0)\right)
\right]\right\}.
\]
Thus $V_m$ is the optimal discounted cost over $m$ stages for the unprojected problem, with zero terminal cost.
If $V_m$ is Lipschitz, the expression inside braces is Lipschitz in $\mu$, uniformly in $\tilde a$, with constant $L_{\tilde f}+\gamma L_{\bar F}\operatorname{Lip}(V_m)$.
Indeed, use the same action and common-noise realization at the two distributions and apply Assumption~\ref{hyp:bdd-smooth-data}.
Taking the infimum over the common action set preserves this bound, so
\[
\operatorname{Lip}(V_{m+1})
\leq L_{\tilde f}+\gamma L_{\bar F}\operatorname{Lip}(V_m).
\]
Consequently, $\operatorname{Lip}(V_m)\leq L_{\tilde f}\sum_{j=0}^{m-1}(\gamma L_{\bar F})^j$.
Bounded costs give the uniform tail estimate
\(
\|V_m-\bar J^*\|_\infty
\leq \frac{\gamma^m C_{\tilde f}}{1-\gamma}.
\)
Passing to the limit in the preceding Lipschitz estimates therefore gives
\[
\operatorname{Lip}(\bar J^*)\leq\frac{L_{\tilde f}}{1-\gamma L_{\bar F}}.
\]
\end{proof}
Assumption~\ref{hyp:covering-time} is the bounded-coverage assumption used in~\cite{EvenDarMansour}.
With a generative simulator that permits queries at arbitrary state-action pairs, it can be satisfied by fixing an ordering of all discretized pairs and looping through it repeatedly.
In this case, $T_{cov} = |\check{\mathfrak{S}} \times \check{\cA}|$.
In practice, it is common to take $\check{a}_{n}$ using a randomized policy, such as an $\epsilon$-greedy policy, which amounts to taking action $\argmin \check{Q}_n(\check{\mu}_n, \cdot)$ with probability $1-\epsilon$, and otherwise to pick an action uniformly at random in $\check{\cA}$.
In standard Q-learning, this is interpreted as a randomized policy.
In our mean-field context, this is a form of common randomization.
The main advantage of this approach is that, as the Q-function is learned, the sequence of actions is mostly concentrated on the optimal choices instead of spending many iterations on irrelevant actions.
However, when the Q-function is not yet well estimated, using an $\epsilon$-greedy policy with a small $\epsilon$ can lead to poor decisions.
This is an illustration of the exploration and exploitation trade-off.
For online trajectory sampling, $\epsilon$-greedy action selection alone does not imply Assumption~\ref{hyp:covering-time}.
A sufficient condition for a probabilistic coverage estimate is that the finite projected MDP is irreducible under uniformly random actions and that the exploration probabilities are bounded below by a constant $\epsilon_0>0$.
Then there exists an integer $L\geq1$ such that, conditional on the past at any starting time, all discretized state-action pairs are visited within the next $L$ updates with probability at least $1/2$.
This gives a high-probability covering bound over a finite run, rather than a deterministic covering bound over all times; see~\cite[Lemma~33 and Corollary~34]{EvenDarMansour}.
The quantitative bounds below are stated under Assumption~\ref{hyp:covering-time}.

We consider projection operators $\proj_{\check{\mathfrak{S}}}: \mathfrak{S} \to \check{\mathfrak{S}}$ and $\proj_{\check{\cA}}: \tilde{\cA} \to \check{\cA}$ such that 
$$
(\proj_{\check{\mathfrak{S}}}(\mu), \proj_{\check{\cA}}(\tilde a) ) := (\check \mu, \check{a})
$$ 
for every $(\mu, \tilde a) \in \mathfrak{S} \times \tilde{\cA}$ where $(\check \mu, \check{a})$ is the closest point (or one of the closest points, in case of ties) in $\check{\mathfrak{S}} \times \check{\cA}$ with respect to $d_{\mathfrak{S}}$ and $d_{\tilde{\cA}}$.
Based on simplex discretizations, this point  satisfies $\| \mu - \check \mu\|_{d_{\mathfrak{S}}} \leq \varepsilon_{\mathfrak{S}}$ and $d_{\tilde{\cA}}( \tilde a , \check{a} ) \leq \varepsilon_{\mathfrak{A}}$.
We set $\beta = (1-\gamma)/2$. %

\begin{theorem}
	\label{th:main-cv-tabular}
	Let $\delta \in (0,1)$ and $0<\varepsilon<\check J_{bound}$ and let us assume Assumptions~\ref{hyp:bdd-smooth-data}--\ref{hyp:covering-time} hold.
	Consider learning rates $(\eta_n)_n$ such that there exists $\kappa \in (1/2,1)$ such that for every $(\check \mu, \check{a}) \in \check{\mathfrak{S}} \times \check{\cA}$, $\eta_n := \eta_n(\check \mu, \check{a}) = 1/ \big ( 1 + C(n, \check \mu, \check{a}) \big)^\kappa$ for each $n \geq 0$, where $C(n, \check \mu, \check{a})$ is the number of times up to $n$ that the pair $(\check \mu, \check{a})$ has been visited in Algorithm~\ref{algo:Qtable-projection}.
	If the number of updates $N_{\mathrm{upd}}$ is of order
		{\small
			\begin{equation}
				\label{eq:LB-Nepi-tabular}
				\Omega\left(
				\left(\frac{(T_{cov})^{1+3 \kappa} \check J_{bound}^2 \, \ln \left(|\check{\mathfrak{S}}| \, |\check{\mathfrak{A}}|^{|S|} \check J_{bound} / (\delta \beta \varepsilon) \right)}{\beta^2 \varepsilon^2}\right)^{\frac{1}{\kappa}}
				\hskip -6pt	+
				\left(\frac{(T_{cov})}{\beta} \ln \left( \frac{\check J_{bound}}{\varepsilon} \right)\right)^{\frac{1}{1-\kappa}} \right),
			\end{equation}
		}
	then with probability $1-\delta$,
	for all $(\mu, \tilde a) \in \mathfrak{S} \times \tilde{\cA}$,
	$$
		\left| \check Q_{N_{\mathrm{upd}}} \Big(\proj_{\check{\mathfrak{S}}} (\mu), \proj_{\check{\cA}}(\tilde{a}) \Big) - \bar Q^*(\mu,  \mu \measprod \tilde a) \right| \le \varepsilon',
	$$
	where
	$ \displaystyle
		\varepsilon'
		=
		\varepsilon
		+ \left( \frac{\gamma}{1 - \gamma} L_{\bar J^*} + L_{\tilde f} + \gamma L_{\bar J^*} L_{\bar F} \right) \varepsilon_{\mathfrak{S}} + \frac{1}{1- \gamma} \left( L_{\tilde f} + \gamma L_{\bar J^*} L_{\bar F} \right) \varepsilon_{\mathfrak{A}}.
	$
\end{theorem}

As an example, if we cover the simplexes with balls of radius $\varepsilon_{\mathfrak{S}}$ and $\varepsilon_{\mathfrak{A}}$ respectively, the number of points in the covering would be of order $\Theta(\varepsilon_{\mathfrak{S}}^{1-|S|})$ and $\Theta(\varepsilon_{\mathfrak{A}}^{1-|A|})$ respectively.
For these coverings, the mixed-control Q-table has
\[
N_{\mathrm{pair}}
:=|\check{\mathfrak S}|\,|\check{\mathfrak A}|^{|S|}
=\Theta\!\left(
\varepsilon_{\mathfrak S}^{-(|S|-1)}
\varepsilon_{\mathfrak A}^{-|S|(|A|-1)}
\right)
\]
entries, for fixed $S$ and $A$ as the mesh radii tend to zero.
In particular, the action-table factor $|\check{\mathfrak A}|^{|S|}$ is exponential in $|S|$ whenever $|\check{\mathfrak A}|>1$; this is a curse of dimensionality of the tabular method.
Under cyclic sampling, $T_{cov}=N_{\mathrm{pair}}$.
Even one complete sweep requires $N_{\mathrm{pair}}$ updates.
For pure level-0 controls, the corresponding table has $|\check{\mathfrak S}|\,|A|^{|S|}$ entries.
The update bound is therefore obtained from~\eqref{eq:LB-Nepi-tabular} by substituting $T_{cov}=N_{\mathrm{pair}}$,
while the error bound $\varepsilon'$ satisfies:
$
	 \varepsilon' \leq \,
	\varepsilon + C (\varepsilon_{\mathfrak{S}} + \varepsilon_{\mathfrak{A}}).
$

\begin{remark}
	\label{ref:holderJ}
	If $\bar J^*$ is only $p-$H\"older continuous, then the proof can be adapted in a straightforward way and the same result holds except that $\varepsilon'$ becomes:
	$$    	    \varepsilon
		+ \frac{\gamma }{1-\gamma} \left( L_{\tilde f} \varepsilon_{\mathfrak{A}} + L_{\bar J^*}(\gamma  (L_{\bar F} \varepsilon_{\mathfrak{A}})^p +  \varepsilon_{\mathfrak{S}}^p ) \right)
		+ \left( L_{\tilde f} (\varepsilon_{\mathfrak{S}} + \varepsilon_{\mathfrak{A}}) + \gamma L_{\bar J^*} \big(L_{\bar F} (\varepsilon_{\mathfrak{S}} + \varepsilon_{\mathfrak{A}}) \big)^p \right).
	$$
\end{remark}

Note that $\varepsilon$ can be chosen as small as desired provided $N_{\mathrm{upd}}$ is large enough.
The second and third terms in the error $\varepsilon'$ are proportional to $\varepsilon_{\mathfrak{S}}$ and $\varepsilon_{\mathfrak{A}}$, a result that is generally unavoidable due to the projection onto the finite sets $\check{\mathfrak{S}}$ and $\check{\mathfrak{A}}$.
However, this error vanishes as $\varepsilon_{\mathfrak{S}} \to 0$ and $\varepsilon_{\mathfrak{A}} \to 0$, i.e., as $\check{\mathfrak{S}}$ and $\check{\mathfrak{A}}$ become increasingly accurate approximations of $\cP(S)$ and $\cP(A)$, respectively.

We prove this result below.
The proof can be summarized in the following three steps:
\textbf{(1)} For $N_{\mathrm{upd}}$ large enough, we have
$
	\check Q_{N_{\mathrm{upd}}} \approx \check Q^*
$
on $\check{\mathfrak{S}} \times \check{\cA}$;
\textbf{(2)}
$
	\check Q^* \approx \tilde Q^*
$
on $\check{\mathfrak{S}} \times \check{\cA}$;
\textbf{(3)} For every $(\mu, \tilde a) \in \mathfrak{S} \times \tilde{\cA}$,
$\tilde Q^*( \proj_{\check{\mathfrak{S}}}(\mu), \proj_{\check{\cA} }(\tilde a ) ) \approx \tilde Q^*(\mu, \tilde a)$.
The first step leverages standard Q-learning convergence results~\citep{EvenDarMansour}, while the other two steps arise from the regularity assumptions and the approximation of $(\mathfrak{S}, \tilde{\cA})$ by $(\check{\mathfrak{S}}, \check{\cA})$.

\begin{proof}[Proof of Theorem~\ref{th:main-cv-tabular}]

	Recall that we denote by $\check J^*$ and $\check Q^*$ respectively the state value function and the state-action value function of the projected MFC problem defined by~\eqref{eq:generic-MFC-fctmeasure-cost-proj}--\eqref{eq:generic-MFC-fctmeasure-dyn-proj}. We first note that, for every $(\mu, \tilde a) \in \mathfrak{S} \times \tilde{\cA}$,
	\begin{align*}
		     & \left| \check Q_{N_{\mathrm{upd}}} \big( \proj_{\check{\mathfrak{S}}}(\mu), \proj_{\check{\cA}}(\tilde a) \big)  - \tilde Q^* \big( \mu, \tilde a \big) \right|
		\\
		 &\hskip 35pt
         \le \left| \check Q_{N_{\mathrm{upd}}} \big( \proj_{\check{\mathfrak{S}}}(\mu), \proj_{\check{\cA}}(\tilde a) \big) - \check Q^* \big( \proj_{\check{\mathfrak{S}}}(\mu), \proj_{\check{\cA}}(\tilde a) \big) \right|
		\\
		     &\hskip 35pt \quad + \left|  \check Q^* \big( \proj_{\check{\mathfrak{S}}}(\mu), \proj_{\check{\cA}}(\tilde a) \big) - \tilde Q^* \big(  \proj_{\check{\mathfrak{S}}}(\mu), \proj_{\check{\cA}}(\tilde a) \big) \right|
		\\
		     & \hskip 35pt\quad + \left| \tilde Q^* \big(  \proj_{\check{\mathfrak{S}}}(\mu), \proj_{\check{\cA}}(\tilde a) \big) - \tilde Q^* \big( \mu, \tilde a \big)  \right|.
	\end{align*}
	We split the proof into three steps, which involve bounding each term on the right-hand side.

	\vskip 6pt\noindent
	\textbf{Step 1.} We first analyze the difference between $\check Q_{N_{\mathrm{upd}}}$ and $\check Q^*$.
	The estimate follows from standard convergence results for Q-learning in finite state-action spaces.
	More precisely, under Assumptions~\ref{hyp:bdd-smooth-data} and~\ref{hyp:covering-time}, with our choice of learning rates, and given that $N_{\mathrm{upd}}$ is of order~\eqref{eq:LB-Nepi-tabular}, we can apply Theorem~4 in~\cite{EvenDarMansour} for asynchronous Q-learning and polynomial learning rates, and we obtain that, with probability at least $1-\delta$,
	$$
		\|\check Q_{N_{\mathrm{upd}}} - \check Q^* \|_{\infty}= \sup_{ (\check \mu, \check{a}) \in \check{\mathfrak{S}} \times \check{\cA}} \Big| \check Q_{N_{\mathrm{upd}}}(\check \mu, \check{a}) - \check Q^*(\check \mu, \check{a}) \Big| \le \varepsilon.
	$$
\textbf{Step 2. } We then turn our attention to the difference between $\check Q^*$ and $\tilde Q^*$.
	The analysis demonstrates that the projection onto $\check{\mathfrak{S}}$ performed at each step does not significantly perturb the value function.
	Recall that for some given common noise $\varepsilon^0$, the operator $\check \Phi^{\varepsilon^0}: \check{\mathfrak{S}} \times \check{\cA} \to \check{\mathfrak{S}}$ is given by $\check \Phi^{\varepsilon^0}( \check \mu, \check{a}) = \proj_{\check{\mathfrak{S}}} \circ \bar F ( \check \mu, \check \mu \measprod \check{a}, \varepsilon^0)$.
	Likewise, we denote the transition dynamics with $\bar F$ by a function $\Phi^{\varepsilon^0}: \mathfrak{S} \times \tilde{\cA} \to \mathfrak{S}$ such that:
	$$
		\Phi^{\varepsilon^0} (\mu, \tilde a) = \bar F( \mu, \mu \measprod \tilde a, \varepsilon^0), \qquad \forall (\mu, \tilde a) \in \mathfrak{S} \times \tilde{\cA}.
	$$
	Let us start by noting that, for every $(\check \mu, \check{a}) \in \check{\mathfrak{S}} \times \check{\cA}$,
	\begin{align*}
\left| \check Q^*(\check \mu, \check{a}) -  \tilde Q^*(\check \mu, \check{a}) \right|
		 & \le  \gamma    \EE \Bigg[ \Big| \check J^*(\check \Phi^{\varepsilon^0}(\check \mu, \check{a})) - \bar J^*( \Phi^{\varepsilon^0}(\check \mu, \check{a}) ) \Big| \Bigg]
		\\
		 & \le   \gamma   \EE \Bigg[  \Big|  \check J^*( \check \Phi^{\varepsilon^0}(\check \mu, \check{a}) ) - \bar J^*( \check \Phi^{\varepsilon^0}(\check \mu, \check{a})) \Big|  + \Big| \bar J^* (\check \Phi^{\varepsilon^0}(\check \mu, \check{a})) - \bar J^*( \Phi^{\varepsilon^0}(\check \mu, \check{a})) \Big|  \Bigg]
		\\
		 & \le \gamma   \EE \Bigg[  \Big|  \inf_{\check{a}' \in \check{\cA}} \check Q^* \left(\check \Phi^{\varepsilon^0}(\check \mu, \check{a}), \check{a}' \right) - \inf_{\tilde a' \in \tilde{\cA} } \tilde Q^* \left(\check \Phi^{\varepsilon^0}(\check \mu, \check{a}), \tilde a' \right) \Big| \Bigg]
		\\
		 & \hskip 125pt +  \gamma L_{\bar J ^*} \EE \left[  \left\Vert \check \Phi^{\varepsilon^0}(\check \mu, \check{a}) -  \Phi^{\varepsilon^0}(\check \mu, \check{a}) \right\Vert_{d_{\mathfrak{S}}} \right],
	\end{align*}
	where the last inequality holds because of the Lipschitz continuity of $\bar J^*$ on $\mathfrak{S}$, see Assumption~\ref{hyp:smooth-V}.
	The second term in the last inequality can be bounded using the simplex discretization properties and Assumption~\ref{hyp:bdd-smooth-data}:
	\begin{align*}
		\EE \left[  \| \check \Phi^{\varepsilon^0}(\check \mu, \check{a}) - \Phi^{\varepsilon^0}( \check \mu, \check{a}) \|_{d_{\mathfrak{S}}}  \right]
		 & = \EE\left[ \| \proj_{\check{\mathfrak{S}}} \circ \bar F(\check \mu, \check \mu \measprod\check{a}, \varepsilon^0) - \bar F(\check \mu, \check \mu \measprod\check{a}, \varepsilon^0)  \|_{d_{\mathfrak{S}}} \right] \leq \varepsilon_{\mathfrak{S}}.
	\end{align*}
	For the first term, let $\check \mu' = \check \Phi^{\varepsilon^0}(\check \mu, \check{a}) \in \check{\mathfrak{S}}$ to simplify the notation, and let us consider $\check{a}^*_1 \in \check{\cA}$ and $\tilde a^*_2 \in \tilde{\cA}$ satisfying:
	$
		\check Q^*(\check \mu', \check{a}^*_1) = \inf_{\check{a}' \in \check{\cA}} \check Q^* \left(\check \mu', \check{a}' \right)$
	and
	$
		\tilde Q^*(\check \mu', \tilde a^*_2) = \inf_{\tilde a' \in \tilde{\cA} } \tilde Q^* \left(\check \mu', \tilde a' \right).
	$
	The existence of $\check{a}^*_1$ follows from the finiteness of $\check{\cA}$, and the existence of $\tilde a^*_2$ follows from compactness of $\tilde{\cA}$ together with Lemma~\ref{le:Qstar_is_lsc}.
	We observe that
	\begin{align*}
		 & \hskip -25pt
         \check Q^*(\check \mu', \check{a}_1^*) - \tilde Q^*( \check \mu', \tilde a_2^*)
		\\
		 & = \Big( \check Q^*(\check \mu', \check{a}_1^*) - \check Q^*( \check \mu', \proj_{\check{\cA}} ( \tilde a_2^*) )  \Big) + \Big( \check Q^*( \check \mu', \proj_{\check{\cA}} ( \tilde a_2^*) ) - \tilde Q^*( \check \mu', \proj_{\check{\cA}}(\tilde a_2^*) ) \Big)
		\\
		 & \qquad\qquad + \Big( \tilde Q^*( \check \mu', \proj_{\check{\cA}}( \tilde a_2^*) ) - \tilde Q^*( \check \mu', \tilde a_2^*)  \Big)
		\\
		 & \leq 0 +  \sup_{ (\check \mu, \check{a}) \in \check{\mathfrak{S}} \times \check{\cA}} \left| (\check Q^* - \tilde Q^*)(\check \mu, \check{a}) \right|
		\\
		 & \qquad\qquad + \left( \tilde f (\check \mu', \proj_{\check{\cA}}(\tilde a^*_2) ) + \gamma \EE_{(\varepsilon^0)'} \left[ \bar J^*( \bar F( \check \mu', \check \mu' \measprod \proj_{\check{\cA}}(\tilde a_2^*), (\varepsilon^0)' ) ) \right]  \right)
		\\
		 & \qquad \qquad - \left( \tilde f (\check \mu', \tilde a^*_2 ) + \gamma \EE_{(\varepsilon^0)'} \left[ \bar J^*( \bar F( \check \mu', \check \mu' \measprod \tilde a_2^*, (\varepsilon^0)' ) ) \right] \right)
		\\
		 & \leq  \| \check Q^* - \tilde Q^* \|_{\infty}  + ( L_{\tilde f} + \gamma L_{\bar J^*} L_{\bar F} ) \varepsilon_{\mathfrak{A}}.
	\end{align*}

	On the other hand,
	\begin{align*}
		\check Q^*( \check \mu', \check{a}_1^*) -  \tilde Q^*(\check \mu', \tilde a_2^*)
		 & = - \Big( \tilde Q^*( \check \mu' ,  \tilde a_2^* ) - \tilde Q^*( \check \mu', \check{a}_1^*) \Big) - \Big( \tilde Q^*( \check \mu' , \check{a}_1^*) - \check Q^*( \check \mu', \check{a}_1^* ) \Big) \\
		 & \geq - \| \check Q^* - \tilde Q^* \|_{\infty}.
	\end{align*}

	Combining the above bounds yields that for every $(\check \mu, \check{a}) \in \check{\mathfrak{S}} \times \check{\cA}$,
	$$
		\left| \check Q^*(\check \mu, \check{a}) - \tilde Q^*(\check \mu, \check{a}) \right| \leq \gamma \left( \| \check Q^* - \tilde Q^* \|_{\infty} + (L_{\tilde f} + \gamma L_{\bar J^*} L_{\bar F} ) \varepsilon_{\mathfrak{A}} \right) + \gamma L_{\bar J^*} \varepsilon_{\mathfrak{S}}.
	$$
	Consequently,
	$$
		\|\check Q^* - \tilde Q^* \|_\infty
		\le
		\frac{\gamma }{1-\gamma} \left( ( L_{\tilde f} + \gamma L_{\bar J^*} L_{\bar F} )\varepsilon_{\mathfrak{A}} + L_{\bar J^*} \varepsilon_{\mathfrak{S}} \right).
	$$
\textbf{Step 3. } Last, we look at the difference between $\tilde Q^*( \proj_{\check{\mathfrak{S}}}(\mu), \proj_{\check{\cA}}(\tilde a )  )$ and $\tilde Q^*(\mu, \tilde a)$.
	For every $\mu \in \mathfrak{S}$ and $\tilde a \in \tilde{\cA}$, letting $\check \mu = \proj_{\check{\mathfrak{S}}} (\mu)$ and $\check{a} = \proj_{\check{\cA}}(\tilde a)$ to simplify the notation, we have $\| \check \mu - \mu \|_{d_{\mathfrak{S}}} \leq \varepsilon_{\mathfrak{S}}$ and $\| \check{a} - \tilde a\|_{d_{\tilde{\cA}}} \leq \varepsilon_{\mathfrak{A}}$.
	We obtain
	\begin{align*}
		 & \left| \tilde Q^* (\check \mu, \check{a}) -  \tilde Q^*(\mu, \tilde a) \right|
		\\
		 & \qquad \le \Bigg|  \tilde f(\check \mu, \check{a})  -  \tilde f(\mu, \tilde a) \Bigg| +  \gamma  \EE \Bigg[  \left| \inf_{\tilde a' \in \tilde{\cA}} \tilde Q^* ( \Phi(\check\mu, \check{a} ), \tilde a') - \inf_{\tilde a' \in \tilde{\cA}}  \tilde Q^*(\Phi(\mu, \tilde a), \tilde a') \right| \Bigg]
		\\
		 & \qquad  \le L_{\tilde f} \left( \| \check \mu - \mu \|_{d_{\mathfrak{S}}} +  \| \check{a} - \tilde a\|_{d_{\tilde{\cA}}} \right)
		+  \gamma  \EE \left[  \left| \bar J^* (\bar F (\check \mu, \check \mu \measprod \check a, \varepsilon^0)) - \bar J^* (\bar F (\mu, \mu \measprod \tilde a, \varepsilon^0)) \right| \right]
		\\
		 & \qquad  \le L_{\tilde f} (\varepsilon_{\mathfrak{S}} + \varepsilon_{\mathfrak{A}} ) + \gamma L_{\bar J^*} \EE \left[ \| \bar F(\check \mu, \check \mu \measprod \check a, \varepsilon^0 ) - \bar F( \mu, \mu \measprod \tilde a, \varepsilon^0) \|_{d_{\mathfrak{S}}} \right]
		\\
		 & \qquad  \le  ( L_{\tilde f} + \gamma L_{\bar J^*} L_{\bar F} )(\varepsilon_{\mathfrak{S}} + \varepsilon_{\mathfrak{A}} ),
	\end{align*}
	where we used the Lipschitz continuity of $\tilde f, \bar J^*, \bar F$ and the assumption on $\check{\mathfrak{S}}$, see Assumptions~\ref{hyp:bdd-smooth-data}, \ref{hyp:smooth-V} and the simplex discretization properties.
\qed\end{proof}

\subsubsection{Q-learning with pure level-0 controls}
The approach described above is designed for cases where we seek optimal actions that may be randomized at the individual level.
Searching in the space $\cP(A)$ comes with a computational cost that is reflected in the bounds through the cardinality of the discrete simplex $\check{\mathfrak{A}}$.
However, in some situations it may be interesting to directly search for actions that are pure at the individual level.
In this paragraph, we restrict our attention to pure level-0 policies, i.e., without idiosyncratic randomization at the agent level.
For such policies, given the recommendation of the central planner (which might depend on a common randomization), each agent chooses their action in a deterministic way, based purely on their own state.
This amounts to removing the random variables $\vartheta_n$ from the probabilistic framework presented in Section~\ref{se:proba-framework}.
For simplicity we will use the same notations for the value functions $\bar{J}^*$, $\bar{Q}^*$ and $\check{Q}$.
The goal is now to approximate $\bar{Q}^*$, which is the Q-function optimized over the set of policies that are pure at the level $0$.
In this case, instead of~\eqref{eq:def-calA-S-PA}, the set of strategy functions is (for simplicity we keep the notation $\tilde{\cA}$):
\begin{equation*}
	\tilde{\cA} := \{ \tilde a: S \to A \} = A^S.
\end{equation*}
In Algorithm~\ref{algo:Qtable-projection}, we replace $\check{a} \in \check{\cA}$ by $\tilde a \in \tilde{\cA}$.

\begin{theorem}
	\label{th:main-cv-tabular-pure}
	Let $\delta \in (0,1)$ and $0<\varepsilon<\check J_{bound}$ and let us assume Assumptions~\ref{hyp:bdd-smooth-data}--\ref{hyp:covering-time} hold.
	Consider learning rates $(\eta_n)_n$ such that there exists $\kappa \in (1/2,1)$ such that for every $(\check \mu, \tilde a) \in \check{\mathfrak{S}} \times \tilde{\cA}$, 
    $$
    \eta_n := \eta_n(\check \mu, \tilde a) = 1/ \big ( 1 + C(n, \check \mu, \tilde a) \big)^\kappa,\qquad  n \geq 0,
    $$
    where $C(n, \check \mu, \tilde a)$ is the number of times up to time $n$ that the pair $(\check \mu, \tilde a)$ has been visited in Algorithm~\ref{algo:Qtable-projection}.
	If the number of updates $N_{\mathrm{upd}}$ is of order
		{\small
			\begin{equation}
				\label{eq:LB-Nepi-tabular-pure}
				\Omega\left(
				\left(\frac{(T_{cov})^{1+3 \kappa} \check J_{bound}^2 \, \ln \left(|\check{\mathfrak{S}}| \, |A|^{|S|} \check J_{bound} / (\delta \beta \varepsilon) \right)}{\beta^2 \varepsilon^2}\right)^{\frac{1}{\kappa}}
				\hskip -6pt +
				\left(\frac{(T_{cov})}{\beta} \ln \left( \frac{\check J_{bound}}{\varepsilon} \right)\right)^{\frac{1}{1-\kappa}} \right),
			\end{equation}
		}
	then with probability $1-\delta$,
	for all $(\mu, \tilde a) \in \mathfrak{S} \times \tilde{\cA}$,
	$$
		\left| \check Q_{N_{\mathrm{upd}}} \Big(\proj_{\check{\mathfrak{S}}} (\mu), \tilde{a} \Big) - \bar Q^*(\mu,  \mu \measprod \delta_{\tilde a}) \right| \le \varepsilon',
	$$
	where
	$ \displaystyle
		\varepsilon'
		=
		\varepsilon
		+ \left( \frac{\gamma}{1 - \gamma} L_{\bar J^*} + L_{\tilde f} + \gamma L_{\bar J^*} L_{\bar F} \right) \varepsilon_{\mathfrak{S}}.
	$
\end{theorem}
A similar result can be obtained if $\bar J^*$ is only H\"older continuous, consistent with Remark~\ref{ref:holderJ}.

\begin{proof}[Proof of Theorem~\ref{th:main-cv-tabular-pure}]

	Recall that we denote by $\check J^*$ and $\check Q^*$ the state value function and the state-action value function of the projected MFC problem defined by~\eqref{eq:generic-MFC-fctmeasure-cost-proj}--\eqref{eq:generic-MFC-fctmeasure-dyn-proj}.	We first note that, for every $(\mu, \tilde a) \in \mathfrak{S} \times \tilde{\cA}$,
	\begin{align*}
		\left| \check Q_{N_{\mathrm{upd}}} \big( \proj_{\check{\mathfrak{S}}}(\mu), \tilde a \big)  - \tilde Q^* \big( \mu, \tilde a \big) \right|
		\leq & \left| \check Q_{N_{\mathrm{upd}}} \big( \proj_{\check{\mathfrak{S}}}(\mu), \tilde a \big) - \check Q^* \big( \proj_{\check{\mathfrak{S}}}(\mu), \tilde a \big) \right|
		\\
		     & \quad + \left|  \check Q^* \big( \proj_{\check{\mathfrak{S}}}(\mu), \tilde a \big) - \tilde Q^* \big(  \proj_{\check{\mathfrak{S}}}(\mu), \tilde a \big) \right|
		\\
		     & \quad + \left| \tilde Q^* \big(  \proj_{\check{\mathfrak{S}}}(\mu), \tilde a \big) - \tilde Q^* \big( \mu, \tilde a \big)  \right|.
	\end{align*}
	We split the proof into three steps to bound separately each term on the right-hand side.

	\vskip 6pt\noindent
	\textbf{Step 1.} We first analyze the difference between $\check Q_{N_{\mathrm{upd}}}$ and $\check Q^*$.
	As before, this step follows from standard convergence results for Q-learning in finite state-action spaces.
	More precisely, under Assumptions~\ref{hyp:bdd-smooth-data} and \ref{hyp:covering-time}, with our choice of learning rates, and given that $N_{\mathrm{upd}}$ is of order~\eqref{eq:LB-Nepi-tabular-pure}, we can apply Theorem~4 in~\cite{EvenDarMansour} for asynchronous Q-learning and polynomial learning rates, and we obtain that, with probability at least $1-\delta$,
	$$
		\|\check Q_{N_{\mathrm{upd}}} - \check Q^* \|_{\infty}= \sup_{ (\check \mu, \tilde a) \in \check{\mathfrak{S}} \times \tilde{\cA}} \Big| \check Q_{N_{\mathrm{upd}}}(\check \mu, \tilde a) - \check Q^*(\check \mu, \tilde a) \Big| \le \varepsilon.
	$$
\textbf{Step 2. } We then turn our attention to the difference between $\check Q^*$ and $\tilde Q^*$.
	The analysis demonstrates that the projection onto $\check{\mathfrak{S}}$ performed at each step does not significantly perturb the value function.
	Recall that for some given common noise $\varepsilon^0$, the operator $\check \Phi^{\varepsilon^0}: \check{\mathfrak{S}} \times \tilde{\cA} \to \check{\mathfrak{S}}$ is given by $\check \Phi^{\varepsilon^0}( \check \mu, \tilde a) = \proj_{\check{\mathfrak{S}}} \circ \bar F ( \check \mu, \check \mu \measprod \tilde a, \varepsilon^0)$.
	Likewise, we denote the transition dynamics with $\bar F$ by a function $\Phi^{\varepsilon^0}: \mathfrak{S} \times \tilde{\cA} \to \mathfrak{S}$ such that:
	$$
		\Phi^{\varepsilon^0} (\mu, \tilde a) = \bar F( \mu, \mu \measprod \tilde a, \varepsilon^0), \qquad \forall (\mu, \tilde a) \in \mathfrak{S} \times \tilde{\cA}.
	$$
	Let us start by noting that, for every $(\check \mu, \tilde a) \in \check{\mathfrak{S}} \times \tilde{\cA}$,
	\begin{align*}
		 & \left| \check Q^*(\check \mu, \tilde a) -  \tilde Q^*(\check \mu, \tilde a) \right|
		\\
		 & \qquad  \le  \gamma    \EE \Bigg[ \Big| \check J^*(\check \Phi^{\varepsilon^0}(\check \mu, \tilde a)) - \bar J^*( \Phi^{\varepsilon^0}(\check \mu, \tilde a) ) \Big| \Bigg]
		\\
		 & \qquad  \le   \gamma   \EE \Bigg[  \Big|  \check J^*( \check \Phi^{\varepsilon^0}(\check \mu, \tilde a) ) - \bar J^*( \check \Phi^{\varepsilon^0}(\check \mu, \tilde a)) \Big|  + \Big| \bar J^* (\check \Phi^{\varepsilon^0}(\check \mu, \tilde a)) - \bar J^*( \Phi^{\varepsilon^0}(\check \mu, \tilde a)) \Big|  \Bigg]
		\\
		 & \qquad  \le \gamma   \EE \Bigg[  \Big|  \inf_{\tilde a' \in \tilde{\cA}} \check Q^* \left(\check \Phi^{\varepsilon^0}(\check \mu, \tilde a), \tilde a' \right) - \inf_{\tilde a' \in \tilde{\cA} } \tilde Q^* \left(\check \Phi^{\varepsilon^0}(\check \mu, \tilde a), \tilde a' \right) \Big| \Bigg]
		\\
		 & \qquad \hskip 75pt  +  \gamma L_{\bar J ^*} \EE \left[  \left\Vert \check \Phi^{\varepsilon^0}(\check \mu, \tilde a) -  \Phi^{\varepsilon^0}(\check \mu, \tilde a) \right\Vert_{d_{\mathfrak{S}}} \right],
	\end{align*}
	where the last inequality holds by Lipschitz continuity of $\bar J^*$ on $\mathfrak{S}$, see Assumption~\ref{hyp:smooth-V}.

	The second term in the last inequality can be bounded using the simplex discretization properties and Assumption~\ref{hyp:bdd-smooth-data}:
	\begin{align*}
		\EE \left[  \| \check \Phi^{\varepsilon^0}(\check \mu, \tilde a) - \Phi^{\varepsilon^0}( \check \mu, \tilde a) \|_{d_{\mathfrak{S}}}  \right]
		 & = \EE_{\varepsilon^0_1}\left[ \| \proj_{\check{\mathfrak{S}}} \circ \bar F(\check \mu, \check \mu \measprod \tilde a, \varepsilon^0) - \bar F(\check \mu, \check \mu \measprod \tilde a, \varepsilon^0)  \|_{d_{\mathfrak{S}}} \right] \leq \varepsilon_{\mathfrak{S}}.
	\end{align*} 

	For the first term, let $\check \mu' = \check \Phi^{\varepsilon^0}(\check \mu, \tilde a) \in \check{\mathfrak{S}}$ to alleviate the notation, and let us consider $\tilde a^*_1 \in \tilde{\cA}$ and $\tilde a^*_2 \in \tilde{\cA}$ satisfying:
	$$
		\check Q^*(\check \mu', \tilde a^*_1) = \inf_{\tilde a' \in \tilde{\cA}} \check Q^* \left(\check \mu', \tilde a' \right)
		\qquad
		\text{ and }
		\qquad
		\tilde Q^*(\check \mu', \tilde a^*_2) = \inf_{\tilde a' \in \tilde{\cA} } \tilde Q^* \left(\check \mu', \tilde a' \right).
	$$
	Since the pure-action set $\tilde{\cA}=A^S$ is finite, the minima defining $\tilde a_1^*$ and $\tilde a_2^*$ are attained.
	We observe that
	\begin{align*}
		 & \check Q^*(\check \mu', \tilde a_1^*) - \tilde Q^*( \check \mu', \tilde a_2^*)
		\\
		 & = \Big( \check Q^*(\check \mu', \tilde a_1^*) - \check Q^*( \check \mu',  \tilde a_2^* )  \Big) + \Big( \check Q^*( \check \mu',  \tilde a_2^* ) - \tilde Q^*( \check \mu', \tilde a_2^* ) \Big)
		\\
		 & \leq 0 +  \sup_{ (\check \mu, \tilde a) \in \check{\mathfrak{S}} \times \tilde{\cA}} \left| (\check Q^* - \tilde Q^*)(\check \mu, \tilde a) \right|                                                  \\
		 & \leq  \| \check Q^* - \tilde Q^* \|_{\infty} .
	\end{align*}
	On the other hand,
	\begin{align*}
		\check Q^*( \check \mu', \tilde a_1^*) -  \tilde Q^*(\check \mu', \tilde a_2^*)
		 & = - \Big( \tilde Q^*( \check \mu' ,  \tilde a_2^* ) - \tilde Q^*( \check \mu', \tilde a_1^*) \Big) - \Big( \tilde Q^*( \check \mu' , \tilde a_1^*) - \check Q^*( \check \mu', \tilde a_1^* ) \Big) \\
		 & \geq - \| \check Q^* - \tilde Q^* \|_{\infty}.
	\end{align*}

	Combining the above bounds yields that for every $(\check \mu, \tilde a) \in \check{\mathfrak{S}} \times \tilde{\cA}$,
	$$
		\left| \check Q^*(\check \mu, \tilde a) - \tilde Q^*(\check \mu, \tilde a) \right| \leq \gamma \left( \| \check Q^* - \tilde Q^* \|_{\infty}   \right) + \gamma L_{\bar J^*} \varepsilon_{\mathfrak{S}}.
	$$
	Consequently,
	$$
		\|\check Q^* - \tilde Q^* \|_\infty
		\le
		\frac{\gamma }{1-\gamma} L_{\bar J^*} \varepsilon_{\mathfrak{S}}.
	$$
\textbf{Step 3. } Last, we look at the difference between $\tilde Q^*( \proj_{\check{\mathfrak{S}}}(\mu), \tilde a  )$ and $\tilde Q^*(\mu, \tilde a)$.
	For every $\mu \in \mathfrak{S}$ and $\tilde a \in \tilde{\cA}$, letting $\check \mu = \proj_{\check{\mathfrak{S}}} (\mu)$  to alleviate the notation, we have $\| \check \mu - \mu \|_{d_{\mathfrak{S}}} \leq \varepsilon_{\mathfrak{S}}$.
	We obtain
	\begin{align*}
		 & \hskip -25pt\left| \tilde Q^* (\check \mu, \tilde a) -  \tilde Q^*(\mu, \tilde a) \right|
		\\
		 & \le \Bigg|  \tilde f(\check \mu, \tilde a)  -  \tilde f(\mu, \tilde a) \Bigg| +  \gamma  \EE \Bigg[  \left| \inf_{\tilde a' \in \tilde{\cA}} \tilde Q^* ( \Phi(\check\mu, \tilde a ), \tilde a') - \inf_{\tilde a' \in \tilde{\cA}}  \tilde Q^*(\Phi(\mu, \tilde a), \tilde a') \right| \Bigg]
		\\
		 & \le L_{\tilde f} \| \check \mu - \mu \|_{d_{\mathfrak{S}}}
		+  \gamma  \EE \left[  \left| \bar J^* (\bar F (\check \mu, \check \mu \measprod \tilde a, \varepsilon^0)) - \bar J^* (\bar F (\mu, \mu \measprod \tilde a, \varepsilon^0)) \right| \right]
		\\
		 & \le L_{\tilde f} \varepsilon_{\mathfrak{S}}  + \gamma L_{\bar J^*} \EE \left[ \| \bar F(\check \mu, \check \mu \measprod \tilde a, \varepsilon^0 ) - \bar F( \mu, \mu \measprod \tilde a, \varepsilon^0) \|_{d_{\mathfrak{S}}} \right]
		\\
		 & \le  ( L_{\tilde f} + \gamma L_{\bar J^*} L_{\bar F} )\varepsilon_{\mathfrak{S}}  ,
	\end{align*}
	where we used the Lipschitz continuity of $\tilde f, \bar J^*, \bar F$ and the assumption on $\check{\mathfrak{S}}$, see Assumptions~\ref{hyp:bdd-smooth-data}, \ref{hyp:smooth-V} and the simplex discretization properties.
\qed\end{proof}

\vskip 6pt
The result of the above theorem provides convergence guarantees for the Q-function.
Let us now derive a consequence in terms of the optimizer.
To this end, we employ the following additional assumption regarding the gap between the values of the best and second-best actions, a standard requirement in approximation algorithms based on tabular Q-functions~\citep{Farahmand,Bellemare}.

\begin{assumption}\label{hyp:action-gap} \textbf{Action gap:} There exists $K_A > 0$ such that:
	$$
		\tilde Q^*( \check \mu, \tilde a) - \inf_{\tilde a' \in \tilde{\cA}} \tilde Q^*( \check \mu, \tilde a' )  \ge K_A, \quad \check \mu \in \check{\mathfrak{S}}, \tilde a \in  \tilde{\cA} \backslash \arg\min_{\tilde a' \in \tilde{\cA}} \tilde Q^*( \check \mu, \tilde a').
	$$
\end{assumption}

For a fixed finite state grid $\check{\mathfrak{S}}$ and the finite pure-action set $\tilde{\cA}=A^S$, the minimum of all strictly positive action gaps is positive whenever there is at least one non-minimizing pair.
If every action is optimal at every grid point, the condition is vacuous.
In particular, the constant $K_A$ may depend on the state grid $\check{\mathfrak{S}}$.
Note that we do not assume existence of a positive lower bound independent of the state grid, and such a bound does not follow from this finite-set argument.

To recover minimizers or approximate minimizers, it will be convenient to work with the following operators.
In general, they are defined on the vector space $\RR^{\texttt{m}}$.
For $\tau>0$ and $x = (x_1,\dots,x_{\texttt{m}}) \in \RR^{\texttt{m}}$, we define $\softmin_\tau : \RR^{\texttt{m}} \to \RR^{\texttt{m}}$ by
$$
	\softmin_\tau(x) = (e^{- \tau x_1}, \dots, e^{- \tau x_{\texttt{m}}}) / \sum_{j} e^{- \tau x_j}.
$$
For $x \in \RR^{\texttt{m}}$, we define $ \argmine: \RR^{\texttt{m}} \to [0,1]^{\texttt{m}}$ by
$$
	\argmine(x) = \left( \mathbf{1}_{i \in \argmin(x)} \right)_{i=1}^{\texttt{m}} / |\argmin(x)|.
$$
where $\argmin(x) = \{ j \in \{1,\ldots, {\texttt{m}} \} \, : \,   x_j = \min\{x_1,\ldots, x_{\texttt{m}} \} \}$.
In the sequel, we use these operators with the dimension ${\texttt{m}} = | \tilde{\cA} | = | A |^{|S|}$.
For any function $q:   \tilde{\cA} \to \RR$, we identify $q$ with the vector $(q(\tilde a))_{\tilde a \in \tilde{\cA}}$.

\begin{corollary}
	\label{coro:action-cv-tabular-pure}
	Assume the same assumptions as in Theorem~\ref{th:main-cv-tabular-pure} hold and, in addition, that Assumption~\ref{hyp:action-gap} holds.
	Let $\check Q_{N_{\mathrm{upd}}}$ be the table returned by Algorithm~\ref{algo:Qtable-projection}, and let $\varepsilon'$ be as in Theorem~\ref{th:main-cv-tabular-pure}.
	Then, with probability at least $1-\delta$, simultaneously for every $\check\mu \in \check{\mathfrak{S}}$,
	\begin{equation*}
		\big\| \softmin_\tau \big( \check Q_{N_{\mathrm{upd}}}( \check\mu, \cdot) \big)	-
		\argmine \big(  \tilde Q^* (\check\mu, \cdot) \big) \big\|_2
		\leq
		\tau \varepsilon' \sqrt{ | \tilde{\cA} | } + 2  e^{- \tau K_A}  |\tilde{\cA}|.
	\end{equation*}
\end{corollary}

The proof is provided below.
The $\argmine$ in the second term is here in case there are several optimal controls.
The $\softmin$ operator regularizes the best action predicted by the estimation $\check Q_{N_{\mathrm{upd}}}$ of the function $\tilde Q^*$.

Notice that, contrary to Theorem~\ref{th:main-cv-tabular}, here $\cP(A)$ is not involved hence the error bound in Theorem~\ref{th:main-cv-tabular-pure} does not depend on any discretization parameter $\varepsilon_{\mathfrak{A}}$.
Theorem~\ref{th:main-cv-tabular-pure} gives $\varepsilon'=\varepsilon+C\varepsilon_{\mathfrak S}$.
Thus, decreasing $\varepsilon$ reduces only the learning contribution to this bound; the discretization term $C\varepsilon_{\mathfrak S}$ remains unchanged when the state grid is fixed.
Along a sequence of state grids with mesh radii $h=\varepsilon_{\mathfrak S}\downarrow0$, let us denote by $\varepsilon'(h)$ the full error bound and by $K_A(h)$ the corresponding gap.
For fixed $S$ and $A$, the bound in Corollary~\ref{coro:action-cv-tabular-pure} tends to zero provided
\[
\tau(h)\varepsilon'(h)\to0,
\qquad
\tau(h)K_A(h)\to\infty.
\]
These conditions can be met if $\varepsilon'(h)/K_A(h)\to0$. For example, one can take $\tau(h)=[K_A(h)\varepsilon'(h)]^{-1/2}$.
Convergence in probability then follows if the failure probabilities $\delta(h)$ also tend to zero.

\begin{proof}[Proof of Corollary~\ref{coro:action-cv-tabular-pure}]
	We use~\cite[Proposition~4]{GaoPavel},  which states that $\softmin_\tau$ is $\tau$-Lipschitz and~\cite[Lemma~7]{guo2019learning} which states that for $(x_i)_{i =1, \dots, {\texttt{m}}}$,
	$$
		\|\softmin_\tau\, (x) - \argmine(x)\|_2 \le 2 {\texttt{m}} e^{- \tau \delta},
	$$
	where $\delta = \inf_{x_j > \inf_i x_i} (x_j-\inf_i x_i)$, and $\delta = \infty$ if all $x_i$ are equal. We can apply this latter result to $\tilde Q^*(\check \mu, \cdot)$ thanks to assumption~\ref{hyp:action-gap}, with ${\texttt{m}} = |\tilde{\cA}|$ and $\delta \ge K_A$.
	Combining this with Theorem~\ref{th:main-cv-tabular-pure},  on the event in that theorem we have, simultaneously for every $\check \mu$,
	\begin{align*}
		 & \bigg\| \softmin_\tau \Big(  \check Q_{N_{\mathrm{upd}}} (\check\mu, \cdot) \Big) - \argmine \Big( \tilde Q^*( \check \mu, \cdot)  \Big) \bigg\|_2
		\\
		 & \qquad \le
		\bigg\| \softmin_\tau \Big(   \check Q_{N_{\mathrm{upd}}} (\check\mu, \cdot) \Big) - \softmin_\tau \Big(  \tilde Q^*(\check \mu, \cdot) \Big) \bigg\|_2
		\\
		 & \qquad \qquad \qquad + \bigg\| \softmin_\tau \Big( \tilde Q^*(\check\mu, \cdot) \Big)  - \argmine \Big( \tilde Q^*(\check \mu, \cdot) \Big) \bigg\|_2
		\\
		 & \qquad \le
		\tau \bigg\|  \check Q_{N_{\mathrm{upd}}} (\check\mu, \cdot) - \tilde Q^* (\check \mu, \cdot) \bigg\|_2 + 2 |\tilde{\cA}| e^{- \tau K_A}
		\\
		 & \qquad  \le \tau  \sqrt{ | \tilde{\cA} | } \sup_{ \tilde a' \in \tilde{\cA}} \left|  \check Q_{N_{\mathrm{upd}}} (\check\mu, \tilde a') - \tilde Q^* (\check \mu, \tilde a') \right| + 2 |\tilde{\cA}| e^{- \tau K_A}
		\\
		 & \qquad \le
		\tau \varepsilon'  \sqrt{ | \tilde{\cA} | } + 2  e^{- \tau K_A} | \tilde{\cA} |.
	\end{align*}
\qed\end{proof}

\subsection{Deep reinforcement learning for MFMDP }\,
\label{sec:DDPG-algo}

 The tabular Q-learning method of the previous section is sufficiently straightforward to allow for a detailed analysis.
However, it cannot be used in practice for large state or action spaces because of the computational cost of discretizing the simplexes.
An alternative is to work directly with continuous spaces, in which case the policies and value functions cannot be represented in a tabular way.
Instead, we can rely on function approximation.
To this end, we now propose methods from deep RL which are more suitable for continuous spaces.
The motivations are twofold.

First, if $S$ and $A$ are finite but we want to learn an optimal policy that is potentially randomized at level-0, the discretization approach proposed in~\S~\ref{se:tabularQ-mixed} has a complexity that increases with the number of points in the discretization of $\cP(A)$, which itself increases exponentially quickly with the cardinality of $A$.
For this reason, it can be interesting to tackle directly $\cP(A)$ as a continuous action space and to use deep RL methods for continuous action space MDPs.

Second, some MFC problems are naturally posed with a continuous state space $S$.
In this case, under mild conditions, the optimal policy is in fact non-randomized not only at the level-1 but even at the level-0 (see, e.g.,~\cite{CarmonaDelarue_book_I,BensoussanFrehseYam} and the references therein).
However, the state of the MFMDP is an element of the infinite dimensional space $\cP(S)$.
As noted in Chapter~\ref{ch:AAP}, this measure-valued state allows the central planner to treat the collective behavior of an infinite population as a single Markovian process on the simplex $\bar S = \cP(S)$.
From a numerical viewpoint, we need two ingredients: (1) a finite-dimensional approximation of the mean field and (2) a parameterized approximation of the value function or the policy taking this finite-dimensional representation of the mean field state as an input.
For the second point, we will again use deep neural networks.
For the first point, for the sake of definiteness, we choose to simply replace $\cP(S)$ by $\cP(\check S)$ where $\check S$ is a discretization of $S$ with a finite number of points.
We assume that, given $\check\mu \in \cP(\check S)$ and a level-1 action $\check a: S \to A$, one can get from the environment a sample of the next state and the associated cost $\tilde f(\check\mu, \check a)$.
The problem thus reduces to an MDP with finite dimensional (but potentially continuous) state and action spaces.
Such MDPs can be solved with a variety of deep RL algorithms.
In the following, we provide numerical illustrations based on the \defi{DDPG}\index[sub]{Deep Deterministic Policy Gradient}\index[sub]{DDPG} proposed in~\cite{lillicrap2015continuous}.
It relies on two neural networks, one for the Q-function (the critic) and one for the policy (the actor).
The core of the algorithm involves alternating updates: the critic is updated by minimizing an empirical square error, while the actor is updated through gradient descent.
To improve exploration, a Gaussian noise $\epsilon^a_{n+1}$ is added to the action prescribed by the actor.
Furthermore, for more stability, target networks are also added.
We refer to~\cite{lillicrap2015continuous} for more details. We write the algorithm below in cost-minimization form: the critic approximates cumulative cost and the actor is updated by descending the estimated cost. The algorithm is summarized in Algorithm~\ref{algo:DDPG}.

\begin{algorithm}[htbp]
	\DontPrintSemicolon
	\KwData{Actor network $\pi_{\theta^{\pi}}$, critic network $Q_{\theta^Q}$, target networks $\pi_{\tilde{\theta}^{\pi}}$ and $Q_{\tilde{\theta}^Q}$, replay buffer $\mathcal{R}$, batch size $N_{\mathrm{batch}}$, exploration noise $\mathcal{N}$, soft update coefficient $\rho_{\mathrm{tar}} \in (0,1)$, discounting factor $\gamma \in (0, 1)$.}
	\KwResult{Trained actor policy $\pi_{\theta^{\pi}}$ mapping states in $\cP(S)$ to actions in $\tilde{\cA}$.}
	\Begin{
	Randomly initialize actor network $\pi_{\theta^{\pi}}$ and critic network $Q_{\theta^Q}$ with weights $\theta^{\pi}$ and $\theta^Q$\;
	Initialize target networks $\tilde{\theta}^{\pi} \leftarrow \theta^{\pi}$ and $\tilde{\theta}^Q \leftarrow \theta^Q$\;
	Initialize replay buffer $\mathcal{R}$\;
	\For{ $\mathrm{episode} = 1, \dots, N_{\mathrm{epi}}$}{
	Initialize a random state $\mu_0 \in \cP(S)$\;
	\For{ $n = 0, \dots, T-1$}{
	\vskip 2pt
	Select action with exploration noise: $\tilde{a}_n = \pi_{\theta^{\pi}}(\mu_n) + \epsilon^a_n$, where $\epsilon^a_n \sim \mathcal{N}$\;
	Execute action $\tilde{a}_n$ in state $\mu_n$, observe cost $\tilde{f}(\mu_n, \tilde{a}_n)$ and new state $\mu_{n+1}$ from the environment\;
	Store transition $(\mu_n, \tilde{a}_n, \tilde{f}(\mu_n, \tilde{a}_n), \mu_{n+1})$ in $\mathcal{R}$\;
	Sample a random minibatch of $N_{\mathrm{batch}}$ transitions $(\mu_i, \tilde{a}_i, \tilde{f}_i, \mu'_{i+1})$ from $\mathcal{R}$\;
	Set target $y_i = \tilde{f}_i + \gamma Q_{\tilde{\theta}^Q} \big( \mu'_{i+1}, \pi_{\tilde{\theta}^{\pi}}(\mu'_{i+1}) \big)$\;
	Update critic network by minimizing the loss: $L(\theta^Q) = \frac{1}{N_{\mathrm{batch}}} \sum_i \big( y_i - Q_{\theta^Q}(\mu_i, \tilde{a}_i) \big)^2$\;
	Update actor policy by descending the sampled cost gradient w.r.t.\ $\theta^{\pi}$:
	$$ \frac{1}{N_{\mathrm{batch}}} \sum_i \nabla_{\tilde{a}} Q_{\theta^Q}(\mu_i, \tilde{a}) \Big|_{\tilde{a} = \pi_{\theta^{\pi}}(\mu_i)} \nabla_{\theta^{\pi}} \pi_{\theta^{\pi}}(\mu_i) $$
	Update target networks: $\tilde{\theta}^Q \leftarrow \rho_{\mathrm{tar}} \theta^Q + (1 - \rho_{\mathrm{tar}}) \tilde{\theta}^Q$ and $\tilde{\theta}^{\pi} \leftarrow \rho_{\mathrm{tar}} \theta^{\pi} + (1 - \rho_{\mathrm{tar}}) \tilde{\theta}^{\pi}$\;
	}
	}
	\KwRet{$\pi_{\theta^{\pi}}$}
	}
	\caption{Mean-Field Deep Deterministic Policy Gradient (MF-DDPG)}
	\label{algo:DDPG}
\end{algorithm}

\section{Numerical Examples}
\label{se:numres}

In this section, we illustrate the above RL algorithms on two different examples.
In each example, the representative agent's state space and action space are finite.
The first example is a cybersecurity model with 4 states, on which we illustrate both the tabular Q-learning algorithm and the DDPG algorithm.
The second example features a larger number of states, leading us to apply only the DDPG algorithm, which offers the advantage of not requiring simplex discretization. 
The code is available at: \url{https://github.com/mlauriere/mfrl-tutorial-code}.

\subsection{Example 1: Cybersecurity model}\,

For a first testbed, we start with a finite state problem.
We revisit the cybersecurity example introduced in~\cite{kolokoltsov2016mean}, now adopting the perspective of a central planner (such as a large company or a state) aiming to protect its system against hacker attacks.
The situation can be phrased as an MFC problem.

In this model, the population consists of a large group of computers which can be either defended (D) or undefended (U), and either infected (I) or susceptible (S) of infection.
Hence the set $S$ has four elements corresponding to the four possible combinations: DI, DS, UI, US.
The action set is $A = \{0,1\}$, where $0$ is interpreted as the fact that the central planner is satisfied with the current protection level (D or U) of the computer under consideration, whereas $1$ indicates a desire to modify this level.
In the latter case, the update occurs at a (fixed) rate $\lambda >0$.
If the controls are pure at level-0, the central planner assigns a single action to each of the four states, which is then applied to all computers at that state.
If the controls are mixed at level-0, the planner chooses a distribution over actions for each state, and each computer in that state then independently selects an action according to the chosen distribution.
When infected, each computer may recover at rate $q_{rec}^D$ or $q_{rec}^U$ depending on whether it is defended or not.
On the other hand, a computer may be infected either directly by a hacker, at rate $v_H q_{inf}^D$ (resp. $v_H q_{inf}^U$) if it is defended (resp. undefended), or by undefended infected computers, at rate $\beta_{UU}\mu(\{UI\})$ (resp. $\beta_{UD}\mu(\{UI\})$) if it is undefended (resp. defended), or by defended infected computers, at rate $\beta_{DU}\mu(\{DI\})$ (resp. $\beta_{DD}\mu(\{DI\})$) if it is undefended (resp. defended).
Here $v_H$ can be interpreted as the attack intensity parameter.

In short, the infinitesimal generator is given by:
\begin{equation}
	\label{eq:cybersecurity-MFC-generator}
	Q^{\mu,a} =
	\begin{pmatrix}
		\dots     & Q^{\mu,a}_{DS \rightarrow DI} & \lambda a & 0
		\\
		q_{rec}^D & \dots                         & 0         & \lambda a
		\\
		\lambda a & 0                             & \dots     & Q^{\mu,a}_{US \rightarrow UI}
		\\
		0         & \lambda a                     & q_{rec}^U & \dots
	\end{pmatrix}
\end{equation}
where
\begin{align*}
	 & Q^{\mu,a}_{DS \rightarrow DI} = v_H q_{inf}^D + \beta_{DD} \mu(\{DI\})  + \beta_{UD} \mu(\{UI\}) ,
	\\
	 & Q^{\mu,a}_{US \rightarrow UI} = v_H q_{inf}^U + \beta_{UU} \mu(\{UI\}) + \beta_{DU} \mu(\{DI\}),
\end{align*}
and all instances of $\dots$ are chosen so that each row sums to zero, i.e., each diagonal entry is the negative of the sum of the off-diagonal rates in its row.
At each time step, the central planner pays a protection cost $k_D>0$ for each defended computer, and a penalty $k_I>0$ for each infected computer.
The instantaneous cost in the MFMDP is thus defined as:
$$
	\bar f(\mu, \bar a) = k_D \mu(\{DI, DS\}) + k_I \mu(\{DI, UI\}), \qquad (\mu, \bar a) \in \bar S \times \bar A.
$$
The optimal control and optimal flow of distributions can be characterized by a forward-backward ODE system which can be obtained in a way similar to what is done in the MFG setting e.g. in ~\cite[\S~7.2.3]{CarmonaDelarue_book_I}.
We will use this solution as a benchmark.

\textbf{Tabular Q-learning.} For the sake of illustration, we present results obtained by tabular Q-learning with simplex discretization as described in \S~\ref{se:tabularQ-discrete}. The state space for the population distribution is $\bar S$, which is identified with the simplex $\mathfrak{S} = \{(\mu^{(i)})_{i=1,\dots,4} \in [0,1]^4\,:\, \sum_{i} \mu^{(i)} = 1\}$. 
In this example, we consider that the actions are pure not only at level 1 (the central planner picks a single control from $\bar A$) but also at level 0 (individual agents do not randomize). 
This means that for each state $i \in \{1,\dots,4\}$, the action $a^{(i)}$ belongs to $\{0, 1\}$. Thus, the set of actions $\bar A$ for the population is the finite set $\{0,1\}^4$ of cardinality 16.
We replace $\mathfrak{S}$ by the finite set:
\begin{equation}
	\label{eq:grid-mu-simplex}
	\check{\mathfrak{S}} = \Big\{ (\mu^{(i)})_{i=1,\dots,4} \in \{0,1/N_m,\dots,1-1/N_m, 1\}^4\,:\, \sum_{i} \mu^{(i)} = 1 \Big\},
\end{equation}
where $N_m+1$ is the number of points in the discretization of each dimension.

This is the running example introduced in Chapter~\ref{ch:prelims}. Here the pure level-0 rule is $p^{\mathrm{sw}}\in\{0,1\}^4$, with $p^{\mathrm{sw}}(x^{(i)})=a^{(i)}$ for $i=1,\dots,4$. For a queried pair $(\check\mu,p^{\mathrm{sw}})$, the simulator returns the cost $c(\check\mu)$ from~\eqref{eq:running-cyber-cost} and the projected successor
\[
\check\mu'=\proj_{\check{\mathfrak S}}\Phi(\check\mu,p^{\mathrm{sw}},\varepsilon^0),
\]
using~\eqref{eq:running-cyber-update}. Algorithm~\ref{algo:Qtable-projection} uses this transition to update the Q-table entry associated with $(\check\mu,p^{\mathrm{sw}})$.

We consider three different training schemes for the Q-function:
\begin{enumerate}
    \item \textbf{SyncFull}: A synchronous approach where the Q-function is updated for all state-action pairs $(\check{\mu}, \check{a}) \in \check{\mathfrak{S}} \times \check{\cA}$ at each iteration using the exact expectation of the next state's value. This serves as a benchmark for the discrete problem.
    \item \textbf{SyncGreedy}: A synchronous approach where at each iteration, for every state $\check{\mu}$, only the Q-value for the current $\epsilon$-greedy action is updated based on a single sample transition.
    \item \textbf{AsyncTraj}: An asynchronous approach where the Q-values are updated along trajectories sampled from the environment, using an $\epsilon$-greedy policy to explore the state space. %
    This is an online trajectory-sampling scheme. (Assumption~\ref{hyp:covering-time} remains an additional hypothesis to check for applying the quantitative bounds to this scheme and has not been verified for these experiments.)
\end{enumerate}

For these simulations, we used the parameters summarized in Table~\ref{tab:params-cyber}. The time step is set to $\Delta t = 0.2$ and the horizon to $T=10.0$ (50 steps per episode). 
The continuous discount rate is set by the convention $\rho=-\log(0.5)/0.1$, so the per-step discount used in these tabular simulations is $\gamma=\exp(-\rho\Delta t)=0.25$.
The discretization mesh size is $N_m=30$ for the case without common noise and $N_m=30$ for the case with common noise.

\begin{table}[htbp]
	\centering
	\begin{tabular}{l c | l c}
		\hline
		\textbf{Parameter}     & \textbf{Value} & \textbf{Parameter}         & \textbf{Value} \\
		\hline
		$\beta_{UU}$           & $0.3$          & $q_{rec}^D$                & $0.5$          \\
		$\beta_{UD}$           & $0.4$          & $q_{rec}^U$                & $0.4$          \\
		$\beta_{DU}$           & $0.3$          & $q_{inf}^D$                & $0.4$          \\
		$\beta_{DD}$           & $0.4$          & $q_{inf}^U$                & $0.3$          \\
		$v_H$                  & $0.6$          & $\lambda$                  & $0.8$          \\
		$k_D$                  & $0.3$          & $k_I$                      & $0.5$          \\
		\hline
		$N_m$ (Mesh size)      & $30$      & $\gamma$ (Discount factor) & $0.25$ \\
		$\Delta t$ (Time step) & $0.2$          & $T$ (Horizon)              & $10$           \\
		\hline
	\end{tabular}
	\caption{Parameters and tabular Q-learning hyperparameters for the cybersecurity model.}
	\label{tab:params-cyber}
\end{table}

\begin{figure}[htbp]
    \centering
        \includegraphics[width=\textwidth]{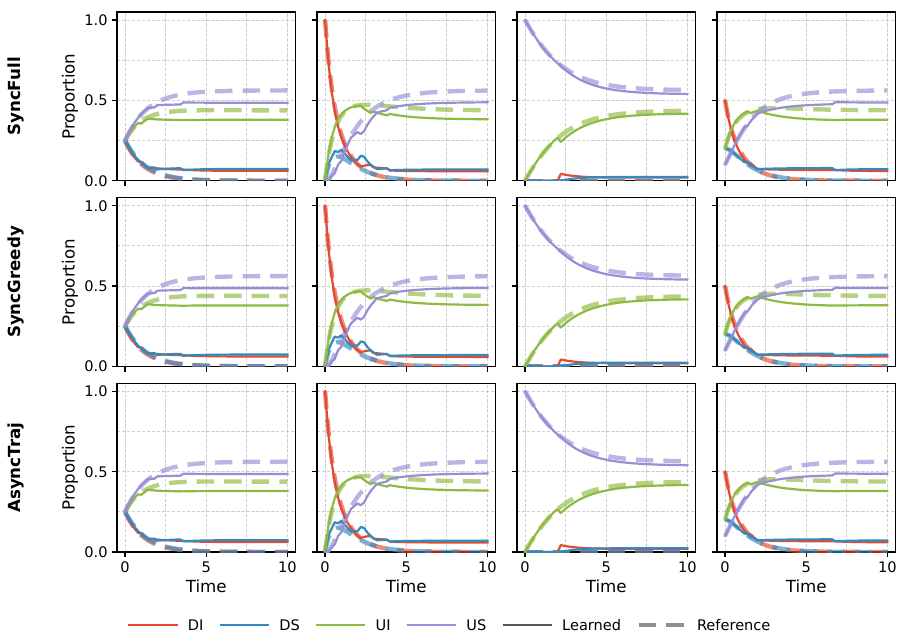}
    \caption{Cybersecurity model (No CN): Evaluation trajectories for the three Q-learning solvers starting from four different initial distributions.}
    \label{fig:cyber-q-eval-nocn}
\end{figure}

\begin{figure}[htbp]
    \centering
    \includegraphics[width=\textwidth]{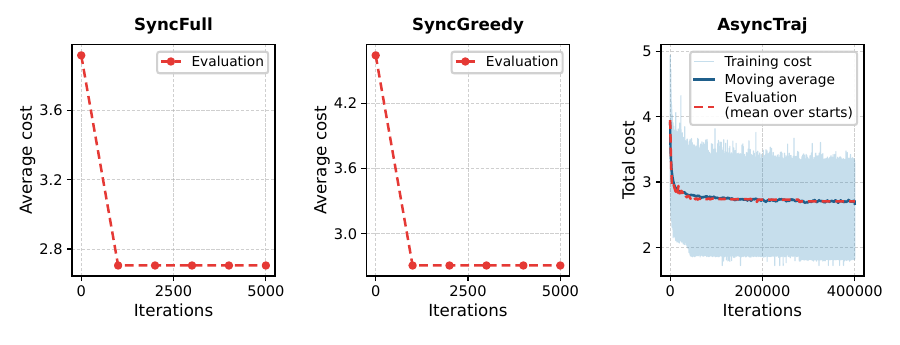}
    \caption{Cybersecurity model (No CN): Training progress for the three solvers. In the AsyncTraj panel, red dashes show evaluation costs, light blue shows training costs, and dark blue shows the plotted 20-iteration moving average after an additional centered average over 5000 iterations, with shorter windows near the endpoints.}
    \label{fig:cyber-q-train-nocn}
\end{figure}

We first present results without common noise (denoted by CN for short). In this case, the MFMDP is deterministic. Figure~\ref{fig:cyber-q-eval-nocn} shows the evaluation trajectories for the three solvers, which all successfully learn to drive the system toward the same stationary configuration although they do not match perfectly the reference ones because of the discretization of the state space. Figure~\ref{fig:cyber-q-train-nocn} displays the training progress in terms of costs.

\begin{figure}[htbp]
    \centering
    \includegraphics[width=\textwidth]{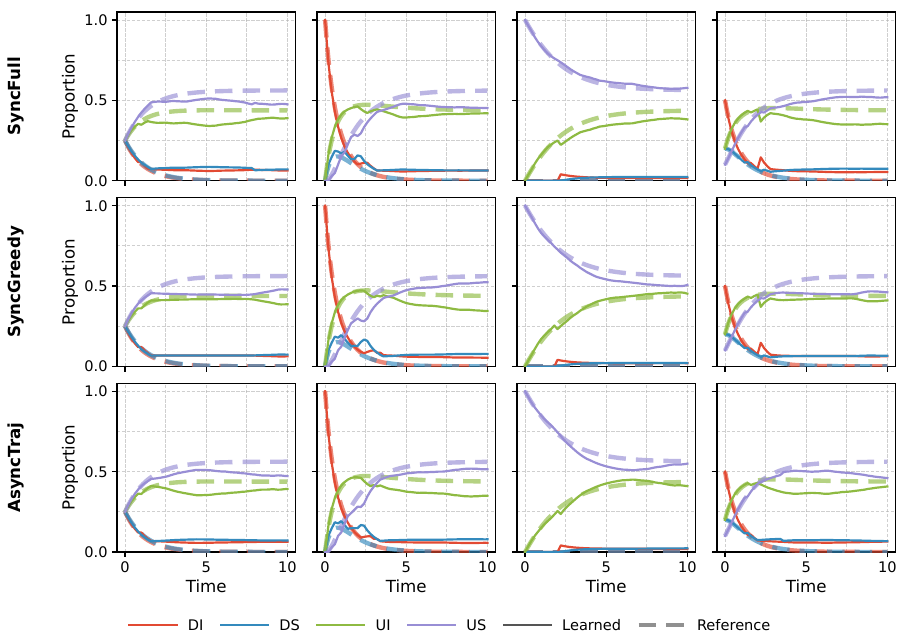}
    \caption{Cybersecurity model (With CN): Evaluation trajectories for the three Q-learning solvers.}
    \label{fig:cyber-q-eval-withcn}
\end{figure}

\begin{figure}[htbp]
    \centering
    \includegraphics[width=\textwidth]{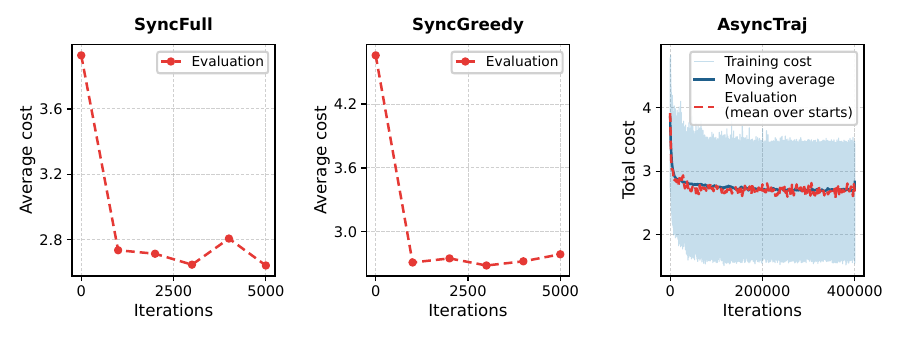}
    \caption{Cybersecurity model (With CN): Training progress for the three solvers. Colours and moving averages are as in Figure~\ref{fig:cyber-q-train-nocn}.}
    \label{fig:cyber-q-train-withcn}
\end{figure}

\begin{figure}[htbp]
    \centering
    \includegraphics[width=\textwidth]{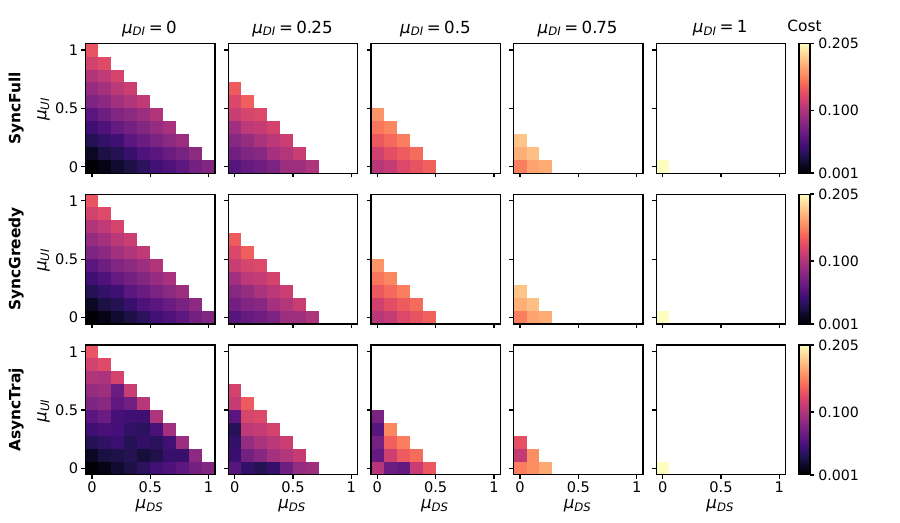}
    \caption{Cybersecurity model (With CN): Value function slices on the simplex for the three solvers compared with the benchmark for the discrete problem; each row corresponds to one method.} 
    \label{fig:cyber-q-simplex-withcn}
\end{figure}

Next, we consider the case with common noise, where the attack intensity $v_H$ oscillates. The results are shown in Figures~\ref{fig:cyber-q-eval-withcn}, \ref{fig:cyber-q-train-withcn}, and \ref{fig:cyber-q-simplex-withcn}. Despite the stochasticity, the solvers manage to learn robust policies that stabilize the system. Figure~\ref{fig:cyber-q-simplex-withcn} specifically shows slices of the learned value function on the simplex, comparing it with the benchmark solution.

\noindent
\textbf{DDPG.} The solution can also be learned without discretizing the simplex $\mathfrak{S}$, by directly considering the value function on the continuous space. We use the DDPG algorithm~\cite{lillicrap2015continuous} to learn a deterministic policy $\pi: \mathfrak{S} \to \bar A$. For this continuous solver, we consider the action space to be the continuous hypercube $\bar A = [0,1]^4$, where $a^{(i)} \in [0,1]$ represents the defense intensity in state $i$. Notice that this is a relaxed action space compared to the tabular setting, where $\bar A = \{0,1\}^4$.

The actor and critic networks both consist of two hidden layers of 32 neurons each, with sigmoid activations. We trained the networks for 500 episodes with a learning rate of $10^{-3}$ for the actor and $2 \times 10^{-3}$ for the critic.
The soft update parameter for the target networks was set to $\rho_{\mathrm{tar}}=0.005$. Here the DDPG implementation uses $\Delta t=0.1$, so the same continuous discount convention $\rho=-\log(0.5)/0.1\approx 6.93$ gives the DDPG per-step discount $\gamma=\exp(-\rho\,0.1)=0.5$.

\begin{figure}[htbp]
    \centering
    \includegraphics[width=\textwidth]{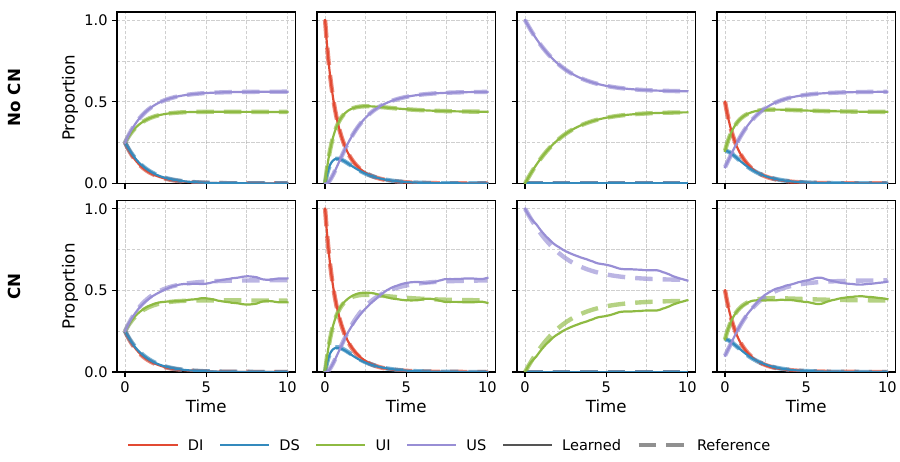}
    \caption{Cybersecurity model: Evaluation trajectories for DDPG starting from four different initial distributions.}
    \label{fig:cyber-ddpg-trajs}
\end{figure}

\begin{figure}[htbp]
    \centering
    \includegraphics[width=\textwidth]{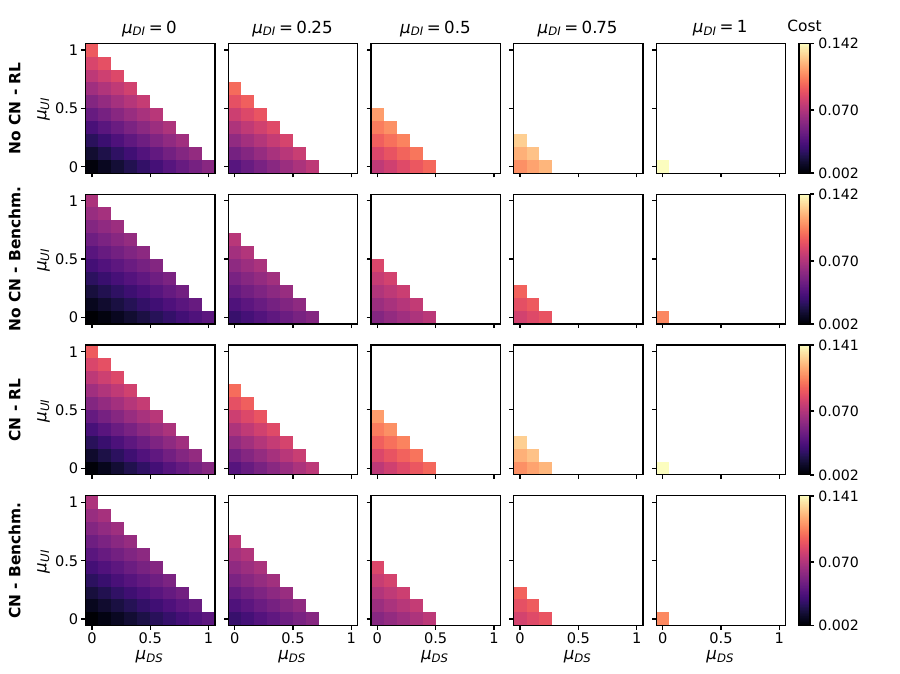}
    \caption{Cybersecurity model: DDPG value function slices on the simplex.}
    \label{fig:cyber-ddpg-simplex}
\end{figure}

\begin{figure}[htbp]
    \centering
    \includegraphics[width=\textwidth]{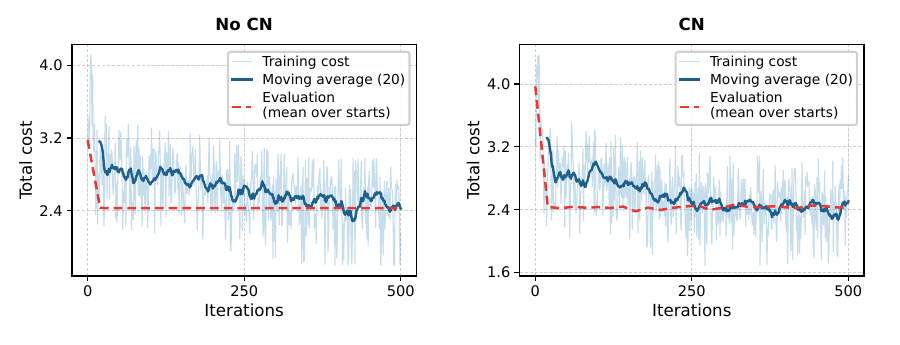}
    \caption{Cybersecurity model: DDPG training progress (cost).}
    \label{fig:cyber-ddpg-costs}
\end{figure}

The results for both deterministic and stochastic environments are displayed in Figure~\ref{fig:cyber-ddpg-trajs} for the evaluation trajectories, Figure~\ref{fig:cyber-ddpg-simplex} for the value function simplex slices, and Figure~\ref{fig:cyber-ddpg-costs} for the training progress in terms of costs. DDPG leads to trajectories that are closer to the reference ones than with tabular Q-learning because the simplex has not been discretized. We also observe that in the presence of common noise, the trajectories remain close to the reference ones, although there are small fluctuations.

\subsection{Example 2: Discrete distribution planning}\,

We now apply the MFC framework to a problem where the objective is to match a target distribution.
We consider a model with parameters summarized in Table~\ref{tab:params-planning}, where the agents' actions correspond to moving left, staying, or moving right.
At the boundaries, attempted moves outside the grid are folded back to the boundary.
The dynamics follow $F(x,a,\mu,e,e^0)=\Pi_X(\Pi_X(x+a)+e^0)$.
Here $\Pi_X$ denotes projection onto the finite grid $X=\{1,\dots,11\}$; the outer projection applies the same boundary convention after the common-noise shift.
A representative agent incurs a one-step cost of $1$ for any movement (left or right), plus the distribution penalty $2\|\mu-\mu_{\texttt{target}}\|_2$.
Thus, $f(x,a,\mu)=\mathbf 1_{\{a\neq0\}}+2\|\mu-\mu_{\texttt{target}}\|_2$.
For this example, the precise target distribution $\mu_{\texttt{target}}$ is given in Table~\ref{tab:params-planning}.
This model does not include idiosyncratic noise.

\begin{table}[htbp]
	\centering
	\begin{tabular}{l c}
		\hline
		\textbf{Parameter}                       & \textbf{Value}                                      \\
		\hline
		$N_{\texttt{states}}$ (Number of states) & $11$                                                \\
		$X$ (State space)                        & $\{1,\dots,11\}$                                    \\
		$A$ (Action space)                       & $\{-1,0,+1\}$                                       \\
		$\mu_{\texttt{target}}$                  & $(0, 0, 0.05, 0.1, 0.2, 0.3, 0.2, 0.1, 0.05, 0, 0)$ \\
		$\epsilon^0_n$ (Common noise)            & $\{-1, 0, 1\}$ with prob. $(0.05, 0.9, 0.05)$       \\
		$T$ (Horizon)                            & $50$                                                \\
		$\Delta t$ (Time step)                   & $1.0$                                               \\
		Episodes                                 & $3000$                                              \\
		NN Architecture                          & $[128, 128, 128]$                                   \\
		\hline
	\end{tabular}
	\caption{Parameters for the discrete distribution planning example.}
	\label{tab:params-planning}
\end{table}
Importantly, except for specific pairs of initial and target distributions, the population generally cannot match the target distribution unless agents are permitted to randomize their actions at the individual level.
We therefore use $\cP(A)^X$ as the level-1 action space.
For a randomized action kernel $\nu\in\cP(A)^X$, the movement-cost term is $\sum_{x\in X}\mu(x)\sum_{a\in A}\mathbf 1_{\{a\neq0\}}\nu_x(a)$.
The resulting continuous action space justifies the application of the DDPG method.

We present results for the DDPG algorithm both with and without common noise. In the stochastic case, the common noise follows the distribution specified in Table~\ref{tab:params-planning}. The networks are trained for 3000 episodes using the hyperparameters detailed in the same table.

\begin{figure}[htbp]
    \centering
    \includegraphics[width=\textwidth]{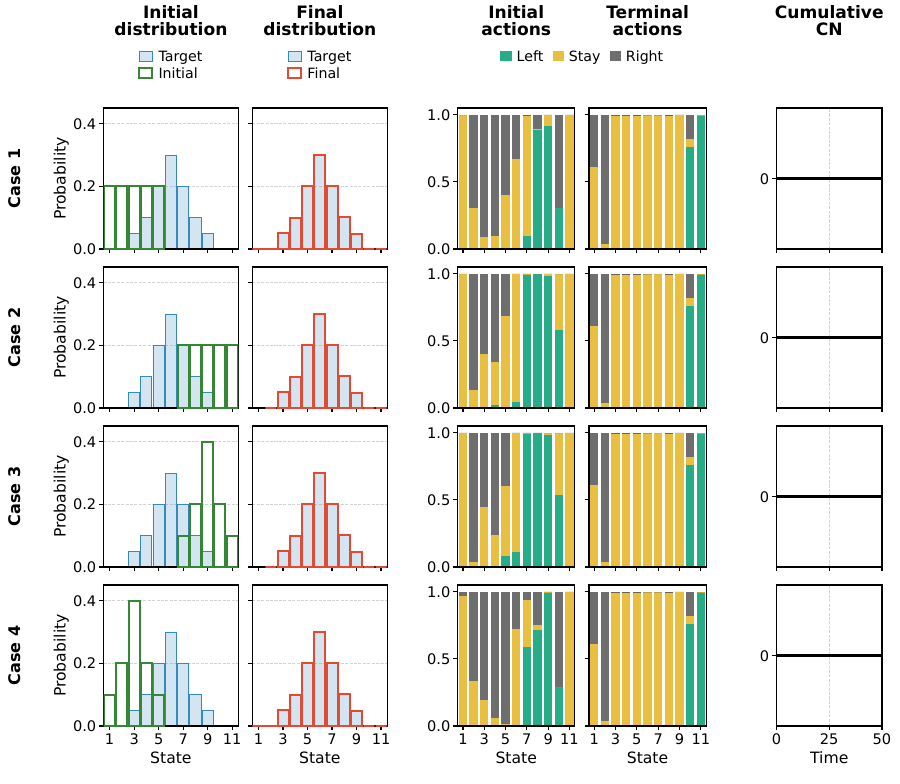}
    \caption{Discrete distribution planning (No CN): Evaluation rows for four initial distributions. Each row shows (from left to right): Initial state, Final state, Initial action, Terminal action, and Cumulative common noise (zero).}
    \label{fig:distribplanning-eval-nocn}
\end{figure}

\begin{figure}[htbp]
    \centering
    \includegraphics[width=\textwidth]{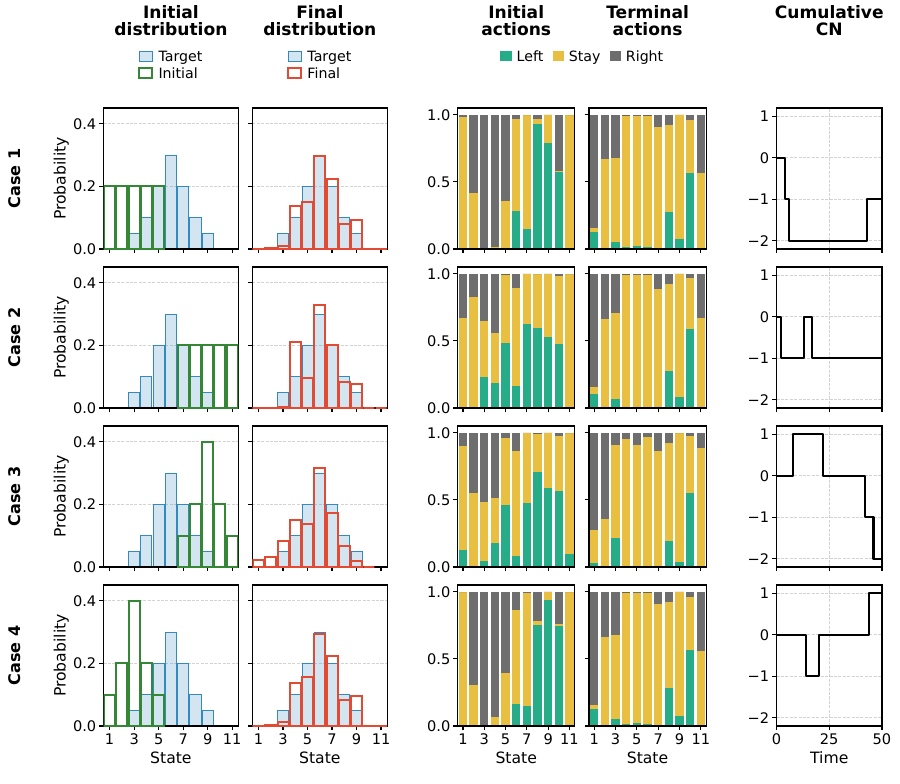}
    \caption{Discrete distribution planning (With CN): Evaluation rows for four initial distributions in a stochastic environment.}
    \label{fig:distribplanning-eval-withcn}
\end{figure}

\begin{figure}[htbp]
    \centering
    \includegraphics[width=\textwidth]{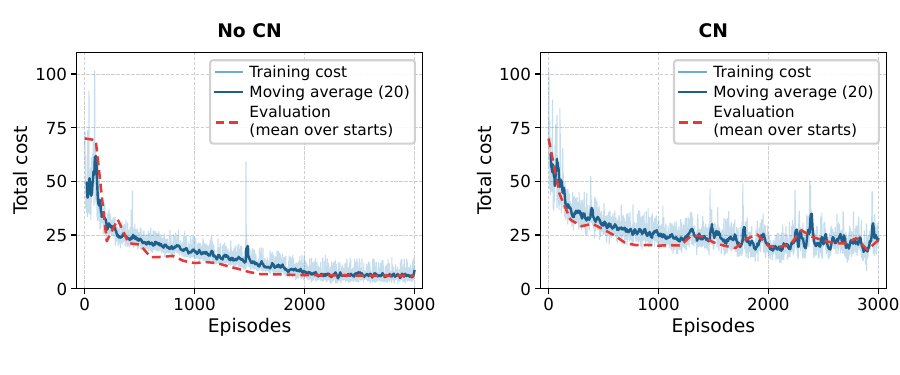}
    \caption{Discrete distribution planning: DDPG training progress (cost). Left: no CN; right: with CN. Pale blue shows training costs, dark blue their 20-episode moving average, and dashed red evaluation costs averaged over five rollouts for each of the four initial distributions. Costs are undiscounted episode totals over 50 time steps.}
    \label{fig:distribplanning-costs}
\end{figure}

Figures~\ref{fig:distribplanning-eval-nocn} and~\ref{fig:distribplanning-eval-withcn} display the evaluation results for four different initial distributions, respectively without and with common noise. Each row corresponds to a single test case and contains five subplots: (1) Initial state distribution vs target, (2) Final state distribution vs target, (3) Initial action distribution across states, (4) Terminal action distribution, and (5) Cumulative common noise (which is zero in the non-stochastic case). Without common noise, the trajectories are deterministic and the final state distribution is very close to the target distribution. With common noise, the trajectories are stochastic and the final state distribution is slightly different from the target distribution. But we observe that in both scenarios, the DDPG agent successfully learns to move the population towards the target distribution. At the initial time, the control predominantly favors movement towards the center of the target distribution, while at the terminal time, the "stay" action becomes dominant. Figure~\ref{fig:distribplanning-costs} shows the training progress in terms of the social cost for both cases.

\section{Notes and Complements}

RL has been applied to solve MFC and MFMDPs using Q-learning as presented in this chapter in~\cite{CarmonaLauriere_AAP}, as well as in~\cite{gu2024meanmarldecentralized,gu2021meanQ}; see also \cite{gu2019dynamicmfc}. The deep RL approach presented in this chapter was first introduced in~\cite{CarmonaLauriere_AAP}. Extensions include applications to mean-field type games, namely, Nash equilibria between a finite number of central planners solving MFC problems: \cite{shao2025reinforcement} studied convergence of Nash Q-learning and application of DDPG to such problems, and an adaptation of the PPO algorithm to MFTGs has been proposed in~\cite{jeloka2025learning}.

Other algorithms include a two-timescale method introduced in~\cite{angiuli2022unified} and whose convergence has been established in~\cite{angiuli2023convergence}. This approach has then been extended to study actor-critic methods for MFC in~\cite{angiuli2023convergence}. Alternative actor-critic methods have been proposed in~\cite{frikha2025actor,PhamWarin}.
\cite{bayraktar2025learning} proposed a learning algorithm based on linear function approximation. \cite{meunier2026model} proposed a method based on an adaptation of the classical policy gradient method. 
In the line of continuous-time RL developed by~\cite{wang2020reinforcement}, \cite{wei2025continuous} studied q-learning for MFC in continuous time. For model-free continuous-time mean-field control from discrete-time transition data, see also~\cite{bayraktar2026mean}. RL methods for partially observable MFC models has been proposed in~\cite{cui2023learningpomfc,cui2024partially}. Related deep-learning methods for MFC with common noise, based on population-dependent controls and neural approximations of the mean field, have been studied in~\cite{dayanikli2024deep}.

The model for cybersecurity considered in this chapter was introduced, in the context of mean field games, by~\cite{kolokoltsov2016mean}; related MFG and MFRL references include \cite{subramanian2019reinforcement,guo2019learning,elie2020convergence,fu2019actorcriticMFG,anahtarci2020qregu,lauriere2022learning}. Other examples of applications of RL for MFC include economics and finance~\cite{angiulia2023reinforcement}, robotic swarm control~\cite{cui2023scalable}, and UAVs~\cite{chen2020mean,gu2025task}.

We conclude by mentioning that while this chapter focuses on a model-free perspective, an alternative approach is to rely on model-based methods such as in~\cite{pasztor2023efficient,jusup2024safe}.

\chapter{Implementation of the LQ Mean Field Model}
\label{ch:numeric_II}

\begin{abstract}
	\emph{
		This chapter establishes global convergence results for policy-gradient methods in the Linear-Quadratic Mean Field Control setting introduced in Chapter~\ref{ch:LQMFC}.
		We analyze three configurations: model-based exact gradients, model-free gradients with an idealized McKean-Vlasov simulator, and model-free gradients with a finite-population simulator.
		We quantify how a finite collection of agents can learn a policy that is approximately optimal for the infinite-population problem.
		Numerical experiments illustrate the behavior of these methods under stochastic dynamics and finite-population simulation.
	}
\end{abstract}

Following the general numerical methods discussed in Chapter~\ref{ch:numeric_I}, this chapter focuses on the Linear-Quadratic Mean Field Control (LQMFC)\index[not]{LQMFC}\index[sub]{linear-quadratic mean field control} setting.
As introduced in Chapter~\ref{ch:LQMFC}, the family of Linear-Quadratic (LQ) models is central to optimal control due to its mathematical tractability and wide-ranging applicability.
In reinforcement learning (RL), these models serve as an important benchmark for analyzing the convergence of iterative algorithms~\cite{fazel2018global,recht2018tour}.

Mean-field type control, or McKean-Vlasov control, extends the LQ framework by incorporating state and action distributions.
While Chapter~\ref{ch:numeric_I} explored tabular and deep RL for general MFMDPs, the specific structure of LQ models allows us to establish global convergence guarantees for Policy Gradient (PG)\index[not]{PG}\index[sub]{policy gradient} methods, even in model-free settings.
Our analysis generalizes the results of Fazel et al.~\cite{fazel2018global} to the mean-field case, proving that the optimal parameters of the feedback control can be effectively learned.

A key feature of our model is the presence of \emph{common noise}, which captures random shocks affecting all agents simultaneously (see Chapter~\ref{ch:LQMFC}).
While it complicates the analysis, it also makes the limiting population distribution stochastic.
In our model-free investigations, we study two types of simulators: an idealized McKean--Vlasov (MKV)\index[not]{MKV}\index[sub]{McKean--Vlasov simulator} simulator in \S~\ref{subsection:modelfree_MKV_simulator}
 and a more realistic finite-population simulator in \S~\ref{subsec:PG-popsimu}.
The idealized MKV setting is of interest in its own right and also serves as an intermediate step toward the population-based method.
We first analyze learning with the MKV simulator, then extend the argument to the population simulator by controlling the additional finite-population approximation error.
Numerical experiments illustrate these methods in \S~\ref{sec:numerics-LQ}. The detailed estimates supporting the convergence proofs are collected in Appendix~\ref{app:lq-pg-proofs}.

\section{Policy gradient methods}
\label{ref:PGconv}

Besides the notation introduced for the LQMFC problem in Section~\ref{se:prelims_LQ} and Chapter~\ref{ch:LQMFC}, it will sometimes be convenient to use the following notation. We will use $u$ for the action. 
The one-step cost function will be denoted by $c: S \times S \times A \times A \to \RR$, with model parameters $(Q, \bar Q, R, \bar R) \in \RR^{d \times d} \times \RR^{d \times d} \times  \RR^{\ell \times \ell} \times \RR^{\ell \times \ell}$ given by 
    \begin{equation}
        \label{fo:lq_one_step_cost}
        c(x, \bar{x}, \ctrl, \bar{\ctrl}) = (x-\bar{x})^\top Q (x-\bar{x}) + \bar{x}^\top (Q + \bar{Q}) \bar{x} + (\ctrl-\bar{\ctrl})^\top R (\ctrl-\bar{\ctrl}) + \bar{\ctrl}^\top (R + \bar{R}) \bar{\ctrl}.
    \end{equation}
Thus, using the notation $f$ introduced in~\eqref{fo:LQ_cost}, we have $f(x, u, \mu) = c(x, \bar\mu^x, u, \bar\mu^a)$, where $\bar\mu^x$ and $\bar\mu^a$ denote the means of $\mu$ for states and actions, respectively. 
We also recall that the system function $F$ is defined in~\eqref{eq:lq_sys_function_F}.

In this section, we search for the optimal control process using PG algorithms.
Theorem~\ref{th:existence_linear_control} states that the optimal control process has feedback coefficients $\theta^* = (K^*, L^*)$.
Inspired by formula~\eqref{fo:optimal_control_linear_in_x_barx}, we limit our search for the optimal control process to the set of controls $\cU_{ad}^\Theta$ that can be parameterized by an element in $\Theta$.
We present model-based and model-free PG algorithms to find the best admissible control parameter $\theta^* \in \Theta$.

\subsection{Reformulation of the MFC problem}
We reformulate the MF problem~\eqref{pb:MFC_L2_admissible_control} with the set of admissible control parameters $\Theta$.
For a parameterized admissible control process $\bctrl^\theta \in \cU_{ad}^\Theta$, we denote the discounted MF cost~\eqref{fo:MKV-discounted_cost} by
\begin{equation}
	\label{eq:def-cost-C-theta}
	C(\theta) = J(\bctrl^\theta).
\end{equation}
Lemma~\ref{le:admissible_of_control_with_feedback_form_on_theta} states that for a control parameter $\theta \in \Theta$, the corresponding parameterized control process $\bctrl^\theta$ belongs to $\cU_{ad}$ (and thus $\bctrl^\theta \in \cU_{ad}^\Theta$), and that the controlled state process $\bX^{\theta} \in \cX$.
For this reformulation, we use Assumption~\ref{ass:optimal-feedback-in-Theta}, which ensures that the optimal feedback coefficients belong to the admissible parameter class.
Proposition~\ref{pr:reformulation_of_MF_problem} below justifies reformulating the search for the optimal control using parameters in $\Theta$.

\begin{proposition}
	\label{pr:reformulation_of_MF_problem}
	Under Assumptions~\ref{ass:finit_cost},~\ref{ass:positivity-qr}, and~\ref{ass:optimal-feedback-in-Theta}, we have
	\begin{equation}
		\label{eq:equivalence_C_and_J_in_MFC}
		\inf_{\theta \in \Theta} C(\theta) = \inf_{\bctrl \in \cU_{ad}} J(\bctrl).
	\end{equation}
\end{proposition}

\begin{proof}
	Let $\theta^* = (K^*, L^*)$ be a control parameter defined in Theorem~\ref{th:existence_linear_control}.
	By Assumption~\ref{ass:optimal-feedback-in-Theta}, $\theta^* \in \Theta$.
	Theorem~\ref{th:existence_linear_control} states that the control process $\bctrl^{\theta^*}$ minimizes the discounted MF cost $J(\bctrl)$ among all admissible controls in $\cU_{ad}$.
	Thus,
	$$
		J(\bctrl^{\theta^*}) = C(\theta^*) \geq \inf_{\theta \in \Theta} C(\theta)
		= \inf_{\bctrl^\theta \in \cU_{ad}^\Theta} J(\bctrl^\theta) \geq \inf_{\bctrl \in \cU_{ad}} J(\bctrl) = J(\bctrl^{\theta^*}).
	$$
	We then conclude the proposition.
\qed\end{proof}

Henceforth, we focus on searching for a minimizer of $C$ in the set $\Theta$ to find the optimal control process for the MFC problem.
We will sometimes adopt a slight abuse of terminology and refer to $\theta \in \Theta$ as a ``control" or say that the state process is ``controlled'' by $\theta$.

In the following, we cast the problem of minimizing the LQMFC cost $C(\theta)$ into two standard LQ problems depending only on the parameters $K$ and $L$, respectively.
For a state process $(X_{n}^{\theta})_{n\geq 0}$ controlled by $\theta = (K, L) \in \Theta$ following dynamics~\eqref{fo:MKV-state} and starting from the initial position $X_0^\theta = \varepsilon_0^0 + \varepsilon_0$, we introduce two auxiliary processes, $\bY = (Y_{n})_{n\geq 0}$ and $\bZ = (Z_{n})_{n\geq 0}$, defined by:
\begin{equation}
	\label{eq:def-Y-Z-fct-X}
	Y_{n} := X^{\theta}_{n} - \bar{X}^{\theta}_{n},
	\qquad\text{and}\qquad
	Z_{n} := \bar{X}^{\theta}_{n},
\end{equation}
for $n \geq 0$, where the dynamics of processes $\bY$ and $\bZ$ satisfy
\begin{align}
	Y_{n+1} & = ({\mathrm A} - {\mathrm B} K) Y_{n} + \varepsilon_{n+1},                    &  & Y_{0} = \varepsilon_0 - \EE[\varepsilon_0];
	\label{eq:dyn_y_theta}
	\\
	Z_{n+1} & = ( \tilde {\mathrm A} - \tilde {\mathrm B}  L) Z_{n}  + \varepsilon^0_{n+1}, &  & Z_{0} = \varepsilon^0_0 + \EE[\varepsilon_0].
	\label{eq:dyn_z_theta}
\end{align}
The processes $\bY$ and $\bZ$ depend on the parameters $K$ and $L$, respectively. The noise processes and their initial states involved in~\eqref{eq:dyn_y_theta} and~\eqref{eq:dyn_z_theta} are independent.
As a result, the processes $\bY$ and $\bZ$ are independent.
We introduce the corresponding cost functions $C_y, C_z: \RR^{\ell \times d} \to \RR$ as:
\begin{equation}
	\label{eq:def_Cy_Cz}
	C_y(K) = \mathbb{E} \big[  \sum_{n \geq 0} \gamma^{n}  (Y_{n})^\top (Q + K^\top R K) Y_{n} \big] , \quad
	C_z(L) = \mathbb{E} \big[  \sum_{n \geq 0} \gamma^{n}   (Z_{n})^\top (\tilde{Q} + L^\top \tilde{R} L) Z_{n}  \big].
\end{equation}
We may sometimes use $\bY^K = (Y_{n}^K)_{n \geq 0}$ and $\bZ^L = (Z_{n}^L)_{n \geq 0}$ to emphasize the dependence on $K$ and $L$ for the processes $\bY$ and $\bZ$.

Let $\Theta_K = \{ K \, | \, \gamma \| {\mathrm A} - {\mathrm B} K \|^2 < 1 \}$ and $\Theta_L = \{ L \, | \, \gamma \| \tilde {\mathrm A} - \tilde {\mathrm B} L \|^2 < 1 \}$.
The following lemma shows the benefit of the re-parametrization \eqref{eq:def-Y-Z-fct-X} inspired by Theorem~\ref{th:existence_linear_control}.

\begin{lemma}
	\label{lem:decompose-cost-reparam}
	For any $\theta = (K,L) \in \Theta$, we have $C(\theta) = C_y(K) + C_z(L),$ and
	\begin{equation}
		\label{eq:reparametrization_of_optimality_LQMF}
		\inf_{ \theta = (K, L) \in \Theta } C(\theta)=\inf_{ K \in \Theta_K} C_y(K) +\inf_{L \in \Theta_L} C_z(L).
	\end{equation}
\end{lemma}

\begin{proof}
	For a parameterized control process $\bctrl^\theta = (\ctrl_{n}^\theta)_{n\geq 0}$ with $\theta = (K, L)$, it can be rewritten as $\ctrl_{n}^\theta = - K Y_{n} - L Z_{n}$ for every $n \geq 0$.
	Thus, the instantaneous cost at any time $n \geq 0$ satisfies $\EE \big[ c(X_{n}^\theta, \bar X_{n}^\theta, \ctrl_{n}^\theta, \bar \ctrl_{n}^\theta) \big] = \EE \big[ (Y_{n})^\top (Q + K^\top R K) Y_{n} + (Z_{n})^\top (\tilde Q + L^\top \tilde R L) Z_{n} \big]$.
	So we have $C(\theta) = C_y(K) + C_z(L)$.
	Moreover, the processes $\bY$ and $\bZ$ are independent because of the independence between the noise sequences $\bvarepsilon$ and $\bvarepsilon^0$, and the independence between the initial perturbations $\varepsilon_0$ and $\varepsilon_0^0$.
	Because $ \Theta_K \times \Theta_L \subset \Theta$, together with the fact that $\bY$ and $\bZ$ are independent processes controlled by $K$ and $L$ respectively, we can then separate the optimization over $\Theta$ into two optimization problems over $\Theta_K$ and $\Theta_L$ for $\bY$ and $\bZ$.
\qed\end{proof}

\subsection{Additional notations}
\label{subsection:PG_additional_notations}

For an admissible control parameter $\theta = (K, L) \in \Theta$, let $\bY$ and $\bZ$ be the two auxiliary processes following dynamics~\eqref{eq:dyn_y_theta}--\eqref{eq:dyn_z_theta}.
The symbols $\Sigma_K^y$ and $\Sigma_L^z$ stand for the infinite discounted sums of variance/covariance matrices given by
\begin{equation}
	\label{eq:variance_matrices_y_and_z}
	\Sigma_K^y := \EE \Big[ \sum \nolimits_{n \geq 0} \gamma^{n} Y_{n} (Y_{n})^\top \Big], \qquad  \Sigma_L^z := \EE \Big[ \sum \nolimits_{n \geq 0} \gamma^{n} Z_{n} (Z_{n})^\top \Big].
\end{equation}
By a slight abuse of terminology, we call the matrices $\Sigma_K^y$ and $\Sigma_L^z$ the variance matrices of the auxiliary processes $\bY$ and $\bZ$.
To alleviate notation, we may omit the superscript in $\Sigma_K^y$ and $\Sigma_L^z$ and use $\Sigma_K, \Sigma_L$ when the context is clear.
With the expression of $C_y(K)$ and $C_z(L)$ in~\eqref{eq:def_Cy_Cz}, we have
\begin{equation}
	\label{eq:expression_cost_sigma}
	C_y(K) = \langle Q + K^\top R K, \, \Sigma_K \rangle_{tr}, \qquad C_z(L) = \langle \tilde Q + L^\top \tilde R L ,\, \Sigma_L  \rangle_{tr}.
\end{equation}
We denote the exact gradients of the cost $C(\theta)$ with respect to $K$ and $L$ by $\nabla_K C(\theta)$ and $\nabla_L C(\theta)$, respectively.
Lemma~\ref{lem:decompose-cost-reparam} implies immediately that
$$
	\nabla_K C(\theta) = \nabla_K C_y(K), \quad \nabla_L C(\theta) = \nabla_L C_z(L).
$$
For the quantitative PG estimates below, we use the following stronger state-cost condition than Assumption~\ref{ass:positivity-qr}.
\begin{assumption}[State-cost coercivity]
\label{ass:pg-state-cost-coercivity}
	The matrices $Q$ and $\tilde Q$ are positive definite.
\end{assumption}
We also make the following non-degeneracy assumption.
\begin{assumption}[Non-degeneracy] We assume that
\label{ass:non-deg}
    \begin{equation}
    \label{eq:non_degenerate_condition}
        \max \Big\{ \lambda_{min} \big(\Sigma_{Y_{0}} \big), \lambda_{min} \big( \Sigma^1 \big) \Big\} > 0,
        \quad
        \max \Big\{ \lambda_{min} \big( \Sigma_{Z_{0}} \big), \lambda_{min} \big(\Sigma^0 \big) \Big\} > 0,
    \end{equation}
    where $Y_{0} = \varepsilon_0 - \EE[ \varepsilon_0 ]$, $Z_{0} = \varepsilon_0^0 + \EE[ \varepsilon_0]$, and 
    $
    \Sigma_{Y_{0}} = \mathbb{E} \left[Y_{0} (Y_{0})^\top \right]
    $,
    $ 
      \Sigma_{Z_{0}} = \mathbb{E} \left[Z_{0} (Z_{0})^\top \right]
    $,
    and where the variance matrices of $(\varepsilon_{n+1}, \varepsilon_{n+1}^0)$ for any $n \geq 0$ in the i.i.d. noise sequences $(\bvarepsilon, \bvarepsilon^0)$ are given by
    $
      \Sigma^1 = \EE[\varepsilon_{n+1} (\varepsilon_{n+1})^\top]
    $,
    $
      \Sigma^0 = \EE[\varepsilon^0_{n+1} (\varepsilon^0_{n+1})^\top]
    $.
\end{assumption}
Under the non-degenerate conditions~\eqref{eq:non_degenerate_condition} in Assumption~\ref{ass:non-deg}, we have
\begin{equation}
	\label{eq:min_eigenvalue_of_variance_matrices}
	\lambda_y^1 := \lambda_{min} \big( \Sigma_{Y_{0}} + \frac{\gamma}{1-\gamma} \Sigma^1 \big) > 0
	, \quad
	\lambda_z^0 := \lambda_{min}\big( \Sigma_{Z_{0}} + \frac{\gamma}{1-\gamma} \Sigma^0 \big) > 0.
\end{equation}
The strict positivity of the above two quantities $\lambda_y^1$ and $\lambda_z^0$ is crucial for establishing a uniform upper bound for the norms of the gradients $\nabla_K C_y(K)$ and $\nabla_L C_z(L)$ within the PG algorithms, and thus they are essential for the global convergence of the PG algorithms.

Since we assume sub-Gaussian distributions for the initial perturbations in the LQ problem, the random variables $Y_{0}$ and $Z_{0}$ are also sub-Gaussian.
We let $C_{init, noise}$ denote an upper bound on the sub-Gaussian norms of $Y_{0}$, $Z_{0}$, and the i.i.d. noise terms $(\varepsilon_{n+1})_{n \geq 0}$ and $(\varepsilon_{n+1}^0)_{n \geq 0}$. That is, for every $n \geq 0$:
\begin{equation}
	\max\big\{ \| Y_{0} \|_{\psi_2}, \| Z_{0} \|_{\psi_2}, \| \varepsilon_{n+1} \|_{\psi_2}, \| \varepsilon_{n+1}^0 \|_{\psi_2} \big\} \leq C_{init, noise}.
\end{equation}

\begin{remark}
	\label{remark:assumption_on_variance}
	In the absence of noise processes, if the initial perturbations $(\varepsilon_0, \varepsilon_0^0)$ are non-degenerate with bounded supports (implying Assumption~\ref{ass:non-deg}), the arguments in~\cite{fazel2018global} can be applied to processes $\bY$ and $\bZ$, respectively, to show the convergence of the PG algorithm using the gradients $(\nabla_K C(\theta), \nabla_L C(\theta) )$.
	In the more general setting considered here, particularly with non-degenerate sub-Gaussian distributions, the strict positivity of $\lambda_y^1$ and $\lambda_z^0$ allows us to derive convergence results without enforcing assumptions on the initial states as in previous work~\cite{fazel2018global}.
	It is important to emphasize that the presence of noise terms poses additional difficulties in proving the convergence of PG algorithms in model-free scenarios.
	The challenges arise not only from the unbounded noise process but also from the quadratic forms involving sub-Gaussian random vectors $Y_{n}$ or $Z_{n}$ whose entries are not independent.
	We cannot simply rely on the classical matrix Bernstein inequality or the Hanson-Wright inequality~\cite{vershynin2018high} for independent random variables; thus, we require more delicate arguments to demonstrate convergence.
\end{remark}

\subsection{Exact PG for MFC}
\label{subsection:exact_PG_for_MFC}

Assuming full knowledge of the model, we compute the optimal parameters using a standard gradient descent algorithm.
With a fixed learning rate $\eta>0$ and initial parameter $\theta^{(0)} = (K_0,L_0)$, we update the parameter from $\theta^{(k)}$ to $\theta^{(k+1)} = (K^{(k+1)}, L^{(k+1)})$ via \defi{exact PG updates}\index[sub]{Exact PG updates}:
\begin{equation}
	\label{eq:one-step-update-exact-pg}
	\left\{
	\begin{array}{rcl}
		K^{(k+1)} & = & K^{(k)} - \eta \nabla_K C_y(K^{(k)}), \\
		L^{(k+1)} & = & L^{(k)} - \eta \nabla_L C_z(L^{(k)}).
	\end{array}
	\right.
\end{equation}
To show the convergence of $C(\theta^{(k)})$ towards the optimal cost $C(\theta^*)$ with the scheme \eqref{eq:one-step-update-exact-pg}, we adapt the gradient descent methods with the Polyak-Lojasiewicz condition~\cite{karimi2016linear, polyak1963gradient, lojasiewicz1963topological} for the nonconvex functions $C_y$ and $C_z$.
Below, we first state the Polyak-Lojasiewicz condition for the auxiliary costs $C_y$ and $C_z$ in Proposition~\ref{pr:PL_cond_Cy_Cz}, followed by a local smoothness condition in the direction of the gradients in Proposition~\ref{pr:local_smoothness_Cy_Cz}.
The proofs of these results can be found in Sections~\ref{subsection:proof_proposition_PL_condition} and \ref{subsection:proof_proposition_local_smoothness}, respectively.

\vskip 6pt
We recall that a differentiable function $f: D \to \RR$ on an open set $D \subseteq \RR^m$ is said to satisfy the \textit{$\nu$-Polyak-Lojasiewicz} ($\nu$-PL for short)\index[not]{PL}\index[sub]{Polyak-Lojasiewicz condition} condition on $D$ for a constant $\nu > 0$ if
\begin{equation}
	\label{eq:PL_condition_definition}
	\nu \big( f(x) - f(x^*) \big) \leq \| \nabla_x f(x) \|^2, \qquad \forall \, x \in D
\end{equation}
where $x^* \in D$ is a global minimum in $D$, i.e., $f(x) \geq f(x^*)$ for all $x \in D$.

\begin{proposition}[PL condition]
	\label{pr:PL_cond_Cy_Cz}
	The cost functions $C_y$ and $C_z$ given in~\eqref{eq:def_Cy_Cz} satisfy the PL condition on $\Theta$: for $\theta = (K, L) \in \Theta$, we have
	\begin{equation}
		\label{eq:PL_cond_Cy_Cz_w_lambda}
		\left\{
		\begin{array}{rcl}
			\nu_{y} \big( C_y(K) - C_y(K^*) \big) & \leq & \| \nabla_K C_y(K) \|_F^2, \\
			\nu_{z} \big( C_z(L) - C_z(L^*) \big) & \leq & \| \nabla_L C_z(L) \|_F^2.
		\end{array}
		\right.
	\end{equation}
	with constants $\nu_y =  4 (\lambda_y^1)^2 \lambda_{min}(R) / \| \Sigma_{K^*} \| $ and $\nu_z =  4 (\lambda_z^0)^2 \lambda_{min}(\tilde R) / \| \Sigma_{L^*} \| $.
\end{proposition}

\begin{proposition}[Local smoothness in exact PG direction]
	\label{pr:local_smoothness_Cy_Cz}
	Consider $\theta=(K,L) \in \Theta$.
	Consider $K'  = K - \eta \nabla_K C_y(K)$ and $L'  = L - \eta \nabla_L C_z(L)$ such that $\theta'=(K', L') \in \Theta$.
	Let $\Delta K = K' - K$ and $\Delta L = L' - L$, then
	\begin{equation}
		\label{eq:local_smoothness_Cy_Cz}
		\left\{
		\begin{array}{rl}
			C_y(K') - C_y(K) & \leq  \big( 1 - \lambda_{var, y} \big) \langle \nabla_K C_y(K),\, \Delta K \rangle_{tr} + \lambda_{hess, y}  \| \Delta K \|_F^2
			\\
			C_z(L') - C_z(L) & \leq   \big( 1 - \lambda_{var, z} \big) \langle \nabla_L C_z(L),\, \Delta L \rangle_{tr} + \lambda_{hess, z}  \| \Delta L \|_F^2
		\end{array}
		\right.
	\end{equation}
	where
	$
		\lambda_{var, y} =  \lambda_{var, y}(K, K') = \| \Sigma_{K'} - \Sigma_{K} \| / \lambda_y^1
	$,
	$
		\lambda_{hess, y} =  \lambda_{hess, y}(K, K') =  \| \Sigma_{K'} \| (\| R \|  + \gamma \|{\mathrm B}\|^2 C_y(K) / \lambda_y^1 )
	$, and $\lambda_{var, z}, \lambda_{hess, z}$ are defined similarly with $\big(\Sigma_{L'}, \Sigma_L, \lambda_z^0, \tilde R, \tilde {\mathrm B}, C_z(L) \big)$.

\end{proposition}

Now, we consider four technical constants $h_{var, y}(C_0)$, $h_{hess, y}(C_0)$, $h_{var, z}(C_0)$, and $h_{hess, z}(C_0)$ defined in Lemma~\ref{lemma:bound_on_coefficent_in_prop_local_smoothness} in Section~\ref{section:proof_of_exact_pg}, where we conduct a detailed perturbation analysis of the control parameters.
These constants are polynomials of $( C_0, 1/\lambda_y^1, 1/\lambda_z^0)$ and other model parameters.
The following Lemma~\ref{lemma:small_learning_rate} states that when the learning rate $\eta$ is small enough, the perturbed control parameter $\theta' = (K', L')$ in Proposition~\ref{pr:local_smoothness_Cy_Cz} is admissible and the coefficients related to local smoothness in equation~\eqref{eq:local_smoothness_Cy_Cz} are bounded.
Its proof is found in Lemma~\ref{lemma:bound_on_coefficent_in_prop_local_smoothness}.

\begin{lemma}[Admissibility]
	\label{lemma:small_learning_rate}
	Consider $\theta = (K, L) \in \Theta$ with $C(\theta) \leq C_0$ for some constant $C_0 \in \RR$.
	If the learning rate $\eta$ is small enough, namely, satisfies
	$$
		\eta \leq \min\{ h_{var, y}(C_0)^{-1} , h_{var, z}(C_0)^{-1} \},
	$$
	and
$$
\|{\mathrm B}\|\,\eta h_{grad,y}(C_0) < 1/\sqrt{\gamma}-\|{\mathrm A}-{\mathrm B}K\|,
\quad
\|\tilde{\mathrm B}\|\,\eta h_{grad,z}(C_0) < 1/\sqrt{\gamma}-\|\tilde{\mathrm A}-\tilde{\mathrm B}L\|,
$$
	then the perturbed control parameter $\theta' = (K', L')$ in Proposition~\ref{pr:local_smoothness_Cy_Cz} with learning rate $\eta$ is admissible in the sense that $\theta ' \in \Theta$.
	Moreover,
	\begin{align*}
		\lambda_{var, y}(K, K') + \eta \lambda_{hess, y}(K, K') & \leq \eta \big( h_{var, y}(C_0) + h_{hess,y}(C_0) \big)
		\\
		\lambda_{var, z}(L, L') + \eta \lambda_{hess, z}(L, L') & \leq \eta \big( h_{var, z}(C_0) + h_{hess, z}(C_0) \big).
	\end{align*}
\end{lemma}

We are now ready to state and prove a global convergence result based on the PL condition and a local smoothness condition on the cost functions.
Consider the sequence $(\theta^{(k)})_{k \geq 0}$ generated by the PG updates~\eqref{eq:one-step-update-exact-pg}.

\begin{theorem}[Exact PG convergence]
	\label{th:exact-CV}
	Let $C_0 = C(\theta^{(0)})$ and $\nu_{pl} = \min\{ \nu_y, \nu_z \} / 2$ .
	Assume that the sublevel set $\{\theta=(K,L)\in\Theta: C(\theta)\le C_0\}$ has positive Euclidean admissibility margins
$$
	m_y(C_0) := \inf_{\theta=(K,L)\in\Theta:\, C(\theta)\le C_0} \big( 1/\sqrt{\gamma}-\|{\mathrm A}-{\mathrm B}K\| \big) > 0,
$$
and
$$
	m_z(C_0) := \inf_{\theta=(K,L)\in\Theta:\, C(\theta)\le C_0} \big( 1/\sqrt{\gamma}-\|\tilde{\mathrm A}-\tilde{\mathrm B}L\| \big) > 0.
$$
	We assume that the learning rate $\eta$ is small enough such that
	\begin{equation}
		\label{eq:condition_learning_rate_exact_pg}
		\begin{split}
			\eta \leq \frac{1}{4} \min \big\{ &h_{var, y}(C_0)^{-1}, h_{hess, y}(C_0)^{-1} ,h_{var, z}(C_0)^{-1} , h_{hess, z}(C_0)^{-1}, 
			\\
			&m_y(C_0)/(\|{\mathrm B}\|h_{grad,y}(C_0)), m_z(C_0)/(\|\tilde{\mathrm B}\|h_{grad,z}(C_0)) \big\}.
		\end{split}
	\end{equation}
	Under Assumptions~\ref{ass:finit_cost},~\ref{ass:positivity-qr},~\ref{ass:optimal-feedback-in-Theta},~\ref{ass:pg-state-cost-coercivity}, and~\ref{ass:non-deg}, for every $k \geq 0$, we have $C(\theta^{(k+1)}) \leq C(\theta^{(k)}) \leq C_0$, and
	\begin{equation}
		\label{eq:linear_convergence_exact_pg}
		C(\theta^{(k+1)}) - C(\theta^*)  \leq \big( 1 - \eta \nu_{pl} \big) \big(  C(\theta^{(k)}) - C(\theta^*) \big).
	\end{equation}
	Moreover, for every $\varepsilon > 0$, when the iteration steps
	$
		k \geq \frac{1}{\eta \nu_{pl} } \log \big( \frac{1}{\varepsilon} \big),
	$
	we have
	\begin{equation}
		\label{eq:epsilon_approx_of_exact_pg}
		C(\theta^{(k)}) - C(\theta^*)  \leq \varepsilon  \big(  C(\theta^{(0)}) - C(\theta^*) \big).
	\end{equation}
\end{theorem}

\begin{proof}

	We assume that at iteration step $j \geq 0$, the control parameter $\theta^{(j)}$ is admissible, and $C(\theta^{(j)}) \leq C_0$.
	By the definitions of $m_y(C_0)$ and $m_z(C_0)$ and by~\eqref{eq:condition_learning_rate_exact_pg}, the margin conditions in Lemma~\ref{lemma:small_learning_rate} hold for $\theta^{(j)}$.
	Lemma~\ref{lemma:small_learning_rate} shows that $\theta^{(j+1)} \in \Theta$ and
	$$
		\lambda_{var, y} + \eta \lambda_{hess, y} \leq \eta \big( h_{var, y}(C_0) + h_{hess, y}(C_0) \big) \leq \frac{1}{2}.
	$$
	Because $Q \succ 0$ and $R \succ 0$, the expression of cost~\eqref{eq:expression_cost_sigma} implies that
	$
		\lambda_y^1 \leq \| \Sigma_{K^*} \| \leq \frac{ C(K^*) }{ \lambda_{min}(Q) } \leq \frac{ C_0 }{ \lambda_{min}(Q)}.
	$
	Together with the definition of $h_{hess, y}(C_0) = \frac{2 C_0}{\lambda_{min}(Q)} \big( \| R \| + \gamma \| {\mathrm B} \|^2 \frac{C_0}{\lambda_y^1} \big)$ in~\eqref{eq:h_hess_y}, and the definition of $\nu_y$ in Proposition~\ref{pr:PL_cond_Cy_Cz}, we have
	$$
		\frac{1}{2} \eta \nu_y = 2 \eta \lambda_{min}(R) \lambda_y^1 \frac{\lambda_y^1}{\| \Sigma_{K^*} \|} \leq  2 \eta \| R \| \frac{ C_0} { \lambda_{min}(Q) } \leq \eta h_{hess, y}(C_0) \leq \frac{1}{4} < 1.
	$$
	Consequently, with the local smoothness condition~\eqref{eq:local_smoothness_Cy_Cz} and the $\nu_y$-PL condition~\eqref{eq:PL_cond_Cy_Cz_w_lambda}, we obtain
	\begin{align*}
		C_y( K^{ (j+1)} ) - C_y ( K^{(j)} ) & \leq - \eta \big( 1 - \lambda_{var, y} - \eta \lambda_{hess, y} \big) \| \nabla_K C_y(K^{(j)})\|_F^2
		\\
		                                    & \leq - \frac{1}{2} \eta \nu_y \big( C_y( K^{ (j)} ) - C_y ( K^* ) \big).
	\end{align*}
	The difference in costs yields directly that $C_y(K^{(j+1)}) \leq C_y(K^{(j)})$, and
	$
		C_y(K^{(j+1)} ) - C_y(K^*) \leq (1 - \eta \nu_{pl} ) \big( C_y( K^{(j)}) - C_y(K^*) \big).
	$
	Applying the same arguments to $C_z$, and adding up the two inequalities for $C_y$ and $C_z$, we obtain equation~\eqref{eq:linear_convergence_exact_pg} for $k=j$ and $C(\theta^{(j+1)}) \leq C(\theta^{(j)}) \leq C_0$.

	Furthermore, by noticing that $ \frac{1}{\eta \nu_{pl}} \log ( \frac{1}{1 - \eta \nu_{pl}}) \geq 1$, when $k \geq \frac{1}{\eta \nu_{pl}} \log( \frac{1}{\varepsilon} )$, we have
	$$
		\log \Big( \frac{ C(\theta^{(0)}) - C(\theta^*) } {C(\theta^{(k)}) - C(\theta^*)  } \Big) \geq k \log( \frac{1}{1 - \eta \nu_{pl}} )  \geq  \log( \frac{1}{\varepsilon} ),
	$$
	which concludes the proof of the theorem.
\qed\end{proof}

\begin{remark}
	\label{remark:importance_of_non_degeneracy_assumption}
	It is worth highlighting that Assumption~\ref{ass:non-deg} on the randomness plays a crucial role in establishing the PL condition and the smoothness condition mentioned in Propositions~\ref{pr:PL_cond_Cy_Cz} and \ref{pr:local_smoothness_Cy_Cz}.
	Notably, it enables convergence even when the system starts from a deterministic initial condition.
\end{remark}

\subsection{Model-free PG for MFC with MKV simulator}
\label{subsection:modelfree_MKV_simulator}

We now turn to a modification of the PG updates~\eqref{eq:one-step-update-exact-pg}, in which the gradient will be approximated in a model-free way.
The new update scheme takes the form of \defi{PG updates with MKV simulator}\index[sub]{PG updates with MKV simulator}:
\begin{equation}
	\label{eq:one-step-update-model-free-MKV}
	\left\{
	\begin{array}{rcl}
		K^{(k+1)} & = & K^{(k)} - \eta \tilde{\nabla}^{T,M,\tau}_K (\theta^{(k)}),
		\\
		L^{(k+1)} & = & L^{(k)} - \eta \tilde{\nabla}^{T,M,\tau}_L (\theta^{(k)}),
	\end{array}
	\right.
\end{equation}
where $\tilde \nabla^{T, M, \tau}(\theta) = \big( \tilde{\nabla}^{T, M,\tau}_K(\theta), \tilde{\nabla}^{T, M,\tau}_L(\theta) \big)$ is defined below, based on samples generated by a simulator.

Let us assume that we have access to the following (stochastic) simulator, called MKV simulator $\cS^T_{MKV}$: given an admissible control parameter $\theta \in \Theta$, $\cS^T_{MKV}(\theta)$ returns a sample of the discounted MF cost (see Section~\ref{subsec:OC-MKV} and equation~\eqref{fo:MKV-discounted_cost}) for the MKV dynamics~\eqref{fo:MKV-state} controlled by the process $(\ctrl_{n}^\theta)_{n \geq 0}$ and truncated at a finite horizon $T$.
In other words, it returns a realization of the truncated MF cost
\begin{equation}
	\label{eq:MKV_simulator_{n}runcated_MF_cost}
	\tilde C^T(\theta) := \sum_{n=0}^{T-1} \gamma^{n} c_{n}(\theta),
\end{equation}
where $c_{n}(\theta) =  c\big(X_{n}^{\theta}, \bar{X}_{n}^{\theta}, \ctrl_{n}^\theta, \bar{\ctrl}_{n}^\theta \big)$ is the instantaneous cost at time $n$ (see equation~\eqref{fo:lq_one_step_cost}).

While the simulator uses the model to generate samples, we aim to develop a model-free technique to compute the optimal control parameter.
Without full knowledge of the model, we use derivative-free methods (see Section~\ref{sec:app-proof-modelfree-MKV-CV}) to estimate the cost gradient.
Here, \emph{model-free} refers to the gradient estimator, which uses perturbed policy parameters and sampled costs without evaluating model coefficients or derivatives of the dynamics.
The MKV simulator remains an idealized oracle: it generates trajectories with the exact conditional-mean process $\bar X_n^\theta$ needed by the dynamics, feedback and costs, and returns the sampled truncated cost $\tilde C^T(\theta)$.
Access to transitions and costs along a single trajectory with an unknown transition kernel does not generally provide this conditional-mean information.
The main idea is to evaluate the function to minimize on a set of perturbed parameters so that one can move towards a weighted average direction that improves the value of the target function~\cite{MR3627456}.
It has been shown that a careful implementation of ES can provide a performance improvement compared with classical RL algorithms based on back-propagation.
Moreover, the method can be easily scaled up and adapted to distributed computer systems \cite{salimans2017evolution,choromanski2018structured}.
We define some notation for the gradient estimation.
Let $\mathbb{B}_\tau\subset \mathbb{R}^{\ell \times d}$ be the ball of radius $\tau$ under the Frobenius norm centered at the origin, and $\mathbb{S}_\tau = \partial \mathbb{B}_\tau$ be its boundary.
The uniform distributions on $\mathbb{B}_\tau$ and $\mathbb{S}_\tau$ are denoted by $\mu_{\mathbb{B}_\tau}$ and $\mu_{\mathbb{S}_\tau}$ respectively.
We should think of $v = (v^{(idy)}, v^{(com)}) \in \SS_\tau \times \SS_\tau$ as a direction for a perturbation of a parameter $\theta = (K, L) \in \Theta$.
For $M$ perturbation directions $\underline{v} = (v_i)_{i=1}^M = \big( (v_i^{(idy)}, v_i^{(com)}) \big)_{i=1}^M$ in which $v_i^{(idy)}$ and $v_i^{(com)}$ are independently sampled from $\mu_{\SS_\tau}$, we denote by $\theta_i = (K + v^{(idy)}_i, L + v^{(com)}_i)$ the control parameter $\theta$ perturbed with $v_i$.
Then, the approximated gradient based on MKV simulator $\cS^T_{MKV}$ and $M$ perturbation directions $\underline{v}$ on $\SS_{\tau} \times \SS_{\tau}$ is defined by
\begin{equation}
	\label{eq:definition_approx_gradient_MKV_simulator}
	\tilde{\nabla}^{T, M,\tau}_K (\theta) = \frac{\ell d}{\tau^2}\frac{1}{M} \sum_{i=1}^M \tilde{C}^T(\theta_i) v^{(idy)}_i,
	\quad
	\tilde{\nabla}^{T, M,\tau}_L (\theta) = \frac{\ell d}{\tau^2}\frac{1}{M} \sum_{i=1}^M \tilde{C}^T(\theta_i) v^{(com)}_i.
\end{equation}
We omit the perturbation directions $\underline v$ in equation~\eqref{eq:definition_approx_gradient_MKV_simulator} for clarity here.
Each call to $\cS^T_{MKV}(\theta_i)$ starts a new trajectory from the prescribed initial law under the queried perturbed policy.
The initial and noise variables for that call are sampled afresh, independently of past calls and of the current perturbation directions, with common and idiosyncratic components generated as specified by the model.
This ability to restart and query the simulator at perturbed policies is part of the access assumption.

\begin{algorithm}
	\caption{Model-free MKV-Based Gradient Estimation}
	\label{algo:MKVestim}
	\begin{algorithmic}
		\STATE {\bfseries Data:} {Parameter $\theta = (K,L)$;  truncation horizon $T$; number of perturbations $M$; perturbation radius $\tau$}
		\STATE {\bfseries Result:} {An approximated gradient $\tilde \nabla^{T, M, \tau}(\theta)$ of $\nabla C(\theta)$}
		\FOR{$i = 1, 2, \dots, M$}
		\STATE Sample $v_{i}^{(idy)}, v_{i}^{(com)}$ i.i.d. $\sim \mu_{\mathbb{S}_\tau}$\;
		\STATE Set $\theta_i = (K_i, L_i) = (K+ v_{i}^{(idy)}, L+ v_{i}^{(com)})$\;
		\STATE Sample $\tilde{C}^T(\theta_i) = \sum_{n=0}^{T-1} \gamma^{n} c_{n}(\theta_i)$ using  $\cS^T_{MKV}(\theta_i)$\;
		\ENDFOR
		\STATE {\bfseries Set } {$\tilde{\nabla}^{T, M, \tau}_K (\theta), \tilde{\nabla}^{T, M, \tau}_L (\theta)$ with equation~\eqref{eq:definition_approx_gradient_MKV_simulator}}
		\STATE {\bfseries Return: }{$\tilde{\nabla}^{T, M, \tau} (\theta) = \big( \tilde{\nabla}^{T, M, \tau}_K (\theta), \tilde{\nabla}^{T, M, \tau}_L (\theta) \big)$}
	\end{algorithmic}
\end{algorithm}

Algorithm~\ref{algo:MKVestim} provides a biased estimator of the true PG $\nabla C(\theta)$.
This method follows the spirit of \cite[Algorithm 1]{fazel2018global}, except that here we have two components playing different roles for the state and the conditional mean.
The approximate gradient $\tilde{\nabla}^{T, M, \tau}(\theta)$ thus depends on three additional hyperparameters, $(T, M, \tau)$: the cost truncation horizon, the number of perturbation directions, and the perturbation radius.
In Proposition~\ref{pr:gradient_approx_with_MKV_simulator} below, we provide conditions on $(T, M, \tau)$ so that Algorithm~\ref{algo:MKVestim} yields a sufficiently good gradient estimator to ensure the convergence of the approximate PG.

\begin{proposition}[gradient approximation with MKV simulator]
	\label{pr:gradient_approx_with_MKV_simulator}
	Consider $\theta \in \Theta$ with $C(\theta) \leq C_0$ for some constant $C_0 \in \RR$.
	Let $\tilde{\varepsilon} > 0$ be a target precision and $\delta_{approx} \in (0,1)$.
	We assume that the parameters $(T, M, \tau)$ in Algorithm~\ref{algo:MKVestim} satisfy
	\begin{align}
		\tau^{-1} & \geq \phi_{pert, radius, MKV}(\tilde{\varepsilon}, C_0),
		\label{eq:model_free_pg_condition_pert_radius}
		\\
		T         & \geq \phi_{trunc, T, MKV}(\tilde{\varepsilon}, \tau, C_0),
		\label{eq:model_free_pg_condition_{n}runcation_T}
		\\
		M         & \geq \phi_{sample, size, MKV} \big(\tilde{\varepsilon}, \tau, T, C_0, \delta_{approx} \big),
		\label{eq:model_free_pg_condition_sample_size}
	\end{align}
	where $\phi_{pert, radius, MKV}$, $\phi_{trunc, T, MKV}$, and $\phi_{sample, size, MKV}$ are defined in Section~\ref{subsection:proof_of_proposition_approx_gradient_w_MKV_simulator} as polynomials of $(d$, $\ell$, $C_0$, $1/\lambda_y^1$, $1/\lambda_z^0$, $C_{init, noise})$ and other model parameters.
	Then, we have
	$$
		\PP \big( \big\| \tilde \nabla^{T, M, \tau}(\theta) - \nabla C(\theta) \big\|  > \tilde{\varepsilon} \big) \leq \delta_{approx},
	$$
	where $\big\| \tilde \nabla^{T, M, \tau}(\theta) - \nabla C(\theta) \big\| = \big\| \tilde \nabla^{T, M, \tau}_K(\theta) - \nabla_K C(\theta) \big\| + \big\| \tilde \nabla^{T, M, \tau}_L(\theta) - \nabla_L C(\theta) \big\|$, and the randomness in $ \tilde \nabla^{T, M, \tau}(\theta)$ comes from the simulator $\cS_{MKV}^T$ and the perturbation directions $\underline v \in (\SS_\tau \times \SS_\tau)^{M}$.
\end{proposition}
The proof of Proposition~\ref{pr:gradient_approx_with_MKV_simulator} is deferred to Section~\ref{subsection:proof_of_proposition_approx_gradient_w_MKV_simulator} after we demonstrate the approximation of the exact PG using derivative-free techniques based on parameter perturbations (Lemma~\ref{lemma:approx_perturbed_gradient_MKV_simulator}), the approximation of the infinite-horizon cost with a truncated finite-horizon cost (Lemma~\ref{lemma:approx_nruncated_policy_gradient_MKV_simulator}), and the approximation of the expected truncated cost using sampled costs from a simulator (Lemma~\ref{lemma:approx_sampled_pg_with_{n}runcated_pg}).
It is worth noting that the bounds on $(\tau, T, M)$ in~\eqref{eq:model_free_pg_condition_pert_radius},~\eqref{eq:model_free_pg_condition_{n}runcation_T}, and~\eqref{eq:model_free_pg_condition_sample_size} are independent of the control parameter $\theta \in \Theta$ as long as $C(\theta) \leq C_0$.

Now, we are ready to state the convergence results in the model-free setting with the MKV simulator $\cS_{MKV}^T$.
We consider a bound for the learning rate:
\begin{align*}
	&\phi_{lrate, MKV}(C_0) 
	\\
	&=\min \Big\{ \min \{ h_{var, y}^{-1}(C_0), h_{var, z}^{-1}(C_0), h_{hess, y}^{-1}(C_0), h_{hess, z}^{-1}(C_0), 
	\\
	& \qquad\qquad\qquad m_y(C_0)/(\|{\mathrm B}\|h_{grad,y}(C_0)), m_z(C_0)/(\|\tilde{\mathrm B}\|h_{grad,z}(C_0)) \} / 4,
	\\
	& \qquad\qquad \min\big\{ h_{small-pert, y}^{-1}(C_0), h_{small-pert, z}^{-1}(C_0) \big\}\, (2 h_{cost}(C_0)) / \nu_{pl}, 
	\\
	&
	 \qquad\qquad m_y(C_0)h_{cost}(C_0)/(\nu_{pl}\|{\mathrm B}\|), m_z(C_0)h_{cost}(C_0)/(\nu_{pl}\|\tilde{\mathrm B}\|)
	  \Big\},
\end{align*}
where the constants $h_{var, y}(C_0)$, $h_{var, z}(C_0)$, $h_{hess, y}(C_0)$, $h_{hess, z}(C_0)$, $h_{small-pert,y}(C_0)$, $h_{small-pert, z}(C_0)$ depending on $C_0$ are defined in Lemma~\ref{lemma:useful_bounds_I} for small perturbation of control parameters, and $h_{cost}(C_0)$ is defined in Lemma~\ref{lemma:perturbation_cost_Cy}.
The margins $m_y(C_0)$ and $m_z(C_0)$ are defined as in Theorem~\ref{th:exact-CV}.
We consider a sequence of control parameters $(\theta^{(k)})_{k \geq 0}$ generated with model-free PG update scheme~\eqref{eq:one-step-update-model-free-MKV} with the approximated gradients $\big(\tilde \nabla^{T, M, \tau} (\theta^{(k)}) \big)_{k \geq 0}$ from Algorithm~\ref{algo:MKVestim} using an MKV simulator $\cS_{MKV}^T$.
For the two model-free convergence theorems below (Theorems~\ref{th:modelfree-MKV-CV} and~\ref{th:modelfree-POP-CV}), fix a deterministic integer $H\geq1$ as the number of updates and a total failure budget $\delta\in(0,1)$, and set $\delta_{approx}=\delta/H$.
For each run, let $\mathcal H_k$ be the sigma-field generated by the initialization and all batches before update $k$, so that the current parameter is $\mathcal H_k$-measurable.
At each update, sample fresh perturbation directions and simulator initial and noise variables, independently of $\mathcal H_k$, with the within-batch independence specified by the model and the estimator.

\begin{theorem}
	\label{th:modelfree-MKV-CV}
	Consider a target precision $0< \varepsilon \leq 1$ and an initial parameter $\theta^{(0)} \in \Theta$ with cost $C(\theta^{(0)}) = C_0$. Assume that the sublevel set $\{\theta\in\Theta:\, C(\theta)\le C_0+1\}$ has positive Euclidean margins $m_y(C_0+1)$ and $m_z(C_0+1)$. Consider a learning rate satisfying
	\begin{equation}
		\label{eq:learning_rate_condition_pg_MKV}
		\eta \leq  \phi_{lrate, MKV}(C_0 + 1).
	\end{equation}
	We choose simulation parameters $\big(\tau, T, M \big)$ in Algorithm~\ref{algo:MKVestim} such that they satisfy equations~\eqref{eq:model_free_pg_condition_pert_radius}, \eqref{eq:model_free_pg_condition_{n}runcation_T}, and \eqref{eq:model_free_pg_condition_sample_size} in Proposition~\ref{pr:gradient_approx_with_MKV_simulator} with cost bound $C_0+1$ and parameters $\tilde{\varepsilon} = \varepsilon  \nu_{pl}/ (2 h_{cost}(C_0 + 1))$ and $\delta_{approx}=\delta/H$.
	In particular, the sample size satisfies
\[
M\geq\phi_{sample, size, MKV}
\big(\tilde\varepsilon,\tau,T,C_0+1,\delta/H\big).
\]
Its perturbation and sampled-cost terms each use failure probability $\delta/(2H)$; for fixed values of the other arguments, the sufficient sample size grows logarithmically in $H/\delta$.
	Then, under our standing assumptions, %
	for every $0\leq k<H$, conditional on $\mathcal H_k$, if $\theta^{(k)}\in\Theta$ and $C(\theta^{(k)})\leq C_0+1$, then we have with probability at least $1-\delta/H$ 
	that $\theta^{(k+1)} \in \Theta$, $C(\theta^{(k+1)}) \leq C_0 + 1$, and
	\begin{equation}
		\label{eq:contraction_pg_MKV}
		C(\theta^{(k+1)}) - C(\theta^*) \leq \big( 1 -  \eta \nu_{pl}/ 2 \big) \max \Big\{ C(\theta^{(k)}) - C(\theta^*) , \, \varepsilon \Big\}.
	\end{equation}
	Moreover, %
	for any $0\leq k\leq H$ satisfying
	$$
		k \geq \frac{2}{\eta \nu_{pl}} \log( \frac{C_0 - C(\theta^*)}{ \varepsilon} ),
	$$
	we have, with probability at least $1-\delta$, that
	\begin{equation}
		C(\theta^{(k)}) - C(\theta^*) \leq \varepsilon.
	\end{equation}
\end{theorem}

\begin{remark}
	The convergence rate for the model-free PG algorithm with an MKV simulator is linear, meaning that for $k=\mathcal{O}(\log(1/\varepsilon))$ we attain an $\varepsilon$-approximation to the optimal cost with probability at least $1-\delta$, provided $k\leq H$.
	With $\delta_{approx}=\delta/H$, the two terms in
\eqref{eq:phi_MKV_sample_size} use
$\delta_{pert}=\delta_{sample}=\delta/(2H)$.
The logarithmic factors in \eqref{eq:phi_pert_size} and
\eqref{eq:phi_sample_size} therefore become
$\log(4(\ell+d)H/\delta)$ and
$\log(8T(\ell+d)H/\delta)$, respectively.
For fixed values of the other arguments, the resulting bound on $M$
depends logarithmically on $H/\delta$.
\end{remark}

\begin{proof}
	Suppose that we are at iteration step $k$, with $0\leq k<H$, of the PG algorithm with MKV simulation $\cS_{MKV}^T(\theta^{(k)})$.
	Let us condition on $\mathcal H_k$, with $\theta^{(k)}\in\Theta$ and $C(\theta^{(k)})\leq C_0+1$.
All probabilities in the one-step argument below are conditional on this history, and all approximation bounds use cost level $C_0+1$.

	We consider first $\theta^{(exact, k)} = (K^{(exact, k)}, L^{(exact, k)})$ generated using the exact PG update scheme~\eqref{eq:one-step-update-exact-pg}: $K^{(exact, k)} = K^{(k)} - \eta \nabla_K C(\theta^{(k)})$, $L^{(exact, k)} = L^{(k)} - \eta \nabla_L C(\theta^{(k)})$.
	Under the assumption on $\eta$ in inequality~\eqref{eq:learning_rate_condition_pg_MKV}, and by the definition of $\phi_{lrate, MKV}(C_0+1)$, the learning rate $\eta$ is then small enough to satisfy condition~\eqref{eq:condition_learning_rate_exact_pg} with $C_0$ replaced by $C_0+1$ in Theorem~\ref{th:exact-CV} for the exact-PG algorithm. This implies that
	$$
		C(\theta^{(exact, k)}) - C(\theta^{*}) \leq ( 1- \eta \nu_{pl}) \big( C(\theta^{(k)}) - C(\theta^*) \big),
	$$
	where $\nu_{pl} = \min\{ \nu_y, \nu_z \} / 2$ is a coefficient related to the PL conditions.
	We also have $C(\theta^{(exact, k)}) \leq C(\theta^{(k)}) \leq C_0 + 1$.

	Now, we consider a control parameter $\theta^{(k+1)} = (K^{(k+1)}, L^{(k+1)})$ generated from the one-step model-free update scheme~\eqref{eq:one-step-update-model-free-MKV}.
	We apply Proposition~\ref{pr:gradient_approx_with_MKV_simulator} on $\theta = \theta^{(k)}$ with $\tilde{\varepsilon} = \kappa \varepsilon$, where $\kappa = \nu_{pl} / \big( 2 h_{cost}(C_0 + 1) \big)$ and with $\delta_{approx}=\delta/H$. 
	If $(T, M, \tau)$ satisfy equations \eqref{eq:model_free_pg_condition_{n}runcation_T}, \eqref{eq:model_free_pg_condition_sample_size}, and \eqref{eq:model_free_pg_condition_pert_radius} for $(\tilde \varepsilon, \delta_{approx})$, we then deduce that with probability at least $1 - \delta_{approx}$,
	$$
		\| \theta^{(k+1)} - \theta^{(exact, k)} \| = \eta  \| \tilde \nabla^{M, T, \tau}(\theta^{(k)}) - \nabla C(\theta^{(k)}) \|  \leq \eta \kappa \varepsilon.
	$$
	When the learning rate $\eta$ is small enough to satisfy~\eqref{eq:learning_rate_condition_pg_MKV}, we have $\eta \kappa. h_{small-pert, y}(C_0 + 1) \leq 1$ and $ \eta \kappa. h_{small-pert, z}(C_0 + 1)  \leq 1$.
	Moreover, by the margin terms in $\phi_{lrate, MKV}(C_0+1)$ and since $\varepsilon \leq 1$,
$$
\|{\mathrm B}\|\,\|K^{(k+1)}-K^{(exact,k)}\| \leq m_y(C_0+1)/2,
\quad
\|\tilde{\mathrm B}\|\,\|L^{(k+1)}-L^{(exact,k)}\| \leq m_z(C_0+1)/2.
$$
	Given that $C(\theta^{(exact, k)}) \leq C(\theta^{(k)}) \leq  C_0 + 1$, the quantity $f_{var}(K^{(exact, k)})$ (resp. $f_{var}(L^{(exact, k)})$) defined in equation~\eqref{eq:f_var_K} is bounded by a constant $h_{small-pert,y}(C_0 + 1)$ (resp. $h_{small-pert, z}(C_0 + 1)$) defined in Lemma~\ref{lemma:useful_bounds_I}.
	As a consequence, with probability at least $1 -\delta_{approx}$, we have
	\begin{align*}
		     & \max \Big\{ \| K^{(k+1)} - K^{(exact, k)} \| f_{var}(K^{(exact, k)}), \,  \| L^{(k+1)} - L^{(exact, k)} \| f_{var}(L^{(exact, k)}) \Big\}
		\\
		\leq & \| \theta^{(k+1)} - \theta^{(exact, k)} \| \max \big\{ f_{var}(K^{(exact, k)}), f_{var}(L^{(exact, k)}) \big\}
		\\
		\leq & \eta \kappa \varepsilon. \max\big\{ h_{small-pert, y}(C_0 + 1), h_{small-pert, z}(C_0 + 1) \big\}
		\\
		\leq & \varepsilon 
        \\
        \leq & 1.
	\end{align*}
	Together with the positive margins at level $C_0+1$, these bounds imply the condition in Lemma~\ref{lemma:stability_of_small_perturbation} for $K^{(k+1)}$ and, analogously, for $L^{(k+1)}$; hence $K^{(k+1)} \in \Theta_{K}$ and $L^{(k+1)} \in \Theta_{L}$.
	Hence, with probability at least $1 - \delta_{approx}$, the control parameter $\theta^{(k+1)}$ is admissible in the sense that $\theta^{(k+1)} \in \Theta$.

	To show the inequality~\eqref{eq:contraction_pg_MKV}, we apply Lemma~\ref{lemma:perturbation_cost_Cy} to the cost function $C(\theta)$ with the small perturbation of parameter $\theta^{(k+1)}$ with respect to $\theta^{(exact, k)}$ mentioned above, and we obtain:
	\begin{align*}
		\Big| C \big(\theta^{(k+1)} \big) - C \big( \theta^{(exact, k)} \big) \Big|
		 & \leq  h_{cost}(C_0 + 1) \big\| \theta^{(k+1)} - \theta^{(exact, k)} \big\|
		\\
		 & \leq h_{cost}(C_0 + 1) \eta \kappa  \varepsilon
		\\
		 & = \eta \nu_{pl} \varepsilon / 2.
	\end{align*}
	Then with probability at least $1 - \delta_{approx}$,
	\begin{align*}
		C(\theta^{(k+1)}) - C(\theta^*)
		 & =  \big( C(\theta^{(exact, k)} )- C(\theta^*) \big) +  \big( C(\theta^{(k+1)}) - C(\theta^{(exact, k)}) \big)
		\\
		 & \leq \big( 1 - \eta \nu_{pl} \big) \big( C(\theta^{(k)}) - C(\theta^*) \big) + (\eta \nu_{pl}/2) \varepsilon
	\end{align*}

	If at iteration step $k$ we have $C(\theta^{(k)}) - C(\theta^*) > \varepsilon$, we deduce that with probability at least $1 - \delta_{approx}$,
	$$
		C(\theta^{(k+1)}) - C(\theta^*) \leq  (1 - \frac{\eta \nu_{pl}}{2} ) \big( C(\theta^{(k)}) - C(\theta^*) \big).
	$$
	In this case, we have a contraction inequality moving from $\theta^{(k)}$ to $\theta^{(k+1)}$.
	Because $C(\theta^*) \leq C(\theta^{(k)}) \leq C_0 + 1$, we have
	$$
		C(\theta^{(k+1)})\leq (1 - \eta \nu_{pl}/2) C(\theta^{(k)}) + \frac{\eta \nu_{pl}}{2} C(\theta^*) \leq C(\theta^{(k)} ) \leq C_0 + 1.
	$$

	Otherwise, if at iteration step $k$, we reach the target precision $\varepsilon$ such that $C(\theta^{(k)}) - C(\theta^*) \leq \varepsilon$, then with probability at least $1 - \delta_{approx}$,
	$$
		C(\theta^{(k+1)}) - C(\theta^*) \leq ( 1 - \eta \nu_{pl} ) \varepsilon + \eta \nu_{pl} /2 \varepsilon \leq \big( 1 - \frac{\eta \nu_{pl}}{2} \big) \varepsilon \leq \varepsilon.
	$$
	In this case, the difference between $C(\theta^{(k+1)})$ and $C(\theta^*)$ are still bounded by $\varepsilon$.
	Because $C(\theta^*) \leq C(\theta^{(0)}) = C_0$ and $\varepsilon \leq 1$, we have
	$$
		C(\theta^{(k+1)}) \leq  C(\theta^*)  + ( 1 - \frac{\eta \nu_{pl}}{2}) \varepsilon \leq C_0 + 1.
	$$

Let $G_j$ be the event that the first $j$ gradient estimates meet their error bound $\tilde\varepsilon=\kappa\varepsilon$, with $G_0=\Omega$.
Define these events recursively, evaluating each new error only on $G_j$.
On $G_j$, the one-step argument gives $\theta^{(j)}\in\Theta$ and $C(\theta^{(j)})\leq C_0+1$.
By fresh sampling and Proposition~\ref{pr:gradient_approx_with_MKV_simulator}, the conditional probability that the next estimate fails its error bound is at most $\delta/H$ on $G_j$.
Hence $\PP(G_j\setminus G_{j+1})\leq\delta/H$, and
\[
\PP(G_H^c)
=\sum_{j=0}^{H-1}\PP(G_j\setminus G_{j+1})
\leq\delta.
\]
On $G_H$, all iterates through step $H$ are admissible and have cost at most $C_0+1$, which justifies the fixed simulation parameters.
Fix $0\leq k\leq H$ satisfying the iteration bound in the theorem; the remaining argument is on $G_H$.
	If, for some $k' < k$, $C(\theta^{(k')}) - C(\theta^*) \leq \varepsilon$, we have $C( \theta^{(j)}) - C(\theta^*) \leq \varepsilon$ for all $j =k', \ldots, k$.
	Otherwise, we have a contraction inequality for all $j=0,\ldots, k-1$ and then
	$$
		C(\theta^{(k)}) - C(\theta^*) \leq (1 - \frac{\eta \nu_{pl}}{2})^k \big(C_0 - C(\theta^*) \big) \leq \varepsilon.
	$$
	This conclusion holds on $G_H$, whose probability is at least $1-\delta$.
\qed\end{proof}

\subsection{Model-free PG for MFC  with population simulator}
\label{subsec:PG-popsimu}

We now turn our attention to a more realistic setting where one does not have access to an oracle for simulating the MKV dynamics, but only to an oracle capable of simulating the evolution of $N$ agents.
We then use the state sample average instead of the theoretical conditional mean, and for the social cost, the empirical average instead of the mean-field cost provided by the MKV simulator.

We rely on the following population simulator $\cS^{T,N}_{pop}$: given a control parameter $\theta$, $\cS^{T,N}_{pop}(\theta)$ returns a sample of the social cost obtained by generating realizations of the $N$-agent state trajectories controlled by $\theta$ and computing the associated cost for the population, see~\eqref{fo:N-multi_state} and~\eqref{fo:N-cost}.
In other words, it returns a realization of
$$
	\tilde{C}^{T,N}(\theta) = \sum_{n=0}^{T-1} \gamma^{n}  \frac{1}{N} \sum_{j=1}^N c \big( X^{(j), \theta}_{n},  \bar{X}^{N, \theta}_{n}, \ctrl^{(j), \theta}_{n}, \bar \ctrl^{N, \theta}_{n} \big)
$$
where $c$ is the one-step cost function~\eqref{fo:lq_one_step_cost}, and $(X_{n}^{(j), \theta}, \ctrl_{n}^{(j), \theta})_{j=1}^N$ are the state-action pairs at time $n$ for the $N$ agents who adopt the same policy function so that
$
	\ctrl_{n}^{(j), \theta} =  -K \big( X_{n}^{(j), \theta} - \bar{X}_{n}^{N, \theta} \big) - L \bar{X}^{N, \theta}_{n}
$
for every $n \geq 0$ and $j=1, \ldots, N$, and where $\bar{X}_{n}^{N, \theta}$ and $\bar{\ctrl}_{n}^{N, \theta}$ are the average state and average control at time $n$.
Let $C^N(\theta) = J^N(\underline{\bU}^\theta)$ in equation~\eqref{fo:social_cost_of_population} when all agents adopt the same control parameter $\theta \in \Theta$.

Notice that this population simulator $\cS^{T,N}_{pop}$ is arguably more realistic, though less powerful, than the previous MKV simulator $\cS^T_{MKV}$.
The former uses only a noisy approximation of the true mean processes, while the MKV simulator generates the exact means.
We emphasize that, in this population simulator, \emph{all} agents adopt the same control and therefore use the \emph{same} perturbed version of the control parameter in Algorithm~\ref{algo:POPestim}.
This design aligns with the idea that the problem corresponds to an optimization problem for a central planner or a group of cooperative agents using the same control rule to minimize the social cost.
In a finite-agent implementation, the empirical state mean must be available through population observations or computed internally by the simulator.
Within each simulated population, the common noise is shared by all agents and the idiosyncratic noises are independent; separate simulator calls use fresh population initializations and noise draws.
The finite-$N$ approximation is additional to horizon truncation and Monte Carlo sampling; see Propositions~\ref{pr:approx_with_social_cost_no_MF_cost} and~\ref{pr:modelfree_population_gradient_approx}.

For the one-step gradient update scheme with population simulator $\cS_{pop}^{T, N}$, we replace the gradient $\tilde{\nabla}^{T, M, \tau}(\theta)$ obtained from Algorithm~\ref{algo:MKVestim} with an MKV simulator by another zero-th order approximation $\tilde{\nabla}^{T, N, M,\tau}(\theta)$ of the gradient based on $\cS_{pop}^{T, N}$ generated by Algorithm~\ref{algo:POPestim}.
The term $\tilde{\nabla}^{T, N, M, \tau}(\theta) = \big( \tilde{\nabla}^{T, N, M, \tau}_K( \theta),  \tilde{\nabla}^{T, N, M, \tau}_L(\theta) \big)$ is called the \emph{sampled population policy gradient} at $\theta \in \Theta$.
It is defined with $M$ perturbation directions $(v_i)_{i=1}^M = (v_i^{(idy)}, v_i^{(com)})_{i=1}^M$ on $\SS_\tau \times \SS_\tau$:
\begin{equation}
	\label{eq:sampled_population_gradient}
	\tilde{\nabla}_K^{T, N, M, \tau}(\theta) =  \frac{\ell d}{\tau^2} \frac{1}{M} \sum_{i=1}^M  \tilde{C}^{T, N}(\theta_i) v_i^{(idy)},
	\quad
	\tilde{\nabla}_L^{T, N, M, \tau}(\theta) = \frac{\ell d}{\tau^2} \frac{1}{M} \sum_{i=1}^M \tilde{C}^{T, N}(\theta_i) v_i^{(com)},
\end{equation}
where $\theta_i = \theta + v_i$ is the perturbed control parameter with $v_i$.
The \defi{PG update with population simulator}\index[sub]{PG update with population simulator} at iteration step $k$ for a parameter $\theta^{(k), pop} = (K^{(k), pop}, L^{(k), pop} )$ is defined by
\begin{equation}
	\label{eq:one-step-update-model-free-PoP}
	\left\{
	\begin{array}{rcl}
		K^{(k+1), pop} & = & K^{(k), pop} - \eta \tilde{\nabla}^{T, N, M,\tau}_K (\theta^{(k), pop}),
		\\
		L^{(k+1), pop} & = & L^{(k), pop} - \eta \tilde{\nabla}^{T, N, M,\tau}_L (\theta^{(k), pop}).
	\end{array}
	\right.
\end{equation}

\begin{algorithm}
	\caption{Model-free Population-Based Gradient Estimation}
	\label{algo:POPestim}
	\begin{algorithmic}
		\STATE {\bfseries Data:} {Parameter $\theta = (K,L)$;  truncation horizon $T$; number of perturbations $M$; perturbation radius $\tau$; number of agents $N$.}
		\STATE {\bfseries Result:} {An approximation of $\nabla C(\theta)$.}
		\FOR{$i = 1, 2, \dots, M$}
		\STATE Sample $v_i^{(idy)}, v_i^{(com)}$ i.i.d. $\sim \mu_{\mathbb{S}_\tau}$\;
		\STATE Set $\theta_i = \big( K+ v_i^{(idy)}, L+ v_i^{(com)} \big)$ \;
		\STATE Sample $\tilde{C}^{T, N}(\theta_i) = \sum_{n=0}^{T-1} \gamma^{n} \frac{1}{N} \sum_{j=1}^N  c( X_{n}^{(j), \theta_i}, \bar{X}_n^{N, \theta_i}, \ctrl_{n}^{(j), \theta_i}, \bar \ctrl_{n}^{N, \theta_i})$ for $N$ agents using  $\cS^{T,N}_{pop}(\theta_i)$ \;
		\ENDFOR
		\STATE {\bfseries Set} {$\tilde{\nabla}_K^{T, N, M, \tau}(\theta)$ and $\tilde{\nabla}_L^{T, N, M, \tau}(\theta)$ with equation~\eqref{eq:sampled_population_gradient}}\;
		\STATE {\bfseries Return: }{$\tilde{\nabla}^{T, N, M, \tau}(\theta) = \Big(\tilde{\nabla}_K^{T, N, M, \tau}(\theta), \tilde{\nabla}_L^{T, N, M, \tau}(\theta)  \Big)$}
	\end{algorithmic}
\end{algorithm}

The main result of this section is Theorem~\ref{th:modelfree-POP-CV}.
It establishes the convergence of the above learning scheme toward the optimal control parameter $\theta^* = (K^*, L^*)$ of the MF problem, using a population simulator $\cS_{pop}^{T, N}$ of the social cost of $N$ homogeneous agents.
Before proving this result, we state Proposition~\ref{pr:approx_with_social_cost_no_MF_cost} and Proposition~\ref{pr:modelfree_population_gradient_approx}, which provide crucial approximation estimates for the social cost by its $N$-agent equivalent, as well as the approximation of the sampled population policy gradient.
Their proofs are given in Section~\ref{subsection:proof_of_approx_modelfree_pop_gradient}.

\begin{proposition}
	\label{pr:approx_with_social_cost_no_MF_cost}
	Consider $\theta \in \Theta$ with $C(\theta) \leq C_0$ for some $C_0 \in \RR$.
	Under our standing assumptions, we have
	$$
		\big| C^N(\theta) - C(\theta) \big| \leq \frac1N\phi_{social-cost, factor}(C_0)
	$$
	where $\phi_{social-cost, factor}(C_0) =  2 d C_0 C_{init, noise}^2 \big( 1/ \lambda_y^1 + 1/ \lambda_z^0 \big) / (1-\gamma)$ is a constant depending on $C_0$.
	Here $C^N(\theta)$ and $C(\theta)$ are expected infinite-horizon costs, so the estimate is deterministic.
\end{proposition}

\begin{proposition}
	\label{pr:modelfree_population_gradient_approx}
	Let $\theta \in \Theta$ with $C(\theta) \leq C_0$ for some $C_0 \in \RR$, choose a target precision $\tilde \varepsilon > 0$ and $\delta_{approx} \in (0,1)$, and let us assume the parameters $(T, N, M, \tau)$ in Algorithm~\ref{algo:POPestim} satisfy
	\begin{align}
		\tau^{-1} & \geq \phi_{pert, radius}(\tilde \varepsilon / 4, C_0)
		\\
		N         & \geq \phi_{agent, size, pop}(\tilde \varepsilon / 4, \tau, C_0)
		\\
		T         & \geq \phi_{trunc, T, pop}(\tilde\varepsilon / 4 , \tau, C_0, N) + 2
		\\
		M         & \geq \max \big\{ \phi_{pert, size}(\tilde\varepsilon / 4, \tau, C_0, \delta_{approx}/2),
		\nonumber                                                                                                            \\
		          & \hspace{60pt} \phi_{sample, size, pop}(\tilde\varepsilon / 4, \tau, T, C_0, \delta_{approx} /2  ) \big\}
	\end{align}
	where $\phi_{pert, radius}$ defined in~\eqref{eq:phi_pert_radius}, and $\phi_{agent, size, pop}$, $\phi_{trunc, T, pop}$, $\phi_{sample, size, pop}$ in Section~\ref{subsection:proof_of_approx_modelfree_pop_gradient} are polynomials in $(d$, $\ell$, $C_0$, $1/ \lambda_{y}^1$, $1/\lambda_z^0$, $C_{init, noise}$, $1/N)$ and other model parameters.
	Then we have
	\begin{equation}
		\PP \big( \| \tilde \nabla^{T,N,M,\tau}(\theta) - \nabla C(\theta) \| > \tilde \varepsilon \big) \leq \delta_{approx}.
	\end{equation}
\end{proposition}

Now, we show that Algorithm~\ref{algo:POPestim} successfully learns the optimal control parameter of the MF problem with a global linear convergence rate in the social cost. To state the theorem, we consider a sequence of control parameters $(\theta^{(k), pop})_{k \geq 0}$ generated by the model-free PG update scheme with the sampled population policy gradient~\eqref{eq:one-step-update-model-free-PoP}.

\begin{theorem}
	\label{th:modelfree-POP-CV}
	We consider an initial control parameter $\theta^{(0), pop} \in \Theta$ with $C(\theta^{(0), pop}) = C_0$ and a target precision $\varepsilon \leq 1$.
	Assume that the sublevel set $\{\theta\in\Theta:\, C(\theta)\le C_0+1\}$ has positive Euclidean margins $m_y(C_0+1)$ and $m_z(C_0+1)$.
	We assume that the number of agents $N$ and the learning rate $\eta$ satisfy
	\begin{align}
		\label{eq:learning_rate_modelfree_pop_simulator}
		\eta & \leq \phi_{lrate, MKV}(C_0 + 1)
		\\
		N    & \geq \phi_{social-cost, factor}(C_0 + 1) / (\rho \varepsilon)
		\label{eq:N_agent_social_cost_approx_condition}
	\end{align}
	where $\rho = (\eta \nu_{pl}) / (16 - 4 \eta \nu_{pl}) \in (0, 1)$.
	Assume also that $(1+2\rho)\varepsilon\leq1$.
	We choose simulation parameters $(T, N, M, \tau)$ in Algorithm~\ref{algo:POPestim} to satisfy conditions in Proposition~\ref{pr:modelfree_population_gradient_approx} with cost bound $C_0+1$ and $\tilde \varepsilon = ( 1 + 2 \rho) \varepsilon \nu_{pl} / \big( 2 h_{cost}(C_0 + 1) \big)$ and $\delta_{approx}=\delta/H$.
	In particular, the sample size satisfies
\[
M\geq\max\Big\{
\phi_{pert, size}(\tilde\varepsilon/4,\tau,C_0+1,\delta/(2H)),
\phi_{sample, size, pop}(\tilde\varepsilon/4,\tau,T,C_0+1,\delta/(2H))
\Big\}.
\]
For fixed values of the other arguments, the sufficient sample size grows logarithmically in $H/\delta$.
	Then, under our standing assumptions, %
	for every $0\leq k<H$, conditional on $\mathcal H_k$, if $\theta^{(k), pop}\in\Theta$ and $C(\theta^{(k), pop})\leq C_0+1$, then we have with probability at least $1-\delta/H$
	 that $\theta^{(k+1), pop} \in \Theta$, $C(\theta^{(k+1), pop})\leq C_0+1$, and
	\begin{equation}
		\label{eq:contraction_inequality_on_C^N_algon_iteration}
		\big| C^N( \theta^{(k+1), pop}) - C^N(\theta^*) \big| \leq \big( 1 - \eta \nu_{pl} / 4 \big) \max \Big\{ \big| C^N( \theta^{(k), pop} ) - C^N(\theta^*) \big|, \, \varepsilon \Big\}.
	\end{equation}
	Moreover, for any $0\leq k\leq H$ satisfying
	$
		k \geq \frac{4}{\eta \nu_{pl}} \log\max\Big\{1,\frac{|C^N(\theta^{(0), pop})-C^N(\theta^*)|}{\varepsilon}\Big\} %
	$
	we have with probability at least $1-\delta$ that
	$$
		\big| C^N(\theta^{(k), pop}) - C^N(\theta^*) \big| \leq \varepsilon.
	$$
\end{theorem}

\begin{proof}
	We use the one-step argument in the proof of Theorem~\ref{th:modelfree-MKV-CV}, with the population gradient estimator.
Fix $0\leq k<H$ and condition on $\mathcal H_k$.
All probabilities in the one-step argument below are conditional on this history, and all approximation bounds use cost level $C_0+1$.
	At iteration step $k=0$, we have by assumption that $\theta^{(0), pop}\in \Theta$ and $C(\theta^{(0), pop}) = C_0 \leq C_0 + 1$.
	Suppose that at iteration step $k \geq 0$, we have $\theta^{(k), pop} \in \Theta$ and $C(\theta^{(k), pop}) \leq C_0 + 1$, that is, the control parameter generated by~\eqref{eq:one-step-update-model-free-PoP} at step $k$ with a population simulator $\cS_{pop}^T$ is admissible.
	We first observe that  Proposition~\ref{pr:approx_with_social_cost_no_MF_cost} with a large $N$ satisfying~\eqref{eq:N_agent_social_cost_approx_condition} implies
	\begin{align*}
		     & C ( \theta^{(k), pop} ) - C (\theta^*)
		\\
		\leq & \big|  C^N(\theta^{(k), pop}) - C^N(\theta^*) \big| + \big| C^{N}(\theta^{(k), pop}) - C(\theta^{(k), pop}) \big| + \big| C(\theta^*) - C^N(\theta^*) \big|
		\\
		\leq & \big| C^N(\theta^{(k), pop}) - C^N(\theta^*) \big| + 2\rho \varepsilon.
	\end{align*}
	Secondly, let $\varepsilon' = ( 1 + 2\rho) \varepsilon$ and $\kappa = \frac{\nu_{pl}}{2 h_{cost}(C_0 + 1)}$.
	If we choose a target precision $\tilde{\varepsilon} = \varepsilon' \kappa$ and a parameter $\delta_{approx}=\delta/H$
	in Proposition~\ref{pr:modelfree_population_gradient_approx} to approximate the true gradient using a population simulator with coefficients $(T, N, M, \tau)$ defined accordingly, and applying the same arguments as in Theorem~\ref{th:modelfree-MKV-CV} under the condition $C(\theta^{(k), pop}) \leq C_0 + 1$ and $\theta^{(k), pop} \in \Theta$, it follows with high probability (at least $1- \delta_{approx}$) that $\theta^{(k+1), pop} \in \Theta$, $C(\theta^{(k+1), pop}) \leq C_0 + 1$, and
	\begin{align*}
		C( \theta^{(k+1), pop} ) - C (\theta^*) & \leq \big( 1 - \frac{\eta \nu_{pl}}{2} \big) \max \Big\{ C ( \theta^{(k), pop} ) - C (\theta^*) , \, (1 + 2 \rho) \varepsilon \Big\}
		\\
		                                        & \leq  \big( 1 - \frac{\eta \nu_{pl}}{2} \big) \Big( \max \Big\{ \big| C^N(\theta^{(k), pop}) - C^N(\theta^*) \big|, \, \varepsilon \Big\} + 2 \rho \varepsilon \Big).
	\end{align*}
	Besides, from the choice of the parameter $\rho$, we have
	$
		2 \rho + 2 ( 1 - \frac{\eta \nu_{pl}}{2} ) \rho  = \eta \nu_{pl}  / 4.
	$
	Therefore, by applying Proposition~\ref{pr:approx_with_social_cost_no_MF_cost} to $\theta^{(k+1), pop} \in \Theta$ with $C(\theta^{(k+1), pop}) \leq C_0 + 1$, we have with probability at least $1 - \delta_{approx}$ that
	\begin{align*}
		     & \big| C^N( \theta^{(k+1), pop} ) - C^N(\theta^*) \big|
		\\
		\leq & 2 \rho \varepsilon + C ( \theta^{(k+1), pop} ) - C(\theta^*)
		\\
		\leq & 2 \rho \varepsilon +  \big( 1 - \frac{\eta \nu_{pl}}{2} \big) \Big( \max \Big\{ \big| C^N(\theta^{(k), pop}) - C^N(\theta^*) \big|, \, \varepsilon \Big\} + 2 \rho \varepsilon \Big)
		\\
		=    & \big( 2 \rho \varepsilon + 2 ( 1 - \frac{\eta \nu_{pl}}{2} ) \rho \varepsilon \big) + \big( 1 - \frac{\eta \nu_{pl}}{2} \big) \max \Big\{ \big| C^N(\theta^{(k), pop}) - C^N(\theta^*) \big|, \, \varepsilon \Big\}
		\\
		\leq & \big( 1 - \eta \nu_{pl} / 4 \big) \max \Big\{ \big| C^N(\theta^{(k), pop}) - C^N(\theta^*) \big|, \, \varepsilon \Big\}.
	\end{align*}
	We thus demonstrate the inequality~\eqref{eq:contraction_inequality_on_C^N_algon_iteration} between the $N$-agent costs from steps $k$ and $k+1$ along Algorithm~\ref{algo:POPestim} using a population simulator $\cS_{pop}^T$ with constant parameters $(\eta, T, N, M, \tau)$.

Set $G_0^{pop}=\Omega$.
For $0\leq j<H$, let $G_{j+1}^{pop}$ be the subset of $G_j^{pop}$ on which the population gradient estimate at step $j$ has error at most $\tilde\varepsilon=\kappa\varepsilon'$.
On $G_j^{pop}$, the one-step argument gives $\theta^{(j),pop}\in\Theta$ and $C(\theta^{(j),pop})\leq C_0+1$, since $\varepsilon'\leq1$.
The conditional failure probability at the next update is at most $\delta/H$ by Proposition~\ref{pr:modelfree_population_gradient_approx}.
The same first-failure calculation gives $\PP(G_H^{pop})\geq1-\delta$.
Let us fix $0\leq k\leq H$ satisfying the iteration bound in the theorem, and let us work on $G_H^{pop}$.
If the social-cost error reaches $\varepsilon$ at some step $k'<k$, inequality~\eqref{eq:contraction_inequality_on_C^N_algon_iteration} keeps it at most $\varepsilon$ through step $k$; otherwise, we have
	$$
		| C^N(\theta^{(k), pop}) - C^N(\theta^*) | \leq (1 - \eta \nu_{pl} / 4 )^k | C^N(\theta^{(0), pop}) - C^N(\theta^*) | \leq \varepsilon
	$$
	This conclusion holds on $G_H^{pop}$, whose probability is at least $1-\delta$.
\qed\end{proof}

\begin{remark}
	In the approximation of the social cost with the MF cost in Proposition~\ref{pr:approx_with_social_cost_no_MF_cost}, the coefficient $\phi_{social-cost, factor}$ is defined using $\lambda_y^1 > 0$ and $\lambda_z^0 > 0$ under Assumption~\ref{ass:non-deg}.
	However, in the $N$-agent dynamics, the sample average of the states at time $n$, $\bar{X}_{n}^N$, would remain stochastic even in the absence of common noise, due to the noise term $\frac{1}{N} \sum_{j=1}^N  \varepsilon^{(j)}_{n+1}$ in the dynamical equation of the sample mean process.
	Hence, we expect that Theorem~\ref{th:modelfree-POP-CV} could be proven without the non-degeneracy of the noise terms in Assumption~\ref{ass:non-deg}.
\end{remark}

\section{Numerical Results}
\label{sec:numerics-LQ}

In this section, we provide numerical illustrations of the RL algorithms developed for the LQMFC problem. We compare the performance of the exact PG (model-based) with the model-free approach using both the MKV simulator and the finite-population simulator.
The code is available at: \url{https://github.com/mlauriere/mfrl-tutorial-code}.
We first examine a scalar model over six population sizes, then consider a coupled two-dimensional example in Section~\ref{subsec:lq-numerics-2d}.

\subsection{Experimental Setup}
We consider a scalar LQMFC model where the dynamics and costs are symmetric across the agent population. The representative agent's objective is to minimize the socially optimal cost by learning the parameters $\theta = (K, L)$, which determine the decentralized control $\ctrl_{n} = -K(X_{n} - \bar{X}_{n}) - L\bar{X}_{n}$. 

The model parameters and hyperparameters used for the simulation are summarized in Table~\ref{tab:lq_mfc_params}. The dynamics are characterized by coefficients $A, \bar{A}, B, \bar{B}$, and the cost function by $Q, \bar{Q}, R, \bar{R}$. The noise in the dynamics is modeled as zero-mean Gaussian with variances $\sigma_\epsilon^2$ and $\sigma_{\epsilon_0}^2$ for the idiosyncratic and common noise, respectively. The common and idiosyncratic initial components are independent and uniformly distributed on $[-1,1]$.

\begin{table}[H]
\centering
\begin{tabular}{|l|c|}
\hline
\textbf{Parameter} & \textbf{Value} \\
\hline
Dynamics coefficients $(A, \bar{A})$ & $(0.5, 0.5)$ \\
Control coefficients $(B, \bar{B})$ & $(0.5, 0.5)$ \\
State cost coefficients $(Q, \bar{Q})$ & $(0.5, 0.5)$ \\
Control cost coefficients $(R, \bar{R})$ & $(0.5, 0.5)$ \\
Discount factor $\gamma$ & $0.9$ \\
Noise variances $(\sigma_\epsilon^2, \sigma_{\epsilon_0}^2)$ & $0.01$ \\
Initial state variances $(\text{Var}(Y_{0}), \text{Var}(Z_{0}))$ & $1/3$ \\
Learning rate base $\eta$ & $0.01$ \\
Perturbation radius $\tau$ & $0.1$ \\
Truncation horizon $T$ & $40$ \\
Number of directions $M$ & $1000$ \\
Total iterations & $800$ \\
Number of seeds & $20$ \\
\hline
\end{tabular}
\caption{Model parameters and hyperparameters for the LQMFC RL experiment.}
\label{tab:lq_mfc_params}
\end{table}

All three methods use the fixed step $\eta=0.01$ and start from $(K_0,L_0)=(0,0.1)$. In the model-free methods, each update averages $M=1000$ independent one-point simulator cost samples, with horizon $T=40$ and independent perturbations $\pm\tau$ in each scalar gain coordinate, where $\tau=0.1$; see~\eqref{eq:definition_approx_gradient_MKV_simulator} and Algorithm~\ref{algo:POPestim}.

We study numerically the three primary configurations analyzed above:
\begin{enumerate}
    \item \textbf{Exact PG}: Uses model-based exact gradients as a baseline.
    \item \textbf{MF-MKV}: Model-free PG using the idealized MKV simulator.
    \item \textbf{MF-Pop ($N$)}: Model-free policy gradient using the finite-population simulator with $N$ agents. We compare $N=2,5,10,20,50,100$ to study the impact of population size on learning.
\end{enumerate}
For the variances chosen here, the mean-field and population costs satisfy
\[
C^N(K,L)=\left(1-\frac1N\right)C_y(K)+\left(1+\frac1N\right)C_z(L),
\]
where $C_y$ and $C_z$ are defined in~\eqref{eq:def_Cy_Cz}. Both weights are positive for $N\ge2$, so the optimal gains $(K^*,L^*)$ are the same for the mean-field and population problems. The population size changes the objective weights and the sampled gradients while preserving these optimal gains. We exclude $N=1$, for which the centered state vanishes and $K$ is not identifiable.

\subsection{Discussion of Results}
\begin{figure}[H]
    \centering
    \includegraphics[width=0.7\textwidth]{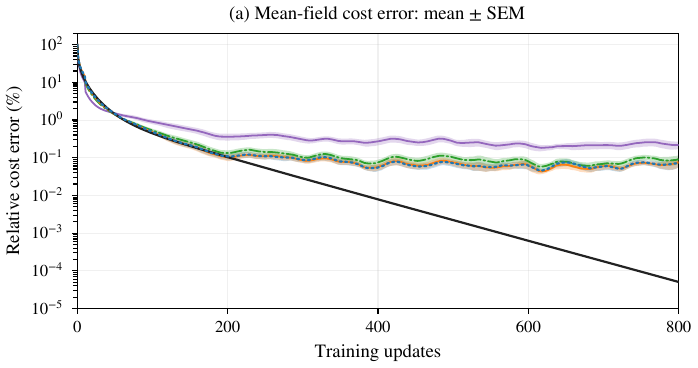}\par\smallskip
    \includegraphics[width=0.4\textwidth]{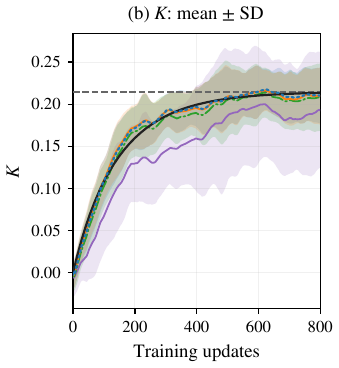}\hfill
    \includegraphics[width=0.4\textwidth]{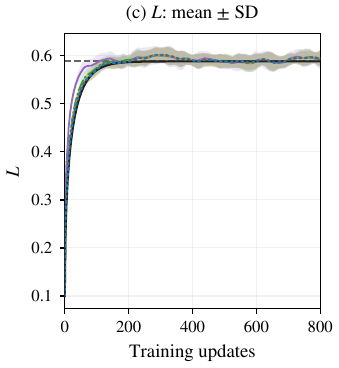}\par\smallskip
    \includegraphics[width=0.7\textwidth]{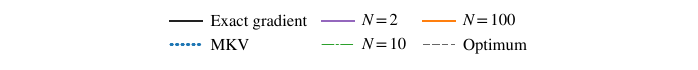}
    \caption{Training with exact mean-field gradients, an MKV simulator, and population simulators ($N=2,10,100$). Means and bands are computed across 20 runs after applying a centered 21-point moving average to each stochastic run. (a) Exact mean-field relative cost error, in percent, with mean $\pm$ SEM. (b,c) Gains $K,L$, with mean $\pm$ one SD; grey dashed lines are the common optimal gains. The exact-gradient reference is deterministic and unsmoothed.}
    \label{fig:lq_costs}
\end{figure}

Each stochastic method is run with 20 independent seeds. Within each replication, common random numbers are used across methods and population sizes to facilitate comparisons. For Figure~\ref{fig:lq_costs}, we apply a centered 21-point moving average to each stochastic run, shortening the window symmetrically near the endpoints. We then compute the mean and standard deviation (SD) across the 20 smoothed runs; the standard error of the mean (SEM) is the SD divided by $\sqrt{20}$. The exact-gradient reference is unsmoothed. The captions specify which quantity is shown.

Figure~\ref{fig:lq_costs} compares the learning trajectories for the exact-gradient method, the MKV simulator, and population simulators with $N=2,10,100$. Panel~(a) reports the mean-field relative cost error
\[
e_{\mathrm{MF}}(\theta)=\frac{C(\theta)-C(\theta^*)}{C(\theta^*)},
\]
expressed in percent. The exact infinite-horizon cost is used only for evaluation of the learned policies. Panels~(b,c) show the gains $K$ and $L$. For $N=2$, learning of $K$ is slower and the early progress in $L$ is faster than for the larger populations. This is consistent with the smaller weight on $C_y$ and larger weight on $C_z$. The stochastic cost errors continue to fluctuate above the exact-gradient reference; MKV and $N=100$ give similar errors.

\begin{figure}[H]
    \centering
    \includegraphics[width=0.9\textwidth]{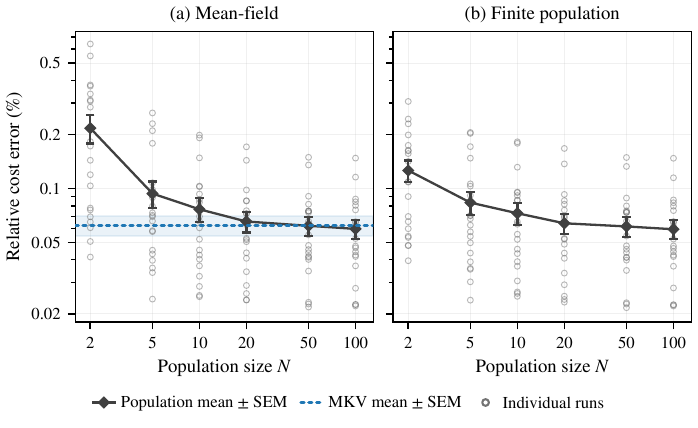}
    \caption{Relative cost errors after training, evaluated in (a) the mean-field problem and (b) each population's own $N$-agent problem. Each run contributes its average unsmoothed error over updates $601$--$800$. Diamonds and error bars show the mean $\pm$ SEM across 20 runs; circles show all individual runs. The blue dotted line and blue band in (a) show the MKV mean $\pm$ SEM.}
    \label{fig:lq_params}
\end{figure}

Figure~\ref{fig:lq_params} shows how the errors depend on $N$ under the same training budget. Panel~(a) uses $e_{\mathrm{MF}}$; panel~(b) uses the relative error for the population's own objective,
\[
e_N(\theta)=\frac{C^N(\theta)-C^N(\theta^*)}{C^N(\theta^*)}.
\]
For each run, the unsmoothed errors are averaged over updates $601$--$800$, then summarized across the 20 runs. The mean-field error decreases from about $0.217\%$ at $N=2$ to $0.0595\%$ at $N=100$, compared with $0.0623\%$ for MKV. The MKV method need not yield the smaller error after a fixed number of stochastic updates. The corresponding own-population errors are $0.126\%$ and $0.0593\%$. Individual runs vary considerably, and the improvement is smaller at large $N$. These errors reflect optimization, gradient sampling, horizon truncation and perturbation effects as well as population size.

\subsection{A Coupled Two-Dimensional Example}
\label{subsec:lq-numerics-2d}

We now take $d=\ell=2$ and learn unrestricted $2\times2$ matrices $K,L$, with eight free entries in total. The coefficients are
\begin{align*}
A&=\begin{pmatrix}0.45&0.20\\-0.10&0.40\end{pmatrix},&
\bar A&=\begin{pmatrix}0.40&-0.05\\0.08&0.45\end{pmatrix},\\
B&=\begin{pmatrix}0.50&0.10\\0.05&0.45\end{pmatrix},&
\bar B&=\begin{pmatrix}0.15&-0.05\\0.03&0.20\end{pmatrix},\\
Q&=\begin{pmatrix}0.65&0.12\\0.12&0.45\end{pmatrix},&
\bar Q&=\begin{pmatrix}0.25&-0.04\\-0.04&0.35\end{pmatrix},\\
R&=\begin{pmatrix}0.55&0.06\\0.06&0.45\end{pmatrix},&
\bar R&=\begin{pmatrix}0.25&-0.03\\-0.03&0.30\end{pmatrix}.
\end{align*}
Both optimal gains have nonzero off-diagonal entries. We retain $\gamma=0.9$. The common and idiosyncratic initial vectors are independent, each with independent coordinates uniformly distributed on $[-1,1]$. The common and idiosyncratic Gaussian innovations are independent and each have covariance $0.01I_2$.

All three methods start from $K_0=0$, $L_0=0.1I_2$ and perform 800 updates with the fixed step $\eta=0.01$. Each stochastic update uses $M=24000$ independent one-point simulator cost samples, horizon $T=40$, and independent uniform perturbations on the Frobenius sphere of radius $\tau=0.15$ for each gain. The estimator normalization is $\ell d/\tau^2=4/\tau^2$ per gain. We compare the exact-gradient reference, MKV, and population simulators with $N=2,10,100$, using 20 independent runs for each stochastic method and common random numbers across methods within each replication.

\begin{figure}[htbp]
    \centering
    \includegraphics[width=0.7\textwidth]{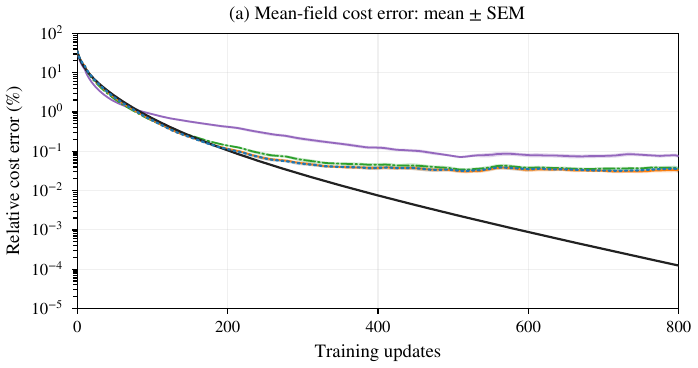}\par\smallskip
    \includegraphics[width=0.4\textwidth]{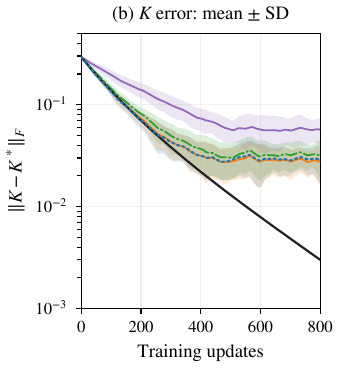}\hfill
    \includegraphics[width=0.4\textwidth]{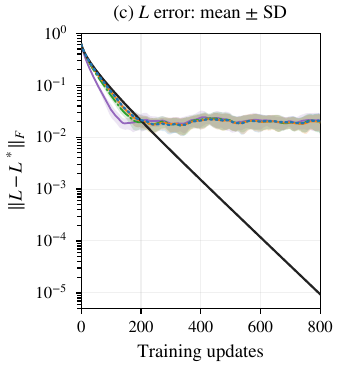}\par\smallskip
    \includegraphics[width=0.7\textwidth]{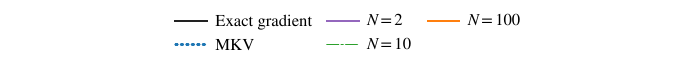}
    \caption{Learning in the coupled two-dimensional example. Means and bands are computed across 20 runs after applying a centered 21-point moving average to each stochastic run, with symmetrically shortened windows near the endpoints. (a) Exact mean-field relative cost error in percent, with mean $\pm$ SEM. (b,c) Frobenius errors $\|K-K^*\|_F$ and $\|L-L^*\|_F$, with mean $\pm$ one SD; norms are computed for each run before smoothing and averaging. The exact-gradient reference is deterministic and unsmoothed.}
    \label{fig:lq-2d}
\end{figure}

Figure~\ref{fig:lq-2d} reports the mean-field relative cost error and the full-matrix gain errors. Averaging unsmoothed cost errors over updates $601$--$800$ within each run and then across the 20 runs gives $0.0339\%$, $0.0781\%$, $0.0370\%$ and $0.0324\%$ for MKV, $N=2$, $N=10$ and $N=100$, respectively. As in the scalar example, $N=2$ learns $K$ more slowly, while $N=100$ closely follows MKV and the late errors in $L$ are similar across simulators. The mean-field and finite-population optimal gains again coincide, so the delay does not indicate a different finite-population optimum. The coefficients, sample count and perturbation radius differ from the scalar example, so these error magnitudes cannot isolate the effect of dimension.

\begin{figure}[htbp]
    \centering
    \includegraphics[width=0.9\textwidth]{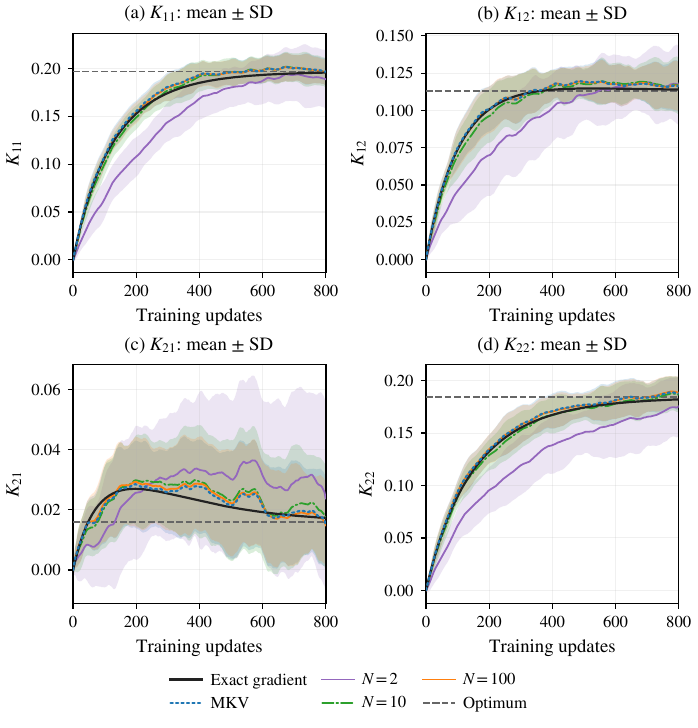}
    \caption{Evolution of the four entries of $K$ in the coupled two-dimensional example. Lines and shading show means $\pm$ one SD across 20 runs after applying the same centered 21-point moving average as in Figure~\ref{fig:lq-2d}. Grey dashed lines show the common mean-field and finite-population optimal entries. The exact-gradient reference is deterministic and unsmoothed.}
    \label{fig:lq-2d-K}
\end{figure}

\begin{figure}[htbp]
    \centering
    \includegraphics[width=0.9\textwidth]{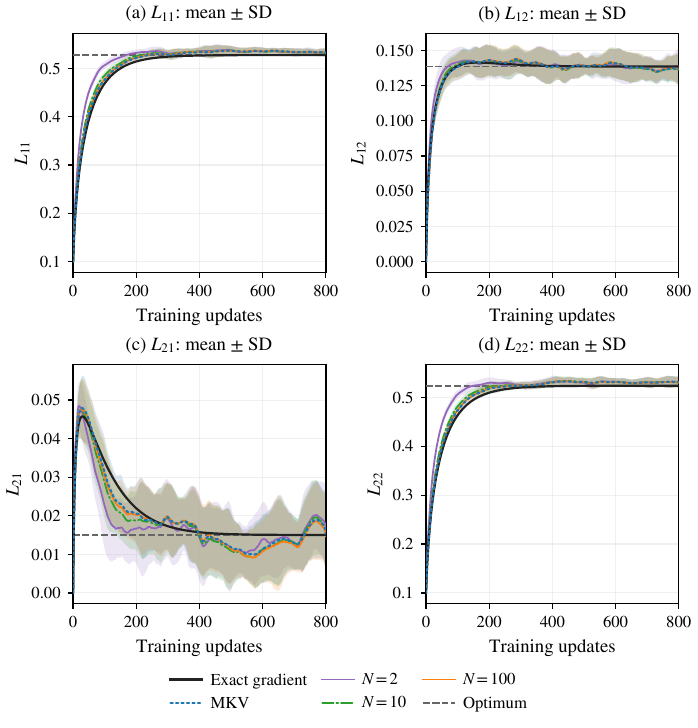}
    \caption{Evolution of the four entries of $L$ in the coupled two-dimensional example. Lines and shading show means $\pm$ one SD across 20 runs after applying the same centered 21-point moving average as in Figure~\ref{fig:lq-2d}. Grey dashed lines show the common mean-field and finite-population optimal entries. The exact-gradient reference is deterministic and unsmoothed.}
    \label{fig:lq-2d-L}
\end{figure}

Figures~\ref{fig:lq-2d-K} and~\ref{fig:lq-2d-L} show all entries of the learned gains. They reveal slower learning and greater variability in $K$ at $N=2$, whereas $N=2$ makes faster early progress in $L$ and the late $L$ trajectories show little separation. The overshoot of $L_{21}$ also occurs along the exact-gradient trajectory, so entrywise learning need not be monotone. Small entries remain noisy despite low cost errors.

\section{Notes and Complements}

The convergence proof in this chapter builds upon the approach proposed by~\cite{fazel2018global} for policy-gradient methods in classical LQR, and extends it by including stochastic dynamics and mean-field interactions.

PG methods have also been applied to linear-quadratic mean field type games, i.e., non-cooperative games between central planners solving MFC problems, see~\cite{carmona2020policyCDC,carmona2021linear,uz2024independent}.
Further related work on LQ mean field learning includes \cite{wang2021globallqmfcg,frikha2024full,delarue2025exploration}.

Other RL algorithms for LQ MFC problems include the two-timescale method of~\cite{angiuli2022unified} and the policy iteration method studied in~\cite{li2024policy}.  \cite{xu2025mean,xu2025robust} proposed an approach based on integral RL for the case of multiplicative noise and robust mean field social control.

\chapter{Conclusions and Perspectives}
\label{ch:conclusion}

\section{From Mean Field Control to Reinforcement Learning}
\label{sec:conclusion-synthesis}

Our point of departure in this monograph has been mean field control (MFC):
we formulate population control problems and develop reinforcement learning
(RL) methods adapted to their structure. The exposition is therefore
organized around control objectives, information structures and distribution
dynamics, rather than around existing multi-agent reinforcement learning
(MARL) algorithms. This offers one route for connecting the two subjects.
Further connections can be developed in both directions: questions about
local information, communication and coordination in MARL can motivate
richer mean field models, while mean field approximations can suggest
policy representations and benchmarks for large cooperative systems.
We hope that this viewpoint will encourage further exchanges between
stochastic control and multi-agent learning.

The finite-agent and mean field formulations describe how a common policy
affects individual states and the distribution of the population.
With common noise, the limiting distribution can remain random, and its
evolution is described through conditional laws. The information available
to the planner and to the agents, including the distinction between common
and private randomization, is therefore part of the control problem.
Chapters~\ref{ch:prelims} and~\ref{ch:AAP} develop these formulations and their
connections. The lifted mean field Markov decision process (MFMDP)
uses the population distribution as its state and a feasible joint
state--action distribution as its action. Under the stated assumptions,
this formulation provides a dynamic programming description of the
population control problem.

The two model classes studied here lead to complementary learning methods.
For the abstract model, Chapter~\ref{ch:numeric_I} builds upon the dynamic programming principle established in Chapter~\ref{ch:AAP} and develops value-based
methods using discretized distribution spaces. It also illustrates a neural actor--critic approach based on deep deterministic policy gradient. The linear quadratic (LQ) model
of Chapter~\ref{ch:LQMFC} has an explicit policy parametrization and hence provides
a setting for a full theoretical analysis of policy-gradient methods in
Chapter~\ref{ch:numeric_II}. The comparison between exact gradients,
mean field simulation and finite-population simulation highlights different
forms of access to the control problem. It also distinguishes the
approximation introduced by a finite population from errors due to
sampling, discretization, time truncation and optimization.
The assumptions and guarantees that we provided for these two classes are distinct, due to the fact that the LQ model has a very specific structure.

These formulations also help locate the book's contribution within the
broader landscape summarized in Table~\ref{tab:learning-formulations}.
In a mean field game, an individual agent optimizes against a given
population flow, and equilibrium additionally requires consistency.
In MFC, the controller accounts for the effect of the policy on that flow.
Networked MARL methods can retain agent identities and exploit interaction
or communication graphs, whereas the mean field models considered here
describe interactions through population distributions.
The suitability of these approaches depends on the objectives,
information and interaction structure of the problem.

\section{Directions for Further Research}
\label{sec:conclusion-directions}

We highlight here a few research directions that we consider important
for extending the models and learning methods studied in this monograph.

\paragraph{Heterogeneous populations.}
Multiple populations or agent types can represent differences
in dynamics, objectives or available information.
Heterogeneous LQ social control was studied in~\cite{huang2012social}, and multi-population control and game formulations were developed in~\cite{bensoussan2018populations}. For an overview of extensions beyond homogeneity and anonymity, including cooperative variants, see~\cite{lauriere2025heterogeneity}. On the learning side, type-dependent mean-action representations in MARL~\cite{subramanian2020multitype} retain a separate aggregate for each type, an idea that can also inform policy representations for heterogeneous MFC. Cooperative MFC approximation and natural policy gradient guarantees for several classes of agents were established in~\cite{mondal2022heterogeneous}.

A central modeling question is which differences must be retained
explicitly and which interactions can still be summarized by distributions.
For learning, this raises questions about the policy representation,
the amount of data needed for each type, and how approximation and
computational costs depend on the number of populations or types.

\paragraph{Partial observation.}
An agent or a planner may observe neither the complete individual state
nor the population distribution.
Maximum-principle methods for partially observed mean field control were developed in~\cite{wang2014partial,buckdahn2017partial}. For general partially observed control, randomized filtering and a Bellman equation on a space of probability measures provide a complementary approach~\cite{bandini2019filtering}. On the learning side, approximate information states have been studied for general partially observed systems~\cite{SubramanianSinhaSeraj2022}, providing methods for learning representations of observation histories with control-performance guarantees. Decentralized partially observable MFC, including a dynamic programming principle, finite-population guarantees and policy gradient methods, is studied in~\cite{cui2023learningpomfc}.

The control formulation must then account for what can be inferred from
observation histories. A central challenge is to identify an information
state sufficient for control and to learn it together with the policy,
while preserving the distinction between common and private information.
Errors in estimating the population distribution must also be related
to the resulting control performance.

\paragraph{Graph-structured interactions.}
Local interactions and communication constraints may require more
information than a single population distribution provides.
Graphon approximations for the control of large linear networks were developed in~\cite{gao2020graphoncontrol}. Related cooperative MARL methods exploit local interactions through localized policies, with computational complexity depending on neighborhood sizes~\cite{QuWiermanLi}. For interactions described by a doubly stochastic matrix,~\cite{mondal2022nonuniform} establishes MFC approximation and learning guarantees under assumptions including affine dependence of rewards on the local mean field. Graphon-based MFC and learning for cooperative MARL are studied in~\cite{hu2023graphoncontrol}, while MFC on sparse graphs using local weak convergence is developed in~\cite{fabian2025sparse}.

One direction is to combine local graph information with aggregate
descriptions of larger groups of agents.
The resulting models should specify which features of the graph are
retained, what information each policy uses, and how the approximation
affects performance. 

\paragraph{Offline learning.}
The model-free methods studied here use access to a simulator to evaluate policies
or state--action choices.
Learning from a fixed collection of trajectories introduces a different
question: which changes of policy and population distribution can be
evaluated reliably from the available data?
In single-agent offline RL, pessimistic value iteration provides data-dependent performance guarantees~\cite{jin2021pessimism}. In mean-field MARL, the method proposed in~\cite{chen2021offlinemeanfield} combines pessimism with mean embeddings to learn from fixed trajectories under assumptions on the function class and data coverage. Related offline MARL methods constrain value estimation to state--action pairs supported by the data to limit extrapolation errors~\cite{yang2021icq}. Extending these approaches to other MFC formulations requires accounting for how changes in the policy affect the population flow.

The coverage of individual states and actions must be considered together
with the population distributions represented in the data.
Useful guarantees will require assumptions connecting this coverage,
the chosen policy class and the intended control objective.

\paragraph{Finite-population performance.}
This monograph mentions a few finite-population approximation results and,
in the LQ setting, learning guarantees under specified assumptions.
Controlled finite-population approximations were developed through fluid limits~\cite{gast2012mean} and controlled McKean--Vlasov limit theory~\cite{lacker2017control}. Quantitative convergence for finite-state control is studied in~\cite{cecchin2021finite}. Common noise and admissible control formulations are addressed in discrete time by~\cite{motte2022mean} and in continuous time by~\cite{djete2022limit}. Further quantitative results include convergence rates with common noise~\cite{cardaliaguet2023rate} and rates under different regularity assumptions~\cite{daudin2024rate}. Under convexity assumptions,~\cite{jackson2025distributed} gives nonasymptotic bounds comparing full-information and distributed cooperative controls.

Learning guarantees connecting MFC algorithms to finite-agent performance include kernel-based Q-learning~\cite{gu2021meanQ} and natural policy gradient for heterogeneous populations~\cite{mondal2022heterogeneous}. Statistical efficiency with model-based function approximation is studied for the cooperative MFC objective in~\cite{huang2024statistics}. Performance bounds separating linear-model approximation from finite-population approximation are obtained in~\cite{bayraktar2025learning}.

A broader objective is to relate the performance of a learned mean field
policy to the optimum of the original finite-agent problem.
Such comparisons must account jointly for population approximation,
statistical error, policy or value-function approximation, and optimization
error. They must also respect the information available to implement the
policy in the finite system, as well as constraints on computation time, memory and energy use.
Extending these guarantees beyond the present model classes would
strengthen the connection between mean field learning and control
of large but finite populations.

\appendix
\chapter{Technical estimates for LQ policy gradient}
\label{app:lq-pg-proofs}

This appendix collects the technical estimates used in Chapter~\ref{ch:numeric_II}. The notation and standing assumptions are those of that chapter.

\section{Technical results in Section~\ref{subsection:exact_PG_for_MFC}}
\label{section:proof_of_exact_pg}

This section collects the technical lemmas used for proving the statements in Section~\ref{subsection:exact_PG_for_MFC}.
These lemmas are adapted from~\cite{fazel2018global} with necessary modifications to cope with the idiosyncratic noise and the common noise processes present in the state dynamics.
To make the paper self-contained, we provide their proofs here.
Because there is an analogy between the processes $\bY$ and $\bZ$ in their definitions and usages, we will state and prove results related to the process $\bY$ and the parameter $K$ only.

\subsection{Cost expressions}
\label{subsection:cost_expression}

For an admissible control $\theta = (K, L) \in \Theta$, consider two $L^2$-discounted integrable processes in $\bY^{K, \xi_y}, \bZ^{L, \xi_z} \in S$ following dynamics~\eqref{eq:dyn_y_theta},~\eqref{eq:dyn_z_theta} with initial values $Y_{0} = Y_{0}^{K, \xi_y} = \xi_y$ and $Z_{0} = Z_{0}^{L, \xi_z}= \xi_z$.
We defined the corresponding value functions $V_y$ and $V_z$ on $\RR^{\ell \times d} \times \RR^d$ by
\begin{equation*}
	V_y(K, \xi_y) = \EE_{\boldsymbol{\varepsilon}} \Big[ \sum_{n \geq 0} \gamma^{n} f(Y_{n}^{K, \xi_y}, K, Q, R) \Big],
	\quad
	V_z(L, \xi_z) = \EE_{\boldsymbol{\varepsilon^0}} \Big[ \sum_{n \geq 0} \gamma^{n} f(Z_{n}^{L, \xi_z}, L, \tilde Q, \tilde R) \Big]
\end{equation*}
where  $f(\xi, K, Q, R) := \xi^\top( Q + K^\top R K ) \xi$ and the expectations $\EE_{\boldsymbol{\varepsilon}}, \EE_{\boldsymbol{\varepsilon^0}}$ are integrations with respect to the idiosyncratic and the common noise processes.

These value functions can be computed simply with the help of solutions $P_y^K$ and $P_z^L$ to the two discrete-time  Lyapunov equations (DLEs):\index[not]{DLE}\index[sub]{discrete-time  Lyapunov equations}
\begin{equation}
	\label{eq:lyapunov_eq_theta}
	\left\{
	\begin{array}{rcl}
		P^y_K & = & Q + K^\top R K + \gamma ({\mathrm A} - {\mathrm B} K)^\top P^y_K ({\mathrm A} - {\mathrm B} K)
		\\
		P^z_L & = & \tilde{Q} + L^\top \tilde{R} L +\gamma ( \tilde{{\mathrm A}}- \tilde {\mathrm B} L)^\top P^z_L ( \tilde{{\mathrm A}} - \tilde {\mathrm B} L).
	\end{array}
	\right.
\end{equation}
The existence of solutions to~\eqref{eq:lyapunov_eq_theta} is deferred to the next Section~\ref{sec:existence_of_DLE}.
Let us also define
\begin{equation}
	\label{eq:value_function_noise_constant}
	\alpha^y_K = \EE \big[ \sum \nolimits_{n \geq 1} \gamma^{n} (\varepsilon_{n})^\top P^y_K \varepsilon_{n} \big],
	\qquad
	\alpha^z_L = \EE \big[ \sum \nolimits_{n \geq 1} \gamma^{n} (\varepsilon^0_{n})^\top P^z_L \varepsilon^0_{n} \big].
\end{equation}
Then, using the dynamics of $\bY^{K, \xi_y}$~\eqref{eq:dyn_y_theta} and the discrete-time Lyapunov equation~\eqref{eq:lyapunov_eq_theta} for $P_K^y$, we have
$
	V_y(K, \xi_y) = \xi_y^\top P_K^y \xi_y + \alpha_K^y.
$
When the initial states are random, e.g., $Y_{0} = X_{0} - \bar{X}_0$, we have
\begin{equation}
	\label{eq:cost_expression_with_PK_y}
	C_y(K) = \EE[ V_y(K, Y_{0}) ] = \langle P_K^y, \,  \Sigma_{Y_{0}} + \frac{\gamma}{1-\gamma} \Sigma^1 \rangle_{tr} = \langle Q + K^\top R K, \, \Sigma_K \rangle_{tr}.
\end{equation}
Because $\lambda_y^1 = \lambda_{min}( \Sigma_{Y_{0}} + \frac{\gamma}{1-\gamma} \Sigma^1 )$, we derive the following bounds:
\begin{equation}
	\label{eq:upper_bound_Riccati_and_Variance_matrix}
	\|P^y_K\| \leq C_y(K) / \lambda_y^1,
	\quad
	\|\Sigma^y_K\| \leq C_y(K) / \lambda_{min}(Q).
\end{equation}
Similar expressions for the cost $C_z(L)$, value $V_z(L, \xi_z)$, and the bounds for $P_L^z$ and $\Sigma_L^z$ follow with $(P_L^z, \Sigma_{Z_{0}}, \Sigma^0, \lambda_{z}^0, \tilde Q, \tilde R)$.

For agent $j$, her control at time $n$ takes the form $\ctrl_{n}^{j} = - K (X_{n}^{j} - \bar{X}_{n}^{N} ) - L \bar{X}_{n}^N$.
And the corresponding auxiliary processes $Y_{n}^{K, j, N} := X_{n}^{j} - \bar{X}^N_{n}$ and $Z_{n}^{L, N} := \bar{X}^N_{n}$ follow the dynamics similar to the ones defined in equations~\eqref{eq:dyn_y_theta}--\eqref{eq:dyn_z_theta}, but by replacing the noise terms with $\varepsilon_{n+1}^{1, j, N} := \varepsilon_{n+1}^{1, (j)} - \frac{1}{N} \sum_{i=1}^N \varepsilon_{n+1}^{1, (i)}$ and $\varepsilon_{n+1}^{0, N} := \varepsilon_{n+1}^0 + \frac{1}{N} \sum_{i=1}^N \varepsilon_{n+1}^{1, (i)}$.
Then, the cost for agent $n$ can be expressed into
By consequence, the social cost of the population given by~\eqref{fo:social_cost_of_population} is bounded by

When these agents adopt the same control parameter $\theta$, with the same arguments used in the homogeneous case, the difference between the auxiliary cost $C^{(j)}_y(K)$ for agent $j$ and the auxiliary cost in the mean-field setting $C_y(K)$ satisfies
The gap between the social cost and the MF cost tends towards 0 when we impose a bounded domain for the control parameters and when $h_{\delta Q} \to 0$ and $N \to \infty$.

\subsection{Gradient expression}
\label{subsection:gradient_expression}
We express the gradient of cost $C_y(K)$ in terms of the variance matrix $\Sigma_K$.
This gradient expression is crucial in the whole paper.
Consequently, the approximation of gradient $\nabla_K C_y(K)$ could hinge on an estimation of variance matrix $\Sigma_K$.
For completeness, we provide the proof here, which is analogous to~\cite[Lemma 1]{fazel2018global}.

\begin{lemma}
	\label{lemma:policy_gradient_expression}
	For $\theta = (K, L) \in \Theta$, let $ E_K = (R + \gamma {\mathrm B}^\top P^y_K {\mathrm B})K - \gamma {\mathrm B}^\top P^y_K {\mathrm A}$, then
	$$
		\nabla_K C(\theta) = \nabla_K C_y(K) = 2 E_K \Sigma_K,
	$$
\end{lemma}

\begin{proof}
	Consider a process $(Y_{n}^{K, \xi_y})_{n \geq 0}$ starting from $Y_{0} = \xi_y$, the value function $V_y(K, \xi_y)$ is rewritten into
	\begin{align*}
		V_y(K, \xi_y) =  \xi_y^\top (Q + K^\top R K) \xi_y  + \gamma \mathbb{E}_{\varepsilon_{n}} \big[ V_y(K, ({\mathrm A} -{\mathrm B} K)\xi_y + \varepsilon_{n} ) \big]
		\label{eq:expression_C_y_K_ytilde}
	\end{align*}
	Noting that $\nabla_{\xi_y} V_y(K, \xi_y) = 2 P^y_K \xi_y$, and $\alpha^y_K$ does not depend on $\xi_y$, we have
	\begin{align*}
		\nabla_K V_y(K, \tilde{Y}) =  \big(2 RK - 2 \gamma {\mathrm B}^\top P^y_K ({\mathrm A}-{\mathrm B} K) \big)\xi_y \xi_y^\top + \gamma \EE \big[ \nabla_K  V_y (K, Y_{n=1}^{k, \xi_y} ) \big].
	\end{align*}
	Continue unrolling the gradient inside the expectation, and note that $C(\theta) = C_y(K) + C_z(L)$ and $C_y(K) = \EE_{\xi_y}[V_y(K, \xi_y)]$, we conclude the result of the lemma.
\qed\end{proof}

\subsection{Cost variation expression}

The following lemma shows the variation in cost $C_y$.
\begin{lemma}
	\label{lemma:difference_in_C_y(K)}
	Consider $\theta'=(K', L') \in \Theta$.
	Let $\Delta K = K' - K$, then
	\begin{equation}
		C_y(K') - C_y(K) =  \langle \Delta K, \, 2 E_K \Sigma_{K'} + (R + \gamma {\mathrm B}^\top P_K^y {\mathrm B} )\Delta K  \Sigma_{K'} \rangle_{tr}.
	\end{equation}
\end{lemma}

\begin{proof}
	Consider a process $(Y_{n}^{K', Y_{0}})_{n \geq 0}$ starting from $Y_{0}$ following dynamics~\eqref{eq:dyn_y_theta} with parameter $K'$, and a sequence of random costs $\big( V_y( K, Y_{n}^{K', Y_{0}} ) \big)_{n \geq 0}$ for parameter $K$ with initial random states $\xi_y = Y_{n}^{K', Y_{0}}$ for $n \geq 0$.
	We consider the advantage quantity at time $n$ given by
	\begin{align*}
		A_y(n, K, Y_{n}^{K', Y_{0}}, K') := & (Y_{n}^{K', Y_{0}})^\top (Q + K'^\top R K') Y_{n}^{K', Y_{0}} + \gamma V_y(K, Y_{n+1}^{K', Y_{0}}) - V_y(K, Y_{n}^{K', Y_{0}})
		\\
		=                               & ( Y_{n}^{K', Y_{0}} )^\top (K' - K)^\top \big( 2 E_K + (R + \gamma {\mathrm B}^\top P_K^y {\mathrm B} ) (K' - K) \big)  Y_{n}^{K', Y_{0}}
	\end{align*}
	where the second equality follows from $V_y(K, \xi_y) = (\xi_y)^\top P_K^y \xi_y + \alpha_K^y$ for any $\xi_y \in \RR^d$.
	Then, the difference between the auxiliary costs $C_y$ at $K$ and $K'$ becomes
	\begin{align*}
		C_y(K') - C_y(K)
		 & = \EE_{Y_{0} \sim X_{0} - \bar{X}_0}\big[ \sum \nolimits_{n\geq 0} \gamma^{n} f( Y_{n}^{K', Y_{0}}, K', Q, R) - V_y(K, Y_{0}) \big]
		\\
		 & = \EE_{Y_{0} \sim X_{0} - \bar{X}_0} \big[ \sum \nolimits_{n \geq 0} \gamma^{n} A_y( n, K, Y_{n}^{K', Y_{0}}, K') \big]
		\\
		 & =  \langle \Delta K, \, 2 E_K \Sigma_{K'} + (R + \gamma {\mathrm B}^\top P_K^y {\mathrm B} )\Delta K  \Sigma_{K'} \rangle_{tr} .
	\end{align*}
\qed\end{proof}

\subsection{Solution to discrete-time Lyapunov equations}
\label{sec:existence_of_DLE}

For any admissible parameters $\theta \in \Theta$, we can construct solutions to the DLEs~\eqref{eq:lyapunov_eq_theta}.

\begin{proposition}
\label{prop:existence_of_lyapunov_equation}
Under Assumptions~\ref{ass:positivity-qr} and~\ref{ass:pg-state-cost-coercivity}, for any admissible control $\theta \in \Theta$, there exists a unique pair of matrices $(P^y_K, P^z_L)$ satisfying (DLEs)~\eqref{eq:lyapunov_eq_theta} and $P^y_K \succ 0$, $P^z_L \succ 0$. 
\end{proposition}

\begin{proof}
    For an admissible parameter $\theta = (K, L) \in \Theta$, the existence of a positive definite symmetric matrix $P_K^y$ (resp. $P^z_L$) as the unique solution to the corresponding Lyapunov equation in~\eqref{eq:lyapunov_eq_theta} is a direct application of~\cite[Theorem 3.2]{bof2018lyapunov} with a Schur matrix $\sqrt{\gamma} ({\mathrm A} - {\mathrm B} K)$ (resp. $\sqrt{\gamma} (\tilde {\mathrm A} - \tilde {\mathrm B} L)$) and a positive definite matrix $Q + K^\top R K$ (resp. $\tilde Q + L^\top \tilde{R} L$). 
    Indeed, we can construct directly the solutions as: 
    \begin{equation*}
    \left\{
    \begin{array}{rcl}
       P_K^y &:=& \sum_{t=0}^\infty \gamma^t \big( ({\mathrm A} - {\mathrm B} K)^\top \big)^t ( Q + K^\top R K ) ({\mathrm A} - {\mathrm B} K)^t 
       \\
       P_L^z &:=& \sum_{t=0}^\infty \gamma^t \big( ( \tilde{{\mathrm A}} - \tilde{{\mathrm B}} L)^\top \big)^t ( \tilde{Q} + L^\top \tilde{R} L ) (\tilde{{\mathrm A}} - \tilde{{\mathrm B}} L)^t.
    \end{array}
    \right.
    \end{equation*}
    We obtain immediately that $P_K^y, P_L^z \succ 0$ since $Q + K^\top R K \succ 0$ and $\tilde Q + L^\top \tilde R L \succ 0$ under Assumptions~\ref{ass:positivity-qr} and~\ref{ass:pg-state-cost-coercivity}.
	
\end{proof}

The following lemma completes the connections between the solution matrices $(P, \bar P)$ of the Riccati equations~\eqref{fo:Riccati} and the solution matrices of (DLEs)~\eqref{eq:lyapunov_eq_theta} with the specific parameters $K^*$ and $L^*$ defined in Theorem~\ref{th:existence_linear_control}. 

\begin{lemma}
We suppose that Assumptions~\ref{ass:positivity-qr} and~\ref{ass:pg-state-cost-coercivity} hold, and we assume that the matrix Riccati equations~\eqref{fo:Riccati} admit solutions $(P, \bar P)$. We define two symmetric matrices $P^{*, y}$ and $P^{*,z}$ in $\RR^{d \times d}$ by
\begin{equation}
\label{eq:riccati_matrix_close_loop}
    P^{*,y} := ( {\mathrm A}^\top P + P^\top {\mathrm A} + 4Q ) / 4, \qquad P^{*,z} := ( \tilde{{\mathrm A}}^\top \bar{P} + \bar{P}^\top \tilde{{\mathrm A}}  + 4 \tilde Q ) / 4.
\end{equation}
We consider the parameters $K^*$ and $L^*$ defined by
$$
    K^* = \frac{1}{2} R^{-1} {\mathrm B}^\top P , \qquad L^* = \frac{1}{2} \tilde{R}^{-1} \tilde{{\mathrm B}}^\top \bar{P}.
$$
Then, we have that $P^{*,y}$ and $P^{*,z}$ satisfy (DLEs)~\eqref{eq:lyapunov_eq_theta} for $K = K^*$ and $L = L^*$.
\end{lemma}

\begin{proof}
    From the Riccati equation~\eqref{fo:Riccati} of $P$, we have
    $$
       2 R K^* =  {\mathrm B}^\top P = \gamma {\mathrm B}^\top ( {\mathrm A}^\top P + 2 Q) ({\mathrm A} - {\mathrm B} K^{*})
    $$
    and then
    \begin{align*}
        {\mathrm A}^\top P + P^\top {\mathrm A} & = \gamma {\mathrm A}^\top \big( {\mathrm A}^\top P + 2 Q \big) ({\mathrm A} - {\mathrm B} K^*) + \gamma ({\mathrm A} - {\mathrm B} K^*)^\top \big(  P^\top {\mathrm A} + 2 Q \big) {\mathrm A} 
        \\ \
        & = \gamma ({\mathrm A} - {\mathrm B} K^*)^\top \big( {\mathrm A}^\top P +  P^\top {\mathrm A} + 4 Q \big) ({\mathrm A} - {\mathrm B} K^*) + 4 (K^*)^\top R K^*.
        \nonumber 
    \end{align*}
    Thus, $P^{*, y} = ( {\mathrm A}^\top P + P^\top {\mathrm A} + 4Q) / 4$ is a solution to the first equation in~\eqref{eq:lyapunov_eq_theta}. Since $Q\succ0$, the stage matrix $Q+(K^*)^\top R K^*$ is positive definite. By applying similar arguments to the Riccati equation of $\bar P$ in~\eqref{fo:Riccati}, we see that $P^{*,z}$ satisfies the second equation in the (DLEs)~\eqref{eq:lyapunov_eq_theta}.

\end{proof}

\subsection{Perturbation analysis}
\label{subsection:perturbation_analysis_of_exact_PG}

\subsubsection{Perturbation on variance}
\label{subsection:perturbation_on_variance}

In this subsection, we show that when $\Delta K = K' - K$ is small enough, the norm between variance matrices is controlled by $\|\Delta K \|$.
To do so, we introduce a function $f_{var}$ on matrix $K$ defined by
\begin{equation}
	\label{eq:f_var_K}
	f_{var}(K) = 4 ( 1 + \| {\mathrm A} - {\mathrm B} K \| ) \| {\mathrm B} \|  \frac{ C_y(K)} { \lambda_y^1 \lambda_{min}(Q) }.
\end{equation}
We also introduce the following operators:
for any $\Sigma \in \cM\cS$,
\begin{equation}
	\label{eq:def-ope-F_K_y}
	\cF_K(\Sigma) = \gamma ({\mathrm A}- {\mathrm B} K) \Sigma ({\mathrm A}-{\mathrm B} K)^\top,
	\,
	\cT_K(\Sigma) = \sum_{n \geq 0} \gamma^{n} ({\mathrm A} - {\mathrm B} K)^{n} \Sigma \big( ({\mathrm A} - {\mathrm B} K)^{n} \big)^\top.
\end{equation}
Consider the operator norms on $\cM\cS$ for $\cF_K$ and $\cT_K$, and denote them by $\vertiii{\cF_K}$ and $\vertiii{\cT_K}$.
Because $\vertiii{\cF_K} = \sup_{\| \Sigma \| = 1} \cF_K(\Sigma) \leq \gamma \| {\mathrm A} - {\mathrm B} K \|^2 < 1$ for a control parameter $K$ in $\theta = (K, L ) \in \Theta$, we then have
\begin{equation}
	\label{eq:cT}
	\cT_K(\Sigma) = \sum \nolimits_{n \geq 0} (\cF_K)^{n}(\Sigma) = ( \bI - \cF_K)^{-1} (\Sigma).
\end{equation}

\begin{lemma} 
	\label{lemma:expression_Sigma_K_with_T_K}
	The variance matrix $\Sigma_K$ satisfies
	\begin{equation}
		\label{eq:link_vairance_y_with_T_K}
		\Sigma_K = \mathcal{T}_K \big(\Sigma_{Y_{0}} +  \tfrac{\gamma}{1- \gamma} \Sigma^1 \big).
	\end{equation}
	Moreover, $\lambda_{min} (\Sigma_K) \geq \lambda_y^1$.
\end{lemma}

\begin{proof}
	Through expanding the variance matrix $\Sigma_K$~\eqref{eq:variance_matrices_y_and_z} with $Y_{n+1}^{K, Y_{0}} = ({\mathrm A} - {\mathrm B} K) Y_{n}^{K, Y_{0}} + \varepsilon_{n+1}$ for $n \geq 0$, we obtain
	\begin{align*}
		\Sigma_K
		 & = \EE [ Y_{0} Y_{0}^\top] + \EE \big[ \sum_{n \ge 1} \gamma^{n} \varepsilon_{n} (\varepsilon_{n})^\top \big] + \EE \big[ \sum_{n \ge 0} \gamma^{n+1} ({\mathrm A} - {\mathrm B} K) Y_{n}^{K, Y_{0}} (Y_{n}^{K, Y_{0}})^\top ({\mathrm A} - {\mathrm B} K)^\top \big]
		\\
		 &
		= \Sigma_{Y_{0}} + \gamma / (1 - \gamma) \Sigma^1 + \cF_K ( \Sigma_K).
	\end{align*}
	Because $\cF_K(\Sigma_K) \succeq 0$, we then conclude the lemma.
\qed\end{proof}

The following lemma is about the bounds with operators $\cF_K$ and $\cT_K$.
\begin{lemma}
	\label{lemma:useful-bounds-ope-T-F}
	Let $\Delta K = K' - K $, we have
	\begin{enumerate}[label=(\roman*)]
		\item   \label{eq:TyK-Ctheta-lambda}
		      $
			      \| \Sigma_K \| / \| \Sigma_{Y_{0}} + \gamma / (1-\gamma) \Sigma^1 \| \leq  \vertiii{\cT_K} \leq \| \Sigma_K \| /  \lambda^1_y.
		      $

		\item \label{eq:perturbation_F_K^y}
		      $
			      \vertiii{\cF_{K'} - \cF_{K} } \leq \gamma \big( 2 \| {\mathrm A} - {\mathrm B} K \| \| {\mathrm B} \| \| \Delta K \| + \| {\mathrm B} \|^2 \| \Delta K \|^2 \big) .
		      $

		\item \label{eq:perturbation_T_K^y_in_T_K^y}
		      For any $\eta < 1$ such that $\vertiii{ \cT_K} \vertiii{ \cF_{K'} - \cF_{K} } \leq \eta$, for any $\Sigma \in \cM\cS$,
		      \begin{equation*}
			      \| (\cT_{K'} - \cT_{K}) (\Sigma) \|
			      \leq \frac{1}{1-\eta} \vertiii{ \cT_K }. \vertiii{\cF_{K'} - \cF_{K} }. \Vert \cT_{K}(\Sigma) \Vert
			      \leq \frac{\eta}{1- \eta} \Vert \cT_K (\Sigma) \Vert.
		      \end{equation*}

	\end{enumerate}
\end{lemma}

\begin{proof} The first inequality in~\ref{eq:TyK-Ctheta-lambda} is by definition of $\vertiii{\cT_K}$ and Lemma~\ref{lemma:expression_Sigma_K_with_T_K}.
	To show the second inequality, we notice that for any $\Sigma, U \in \cM\cS$ and $\Sigma$ invertible, we have
	\begin{align*}
		\| \cT_K(U) \|
		 & = \sup_{\| \xi \| = 1}\sum \nolimits_{n \geq 0} \langle (({\mathrm A}-{\mathrm B} K)^{n})^\top \xi \xi^\top({\mathrm A} - {\mathrm B} K)^{n}, \,  X  \rangle_{tr}
		\\
		 & \leq \sup_{\| \xi \| = 1}  \sum \nolimits_{n \geq 0} \gamma^{n} \langle (({\mathrm A}-{\mathrm B} K)^{n})^\top \xi \xi^\top({\mathrm A} - {\mathrm B} K)^{n} \Sigma^{1/2}, \,  \Sigma^{1/2}  \rangle_{tr} \| \Sigma^{-1/2} U \Sigma^{-1/2} \|
		\\
		 & = \| \Sigma^{-1/2} U \Sigma^{-1/2} \| \| \cT_K(\Sigma) \|.
	\end{align*}
	Then $\vertiii{\cT_K} = \sup_{\| U \| = 1} \| \cT_K(U) \| \leq \| \cT_K \big( \Sigma) \|  \| \Sigma^{-1/2} \|^2   \leq \| \Sigma_K \| / \lambda_y^1$ where the last inequality is due to Lemma~\ref{lemma:expression_Sigma_K_with_T_K} with $\Sigma = \Sigma_{Y_{0}} + \gamma / (1-\gamma) \Sigma^1$.
	The inequality in~\ref{eq:perturbation_F_K^y} is derived from the definitions of $\cF_K$ and $\cF_{K'}$.
	The inequality in~\ref{eq:perturbation_T_K^y_in_T_K^y} is adapted from~\cite[Lemma 20]{fazel2018global}, and we generalize it with a coefficient $\eta < 1$.
\qed\end{proof}

Now, we are ready to show a perturbation result for the variance matrix.

\begin{lemma}
	\label{lemma:perturbation_variance_Cy}
	Consider $\theta, \theta' \in \Theta$.
	If $ f_{var}(K) \| K' - K \| \leq 1$, then
	\begin{equation}
		\label{eq:perturbation_variance}
		\| \Sigma_{K'} - \Sigma_{K} \| \leq  f_{var}(K) \| \Sigma_K \| \| K' - K \|  \leq \| \Sigma_{K} \|
	\end{equation}
\end{lemma}

\begin{proof}
	Let $\Delta K = K' - K$.
	Because $\|{\mathrm B}\| \|\Delta K \| \leq \|{\mathrm B} \| f_{var}(K)^{-1} \leq 1/4$, so by Lemma~\ref{lemma:useful-bounds-ope-T-F}\ref{eq:TyK-Ctheta-lambda}\ref{eq:perturbation_F_K^y}, we have
	$$
		\vertiii{\cT_K}. \vertiii{\cF_{K'} - \cF_K} \leq f_{var}(K) \|\Delta K \| \gamma /2 \leq 1/2.
	$$
	Then, by Lemma~\ref{lemma:useful-bounds-ope-T-F}~\ref{eq:perturbation_T_K^y_in_T_K^y} with $\eta = 1/2$, we obtain that
	\begin{align*}
		\| \Sigma_{K'} - \Sigma_K \| = \| ( \cT_{K'} - \cT_K) (\Sigma_{Y_{0}} + \tfrac{\gamma}{1- \gamma} \Sigma^1 ) \|
		 & \leq 2 \vertiii{\cT_K} \vertiii{\cF_{K'} - \cF_K} \| \Sigma_K \|
		\\
		 & \leq  f_{var}(K) \|\Delta K \| \| \Sigma_K \|
	\end{align*}
\qed\end{proof}

\subsubsection{Perturbation on the DLE solution matrix}
\label{subsection:pertubation_on_DLE}

Consider the solutions $P_K^y$ and $P_{K'}^y$ of the first DLE equation in~\eqref{eq:lyapunov_eq_theta} for parameters $K$ and $K'$.
We introduce a function $f_{riccati}$ on the matrix $K$ defined by
\begin{equation}
	\label{eq:f_riccati_K}
	f_{riccati}(K, \Delta K) = \frac{ C_y(K) }{\lambda_y^1} \bigg( \gamma f_{var}(K) + \frac{4\|R\|}{\lambda_{min}(Q) } \| K \| + \frac{ 2 \| R \|}{ \lambda_{min}(Q)} \| \Delta K \| \bigg)
\end{equation}

\begin{lemma}
	\label{lemma:perturbation_riccati_Cy}
	Consider $\theta, \theta' \in \Theta$.
	If $f_{var}(K) \| K' - K \| \leq 1$, then
	\begin{equation}
		\label{eq:perturbation_DLE_solution}
		\| P_{K'}^y - P_K^y \|  \leq  f_{riccati}(K, \Delta K ) \| K' - K \|
	\end{equation}
\end{lemma}
\begin{proof} From the DLE~\eqref{eq:lyapunov_eq_theta} for $P_K^y$, we have $P_K^y = Q + K^\top R K + \cF_K(P_K^y)$, which implies that $P_K^y = \cT_K( Q + K^\top R K )$.
	Then by Lemma~\ref{lemma:useful-bounds-ope-T-F}\ref{eq:perturbation_T_K^y_in_T_K^y} with $\eta = 1/2$, we have
	\begin{align*}
		\| P_{K'}^y - P_K^y \| & = \| (\cT_{K'} - \cT_K)(Q + K^\top R K )  + \cT_{K'} \big( (K')^\top R K' - K^\top R K\big) \|
		\\
		                       & \leq  f_{var}(K, \Delta K)  \| \Delta K \|. \| P_K^y \| +  2 \vertiii{\cT_{K'}}. \| R \|  ( 2 \| K \|\|\Delta K \| + \| \Delta K \|^2 )
	\end{align*}
	Lemma~\ref{lemma:useful-bounds-ope-T-F}~\ref{eq:TyK-Ctheta-lambda} and Lemma~\ref{lemma:perturbation_variance_Cy} implies $\vertiii{\cT_{K'}} \leq \| \Sigma_{K'} \|/ \lambda_y^1 \leq 2 \| \Sigma_K \| / \lambda_y^1$.
	Together with inequalities~\eqref{eq:upper_bound_Riccati_and_Variance_matrix} for $\| P_K^y \|$ and $\| \Sigma_K \|$, we conclude the lemma.
\qed\end{proof}

\noindent \textbf{Perturbation on gradient}
\label{subsection:perturbation_on_gradient}

We show in this subsection that the difference between gradients $\nabla_K C_y$ at $K$ and $K'$ is controlled by $\| K' - K \|$.
To do so, we consider a function
\begin{equation}
	\label{eq:f_grad_K}
	f_{grad}(K) = \Big( 6 \| R \| +  \frac{ 4 \| {\mathrm B}\|^2 C_y(K)}{\lambda_y^1}  + 6 f_{var}(K) \| R \| \| K \| + \frac{3}{2} ( f_{var}(K))^2 \lambda_{min}(Q) \Big)\frac{C_y(K)}{\lambda_{min}(Q)}.
\end{equation}

\begin{lemma}
	\label{lemma:perturbation_gradient_Cy}
	Consider $\theta, \theta' \in \Theta$.
	If $f_{var}(K) \| K' - K \| \leq 1$, then
	\begin{equation}
		\label{eq:perturbation_gradient}
		\| \nabla_K C_y(K') - \nabla_K C_y(K) \| \leq f_{grad}(K) \| K' - K \|.
	\end{equation}
\end{lemma}

\begin{proof}
	By Lemma~\ref{lemma:policy_gradient_expression} on gradient expression $\nabla_K C(K) = 2 E_K \Sigma_K$, we have
	$$
		\| \nabla_K C(K') - \nabla_K C(K) \| \leq 2 \| E_{K'} - E_K \|. \|\Sigma_{K'} \| + 2 \| E_{K} \|. \| \Sigma_{K'} - \Sigma_K \|
	$$
	We bound each term in the above inequality to derive the function $f_{grad}$.
	First, with inequality~\eqref{eq:upper_bound_Riccati_and_Variance_matrix} on $\| P_K^y \|$, we notice that
	$$
		\| E_K \| = \| RK - \gamma {\mathrm B}^\top P_K^y ({\mathrm A} - {\mathrm B} K) \| \leq \| R \| \, \| K \| + \gamma f_{var}(K) \lambda_{min}(Q) / 4.
	$$
	Next, with Lemma~\ref{lemma:perturbation_riccati_Cy} and $\| {\mathrm B} \| \| \Delta K \| \leq f_{var}(K) \| \Delta K \| \leq 1$, we have
	\begin{align*}
		\| E_{K'} - E_K \| & = \| R \Delta K - \gamma {\mathrm B}^\top ( P_{K'}^y - P_K^y) ( {\mathrm A} - {\mathrm B} K - {\mathrm B} \Delta K ) + \gamma {\mathrm B}^\top P_K^y {\mathrm B} \Delta K \|
		\\
		                   & \leq \big( \| R \| + \gamma \| {\mathrm B} \|  ( 1 + \| {\mathrm A} - {\mathrm B} K \| ) f_{riccati}(K, \Delta K ) + \gamma \| {\mathrm B} \|^2 \| P_K^y \| \big) \| \Delta K \|.
	\end{align*}
	By expanding the term $f_{riccati}(K, \Delta K)$ in~\eqref{eq:f_riccati_K}, we obtain
	\begin{align*}
		  & \gamma \| {\mathrm B} \| ( 1 + \| {\mathrm A} - {\mathrm B} K\| ) f_{riccati}(K, \Delta K)
		\\
		= & \frac{\gamma^2}{4} ( f_{var}(K) )^2 \lambda_{min}(Q) + \gamma f_{var}(K) \| R \| \| K \| + \frac{\gamma \| R \|}{2} f_{var}(K) \, \| \Delta K \|.
	\end{align*}
	Together with Lemma~\ref{lemma:perturbation_variance_Cy} on $\| \Sigma_{K'} - \Sigma_K \|$, we conclude that
	\begin{align*}
		\| \nabla_K C(K') - \nabla_K C(K) \| \leq & 2 \Big(  \|R\|  + \gamma^2 ( f_{var}(K) )^2 \lambda_{min}(Q) / 4 + \gamma f_{var}(K) \| R \| \| K \|
		\\
		                                          & \hspace{10pt} + \gamma \| R \| /2 + \gamma \| {\mathrm B} \|^2 C_y(K) / \lambda_y^1 \Big) \| \Delta K \| \big( 2 \| \Sigma_K \| \big)
		\\
		                                          & + 2 \Big( \| R \| \| K \| + \gamma f_{var}(K) \lambda_{min}(Q) / 4 \Big)  f_{var}(K) \| \Sigma_K\| \| \Delta K \|
		\\
		\leq                                      & f_{grad}(K) \, \| \Delta K \|.
	\end{align*}
\qed\end{proof}

\subsection{Admissibility of perturbed control parameter}
\label{subsection:stability_of_small_perturbation}
This subsection is to show that when the perturbations $\| \Delta K \| = \| K' - K \|$ and $\| \Delta L \| = \| L' - L \|$ are small enough, the perturbed parameter $\theta' = (K', L')$ is still admissible, i.e. $\theta' \in \Theta$.

\begin{lemma}
	\label{lemma:stability_of_small_perturbation}
	Consider $K \in \RR^{\ell \times d}$ such that $\gamma \| {\mathrm A} - {\mathrm B} K \|^2 < 1$.
	If a matrix $K' \in \RR^{\ell \times d}$ satisfies $\|{\mathrm B}\|\,\| K' - K \| < 1/\sqrt{\gamma} - \|{\mathrm A}-{\mathrm B}K\|$, then
	$$
		\gamma \| {\mathrm A} - {\mathrm B} K' \|^2 < 1.
	$$
\end{lemma}
\begin{proof}
	By the triangle inequality,
	$$
		\|{\mathrm A}-{\mathrm B}K'\| \leq \|{\mathrm A}-{\mathrm B}K\| + \|{\mathrm B}\|\,\|K'-K\| < 1/\sqrt{\gamma}.
	$$
	This gives $\gamma\|{\mathrm A}-{\mathrm B}K'\|^2<1$.
\qed\end{proof}

\subsection{Proof of Proposition~\ref{pr:PL_cond_Cy_Cz} for PL condition}
\label{subsection:proof_proposition_PL_condition}

\begin{proof} Similar to~\cite[Lemma 11]{fazel2018global}, we first notice that for real values $a, b, c \in \RR$ with $c \neq 0$, $2ab + a^2 c = (a + b/c)^2 c - b^2 / c \geq - b^2 / c$, and this inequality still holds for matrices with proper dimensions.
	For $n \geq 0$, let $a = (K' - K) Y_{n}^{K', Y_{0}}$, $b=E_K Y_{n}^{K', Y_{0}} $, $c=R + \gamma {\mathrm B}^\top P_K^y {\mathrm B}$, the advantage quantity $A_y(n, K, Y_{n}^{K', Y_{0}}, K')$ in Lemma~\ref{lemma:difference_in_C_y(K)} satisfies
	$$
		A_y(n,K, Y_{n}^{K', Y_{0}}, K')  \geq - ( E_K Y_{n}^{K', Y_{0}})^\top (R + \gamma {\mathrm B}^\top P_K^y {\mathrm B})^{-1} E_K Y_{n}^{K', Y_{0}}.
	$$
	Then, together with Lemma~\ref{lemma:policy_gradient_expression} for the gradient expression, we deduce that
	\begin{align*}
		C_y(K) - C_y(K^*)
		 & = - \EE \Big[ \sum_{n \geq 0} \gamma^{n} A_y(t, K, Y_{n}^{K^*, Y_{0}}, K^*) \Big]
		\\
		 & \leq \| (R + \gamma {\mathrm B}^\top P_K^y {\mathrm B})^{-1} \| \| \Sigma_{K^*} \|  \| E_K \|_F^2
		\\
		 & \leq  \frac{\Vert \Sigma_{K^*} \Vert}{4 \lambda_{min}(R) \lambda_{min}(\Sigma_K)^2} \| \nabla_K C_y(K) \|_F^2.
	\end{align*}
	Because $\lambda_y^1 \leq \lambda_{min}( \Sigma_K)$ (see Lemma~\ref{lemma:expression_Sigma_K_with_T_K}), let $\nu_y = 4 \lambda_{min}(R). (\lambda_y^1)^2  / \| \Sigma_{K^*} \| $, we conclude that $C_y$ satisfies the $\nu_y-$PL condition.
\qed\end{proof}

We can also derive a lower bound for $C_y(K) - C_y(K^*)$.
Consider $K' = K - (R + {\mathrm B}^\top P_K^y {\mathrm B})^{-1} E_K$ so that the advantage quantities in Lemma~\ref{lemma:difference_in_C_y(K)} attain the lower bounds.
Consequently, with $R \succ 0$ and $P_K^y \succ 0$ (Assumption~\ref{ass:positivity-qr} and Proposition~\ref{prop:existence_of_lyapunov_equation}), we have
\begin{align}
	C_y(K) - C_y(K^*)
	 & \geq C_y(K) - C_y(K')
	\\
	 & = \langle (R + \gamma {\mathrm B}^\top P_K^y {\mathrm B})^{-1} E_K, \,  E_K \Sigma_{K'}  \rangle_{tr}
	\nonumber                                                                                                \\
	 & \geq  \lambda_y^1 \, \| E_K \|_F^2  / \| R + \gamma {\mathrm B}^\top P_K^y {\mathrm B} \|.
	\label{eq:lower_bound_for_Cy_minus_Copt}
\end{align}

\subsection{Proof of Proposition~\ref{pr:local_smoothness_Cy_Cz} for local smoothness}
\label{subsection:proof_proposition_local_smoothness}

\begin{proof}
	Consider $\Delta K = K' - K = - \eta \nabla_K C_y(K)$ with $\eta > 0$.
	By Lemma~\ref{lemma:policy_gradient_expression} on the gradient $\nabla_K C_y(K)$, we have
	\begin{align*}
		\langle \Delta K , 2 E_K \Sigma_{K'} \rangle_{tr}
		 & =  \langle \Delta K,  \nabla_K C_y(K) + 2 E_K (\Sigma_{K'} - \Sigma_K) \rangle_{tr}
		\\
		 & = \langle \Delta K,  \nabla_K C_y(K) \rangle_{tr} -  \langle 2 \eta E_K \Sigma_K, 2E_K (\Sigma_{K'} - \Sigma_K) \rangle_{tr}
		\\
		 & \leq \langle \Delta K,  \nabla_K C_y(K) \rangle_{tr} + 4\eta \big| \langle (E_K)^\top E_K \Sigma_K, \Sigma_{K'} - \Sigma_K \rangle_{tr} \big|
		\\
		 & \leq \langle \Delta K,  \nabla_K C_y(K) \rangle_{tr} + 4\eta \langle (E_K)^\top E_K \Sigma_K, \Sigma_K \rangle_{tr} \frac{\| \Sigma_{K'} - \Sigma_{K} \|}{\lambda_{min}(\Sigma_K)}
		\\
		 & \leq \big( 1 - \lambda_{var, y}(K, K') \big)  \langle \Delta K, \nabla_K C_y(K) \rangle_{tr},
	\end{align*}
	where $\lambda_{var, y}(K, K') = \| \Sigma_{K'} - \Sigma_K \| / \lambda_y^1$.
	Combining this with Lemma~\ref{lemma:difference_in_C_y(K)}, we obtain the local smoothness inequality~\eqref{eq:local_smoothness_Cy_Cz} for $C_y$ in Proposition~\ref{pr:local_smoothness_Cy_Cz}.
\qed\end{proof}

\subsection{Bounds in perturbation analysis}
\label{subsection:bound_on_local_smoothness_coefficient}
To demonstrate the convergence results for the exact PG descent algorithm in Theorem~\ref{th:exact-CV}, we first establish several bounds on the coefficients $\lambda_{var, y}$ and $\lambda_{hess,y}$ introduced in Proposition~\ref{pr:local_smoothness_Cy_Cz}.

\begin{lemma}
	\label{lemma:useful_bounds_I}
	Consider $\theta \in \Theta$ and $C(\theta) \leq C_0$ for some constant $C_0 \in \RR$.
	We define the following constants:
	\begin{align}
		 & h_{small-pert, y}(C_0) = 4 \|{\mathrm B} \| \big( 1 + 1 / \sqrt{\gamma} \big) \frac{C_0}{\lambda_y^1 \lambda_{min}(Q)}
		\label{eq:bound_f_var_K}
		\\
		 & h_{grad,y}(C_0)  = 2 \frac{C_0}{\lambda_{min}(Q)} \sqrt{\Big( \| R \| + \gamma \| {\mathrm B} \|^2 \frac{C_0}{ \lambda_y^1} \Big) \frac{C_0}{\lambda_y^1}}
		\label{eq:bound_grad_Cy_K}
		\\
		 & h_{K}(C_0)  = \frac{1}{\lambda_{min}(R)} \Big(   \sqrt{\Big( \| R \| + \gamma \| {\mathrm B} \|^2 \frac{C_0}{ \lambda_y^1} \Big) \frac{C_0}{\lambda_y^1}} + \gamma \|{\mathrm B} \| \| {\mathrm A} \| \frac{C_0}{\lambda_y^1} \Big)
		\label{eq:bound_K}
	\end{align}
	Then the following inequalities hold:
	\begin{equation}
		\label{eq:inequality_useful_bounds_I}
		f_{var}(K) \leq h_{small-pert,y}(C_0),
		\quad \| \nabla_K C_y(K) \| \leq h_{grad,y}(C_0), \quad \| K \| \leq h_{K}(C_0).
	\end{equation}
	Similar bounds $h_{small-pert, z}(C_0)$, $h_{grad, z}(C_0)$, and $h_{L}(C_0)$ can be defined by replacing $({\mathrm B}, Q, R, \lambda_y^1)$ with $(\tilde {\mathrm B}, \tilde Q, \tilde R, \lambda_z^0)$ in~\eqref{eq:bound_f_var_K},~\eqref{eq:bound_grad_Cy_K}, and~\eqref{eq:bound_K}.
\end{lemma}
\begin{proof}
	By definition of $\theta \in \Theta$, we have $\| {\mathrm A} - {\mathrm B} K \| \leq 1 / \sqrt{\gamma}$, then equation~\eqref{eq:f_var_K} for $f_{var}(K)$ implies the first inequality in~\eqref{eq:inequality_useful_bounds_I}.
	From equation~\eqref{eq:lower_bound_for_Cy_minus_Copt}, we have $\| E_K \|^2 \leq \| E_K \|_F^2 \leq ( \|R\| + \gamma \|{\mathrm B}\|^2 \| P_K^y \|) C_y(K) / \lambda_y^1$.
	Together with the bounds on $\| P_K^y \|$ and $\| \Sigma_K \|$ in~\eqref{eq:upper_bound_Riccati_and_Variance_matrix} and the gradient expression from Lemma~\ref{lemma:policy_gradient_expression}, we derive the second inequality in~\eqref{eq:inequality_useful_bounds_I}.
	Moreover, we notice that
	$$
		\| K \| = \| (R + \gamma {\mathrm B}^\top P_K^y {\mathrm B})^{-1} ( R + \gamma {\mathrm B}^\top P_K^y {\mathrm B}) K \| \leq \frac{1}{\lambda_{min}(R)} ( \| E_K \| + \gamma \| {\mathrm B}^\top P_K^y {\mathrm A} \| ),
	$$
	we obtain the third inequality in~\eqref{eq:inequality_useful_bounds_I}.
\qed\end{proof}

\begin{lemma}
	\label{lemma:bound_on_coefficent_in_prop_local_smoothness}
	Consider $\theta = (K, L) \in \Theta$ such that $C_y(K) \leq C(\theta) \leq C_0$ with a constant $C_0 \in \RR$.
	Consider the following two constants in $\RR$ defined by
	\begin{align}
		h_{var, y}(C_0)  & = 8 \frac{ \big( 1 + 1 / \sqrt{\gamma} \big) \|{\mathrm B} \| ( C_0 )^3 }{ (\lambda_y^1)^2 (\lambda_{min}(Q) )^3} \sqrt{\Big(\| R \| + \gamma \| {\mathrm B} \|^2 \frac{C_0}{\lambda_y^1} \Big) \frac{C_0}{\lambda_y^1}},
		\label{eq:h_var_y}
		\\
		h_{hess, y}(C_0) & = 2 \frac{C_0}{\lambda_{min}(Q)}\Big( \| R \| + \gamma \| {\mathrm B} \|^2 \frac{C_0}{\lambda_y^1} \Big).
		\label{eq:h_hess_y}
	\end{align}
	Then
	$$
		f_{var}(K) \| \nabla_K C_y(K) \|\frac{ \| \Sigma_K \| }{\lambda_y^1} \leq h_{var, y}(C_0).
	$$
	Moreover, for a matrix $K' \in \RR^{\ell \times d}$ defined by $K' = K - \eta \nabla_K C_y(K)$ with some $\eta > 0$, if $\eta \leq h_{var, y}^{-1}(C_0)$ and $\|{\mathrm B}\|\,\eta h_{grad,y}(C_0) < 1/\sqrt{\gamma}-\|{\mathrm A}-{\mathrm B}K\|$, then
	\begin{align*}
		 & \| K' - K \| f_{var}(K) \leq 1,
		\\
		 & \gamma \| A - B K' \|^2 < 1,
		\\
		 & \lambda_{var, y}(K, K') := \| \Sigma_{K'} - \Sigma_K \| / \lambda_y^1  \leq \eta h_{var, y}(C_0),
		\\
		 & \lambda_{hess, y}(K, K') := \| \Sigma_{K'} \| \big( \|R \| + \gamma \|B \|^2 C_y(K) / \lambda_y^1 \big) \leq h_{hess, y}(C_0).
	\end{align*}
\end{lemma}
\begin{proof}
	From inequality~\eqref{eq:bound_f_var_K} and~\eqref{eq:bound_grad_Cy_K} and $ 1 \leq \frac{ \| \Sigma_K \| }{ \lambda_y^1} \leq \frac{ C_0 }{ (\lambda_y^1 \lambda_{min}(Q)) }$, we have
	$$
		f_{var}(K) \| \nabla_K C_y(K) \| \frac{ \| \Sigma_K \| }{\lambda_y^1}  \leq h_{small-pert,y}(C_0) h_{grad,y}(C_0) \frac{C_0}{\lambda_y^1 \lambda_{min}(Q)} = h_{var, y}(C_0).
	$$
	For the $K'$ considered in the lemma, we have
	$$
		\| K' - K \| f_{var}(K) = \eta f_{var}(K) \| \nabla_K C_y(K) \| \leq \eta h_{var, y}(C_0) \leq 1.
	$$
	Since $\|K'-K\| \leq \eta h_{grad,y}(C_0)$, the margin condition and Lemma~\ref{lemma:stability_of_small_perturbation} imply that $\gamma \| A - B K' \|^2 < 1$. Lemma~\ref{lemma:perturbation_variance_Cy} then implies
	$
		\| \Sigma_{K'} - \Sigma_K \| \leq f_{var}(K) \| \Sigma_K \| \| \Delta K \| \leq \| \Sigma_K \|.
	$
	As a result,
	$$
		\lambda_{var, y}(K, K') \leq  \eta  f_{var}(K) \| \nabla_K C_y(K) \| \frac{ \| \Sigma_K \|}{\lambda_y^1} \leq \eta h_{var, y}(C_0)
	$$
	and
	$$
		\lambda_{hess, y}(K, K') \leq \big( \| \Sigma_K \| + \| \Sigma_{K'} - \Sigma_K \| \big) \big( \|R \| + \gamma \|B \|^2 C_y(K) / \lambda_y^1 \big) \leq h_{hess, y}(C_0).
	$$
\qed\end{proof}

\subsubsection{Perturbation on cost function}
\label{section:perturbation_on_cost_function}

\begin{lemma}
	\label{lemma:perturbation_cost_Cy}
	Consider $\theta \in \Theta$ and $C(\theta) \leq C_0$ with some $C_0 \in \RR$.
	Let
	\begin{align}
		\label{eq:h_raccati}
		h_{riccati, y}(C_0) & = \frac{C_0}{\lambda_y^1} \Big( \gamma h_{small-pert,y}(C_0) + \frac{4 \| R \|}{\lambda_{min}(Q)} h_{K}(C_0) + \frac{2 \| R \| }{\| {\mathrm B} \| \lambda_{min}(Q) } \Big)
		\\
		h_{cost, y}(C_0)    & = \frac{2d. h_{riccati, y}(C_0) C_{init, noise}^2}{1- \gamma}
		\label{eq:h_cost_y}
		\\
		h_{f\_grad\_K}(C_0) & = \frac{C_0}{\lambda_{min}(Q)}\big( 6 \| R \| ( 1+ h_{small-pert, y}(C_0) h_K(C_0) )
		\nonumber                                                                                                                                                                                         \\
		                    & \hspace{60pt} + 4 \| {\mathrm B} \|^2 C_0 / \lambda_y^1 + 3h_{small-pert,y}^2(C_0) \lambda_{min}(Q) / 2 \big).
		\label{eq:h_f_grad_K}
	\end{align}
	If $f_{var}(K) \| K' - K \| \leq 1$, then
	\begin{align*}
		| C_y(K') - C_y(K) |                     & \leq h_{cost, y}(C_0) \| K' - K \|,
		\\
		\| \nabla_K C_y(K') - \nabla_K C_y(K) \| & \leq h_{f\_grad\_K}(C_0) \| K' - K \|.
	\end{align*}
	Similar constants $h_{riccati, z}(C_0)$, $h_{cost, z}(C_0)$, $h_{f\_grad\_L}(C_0)$ are defined with $(\lambda_z^0, \tilde {\mathrm B}, \tilde Q, \tilde R)$ and with $h_{small-pert, z}(C_0), h_L(C_0)$ from Lemma~\ref{lemma:useful_bounds_I}.
	Let $h_{cost}(C_0) = h_{cost, y}(C_0) + h_{cost, z}(C_0)$ and $\| \theta' - \theta \| = \| K' - K \| + \| L' - L \|$, then
	\begin{equation}
		\label{eq:perturbation_cost}
		| C(\theta') - C(\theta) | \leq h_{cost}(C_0) \, \| \theta' - \theta \|.
	\end{equation}
\end{lemma}

\begin{proof}
	From the cost expression~\eqref{eq:cost_expression_with_PK_y} $C_y(K) = \langle P_K^y, \Sigma_{Y_{0}} + \frac{\gamma}{1-\gamma} \Sigma^1 \rangle_{tr}$, and from Lemma~\ref{lemma:perturbation_riccati_Cy} on the perturbation between $\| P_{K'}^y - P_K^y \|$, we have
	\begin{align*}
		| C_y(K') - C_y(K) | & =  \Big| \langle P_{K'}^y - P_{K}^y,  \Sigma_{Y_{0}} + \gamma \Sigma^1 / (1-\gamma) \rangle_{tr} \Big|
		\\
		                     & \leq h_{riccati, y}(C_0) \| \Delta K \|  ( \| Y_{0} \|^2_{L^2(\RR^d)} + \gamma \| \varepsilon_{n} \|^2_{L^2(\RR^d)} / ( 1- \gamma) ).
	\end{align*}
	Then, by the fact that $\| \xi \|_{L^2(\RR^d)} \leq \sqrt{2d} \| \xi \|_{\psi_2} $ for a sub-Gaussian random vector $\xi \in \RR^d$, and $\| Y_{0} \|_{\psi_2}, \| \varepsilon_{n} \|_{\psi_2} \leq C_{init, noise}$, we obtain inequality $| C_y(K') - C_y(K) | \leq h_{cost, y}(C_0) \| \Delta K \|$.
	By Lemma~\ref{lemma:perturbation_gradient_Cy} and the definition of $f_{grad}(K)$ in~\eqref{eq:f_grad_K}, we obtain inequality $\| \nabla_K C_y(K') - \nabla_K C_y(K) \| \leq h_{f\_grad\_K}(C_0) \| \Delta K \|$.
\qed\end{proof}

\section{Technical lemma for Section~\ref{subsection:modelfree_MKV_simulator}}
\label{sec:app-proof-modelfree-MKV-CV}

\subsection{Intuition for zeroth-order optimization}

Derivative-free optimization (see e.g.~\cite{MR2487816,MR3627456}) tries to find an optimizer of a function $f: x \mapsto f(x)$ using only pointwise values of $f$, without having access to its gradient.
The basic idea is to express the gradient of $f$ in a point $x$ by using values of $f(\cdot)$ at a Gaussian perturbation region around point $x$.
This can be achieved by introducing an approximation of $f$ using a Gaussian smoothing.
Formally, if we define the perturbed function for some $\sigma >0$ as
$
	f_{\sigma^2}(x) = \mathbb{E}[f(x+ \varepsilon)],
$ where
$\varepsilon \sim \mathcal{N}(0,\sigma^2 I)$,
then its gradient can be written as
$
	\frac{d}{dx}f_{\sigma^2}(x) = \frac{1}{\sigma^2} \mathbb{E}[f(x+\varepsilon) \varepsilon].
$
We can replace the Gaussian smoothing term $\varepsilon$ by a uniform random variable with bounded norms.

Consider $U = (U^{(idy)}, U^{(com)})$ with $U^{(idy)},U^{(com)}$ being two independent random matrices in $\RR^{\ell \times d}$  uniformly distributed on the ball with radius $\tau$: $\mathbb{B}_{\tau} = \{ W \in \RR^{\ell \times d} : \| W \|_F \leq \tau \}$.
The unit ball under Frobenius norm for matrices in $\RR^{\ell \times d}$ can be identified as a unit ball for vectors in $\RR^{\ell d}$.
Let $\mathbb{S}_\tau = \{ W \in \RR^{\ell \times d}: \| W \|_F = \tau \}$ be the sphere of $\mathbb{B}_\tau$ of dimension $\ell d -1$.
The ratio between the sphere area of a ball with radius $\tau$ in $\RR^{\ell d}$  to its volume equals $(\ell d) / \tau$~\cite{MR2298287}.
Furthermore, the spectral norm of $W \in \mathbb{B}_\tau$ satisfies that $\| W \| \leq \| W \|_F = \tau$.

For any $\theta = (K,L)$ and $\tau > 0$, we introduce the following smoothed version $C_\tau$ of $C$, defined by:
\begin{equation}
	\label{eq:formula_perturbation_definition_C}
	C_\tau(\theta) = \mathbb{E}_U[C(\theta + U)] = \mathbb{E}_{(U^{(idy)},U^{(com)})}[C_y(K+U^{(idy)}) + C_z(L+U^{(com)})],
\end{equation}
where $\mathbb{E}_U$ denotes the expectation over $U$, and the costs $C_y(K+U^{(idy)})$ and $C_z(L + U^{(com)})$ are evaluated using the auxiliary processes $\bY^{K + U^{(idy)}}$ and $\bZ^{L + U^{(com)}}$ as defined in equations~\eqref{eq:def_Cy_Cz}.
Let $V^{(idy)} \sim \tau U^{(idy)} / \| U^{(idy)} \|_F$ and $V^{(com)} \sim \tau U^{(com)} / \| U^{(com)} \|_F$ be two random matrices independently and uniformly distributed on the sphere $\SS_\tau$; then $V = (V^{(idy)}, V^{(com)}) \in \SS_\tau \times \SS_\tau$.
From~\cite[Lemma 2.1]{MR2298287} or~\cite[Lemma 26]{fazel2018global}, we get the following expressions for the gradients of $C_\tau(\theta)$:
\begin{align}
	\label{eq:gradCK-zero-order}
	 \nabla_K C_\tau(\theta) & = \frac{\ell d}{\tau^2} \mathbb{E}_{V^{(idy)}} \big[ C_y(K + V^{(idy)}) V^{(idy)} \big] = \frac{\ell d}{\tau^2}\mathbb{E}_V \big[C(\theta + V) V^{(idy)} \big],
	\\
	\label{eq:gradCL-zero-order}
	 \nabla_L C_\tau(\theta) & = \frac{\ell d}{\tau^2} \mathbb{E}_{V^{(com)}} \big[ C_z(L + V^{(com)}) V^{(com)} \big] = \frac{\ell d}{\tau^2} \mathbb{E}_V \big[C(\theta+V) V^{(com)} \big].
\end{align}
The last equalities in~\eqref{eq:gradCK-zero-order} and~\eqref{eq:gradCL-zero-order} are justified by Lemma~\ref{lem:decompose-cost-reparam} and the fact that $V^{(idy)}$ and $V^{(com)}$ are independent zero-mean random variables.

To show that the gradient estimator from Algorithm~\ref{algo:MKVestim} provides a good approximation of the exact PG, we introduce two more gradient estimators before we dive into the perturbation analysis in this section.
Let $M$ be the number of perturbation directions, $\tau > 0$ be the perturbation radius, and $T$ be the truncation horizon for the infinite discounted MF cost. The following four notations of gradient will be used in the next subsections:
\begin{enumerate}[label=(\roman*)]
	\item \textbf{[Policy gradient]} $\nabla C(\theta) = \big(\nabla_K C(\theta), \nabla_L C(\theta) \big)$ represents the exact gradient of cost $C(\theta)$ with respect to $\theta = (K, L)$.

	\item \textbf{[Perturbed policy gradient]}	$\hat{\nabla}^{M,\tau}(\theta, \underline v)$ defined in~\eqref{eq:def_perturbed_PG_MKV} below is an approximation of exact PG based on the zeroth-order approximation technique with $M$ perturbation directions $\underline v = (v_i)_{i=1}^M$ where $v_i \sim \mu_{\mathbb{S}_\tau} \otimes \mu_{\mathbb{S}_\tau}$.

	\item \textbf{[Truncated policy gradient]} $\hat{\nabla}^{T, M, \tau}(\theta, \underline v)$ defined in~\eqref{eq:def_{n}runcated_gradient_MKV} below is an approximation of perturbed policy gradient using a cost truncated to a finite horizon $T$.

	\item \textbf{[Sampled policy gradient]} $\tilde{\nabla}^{T, M,\tau} \big(\theta, \underline v \big)$ defined in equation~\eqref{eq:definition_approx_gradient_MKV_simulator} is the output of Algorithm~\ref{algo:MKVestim}.
\end{enumerate}

Before detailing the exact bounds, we sketch the asymptotic requirements to achieve an $\varepsilon$-approximation of the exact PG. We will show that the perturbation radius must scale as $\tau = \mathcal{O}(\varepsilon)$, the truncation horizon as $T = \mathcal{O}(\log(1/\varepsilon))$, and the sample size as $M = \tilde{\mathcal{O}}(\varepsilon^{-4})$. The following sections make these dependencies explicit through threshold functions $\phi_{pert, radius}, \phi_{trunc, T},$ and $\phi_{sample, size}$.

Our proof strategy bounds the total approximation error via the triangle inequality:
\begin{align*}
	\| \tilde{\nabla}^{T, M, \tau} - \nabla C \| & \leq \underbrace{\| \tilde{\nabla}^{T, M, \tau} - \hat{\nabla}^{T, M, \tau} \|}_{\text{Sampling Error}} 
	+ \underbrace{\| \hat{\nabla}^{T, M, \tau} - \hat{\nabla}^{M, \tau} \|}_{\text{Truncation Error}}
	+ \underbrace{\| \hat{\nabla}^{M, \tau} - \nabla C \|}_{\text{Zeroth-order Error}}.
\end{align*}
The following subsections bound each of these distinct sources of error.

\subsection{Approximation with smoothed cost}
\label{subsection:approximation_with_smoothed_cost}
In this subsection, we show an approximation result of the exact PG $\nabla C(\theta)$ using the perturbed policy gradient $\hat{\nabla}^{M, \tau}(\theta, \underline v)$ defined in~\eqref{eq:def_perturbed_PG_MKV} below.
To start, we consider the gradient of a smooth cost $C_\tau(\theta)$ defined in~\eqref{eq:formula_perturbation_definition_C} with respect to $\theta = (K, L)$ as $\nabla C_\tau(\theta) = \big( \nabla_K C_\tau(\theta) , \nabla_L C_\tau(\theta) \big)$.
Equations~\eqref{eq:gradCK-zero-order} and~\eqref{eq:gradCL-zero-order} provide the expression for $\nabla C_\tau(\theta)$ without taking the gradients.
Now, we approximate the expectations $\EE_V$ therein using $M$ independent samples $(v_i^{(idy)})_{i=1}^M$ of $V^{(idy)}$ and $M$ independent samples $(v_i^{(com)})_{i=1}^M$ of $V^{(com)}$ respectively.
Let $\underline v = (v_i)_{i=1}^M$ with $v_i = (v_{i}^{(idy)}, v_{i}^{(com)}) \in \RR^{\ell \times d} \times \RR^{\ell \times d}$ and $v_i \sim \mu_{\mathbb{S}_\tau} \otimes \mu_{\mathbb{S}_\tau}$.
We define the ``perturbed policy gradient" in $K$ and in $L$ by
\begin{equation}
	\label{eq:def_perturbed_PG_MKV}
	\hat{\nabla}_K^{M, \tau}(\theta, \underline v) = \frac{\ell d}{\tau^2} \frac{1}{M} \sum_{i=1}^M C(\theta + v_i) v_i^{(idy)},
	\quad
	\hat{\nabla}_L^{M, \tau}(\theta, \underline v) = \frac{\ell d}{\tau^2} \frac{1}{M} \sum_{i=1}^M C(\theta + v_i) v_i^{(com)}.
\end{equation}
Let $\hat{\nabla}^{M, \tau}(\theta, \underline v) = \big( \hat{\nabla}_K^{M, \tau}(\theta, \underline v),  \hat{\nabla}_L^{M, \tau}(\theta, \underline v) \big)$.
In~\eqref{eq:def_perturbed_PG_MKV}, we consider $\tau$ small enough so that by the stability Lemma~\ref{lemma:stability_of_small_perturbation}, we have $\theta + v_i \in \Theta$ for all $i=1, \ldots, M$, and thus $\hat{\nabla}^{M, \tau}(\theta, \underline v)$ is well-defined.
It is clear that the randomness of the perturbed policy gradient comes from the random directions $\underline v \in (\mathbb{S}_{\tau} \times \mathbb{S}_\tau )^{M}$.

\begin{lemma}
	\label{lemma:approx_perturbed_policy_gradient_K}
	Consider $\theta \in \Theta$ with $C(\theta) \leq C_0$ for some constant $C_0 \in \RR$.
	For any $\varepsilon > 0$ and $\delta_{pert, grad} \in (0,1)$, if
	\begin{align}
		\tau & \leq \min\{ C_0 h_{cost}^{-1}(C_0), h_{small-pert, y}^{-1}(C_0), \varepsilon h_{f\_grad\_K}(C_0)^{-1} / 2 \}
		\label{eq:condition_{n}au_K}
	\end{align}
	where $h_{cost}(C_0), h_{small-pert, y}(C_0), h_{f\_grad\_K}(C_0)$ are defined in Lemma~\ref{lemma:useful_bounds_I} and Lemma~\ref{lemma:perturbation_cost_Cy}, then, we have
	$$
		\PP_{\underline v} \left( \| \hat\nabla^{M, \tau}_K(\theta, \underline v) - \nabla_K C(\theta) \| > \varepsilon \right) \leq (\ell + d) \exp \Big( \frac{-M \varepsilon^2 \tau^2 } { 64 \ell^2 d^2 C_0^2 + 16 \ell d C_0 \varepsilon / 3 } \Big).
	$$
\end{lemma}

\begin{proof}
	For $i=1, \ldots, M$, let $\theta_i = \theta + v_i$ and $W_i = \frac{\ell d}{\tau^2 M} \big( C(\theta_i) \tilde{v}_i^{(idy)} - \EE_{v_i}[ C(\theta + v_i) \tilde{v}_i^{(idy)} ] \big)$ where $\tilde v_i^{(idy)}$ is the dilatation matrix of $v_i^{(idy)}$ and thus $\| \tilde v_i^{(idy)} \| = \| v_i^{(idy)} \| \leq \| v_i^{(idy)} \|_F = \tau$.
	For $\tau$ small enough such that $\tau \leq C_0 h_{cost}(C_0)^{-1}$ where $h_{cost}(C_0)$ defined in Lemma~\ref{lemma:perturbation_cost_Cy}, we have $| C(\theta_i) - C(\theta) | \leq h_{cost, y}(C_0) \| v_i^{(idy)} \| + h_{cost, z}(C_0) \| v_i^{(com)} \| \leq  h_{cost}(C_0) \tau \leq C_0$ and so $\| C(\theta_i) \tilde{v}_i^{(idy)} - \EE_{v_i}[C(\theta_i) \tilde{v}_i^{(idy)}] \| \leq 4C_0 \tau$.
	Thus, for $i=1, \ldots, M$,
	$$
		\| W_i \| \leq \frac{4 \ell d C_0 }{M \tau} =: \lambda_{W},
		\qquad
		\Big \| \sum_{i=1}^M \EE[ W_i^2 ] \Big\| \leq M \sup_{i=1,\ldots, M} \EE[ \| W_i \|^2 ] \leq \frac{8 \ell^2 d^2 C_0^2}{M \tau^2} = :\sigma^2_W.
	$$
	By matrix Bernstein's inequality with bounded random matrices~\cite[Theorem 6.1]{MR2946459}, and by gradient expressions~\eqref{eq:gradCK-zero-order} for $\nabla_K C_\tau(\theta)$, we obtain that for any $\tilde \varepsilon > 0$,
	\begin{align*}
		\PP \big( \| \hat{\nabla}_K^{M, \tau}(\theta, \underline v) - \nabla_K C_\tau(\theta) \| > \tilde{\varepsilon} \big) & =
		\PP \big( \big\| \sum_{i=1}^M W_i \big\|  > \tilde{\varepsilon} \big) \leq (\ell + d) \exp \Big(\frac{- \tilde{\varepsilon}^2 / 2}{ \sigma^2_W + \tilde{\varepsilon} \lambda_W / 3}  \Big).
	\end{align*}
	Moreover, when $\| (K + U^{(idy)}) - K \| = \| U^{(idy)} \| \leq \| U^{(idy)} \|_F = \tau \leq h_{small-pert,y}^{-1}(C_0)$, Lemma~\ref{lemma:perturbation_gradient_Cy} on the perturbation of gradients implies that
	\begin{align*}
		\| \nabla_K C_\tau(\theta) - \nabla_K C(\theta) \| & = \| \EE_{U}\big[ \nabla_K C_y( K+ U^{(idy)}) - \nabla_K C_y(K) \big] \|
		\\
		                                                   & \leq \EE_U[ f_{grad}(K) \| U \| ] \leq h_{f\_grad\_K}(C_0) \tau,
	\end{align*}
	where $h_{f\_grad\_K}(C_0)$ is defined in~\eqref{eq:h_f_grad_K}.
	Hence, with high probability, we have
	\begin{align*}
		\| \hat\nabla^{M, \tau}_K( \theta, \underline v) - \nabla_K C(\theta) \| & \leq
		\|  \hat\nabla^{M, \tau}_K(\theta, \underline v) - \nabla_K C_\tau(\theta) \| + \| \nabla_K C_\tau(\theta) - \nabla_K C(\theta) \|
		\\
		                                                                         & \leq  \tilde{\varepsilon} + h_{f\_grad\_K}(C_0) \tau.
	\end{align*}
	By choosing $\tilde \varepsilon = \varepsilon /2$, and $\tau$ satisfies conditions~\eqref{eq:condition_{n}au_K}, we have $ \tilde{\varepsilon} + h_{f\_grad\_K}(C_0) \tau \leq \varepsilon$ and we conclude the results.
\qed\end{proof}

Based on the conditions in Lemma~\ref{lemma:approx_perturbed_policy_gradient_K}, we define the following two functions for the perturbation radius and the number of perturbation directions:
\begin{align}
	 & \phi_{pert, radius}(\varepsilon, C_0) = \max \big\{ h_{cost}(C_0) / C_0,  h_{small-pert, y}(C_0), h_{small-pert, z}(C_0),
	\nonumber                                                                                                                                                                          \\
	 & \hspace{120pt}  4 h_{f\_grad\_K}(C_0) / \varepsilon, 4 h_{f\_grad\_L}(C_0) / \varepsilon \big\}
	\label{eq:phi_pert_radius}
	\\
	 & \phi_{pert, size}(\varepsilon, \tau, C_0, \delta_{pert} ) =  \log \Big( \frac{ 2(\ell + d) }{\delta_{pert}} \Big) \Big( \frac{16 \ell d  C_0 }{ \tau \varepsilon } + 1 \Big)^2
	\label{eq:phi_pert_size}
\end{align}
where $h_{cost}(C_0)$, $h_{small-pert, y}(C_0)$, $h_{small-pert, z}(C_0)$, $h_{f\_grad\_K}(C_0)$, $h_{f\_grad\_L}(C_0)$ are defined in Lemma~\ref{lemma:useful_bounds_I} and Lemma~\ref{lemma:perturbation_cost_Cy}.

\begin{lemma}
	\label{lemma:approx_perturbed_gradient_MKV_simulator}
	Consider $\theta \in \Theta$ with $C(\theta) \leq C_0$ for some constant $C_0 \in \RR$.
	For a target precision $\varepsilon > 0$, if
	\begin{equation}
		\label{eq:condition_{n}au_M_approx_perturbed_gradient}
		\tau^{-1} \geq \phi_{pert, radius}(\varepsilon, C_0)
		\quad \text{and} \quad
		M \geq \phi_{pert, size}(\varepsilon, \tau, C_0, \delta_{pert}),
	\end{equation}
	then
	$$
		\PP_{\underline v}\big( \| \hat \nabla^{M, \tau}(\theta, \underline v) - \nabla C(\theta) \| > \varepsilon \big) \leq \delta_{pert},
	$$
	where $\| \hat \nabla^{M, \tau}(\theta, \underline v) - \nabla C(\theta) \|  = \| \hat \nabla^{M, \tau}_K(\theta, \underline v) - \nabla_K C(\theta) \| + \| \hat \nabla^{M, \tau}_L(\theta, \underline v) - \nabla_L C(\theta) \|$.
\end{lemma}

\begin{proof} We apply Lemma~\ref{lemma:approx_perturbed_policy_gradient_K} with target precision $\varepsilon / 2$, then for $\tau$ small enough and $M$ large enough satisfying condition~\eqref{eq:condition_{n}au_M_approx_perturbed_gradient}, we have
	$
		\PP_{\underline v}\big( \| \hat \nabla^{M, \tau}_K(\theta, \underline v) - \nabla_K C(\theta) \| > \varepsilon / 2 \big) \leq \delta_{pert}.
	$
	Together with the similar concentration inequality for $\| \hat \nabla^{M, \tau}_L(\theta, \underline v) - \nabla_L C(\theta) \|$, we conclude the lemma.
\qed\end{proof}

\subsection{Approximation with truncated cost \label{subsection:approximation_C_with_C_T}}
In this subsection, we first provide an approximation of the MF cost $C(\theta)$ with a truncated cost, then we show a result for the approximation of truncated policy gradient $\hat\nabla^{T, M, \tau}$ defined later in the subsection.
For $T \geq 0$, let us define the truncated variances and the truncated costs by:
\begin{align}
	\label{eq:truncation_Sigma_y_Sigma_z_T}
	\Sigma_{K}^{y,T} & = \mathbb{E}\big[\sum_{n=0}^{T-1} \gamma^{n} (Y_{n}^{K}) (Y_{n}^K)^\top \big],
	\quad
	\Sigma_{L}^{z,T}  = \mathbb{E} \big[ \sum_{n=0}^{T-1} \gamma^{n} (Z_{n}^L) (Z_{n}^L)^\top \big],
	\\
	\label{eq:truncation_Cy_Cz_T}
	C^T_y(K)         & = \mathbb{E} \big[  \sum_{n=0}^{T-1}  \gamma^{n} f(Y_{n}^{K}, K, Q, R) \big] ,
	\quad
	C^T_z(L) = \mathbb{E} \big[  \sum_{n=0}^{T-1}  \gamma^{n} f(Z_{n}^{L}, L, \tilde{Q}, \tilde{R} ) \big],
\end{align}
where $f(\xi, \phi, q, r) = \xi^{\top} (q + \phi^\top r \phi)\xi$.
We define $C^T(\theta) = \EE[ \sum_{n=0}^{T-1} \gamma^{n} c_{n}(\theta) ]$ where $c_{n}(\theta) = c(X_{n}^\theta, \bar{X}_{n}^\theta, \ctrl_{n}^\theta, \bar \ctrl_{n}^\theta)$ is the instantaneous cost at time $n$.
With the help of auxiliary process $(\bY^K, \bZ^L)$ from~\eqref{eq:dyn_y_theta} and~\eqref{eq:dyn_z_theta}, we have $C^T(\theta) = C_y^T(K) + C_z^T(L)$.

For $\theta = (K, L) \in \Theta$ with $C(\theta) \leq C_0$, assume that the sublevel set has the positive Euclidean margins $m_y(C_0)$ and $m_z(C_0)$ defined in Theorem~\ref{th:exact-CV}. Then
$$
	\gamma \| {\mathrm A} - {\mathrm B} K \|^2 \leq  \gamma (1/\sqrt{\gamma}-m_y(C_0))^2, \quad \gamma \| \tilde {\mathrm A} - \tilde {\mathrm B} L \|^2 \leq \gamma (1/\sqrt{\gamma}-m_z(C_0))^2.
$$
Let us define a constant $\gamma_{pert}$ and a function $\phi_{truncation}(\varepsilon)$ defining the lower bound on the truncation time horizon:
\begin{align}
	 & \gamma_{pert} = \max \big\{ \gamma, \gamma (1/\sqrt{\gamma}-m_y(C_0))^2,\, \gamma (1/\sqrt{\gamma}-m_z(C_0))^2 \big\},
	\label{eq:gamma_pert}
	\\
	 & \phi_{truncation}(\varepsilon) = \Big( \frac{\log( 2d C_{init, noise}^2) - \log(\varepsilon (1 - \gamma_{pert})^2 ) + 1 }{\log(1/\gamma_{pert})} \Big)^2.
	\label{eq:phi_{n}runcation}
\end{align}
We have the following result between the truncated cost $C^T(\theta)$ and the MF cost $C(\theta)$.

\begin{lemma}
	\label{lemma:approx_nruncated_cost}
	Consider $\theta \in \Theta$ with $C(\theta) \leq C_0$.
	For any $\varepsilon > 0$ and $T \geq 2$, if $T \geq \phi_{truncation}(\varepsilon)$, we have
	\begin{equation}
		C(\theta) - C^T(\theta)  \leq \varepsilon \, h_{T}(C_0),
	\end{equation}
	where $h_T(C_0) = d.(\| Q \| + \| \tilde Q\| + \|R \| h_K(C_0)^2 + \| \tilde R \| h_L(C_0)^2 )$ with $h_K(C_0)$, $h_L(C_0)$ defined in Lemma~\ref{lemma:useful_bounds_I}.
\end{lemma}

\begin{proof}
	We first show that if $T \geq \max\{ 2, \phi_{truncation}(\varepsilon)\}$, then $ \big\| \Sigma_K^y - \Sigma_K^{y,T} \big\| \leq \varepsilon$.
	We notice that
	$
		\mathbb{E}\left[ \gamma^{n} y^K_{n} (y^K_{n})^\top \right] = (\cF_K^y)^{n} \left( \Sigma_{Y_{0}} \right) + \sum_{s=0}^{n-1} (\cF_K^y)^s \left( \gamma^{n-s} \Sigma^1 \right).
	$
	So, by Lemma~\ref{lemma:expression_Sigma_K_with_T_K} and the definition of $\cF_K$, $\cT_K$ in~\eqref{eq:def-ope-F_K_y}, we have
	\begin{align*}
		\Sigma_K^y - \Sigma_K^{y, T} &= \Sigma_K^y - \mathbb{E} \big[ \sum_{n=0}^{T-1} \gamma^{n} Y_{n}^K (Y_{n}^K)^\top \big]
		\\
		&= \cT_K \big(\Sigma_{Y_{0}} + \frac{\gamma}{1-\gamma}\Sigma^1 \big) -  \Big[ \sum_{n=0}^{T-1} (\cF_K)^{n} (\Sigma_{Y_{0}})  + \sum_{n=1}^{T-1} \Big(\sum_{s=0}^{n-1} (\cF_K)^s (\gamma^{n-s} \Sigma^1) \Big) \Big]
		\\
		&= \frac{\gamma}{1-\gamma}\cT_K( \Sigma^1 ) +  \sum_{n=T}^\infty (\cF_K)^{n}(\Sigma_{Y_{0}})  - \sum_{s=1}^{T-1}  \gamma^s \Big( \cT_K(\Sigma^1) - \sum_{n=T-s}^\infty (\cF_K)^{n}(\Sigma^1) \Big)
		\\
		=                              & \frac{\gamma^T}{1- \gamma} \cT_K \left( \Sigma^1 \right) + (\cF_K)^T \left( \cT_K(\Sigma_{Y_{0}})  \right)
		+ \sum_{s=1}^{T-1} \gamma^s \left(\cF_K \right)^{T-s} \left( \cT_K (\Sigma^1) \right).
	\end{align*}
	Because $\gamma \leq \gamma_{pert}$, $\vertiii{\cF_K} = \gamma \| {\mathrm A} - {\mathrm B} K \|^2 \leq \gamma_{pert}$, $\vertiii{ \cT_K } \leq (1- \gamma_{pert})^{-1}$, $\| \Sigma_{Y_{0}} \| \leq \EE[\| Y_{0}\|^2] \leq  \EE[  2 d \| Y_{0} \|_{\psi_2}^2 ] \leq 2d C_{init, noise}^2$, and $\| \Sigma^1 \| \leq 2d C_{init, noise}^2$, we obtain that
	\begin{align*}
		\big\| \Sigma_K^y - \Sigma_K^{y,T} \big\|
		 & \leq  \Big( \frac{\gamma_{pert}^T}{1- \gamma_{pert}} + \gamma_{pert}^T + (T-1) \gamma_{pert}^T \Big)  \frac{2d C_{init, noise}^2}{1- \gamma_{pert}}
		\\
		 & \leq  2d \gamma_{pert}^T (T +1) C_{init, noise}^2 / (1- \gamma_{pert})^2 .
	\end{align*}
	For $T \geq 2$ and $T \geq \phi_{truncation}(\varepsilon)$, we have
	$$
		T \log( 1/ \gamma_{pert}) - \log(T + 1) \geq \sqrt{T} \log(1/\gamma_{pert}) - 1 \geq \log \Big( \frac{2d C_{init, noise}^2}{\varepsilon ( 1- \gamma_{pert})^2 } \Big),
	$$
	which implies directly that $\| \Sigma_K^y - \Sigma_K^{y, T} \| \leq \varepsilon$.
	Similar arguments imply that $\| \Sigma_L^z - \Sigma_{L}^{z, T} \| \leq \varepsilon$.
	From equation~\eqref{eq:cost_expression_with_PK_y} and the definition of $C_y^T(K)$, $C_z^T(L)$ in~\eqref{eq:truncation_Cy_Cz_T}, we have
	\begin{align*}
		| C(\theta) - C^T(\theta) | & = \Big| \langle Q + K^\top R K,\, \Sigma_K^y - \Sigma_K^{y, T} \rangle_{tr} + \langle \tilde Q + L^\top \tilde R L,\, \Sigma_L^z - \Sigma_L^{z, T} \rangle_{tr} \Big|
		\\
		                            & \leq \varepsilon. h_T(C_0).
	\end{align*}
\qed\end{proof}

We define the truncated policy gradients in $K$ and in $L$ at time $T$ with perturbation size, perturbation radius, and perturbation directions $(M, \tau, \underline v)$ by
\begin{equation}
	\label{eq:def_{n}runcated_gradient_MKV}
	\hat{\nabla}_K^{T, M, \tau}(\theta, \underline v) = \frac{\ell d}{M \tau^2} \sum_{i=1}^M C^T(\theta + v_i) v_i^{(idy)},
	\quad
	\hat{\nabla}_L^{T, M, \tau}(\theta, \underline v) = \frac{\ell d}{M \tau^2} \sum_{i=1}^M C^T(\theta + v_i) v_i^{(com)}.
\end{equation}
Let $\hat{\nabla}^{T, M, \tau}(\theta, \underline v) = \big( \hat{\nabla}_K^{T, M, \tau}(\theta, \underline v), \,  \hat{\nabla}_L^{T, M, \tau}(\theta, \underline v)  \big)$.
Here, we assume that the perturbation radius $\tau$ is small enough so that $\theta + v_i \in \Theta$ for all $i=1, \ldots, M$, and thus the truncated policy gradient $\hat{\nabla}^{T, M, \tau}(\theta, \underline v)$ is well-defined.

To show the approximation of perturbed policy gradient $\hat{\nabla}^{M, \tau}$ with a truncated policy gradient $\hat{\nabla}^{T, M, \tau}$,  based on Lemma~\ref{lemma:approx_nruncated_cost}, we define
\begin{align}
	h_T(C_0)                                & = d.(\| Q \| + \| \tilde Q\| + \|R \| h_K(C_0)^2 + \| \tilde R \| h_L(C_0)^2 ),
	\label{eq:h_T}
	\\
	\phi_{trunc, T}(\varepsilon, \tau, C_0) & = \phi_{truncation}\big( \frac{\tau  \varepsilon }{2 \ell d. h_T(C_0)} \big)
	\nonumber
	\\
	                                        & = \Big( \frac{ \log( 4 \ell d^2 C_{init, noise}^2 h_T(C_0) ) - \log( \varepsilon \tau (1 - \gamma_{pert})^2 ) + 1 }{\log(1/\gamma_{pert})} \Big)^2.
	\label{eq:phi_{n}runc_T}
\end{align}
where $h_{K}(C_0)$ and $h_L(C_0)$ are defined in Lemma~\ref{lemma:useful_bounds_I}, and $\gamma_{pert}$, $\phi_{truncation}$ are defined in equations~\eqref{eq:gamma_pert} and~\eqref{eq:phi_{n}runcation}.

\begin{lemma}
	\label{lemma:approx_nruncated_policy_gradient_MKV_simulator}
	Consider $\theta \in \Theta$ with $C(\theta) \leq C_0$ for some constant $C_0 \in \RR$.
	For a given $\varepsilon > 0$, if the truncation horizon $T \geq \max \{2, \phi_{trunc, T}(\varepsilon, \tau, C_0) \}$, then we have,
	\begin{equation}
		\big\| \hat{\nabla}^{T, M, \tau}(\theta, \underline v) -  \hat{\nabla}^{M, \tau}(\theta, \underline v) \big\| \leq \varepsilon
	\end{equation}
	where $ \big\| \hat{\nabla}^{T, M, \tau}(\theta, \underline v) -  \hat{\nabla}^{M, \tau}(\theta, \underline v)  \big\| =  \big\| \hat{\nabla}_K^{T, M, \tau}(\theta, \underline v) -  \hat{\nabla}_K^{M, \tau}(\theta, \underline v)  \big\| +  \big\| \hat{\nabla}_L^{T, M, \tau}(\theta, \underline v) -  \hat{\nabla}_L^{M, \tau}(\theta, \underline v)  \big\|$.
\end{lemma}

\begin{proof}
	Let $\theta_i = \theta + v_i$ for $i=1,\ldots, M$.
	From the definitions~\eqref{eq:def_perturbed_PG_MKV} and~\eqref{eq:def_{n}runcated_gradient_MKV}, Lemma~\ref{lemma:approx_nruncated_cost} implies that
	$$
		\| \hat{\nabla}_K^{T, M, \tau}(\theta, \underline v) -  \hat{\nabla}_K^{M, \tau}(\theta, \underline v)  \big\|
		\leq  \frac{\ell d}{M \tau^2} \sum_{i=1}^M | C(\theta_i) - C^T(\theta_i) |. \| v_i^{(idy)} \| \leq  \frac{\varepsilon}{2}.
	$$
	We then conclude the result with similar inequalities for gradients in $L$.
\qed\end{proof}

\subsection{Sampled costs with MKV simulator}
\label{subsection:approx_sampled_costs_with_MKV_simulator}
In this section, we will first show a matrix concentration inequality (Lemma~\ref{lemma:concentration_ineq_for_yt_K}) for a quadratic form $(Y_{n}^K)^\top ( Q + K^\top R K) Y_{n}^K$ where $Y_{n}^K \in \RR^{d}$ is a term of the auxiliary process $\bY^{K}$ at time $n$, specifically a sub-Gaussian random vector with dependent coordinate entries.
These quadratic forms from the auxiliary process $\bY^K$ as well as those from $\bZ^L$ are closely related to the sample costs simulated from an MKV simulator $\cS_{MKV}^T$.
Based on this result, we will quantify the approximation error between the sampled policy gradients as outputs of Algorithm~\ref{algo:MKVestim} and the truncated policy gradient defined in equation~\eqref{eq:def_{n}runcated_gradient_MKV}.
To start, we define a constant
\begin{equation}
	\label{eq:h_psi_2, y}
	h_{\psi_2, y}(\tau, C_0) = \frac{2 \tau d^2 C_{init, noise}^2  \Big( \| Q \| + \|R\| ( \tau + h_K(C_0) )^2 \Big) }{(1 - \sqrt{\gamma})^2}.
\end{equation}
A similar constant $h_{\psi_2, z}(\tau, C_0)$ is defined with $(\tilde Q, \tilde R, h_L(C_0))$ in the above equation.
To abbreviate the notation, we also consider a real-value function $f$ on $\RR^{d} \times \RR^{\ell \times d} \times \RR^{d \times d} \times \RR^{\ell \times \ell}$ given by $f(\xi, \varphi, q, r) \mapsto \xi^\top ( q + \varphi^\top r \varphi ) \xi$.

\begin{lemma}
	\label{lemma:concentration_ineq_for_yt_K}
	Consider $\theta = (K, L) \in \Theta$ with $C(\theta) \leq C_0$ for some constant $C_0 \in \RR$.
	Let $\bY^K = (Y_{n}^K)_{n \geq 0}$ be the process following dynamics~\eqref{eq:dyn_y_theta}.
	For any $n \geq 0$, consider a quadratic form  of $Y_{n}^K$ given by $\zeta_{n}^{K} = f(Y_{n}^K, K, Q, R) = (Y_{n}^K)^\top (Q + K^\top R K ) Y_{n}^K$.
	Then, the random vector $Y_{n}^K$ is sub-Gaussian in $\RR^d$, and the real-valued random variable $\zeta_{n}^K$ is sub-exponential.
	Their norms are bounded by
	\begin{equation}
		\label{eq:bound_on_Y_{n}_K_and_quadratic_form}
		\| Y_{n}^K \|_{\psi_2} \leq \gamma^{-n/2} d C_{init, noise} / (1 - \sqrt{\gamma}),
		\qquad
		\| \zeta_{n}^K \|_{\psi_1} \leq Tr( Q + K^\top R K) \| Y_{n}^{K} \|_{\psi_2}^2.
	\end{equation}
	Moreover, for $M$ perturbed control parameters $K_i = K + v_i^{(idy)}$ with perturbation directions $(v_i^{(idy)})_{i=1, \ldots, M}$ on sphere $\SS_{\tau}$ independently sampled from $\mu_{\SS_\tau}$.
	Let $\zeta_{n}^{K_i} = f(Y_{n}^{K_i}, K_i, Q, R)$ for all $n \geq 0$ and $i = 1, \ldots, M$.
	Then we have, for any $\varepsilon > 0$,
	\begin{equation}
		\label{eq:concentration_ineq_for_quadratic_yt_K}
		\PP \Big( \Big\| \frac{\gamma^{n}}{M} \sum_{i=1}^M ( \zeta_{n}^{K_i} -  \EE[ \zeta_{n}^{K_i} ] ) v_i^{(idy)} \Big\| \geq \varepsilon \Big) \leq (\ell + d)  \exp \Big( \frac{-M \varepsilon^2 }{ 4 h_{\psi_2, K}^2(\tau, C_0) + \varepsilon h_{\psi_2, K}(\tau, C_0) } \Big)
	\end{equation}
	where the expectation $\EE$ and the probability $\PP$ are with respect to the randomness in $(Y_{n}^{K_i})_{i=1,\ldots, M}$, but not the randomness in perturbation directions $(v^{(idy)}_i)_{i=1, \ldots, M}$.
\end{lemma}

\begin{proof}
	We first show the two bounds on $Y_{n}^K$ and $\zeta_{n}^K$ in equation~\eqref{eq:bound_on_Y_{n}_K_and_quadratic_form}.
	From the dynamics~\eqref{eq:dyn_y_theta}, we have $Y_{n}^K = ({\mathrm A} - {\mathrm B} K)^{n} Y_{0} + \sum_{s=1}^{n} ( {\mathrm A} - {\mathrm B} K)^{n-s} \varepsilon_s$ where $Y_{0} = \varepsilon_0 - \EE[ \varepsilon_0 ]$.
	Because $(\varepsilon_s)_{s \geq 0}$ are all sub-Gaussian random variables, then $Y_{n}^K$ is also sub-Gaussian.
	Let $\xi_0 = Y_{0}$, $\xi_{s+1} = \varepsilon_{s+1}$ and $H_s = ({\mathrm A} - {\mathrm B} K)^{n-s}$ for $s=0, \ldots, n$.
	Let $\xi_{s,j} \in \RR$ and $H_{s, j} \in \RR^d$ be the $j$-th entry of $\xi_s$ and the $j-$th column of $H_s$, for all $j=1, \ldots, d$.
	Then $\| \xi_{s, j} \|_{\psi_2} \leq \| \xi_s \|_{\psi_2} \leq C_{init, noise} $, $\| H_{s, j} \| \leq \| H_s \|$, and
	$
		Y_{n}^K = \sum_{s=0}^{n} H_s \xi_s.
	$
	We deduce that
	\begin{align*}
		\gamma^{n} \| Y_{n}^{K}\|_{\psi_2}^2
		\leq \gamma^{n} \Big( \sum_{s=0}^{n} \sum_{j=1}^d \| H_{s, j} \xi_{s, j} \|_{\psi_2} \Big)^2
		 & \leq \gamma^{n} \Big( \sum_{s=0}^{n}  \sum_{j=1}^d \| H_{s, j} \|_2. \| \xi_{s, j} \|_{\psi_2} \Big)^2
		\\
		 & \leq   d^2 C_{init, noise}^2 \Big( \sum_{s=0}^{n} \gamma^{s/2} ( \gamma^{(n-s)/2} \| H_{s} \| ) \Big)^2
		\\
		 & \leq d^2 C_{init, noise}^2 / (1 - \sqrt{\gamma})^2.
	\end{align*}
	Furthermore, from Lemma~\ref{lemma:subexponential_quadratic_form_bound}, we know that the quadratic form $\zeta_{n}^K = f(Y_{n}^K, K, Q, R)$ is a sub-exponential random variable and $\| \zeta_{n}^K \|_{\psi_1} \leq Tr(Q + K^\top R K) \| Y_{n}^K \|_{\psi_2}^2$.

	Now, we will show the concentration inequality~\eqref{eq:concentration_ineq_for_quadratic_yt_K} with small perturbation directions $(v_i^{(idy)})_{i=1}^M$.
	Let $(\Omega^{(i)}, \cF^{(i)}, \PP^{(i)})_{i=1}^M$ be M independent copies of the probability space $(\Omega, \cF, \PP)$, and let $(\tilde \Omega, \tilde \cF, \tilde \PP)$ with $\tilde \Omega = \Omega^{(1)} \times \ldots \times \Omega^{(M)}$, $\tilde \cF = \cF^{(1)} \otimes \ldots \otimes \cF^{(M)}$, and $\tilde \PP = \PP^{(1)} \otimes \ldots \otimes \PP^{(M)}$ be the extended product probability space.
	Consider a small enough perturbation radius $\tau$, so that the perturbed control parameter $K_i = K + v_i^{(idy)}$ satisfies $\gamma \| {\mathrm A} - {\mathrm B} K_i \|^2 < 1$ for all $i = 1, \ldots, M$.
	For each $K_i$, we define an auxiliary process $\bY^{K_i} = (Y_{n}^{K_i})_{n \geq 0}$ following dynamics~\eqref{eq:dyn_y_theta}, but with an independent copy $\boldsymbol{\varepsilon^{i}}$ of the noise process $\boldsymbol{\varepsilon}$ and an independent copy of the initial status $Y_{0}^{K_i} \sim Y_{0}$.
	Consequently, for each time $n \geq 0$, $(\bY^{K_1}_{n}, \ldots, \bY^{K_M}_{n})$ are $M$ i.i.d. random vectors in $\RR^{d}$.
	For a chosen time $n \geq 0$ and $i=1, \ldots, M$, consider the quadratic form of $Y_{n}^{K_i}$ given by $\zeta_{n}^{K_i} = f(Y_{n}^{K_i}, K_i, Q, R)$ and its unbiased term $W_{n}^i = \zeta_{n}^{K_i} - \EE[ \zeta_{n}^{K_i}]$.

	Let $\tilde v_i^{(idy)} \in \RR^{ (\ell + d) \times (\ell + d) }$ be the dilatation matrix of $v_i^{(idy)} \in \RR^{\ell \times d}$. The dilatation trick (detailed in Lemma~\ref{lemma:dilatation_matrix}) maps an asymmetric random matrix into a symmetric space where the Matrix Bernstein inequality can be applied.
	For sub-exponential random variables, we have
	$
		\| \zeta_{n}^{K_i} \|_{L^p}^p \leq 2 p!  \| \zeta_{n}^{K_i} \|_{\psi_1}^p
	$
	so that
	$$
		\EE\big[ ( W_{n}^i \tilde{v}_i^{(idy)})^p \big] = \tau^p \big\| W_{n}^i \big\|_{L^p}^p \big( \tilde{v}_i^{(idy)} / \tau \big)^p \preceq \tau^p 2^p ( 2 p! \| \zeta_{n}^{K_i} \|_{\psi_1}^p ) \mathbf{I}_{(\ell + d)}
	$$
	where $\mathbf{I}_{(\ell + d)}$ is an identity matrix on $\RR^{(\ell + d) \times (\ell + d)}$.
	Let us define two time-dependent coefficients:
	$$
		r_{n} = \gamma^{-n} h_{\psi_2}(\tau, C_0) ,
		\quad
		\sigma_{n} = 2  \gamma^{-n} \sqrt{M} h_{\psi_2}(\tau, C_0).
	$$
	We observe that for all $i=1, \ldots, M$, 
    $2 \tau \| \zeta_{n}^{K_i} \|_{\psi_1} \leq \gamma^{-n} h_{\psi_2}(\tau, C_0)$, so we have $\EE\big[ ( W_{n}^i \tilde{v}_i^{(idy)})^p \big] \preceq ( p! / 2 ) r_{n}^{p-2} \big( 4 \tau \| \zeta_{n}^{K_i} \|_{\psi_1} \mathbf{I}_{(\ell + d)} \big)^2$, and $\| \sum_{i=1}^M \big( 4 \tau \| \zeta_{n}^{K_i} \|_{\psi_1} \mathbf{I}_{(\ell + d)} \big)^2 \| \leq \sigma_{n}^2$.
	Thus, by the Matrix Bernstein inequality for sub-exponential case~\cite[Theorem 6.2]{MR2946459}, we have
	\begin{align*}
		\PP \Big( \Big\| \frac{\gamma^{n}}{M} \sum_{i=1}^M ( \zeta_{n}^{K_i} -  \EE[ \zeta_{n}^{K_i} ] ) v_i^{(idy)} \Big\| \geq \varepsilon \Big) & = \PP \Big( \big\| \sum_{i=1}^M W_{n}^i \tilde{v}_i^{(idy)} \big\| \geq \frac{\varepsilon M}{\gamma^{n}} \Big)
		\\
		                                                                                                                                       & \leq (\ell + d) \exp \Big( \frac{- (\varepsilon M \gamma^{-n} )^2 }{ \sigma_{n}^2 + r_{n} (\varepsilon M \gamma^{-n}) }\Big).
	\end{align*}
	Replacing the coefficient $r_{n}$ and $\sigma_{n}$ with $h_{\psi_2, y}(\tau, C_0)$, we obtain the concentration inequality~\eqref{eq:concentration_ineq_for_quadratic_yt_K}.
\qed\end{proof}

\begin{remark}
	It should be noted that the concentration inequality~\eqref{eq:concentration_ineq_for_quadratic_yt_K} is related to a quadratic form of a sub-Gaussian random vector with dependent coordinate entries, so we cannot apply the Hanson-Wright concentration inequality~\cite{vershynin2018high} or the classical Bernstein's inequality with bounded variables as considered in~\cite{fazel2018global}.
	Lemma~\ref{lemma:concentration_ineq_for_yt_K} constructs a tail bound with the help of the discount coefficient $\gamma$.
\end{remark}

Now, we show the approximation between the sampled policy gradients from Algorithm~\ref{algo:MKVestim} and truncated policy gradients defined in~\eqref{eq:def_{n}runcated_gradient_MKV}.
We recall that the sampled policy gradients at $\theta = (K, L)$ with $M$ perturbation directions $\underline v = ( (v_i^{(idy)}, v_i^{(com)} ) )_{i=1}^M$ on $\SS_{\tau} \times \SS_{\tau}$ with radius $\tau > 0$ and from an MKV simulator $\cS_{MKV}^{T}$ is given by
\begin{equation}
	\label{eq:def_sampeld_policy_gradient_MKV_restate}
	\tilde{\nabla}^{T, M,\tau}_K (\theta, \underline v) = \frac{\ell d}{\tau^2}\frac{1}{M} \sum_{i=1}^M \tilde{C}^T(\theta_i) v^{(idy)}_i,
	\quad
	\tilde{\nabla}^{T, M,\tau}_L (\theta, \underline v) = \frac{\ell d}{\tau^2}\frac{1}{M} \sum_{i=1}^M \tilde{C}^T(\theta_i) v^{(com)}_i,
\end{equation}
where $\theta_i = (K + v_i^{(idy)}, L + v_i^{(com)})$ for $i=1, \ldots, M$.
We also define the following constants to bound the perturbation size $M$ in the sampled costs:
\begin{align}
	\phi_{sample,size}(\varepsilon, \tau, T, C_0, \delta_{sample}) = \log \Big( \frac{4T (\ell + d)}{\delta_{sample}} \Big) \Big( \frac{8 \ell d T}{\varepsilon \tau^2 }   h_{\psi_2}(\tau, C_0)  + \frac{1}{4} \Big)^2
	\label{eq:phi_sample_size}
\end{align}
where $h_{\psi_2}(\tau, C_0) = \max \big\{ h_{\psi_2, y}(\tau, C_0),\, h_{\psi_2, z}(\tau, C_0) \big\}$ with $h_{\psi_2, y}(\tau, C_0)$ and $h_{\psi_2, z}(\tau, C_0)$ defined in~\eqref{eq:h_psi_2, y}.

\begin{remark}
	We note that the sample complexity bounds, such as $\phi_{sample,size}$ above and the truncation horizon $\phi_{truncation}$, are highly sensitive to the discount factor $\gamma$. In particular, as $\gamma \to 1$, the terms $(1 - \sqrt{\gamma})^2$ and $(1 - \gamma_{pert})^2$ in the denominators cause these bounds to grow significantly. This reflects a known challenge in infinite-horizon reinforcement learning, where the effective horizon diverges as the discount factor approaches $1$.
\end{remark}

\begin{lemma}
	\label{lemma:approx_sampled_pg_with_{n}runcated_pg}
	Consider $\theta \in \Theta$ with $C(\theta) \leq C_0$ for some $C_0 \in \RR$.
	If the perturbation size $M$ satisfies
	$$
		M \geq \phi_{sample, size}(\varepsilon, \tau, T, C_0, \delta_{sample}),
	$$
	then for a given sequence of $2M$ perturbation directions $\underline v$ independently and uniformly sampled from $\SS_\tau$, we have, for any $\varepsilon > 0$,
	\begin{equation}
		\label{eq:approx_sampled_pg_with_{n}runcated_pg}
		\PP \big( \| \tilde \nabla^{T, M, \tau}(\theta, \underline v) - \hat{\nabla}^{T, M, \tau}(\theta, \underline v) \| \geq \varepsilon \big) \leq \delta_{sample}.
	\end{equation}
\end{lemma}

\begin{proof}
	Let $\theta_i = \theta + v_i = ( K + v_i^{(idy)}, L + v_i^{(com)}) = (K_i, L_i)$ with perturbation direction $v_i = (v_i^{(idy)}, v_i^{(com)})$, for $i=1, \ldots, M$.
	Consider the auxiliary processes $(\bY^{K_i}, \bZ^{L_i})$ and for every $n\geq 0$ we define $\zeta_{n}^{K_i} = f(Y_{n}^{K_i}, K_i, Q, R)$ and $\tilde \zeta_{n}^{L_i} = f(Z_{n}^{L_i}, L_i, \tilde Q, \tilde R)$.
	Then
	\begin{align*}
		     & \|  \tilde \nabla^{T, M, \tau}_K(\theta, \underline v) - \hat{\nabla}^{T, M, \tau}_K(\theta, \underline v)  \|
		\\
		=    & \Big\| \frac{\ell d}{M \tau^2} \sum_{i=1}^M \sum_{n=0}^{T-1} \gamma^{n} \big[ \big( \zeta_{n}^{K_i} - \EE[\zeta_{n}^{K_i}] \big) + \big( \tilde{\zeta}_{n}^{L_i} - \EE[ \tilde{\zeta}_{n}^{L_i} ] \big) \big] v_i^{(idy)} \Big\|
		\\
		\leq & \frac{\ell d}{\tau^2} \sum_{n=0}^{T-1} \Big( \Big\| \frac{\gamma^{n}}{M} \sum_{i=1}^M  \big( \zeta_{n}^{K_i} - \EE[\zeta_{n}^{K_i}] \big) v_i^{(idy)} \Big\| + \Big\| \frac{\gamma^{n}}{M} \sum_{i=1}^M  \big( \tilde{\zeta}_{n}^{L_i} - \EE[ \tilde{\zeta}_{n}^{L_i} ] \big) v_i^{(idy)} \Big\| \Big)
		\\
		\leq & \frac{\ell d}{\tau^2}. 2T \tilde{\varepsilon}
	\end{align*}
	with probability at least $1 - 2T (\ell + d) \exp\big( \frac{- M \tilde{\varepsilon}^2 }{4 h_{\psi_2}^2(\tau, C_0) + \tilde{\varepsilon} h_{\psi_2}(\tau, C_0) } \big)$.
	Similar inequality holds for $\|  \tilde \nabla^{T, M, \tau}_L (\theta, \underline v) - \hat{\nabla}^{T, M, \tau}_L(\theta, \underline v)  \| $.
	Then, by choosing $\tilde{\varepsilon}$ and $M$ satisfy that
	$$
		\Big( \frac{4 \ell d T}{\tau^2} \Big) \tilde{\varepsilon} = \varepsilon,
		\quad
		4 T (\ell + d) \exp \Big( \frac{-M \tilde{\varepsilon}^2 }{4 h_{\psi_2}^2(\tau, C_0)  + \tilde{\varepsilon} h_{\psi_2}(\tau, C_0) } \Big) \leq \delta_{sample},
	$$
	or letting $M \geq \phi_{sample, size}(\varepsilon, \tau, T, C_0, \delta_{sample})$, we conclude that with probability at least $1 - \delta_{sample}$, $\| \tilde \nabla^{T, M, \tau}(\theta, \underline v) - \hat{\nabla}^{T, M, \tau}(\theta, \underline v) \| \leq \varepsilon$.
\qed\end{proof}

\subsection{Proof of Proposition~\ref{pr:gradient_approx_with_MKV_simulator}}
\label{subsection:proof_of_proposition_approx_gradient_w_MKV_simulator}
In this subsection, we show that when the truncation horizon $T$ is large enough for an MKV simulator $\cS_{MKV}^T$, the perturbation radius $\tau$ is small enough, and the perturbation size $M$ is large enough in Algorithm~\ref{algo:MKVestim}, the output of the approximated policy gradient $\tilde \nabla^{T, M, \tau}$ from the algorithm is close to the exact PG $\nabla C$ with high probability.
To do so, we define
\begin{align}
	\phi_{pert, radius, MKV}( \varepsilon, C_0)
	 & = \phi_{pert, radius}( \varepsilon / 3, C_0)
	\label{eq:phi_MKV_pert_radius}
	\\
	\phi_{trunc, T, MKV}(\varepsilon, \tau,  C_0)
	 & = \phi_{trunc, T}( \varepsilon / 3, \tau, C_0) + 2
	\label{eq:phi_MKV_{n}runc_T}
	\\
	\phi_{sample, size, MKV}(\varepsilon, \tau, T, C_0, \delta_{approx})
	 & = \phi_{pert, size} \big(\varepsilon / 3, \tau, C_0, \delta_{approx} / 2 \big)
	\nonumber                                                                                           \\
	 & \hspace{10pt} + \phi_{sample, size} \big(\varepsilon / 3, \tau, T, C_0, \delta_{approx} / 2\big)
	\label{eq:phi_MKV_sample_size}
\end{align}
where $\phi_{pert, radius}$~\eqref{eq:phi_pert_radius}, $\phi_{pert, size}$~\eqref{eq:phi_pert_size}, $\phi_{trunc, T}$~\eqref{eq:phi_{n}runc_T}, $\phi_{sample, size}$~\eqref{eq:phi_sample_size} are polynomial in $(d, \ell, C_0, (\lambda_y^1)^{-1}, (\lambda_z^0)^{-1}, C_{init, noise}, \varepsilon^{-1}, \tau^{-1})$ and model parameters defined in sections~\ref{subsection:approximation_with_smoothed_cost},~\ref{subsection:approximation_C_with_C_T}, and~\ref{subsection:approx_sampled_costs_with_MKV_simulator}.\\

Let $\varepsilon > 0$ be a target precision.
We assume that

\begin{proof}
	From Lemma~\ref{lemma:approx_perturbed_gradient_MKV_simulator}, Lemma~\ref{lemma:approx_nruncated_policy_gradient_MKV_simulator}, and Lemma~\ref{lemma:approx_sampled_pg_with_{n}runcated_pg}, we have
	\begin{align*}
		     & \PP \big( \big\| \tilde \nabla^{T, M, \tau}(\theta, \underline v) - \nabla C(\theta) \big\|  > \varepsilon \big)
		\\
		\leq & \PP \big( \big\| \hat{\nabla}^{M, \tau}(\theta, \underline v) - \nabla C(\theta) \big\| + \big\| \hat{\nabla}^{T, M, \tau}(\theta, \underline v) - \hat{\nabla}^{M, \tau}(\theta, \underline v) \big\|
		\\
		     & \hspace{10pt} + \big\| \tilde{\nabla}^{T,M, \tau}(\theta, \underline v) - \hat{\nabla}^{T, M, \tau}(\theta, \underline v) \big\| > \varepsilon \big)
		\\
		\leq & \PP \big(  \big\| \hat{\nabla}^{M, \tau}(\theta, \underline v) - \nabla C(\theta) \big\| > \varepsilon / 3 \big) + \PP\big( \big\| \hat{\nabla}^{T, M, \tau}(\theta, \underline v) - \hat{\nabla}^{M, \tau}(\theta, \underline v) \big\|  > \varepsilon / 3 \big)
		\\
		     & \hspace{10pt} + \PP \big( \big\| \tilde{\nabla}^{T,M, \tau}(\theta, \underline v) - \hat{\nabla}^{T, M, \tau}(\theta, \underline v) \big\| > \varepsilon / 3 \big)
		\\
		\leq & \delta_{approx} / 2 + 0 + \delta_{approx} / 2
		\\
		=    & \delta_{approx}.
	\end{align*}
    \,
\qed\end{proof}

\subsection{Additional lemmas for sub-Gaussian random vectors}
\label{subsection:additional_lemma_subgaussian}

The following lemma is similar to~\cite[Proposition 2.4]{zajkowski2020bounds}.
\begin{lemma}
	\label{lemma:subexponential_quadratic_form_bound}
	Consider a sub-Gaussian random vector $\xi \in \RR^{d}$ and a PSD matrix $X \succeq 0$ in $\RR^{d \times d}$. The quadratic form $\xi^\top X \xi$ is a sub-exponential random variable, and
	\begin{equation}
		\| \xi^\top X \xi \|_{\psi_1} \leq Tr(X) \| \xi \|_{\psi_2}^2.
	\end{equation}
\end{lemma}
\begin{proof}
	Because $X \succeq 0$, let $X = U^\top \Sigma U$ be the eigenvalue decomposition of $X$ with $\Sigma = \diago(\lambda_{1}, \lambda_{2}, \ldots, \lambda_{d})$, where $\lambda_{i=1,\ldots,d}$ are the eigenvalues of $X$.
	Let $\zeta = U \xi = [\zeta_1, \ldots, \zeta_d]^\top$, then the $d$ entries $\zeta_{i=1,\ldots d}$ are sub-Gaussian random variables, and $\zeta$ is also a sub-Gaussian random vector such that $\| \zeta \|_{\psi_2} = \sup_{ s \in \mathbb{S}^{d-1}}\| \langle U \xi, s\rangle\|_{\psi_2} = \| \xi \|_{\psi_2}$.
	Consequently, we have $\xi^\top X \xi = \zeta^\top \Sigma \zeta =  \sum_{i=1}^d \lambda_i \zeta_i^2 $, and~\cite[Lemma 2.7.7]{vershynin2018high} implies that $\xi^\top X \xi$ is sub-exponential.
	Moreover,
	\begin{align*}
		\| \xi^\top X \xi \|_{\psi_1} = \Big\| \sum_{i=1}^d \lambda_i \zeta_i^2 \Big\|_{\psi_1} \leq \sum_{i=1}^d\lambda_i  \| \zeta_i^2 \|_{\psi_1} = \sum_{i=1}^d \lambda_i  \| \zeta_i \|_{\psi_2}^2 \leq \big( \sum_{i=1}^d \lambda_i \big) \| \zeta \|_{\psi_2}^2,
	\end{align*}
	where the second equality is justified by the definitions of $\| \cdot \|_{\psi_1}$ and $\| \cdot \|_{\psi_2}$ norms (see~\cite[Lemma 2.7.6]{vershynin2018high}), and the last inequality is due to $\| \zeta_i \|_{\psi_2} = \| \langle \zeta, e_i \rangle \|_{\psi_2} \leq \| \zeta \|_{\psi_2}$ with $e_i \in \mathbb{S}^{d-1}$ being the base vector for dimension $i$.
	We conclude the lemma with $Tr(X) = \sum_{i=1}^d \lambda_{i}$.
\qed\end{proof}

The following lemma can be found in~\cite[Proposition 2.7.1]{vershynin2018high} to connect the $L^p$ norms with the sub-exponential norm and the sub-Gaussian norm.
\begin{lemma}
	\label{lemma:relationship_between_Lp_and_subsexp}
	For a sub-Gaussian random vector $\xi \in \RR^{d}$ and a sub-exponential random variable $\zeta \in \RR$, we have
	$$
		\| \xi \|_{L^2(\RR^d)} \leq \sqrt{2 d} \| \xi \|_{\psi_2}, \qquad \| \zeta \|_{L^2(\RR)} \leq 2 \| \zeta \|_{\psi_1}
	$$
	Moreover, for any $p \geq 2$, we have
	$$
		\| \zeta \|_{L^p(\RR)} \leq (2 p !)^{1/p} \| \zeta \|_{\psi_1} \leq p \| \zeta \|_{\psi_1}.
	$$
\end{lemma}

The following lemma is about the dilatation technique used for matrix concentration inequalities with non-symmetric matrices; see~\cite{MR2946459}.
\begin{lemma}
	\label{lemma:dilatation_matrix}
	For a matrix $v \in \RR^{\ell \times d}$, its dilatation matrix $\tilde v$ is defined as
	$$
		\tilde{v} = \left[ \begin{array}{cc} 0 & v \\ v^\top & 0 \end{array} \right] \in \RR^{(\ell + d) \times (\ell + d)}.
	$$
	We have the following properties for the dilatation matrix: $\tilde v$ is symmetric and $ \| \tilde v \| = \| v \| $.
\end{lemma}

\section{Technical proofs for Section~\ref{subsec:PG-popsimu}}
\label{sec:app-proof-modelfree-POP-CV}

\subsection{Auxiliary costs for \texorpdfstring{$N$}{N}-agent problem}

To show the model-free population-based gradient estimator $\tilde{\nabla}^{T,N,M, \tau}$ from Algorithm~\ref{algo:POPestim} provides a good approximation for the exact PG $\nabla C$, similar to the approach in Section~\ref{sec:app-proof-modelfree-MKV-CV} for the MKV simulator, we define the \emph{perturbed population policy gradient} $\hat{\nabla}^{N, M, \tau}$ and the \emph{truncated population policy gradient} $\hat{\nabla}^{T, N, M, \tau}$ by replacing the perturbed MF cost $C(\theta + v_i)$ in~\eqref{eq:def_perturbed_PG_MKV} and the truncated perturbed cost $C^T(\theta + v_i)$ in~\eqref{eq:def_{n}runcated_gradient_MKV} with the corresponding $N$-agent social costs adapted from equation~\eqref{fo:social_cost_of_population}.

Let $[N] = \{1, \ldots, N\}$.
For $\theta = (K, L) \in \Theta$, we introduce auxiliary processes $(\bY^{K, j, N})_{j \in [N]}$ and $\bZ^{L, N}$ for the $N$-agents by $ Y_{n}^{K, j, N} = X_{n}^{(j), \theta} - \bar X_{n}^{N, \theta}$ and $Z_{n}^{L, N} = \bar X_{n}^{N, \theta} = \frac{1}{N} \sum_{j=1}^N X_{n}^{(j), \theta}$ for all $n \geq 0$ and $j \in [N]$, where $(X_{n}^{(j), \theta})_{n\geq 0}$ is the state process for agent $j$ following dynamics~\eqref{fo:N-multi_state} with control at time $n$ given by $\ctrl_{n}^{(j), \theta} = -  K (X_{n}^{(j), \theta} - \bar{X}_{n}^{N, \theta}) - L \bar{X}_{n}^{N, \theta}$.
The dynamics of $(\bY^{K, j, N})_{j \in [N]}$ and $\bZ^{L, N}$ then follow
\begin{align}
	 & Y_{n+1}^{K, j, N} = ( {\mathrm A} - {\mathrm B} K ) Y_{n}^{K, j, N} + \tilde{\varepsilon}_{n+1}^{(j)}, \quad                                     &  & Y_{0}^{j, N} = \varepsilon_0^{(j)} - \tfrac{1}{N}\sum \nolimits_{i \in [N]} \varepsilon_0^{(i)},
	\label{eq:def_dyn_y_K_n_N}
	\\
	 & Z_{n+1}^{L, N}  = ( \tilde {\mathrm A} - \tilde {\mathrm B} L ) Z_{n}^{L, N} + \varepsilon_{n+1}^0 + \overline{{\varepsilon}_{n+1}^{N}}, \quad &  & Z_{0}^{N} = \varepsilon_0^0 +  \tfrac{1}{N}\sum \nolimits_{i \in [N]} \varepsilon_0^{(i)},
	\label{eq:def_z_L_N}
\end{align}
where $\overline{{\varepsilon}_{n+1}^{N}} = \frac{1}{N} \sum_{j=1}^N \varepsilon_{n+1}^{(j)}$
and $\tilde{\varepsilon}_{n+1}^{(j)} = \varepsilon_{n+1}^{(j)} - \overline{{\varepsilon}_{n+1}^{N}}$ for $n \geq 0$.
We define the variance matrices $\Sigma_K^{y, j, N}$ and $\Sigma_L^{z, N}$ by replacing $Y_{n}$ and $Z_{n}$ in~\eqref{eq:variance_matrices_y_and_z} with $Y_{n}^{K, j, N}$ and $Z_{n}^{L, N}$, and define $\Sigma_K^{y, N} = \frac{1}{N} \sum_{j=1}^N  \Sigma_K^{y, j, N}$.
Similar modifications are applied to the initial perturbation variances $\Sigma_{Y_{0}}^N = \frac{1}{N} \sum_{j=1}^N  \EE[ Y_{0}^{j, N} (Y_{0}^{j,N})^\top ] $ and $\Sigma_{Z_{0}}^N = \EE[ Z_{0}^N (Z_{0}^N)^\top ]$, and for the noise variances
$
	\Sigma^{1, N} = \frac{1}{N} \sum_{j=1}^N  \EE\big[ \tilde{\varepsilon}_{n}^{(j)} \big( \tilde{\varepsilon}_{n}^{(j)} \big)^\top \big]$ and $\Sigma^{0, N} = \EE\big[ \big(  \varepsilon_{n}^0 + \overline{{\varepsilon}_{n}^{N}} \big) \big(  \varepsilon_{n}^0 + \overline{{\varepsilon}_{n}^{N}} \big)^\top \big] $.
For the auxiliary costs in the $N-$agent problem, we adjust equation~\eqref{eq:def_Cy_Cz} with $(\bY^{K, j, N})_{j\in [N]}$ and $\bZ^{L, N}$ so that
\begin{align*}
	C_y^{N}(K) 
    &= \frac{1}{N} \sum \nolimits_{j \in [N]} \EE \big[ \sum \nolimits_{n\geq 0} \gamma^{n} (Y_{n}^{K, j, N})^\top (Q + K^\top R K ) Y_{n}^{K, j, N} \big] 
    = \langle Q + K^\top R K, \Sigma_{K}^{y, N} \rangle_{tr}
	\\
	C_z^{N}(L) 
    &= \EE \big[ \sum \nolimits_{n\geq 0} \gamma^{n} (Z_{n}^{L, N})^\top ( \tilde{Q} + L^\top \tilde{R} L ) Z_{n}^{L, N} \big]
    = \langle \tilde{Q} + L^\top \tilde{R} L, \Sigma_{L}^{z, N} \rangle_{tr}.
\end{align*}
For $\theta = (K, L) \in \Theta$, $C^N(\theta) = C_y^N(K) + C_z^N(L)$.
Because all agents adopt the same control parameter $\theta$, the $N$-agent auxiliary costs $C_y^N(K)$ and $C_z^{N}(L)$ can then be expressed with solution matrices $P_K^y$ and $P_L^z$ to the DLEs~\eqref{eq:lyapunov_eq_theta}:
\begin{equation}
	\label{eq:N_agent_Cy_Cz_expression_PK_PL}
	C_y^{N}(K) = \langle P_K^y, \, \Sigma_{Y_{0}}^N + \frac{\gamma}{1 - \gamma} \Sigma^{1, N} \rangle_{tr},
	\quad
	C_z^N(L) = \langle P_L^z, \, \Sigma_{Z_{0}}^N + \frac{\gamma}{1 - \gamma} \Sigma^{0, N} \rangle_{tr}.
\end{equation}
The truncated versions at horizon $T$ for the variances and the costs $\Sigma_K^{y, N, T}$, $\Sigma_{L}^{z, N, T}$, $C^{N, T}_y(K)$, $C^{N, T}_z(L)$, $C^{N, T}(\theta)$, and the sampled truncated costs with population simulator $\tilde C_y^{N, T}(K),$ $\tilde C_z^{N, T}(L),$ $\tilde C^{N, T}(\theta)$ are defined accordingly.

By definition of the matrices $\Sigma^{1,N}$, $\Sigma^{0, N}$, $\Sigma_{Y_{0}}^{N}$ and $\Sigma_{Z_{0}}^N$, we have
\begin{align}
	\Sigma^{1, N}    & = (1 - \tfrac{1}{N} ) \Sigma^1,    & \Sigma^{0, N}    & = \Sigma^0 + \tfrac{1}{N}\Sigma^1,
	\label{eq:variance_noise_N_agent}
	\\
	\Sigma_{Y_{0}}^{N} & = ( 1 - \tfrac{1}{N})\Sigma_{Y_{0}}, & \Sigma_{Z_{0}}^{N} & =  \Sigma_{Z_{0}} + \tfrac{1}{N}\Sigma_{Y_{0}}.
	\label{eq:variance_init_N_agent}
\end{align}
Moreover, because the sub-Gaussian norms of $Y_{0}$, $Z_{0}$ and the noise processes are bounded by $C_{init,noise}$, we have for all $n \geq 0$ and $j \in [N]$,
\begin{equation}
	\label{eq:bound_phi2_y0_z0_N_agent}
	\max \Big\{ \| Y_{0}^{j, N} \|_{\psi_2}, \| Z_{0}^{N} \|_{\psi_2},
	\| \tilde{\varepsilon}_{n+1}^{(j)} \|_{\psi_2}, \big\| \varepsilon_{n+1}^0 + \overline{{\varepsilon}_{n+1}^{N}} \big\|_{\psi_2}  \Big\}  \leq 2 C_{init, noise}
\end{equation}
We define the minimum eigenvalues associated with the initial and noise variance matrices:
\begin{equation}
	\lambda_y^{1,N} = \lambda_{min}\big( \Sigma_{Y_{0}}^N + \gamma \Sigma^{1, N } / (1 - \gamma) \big), \quad \lambda_z^{0,N} = \lambda_{min}\big( \Sigma_{Z_{0}}^N + \gamma \Sigma^{0, N} / ( 1- \gamma) \big).
\end{equation}

\subsection{Proof of Proposition~\ref{pr:approx_with_social_cost_no_MF_cost} and~\ref{pr:modelfree_population_gradient_approx}}
\label{subsection:proof_of_approx_modelfree_pop_gradient}

\noindent \textbf{Proof of Proposition~\ref{pr:approx_with_social_cost_no_MF_cost}}
\begin{proof}
	From the expressions~\eqref{eq:N_agent_Cy_Cz_expression_PK_PL} for auxiliary costs $C_y^N$ and $C_z^N$ for the $N$-agent problem, and the expression~\eqref{eq:cost_expression_with_PK_y} for $C_y$ and $C_z$ for the MF problem, we have
	\begin{align*}
		| C^N(\theta) - C(\theta) |
		 & = \Big| \langle P_L^z - P_K^y, \, \frac{1}{N} \big( \Sigma_{Y_{0}} + \frac{\gamma}{1- \gamma} \Sigma^1 \big) \rangle_{tr} \Big|
		\\
		 & \leq \frac{1}{N} ( \|P_K^y\| + \|P_L^z \|) \frac{ \max \big\{ Tr(\Sigma_{Y_{0}}), Tr(\Sigma^1) \big\} }{1-\gamma}
	\end{align*}
	Together with $Tr(\Sigma_{Y_{0}}) = \EE[ \| Y_{0} \|^2 ] \leq 2d C_{init,noise}^2$ and $Tr(\Sigma^1) \leq 2 d C_{init, noise}^2$, and the bounds~\eqref{eq:upper_bound_Riccati_and_Variance_matrix}, we conclude the result by setting
	$
		\phi_{social-cost, factor}(C_0) = \frac{2 d C_{init,noise}^2 C_0}{(1 - \gamma)} \Big( \frac{1}{\lambda_y^1} + \frac{1}{\lambda_z^0} \Big).
	$
\qed\end{proof}

We show in the following that the output of Algorithm~\ref{algo:POPestim}, $\tilde \nabla^{T, N, M, \tau}$, is close enough to the exact PG $\nabla C$.
We use the same arguments as in Section~\ref{sec:app-proof-modelfree-MKV-CV}.
In the following, we specify only the differences in the corresponding proof arguments when we use a social cost instead of a mean field cost. To begin with, we define a few bounds for the parameters $(T, N, M)$.
\begin{align}
	 & \phi_{agent, size, pop}(\varepsilon, \tau, C_0) =  \big( \ell d / (\varepsilon \tau) \big) \phi_{social-cost, factor}(C_0)
	\nonumber                                                                                                                                                                              \\
	 & = \frac{2 \ell d^2 C_{init, noise}^2 C_0}{\tau \varepsilon ( 1- \gamma)}  \Big( \frac{1}{\lambda_y^1} + \frac{1}{\lambda_z^0} \Big)
	\label{eq:phi_agent_size_N}
	\\
	 & \phi_{sample,size, pop}(\varepsilon, \tau, T, C_0, \delta_{sample}) = \phi_{sample, size}(\varepsilon / 4, \tau, T, C_0, \delta_{sample})
	\nonumber                                                                                                                                                                              \\
	 & =  \log \Big( \frac{4T (\ell + d)}{\delta_{sample}} \Big) \Big(  \frac{32 \ell d T}{\varepsilon \tau^2}  (  h_{\psi_2}(\tau, C_0) )  + \frac{1}{4} \Big)^2
	\label{eq:phi_sample_size_N-agent}
	\\
	 & \phi_{trunc, T, pop}(\varepsilon, \tau, C_0, N) = \phi_{trunc, T}(\varepsilon / (1 + 1/N), \tau, C_0)
	\nonumber                                                                                                                                                                              \\
	 & = \Big( \frac{1}{\log(1 / \gamma_{pert})} \Big)^2 \Big[ \log \Big( \frac{4 \ell d C_{init, noise}^2 h_T(C_0)}{\varepsilon \tau ( 1- \gamma_{pert})^2} (1 + \tfrac{1}{N}) \Big) + 1 \Big]^2
	\label{eq:phi_{n}runc_T_N_agent}
\end{align}
where $\phi_{social-cost, factor}(C_0)$ is defined in Proposition~\ref{pr:approx_with_social_cost_no_MF_cost} , and $\phi_{sample, size}$, $\phi_{trunc, T}$ are defined in equations~\eqref{eq:phi_sample_size} and~\eqref{eq:phi_{n}runc_T}.

\begin{lemma}
	\label{lemma:approx_gradient_pop_simulator}
	Consider $\theta \in \Theta$ such that $C(\theta) \leq C_0$ for some $C_0 \in \RR$.
	For a target precision $\varepsilon > 0$ and $\delta_{approx} \in (0,1)$, if the following conditions are satisfied:
	\begin{align*}
		\tau^{-1} & \geq \phi_{pert, radius}(\varepsilon, C_0)
		\\
		N         & \geq \phi_{agent, size, pop}(\varepsilon, \tau, C_0)
		\\
		T         & \geq \phi_{trunc, T, pop}(\varepsilon, \tau, C_0, N)
		\\
		M         & \geq \max \big\{ \phi_{pert, size}(\varepsilon, \tau, C_0, \delta_{approx}),  \phi_{sample, size, pop}( \varepsilon, \tau, T, C_0, \delta_{approx})  \big\},
	\end{align*}
	then, using $M$ perturbation directions $\underline v = (v_i^{(idy)}, v_i^{(com)})_{i=1}^M$ sampled from $\SS_\tau \times \SS_\tau$, the following inequalities hold:
	\begin{align}
		 & \PP \big( \| \tilde{\nabla}^{T, N, M, \tau}(\theta, \underline v) - \hat{\nabla}^{T, N, M, \tau}(\theta, \underline v) \| > \varepsilon \big) \leq \delta_{approx},
		\label{eq:ineq_approx_sample_gradient_pop}
		\\
		 & \| \hat{\nabla}^{T, N, M, \tau}(\theta, \underline v) - \hat{\nabla}^{N, M, \tau}(\theta, \underline v) \|  \leq \varepsilon,
		\label{eq:ineq_approx_nruncated_gradient_pop}
		\\
		 & \| \hat{\nabla}^{N, M, \tau}(\theta, \underline v) - \hat{\nabla}^{M, \tau} (\theta, \underline v) \|  \leq \varepsilon.
		\label{eq:ineq_approx_perturbed_gradient_pop}
	\end{align}
\end{lemma}

\begin{proof}
	We first prove inequality~\eqref{eq:ineq_approx_sample_gradient_pop}.
	For $\theta = (K, L) \in \Theta$ and $K_i = K + v_i^{(idy)}$, using the bounds in~\eqref{eq:bound_phi2_y0_z0_N_agent} and Lemma~\ref{lemma:concentration_ineq_for_yt_K}, we deduce that $\gamma^{n} \| Y_{n}^{K_i, j, N} \|_{\psi_2}^2 \leq \frac{4 d^2 C_{init, noise}^2}{( 1- \sqrt{\gamma})^2}$.
	Defining $\zeta_{n}^{K_i, j, N} = (Y_{n}^{K_i, j, N})^\top ( Q + K_i^\top R K_i) Y_{n}^{K_i, j, N}$ as in Lemma~\ref{lemma:concentration_ineq_for_yt_K}, we obtain:
	\begin{align*}
		     & \PP \Big( \big\| \frac{\gamma^{n}}{M} \sum_{i=1}^M \big( \frac{1}{N} \sum_{j=1}^N  ( \zeta_{n}^{K_i, j, N} - \EE[\zeta_{n}^{K_i, j, N} ] \big) v_i^{(idy)}  \big\| \geq \varepsilon \Big)
		\\
		\leq & (\ell + d ) \exp \Big( (- M \varepsilon^2 ) / (64 h_{\psi_2, K}^2(\tau, C_0) + 8 \varepsilon h_{\psi_2, K}(\tau, C_0) ) \Big).
	\end{align*}
	Consequently, following the reasoning in Lemma~\ref{lemma:approx_sampled_pg_with_{n}runcated_pg}, inequality~\eqref{eq:ineq_approx_sample_gradient_pop} holds provided $M \geq \phi_{sample, size, pop}( \varepsilon, \tau, T, C_0, \delta_{approx})$.

	Next, we establish inequality~\eqref{eq:ineq_approx_nruncated_gradient_pop} for the truncated gradient with social cost.
	We observe a similar bound as in Lemma~\ref{lemma:approx_nruncated_cost} regarding:
	\begin{align*}
		\big\| \Sigma_K^{y, N} - \Sigma_K^{y,T, N} \big\|
		 & \leq  \Big( \frac{\gamma_{pert}^T}{1- \gamma_{pert}} + \gamma_{pert}^T + (T-1) \gamma_{pert}^T \Big)  \frac{\max\big\{  \| \Sigma_{Y_{0}}^N \|, \| \Sigma^{1, N} \| \big\} }{1- \gamma_{pert}}
		\\
		 & \leq ( 1 + \tfrac{1}{N})  2d \gamma_{pert}^T (T +1) C_{init, noise}^2  / (1- \gamma_{pert})^2.
	\end{align*}
	Then, following the same steps as in Lemma~\ref{lemma:approx_nruncated_policy_gradient_MKV_simulator}, we obtain inequality~\eqref{eq:ineq_approx_nruncated_gradient_pop}.

	Finally, for inequality~\eqref{eq:ineq_approx_perturbed_gradient_pop} concerning the perturbed policy gradient with social cost, we note that:
	$$
		\| \hat{\nabla}^{N, M, \tau}(\theta, \underline v) - \hat{\nabla}^{M, \tau} (\theta, \underline v) \|
		\leq \frac{\ell d}{M \tau} \sum_{i=1}^M | C^{N}(\theta_i) - C(\theta_i) | \| v_i^{(idy)} \|.
	$$
	Applying Proposition~\ref{pr:approx_with_social_cost_no_MF_cost} completes the proof.
\qed\end{proof}

\noindent \textbf{Proof of Proposition~\ref{pr:modelfree_population_gradient_approx}}
\begin{proof}
	The result is an immediate consequence of Lemma~\ref{lemma:approx_gradient_pop_simulator} and Lemma~\ref{lemma:approx_perturbed_gradient_MKV_simulator} with the following inequality:
	\begin{align*}
		     & \big\| \tilde \nabla^{T, N, M, \tau}(\theta, \underline v) - \nabla C(\theta) \big\|
		\\
		\leq & \| \tilde{\nabla}^{N, T, M, \tau}(\theta, \underline v) - \hat{\nabla}^{N, T, M, \tau}(\theta, \underline v) \|
		+
		\big\| \hat{\nabla}^{N, T, M, \tau}(\theta, \underline v) - \hat{\nabla}^{N, M, \tau}(\theta, \underline v) \big\|
		\\
		     & + \big\| \hat{\nabla}^{N, M, \tau}(\theta, \underline v) - \hat{\nabla}^{M, \tau}(\theta, \underline v) \big\|
		+ \big\| \hat{\nabla}^{M, \tau}(\theta, \underline v) - \nabla C(\theta) \|.
	\end{align*}
\qed\end{proof}

\backmatter

\cleardoublepage
\phantomsection
\addcontentsline{toc}{chapter}{Bibliography}
\bibliographystyle{plainnat}
\bibliography{MFRL_biblio,MFRL_biblio_complem}

\begin{thebibliography}{178}
\providecommand{\natexlab}[1]{#1}
\providecommand{\url}[1]{\texttt{#1}}
\expandafter\ifx\csname urlstyle\endcsname\relax
  \providecommand{\doi}[1]{doi: #1}\else
  \providecommand{\doi}{doi: \begingroup \urlstyle{rm}\Url}\fi

\bibitem[Achdou and Lasry(2019)]{achdoulasry2019meancrowd}
Yves Achdou and Jean-Michel Lasry.
\newblock Mean field games for modeling crowd motion.
\newblock In \emph{Contributions to Partial Differential Equations and
  Applications}, volume~47 of \emph{Computational Methods in Applied Sciences},
  pages 17--42. Springer, 2019.
\newblock \doi{10.1007/978-3-319-78325-3_4}.

\bibitem[Achdou et~al.(2014)Achdou, Buera, Lasry, Lions, and
  Moll]{achdou2014pde}
Yves Achdou, Francisco Buera, Jean-Michel Lasry, Pierre-Louis Lions, and
  Benjamin Moll.
\newblock {PDE} models in macroeconomics.
\newblock \emph{Philosophical Transactions of the Royal Society A}, 2014.

\bibitem[Achdou et~al.(2017)Achdou, Han, Lasry, Lions, and
  Moll]{achdou2017income}
Yves Achdou, Jiequn Han, Jean-Michel Lasry, Pierre-Louis Lions, and Benjamin
  Moll.
\newblock Income and wealth distribution in macroeconomics: a continuous-time
  approach.
\newblock Technical Report 23732, National Bureau of Economic Research, 2017.

\bibitem[Alasseur et~al.(2020)Alasseur, Taher, and
  Matoussi]{alasseur2020extended}
Cl{\'e}mence Alasseur, Imen~Ben Taher, and Anis Matoussi.
\newblock An extended mean field game for storage in smart grids.
\newblock \emph{Journal of Optimization Theory and Applications}, 184\penalty0
  (2), 2020.

\bibitem[Anahtarci et~al.(2023)Anahtarci, Kariksiz, and
  Saldi]{anahtarci2020qregu}
Berkay Anahtarci, Can~Deha Kariksiz, and Naci Saldi.
\newblock {Q}-learning in regularized mean-field games.
\newblock \emph{Dynamic Games and Applications}, 13:\penalty0 89--117, 2023.

\bibitem[Andersson and Djehiche(2011)]{AnderssonDjehiche2011}
Daniel Andersson and Boualem Djehiche.
\newblock A maximum principle for {SDEs} of mean-field type.
\newblock \emph{Applied Mathematics \& Optimization}, 63:\penalty0 341--356,
  2011.
\newblock \doi{10.1007/s00245-010-9123-8}.

\bibitem[Angiuli et~al.(2022)Angiuli, Fouque, and
  Lauri{\`e}re]{angiuli2022unified}
Andrea Angiuli, Jean-Pierre Fouque, and Mathieu Lauri{\`e}re.
\newblock Unified reinforcement {Q}-learning for mean field game and control
  problems.
\newblock \emph{Mathematics of Control, Signals, and Systems}, 34\penalty0
  (2):\penalty0 217--271, 2022.

\bibitem[Angiuli et~al.(2023{\natexlab{a}})Angiuli, Fouque, and
  Lauri{\`e}re]{angiulia2023reinforcement}
Andrea Angiuli, Jean-Pierre Fouque, and Mathieu Lauri{\`e}re.
\newblock Reinforcement learning for mean field games, with applications to
  economics.
\newblock \emph{Machine Learning and Data Sciences for Financial Markets: A
  Guide to Contemporary Practices}, page 393, 2023{\natexlab{a}}.

\bibitem[Angiuli et~al.(2023{\natexlab{b}})Angiuli, Fouque, Lauriere, and
  Zhang]{angiuli2023convergence}
Andrea Angiuli, Jean-Pierre Fouque, Mathieu Lauriere, and Mengrui Zhang.
\newblock Convergence of multi-scale reinforcement {Q}-learning algorithms for
  mean field game and control problems.
\newblock \emph{arXiv preprint arXiv:2312.06659}, 2023{\natexlab{b}}.

\bibitem[Aurell and Djehiche(2019)]{aurell2019modeling}
Alexander Aurell and Boualem Djehiche.
\newblock Modeling tagged pedestrian motion: a mean-field type game approach.
\newblock \emph{Transportation Research Part B: Methodological}, 121, 2019.

\bibitem[Aurell et~al.(2022)Aurell, Carmona, Dayanikli, and
  Lauri{\`e}re]{aurell2022optimalincentives}
Alexander Aurell, Ren{\'e} Carmona, Gokce Dayanikli, and Mathieu Lauri{\`e}re.
\newblock Optimal incentives to mitigate epidemics: a {S}tackelberg mean field
  game approach.
\newblock \emph{SIAM Journal on Control and Optimization}, 60\penalty0
  (2):\penalty0 S294--S322, 2022.

\bibitem[Bagagiolo and Bauso(2014)]{bagagiolo2014mean}
Fabio Bagagiolo and Dario Bauso.
\newblock Mean-field games and dynamic demand management in power grids.
\newblock \emph{Dynamic Games and Applications}, 4\penalty0 (2), 2014.

\bibitem[Bandini et~al.(2019)Bandini, Cosso, Fuhrman, and
  Pham]{bandini2019filtering}
Elena Bandini, Andrea Cosso, Marco Fuhrman, and Huy{\^e}n Pham.
\newblock Randomized filtering and bellman equation in wasserstein space for
  partial observation control problem.
\newblock \emph{Stochastic Processes and their Applications}, 129\penalty0
  (2):\penalty0 674--711, 2019.
\newblock \doi{10.1016/j.spa.2018.03.014}.

\bibitem[Bardi and Cardaliaguet(2021)]{BardiCardaliaguet}
M.~Bardi and P.~Cardaliaguet.
\newblock Convergence of some mean field games systems to aggregation and
  flocking models.
\newblock \emph{Nonlinear Analysis}, 2021.

\bibitem[Ba{\c{s}}ar et~al.(2026)Ba{\c{s}}ar, Djehiche, and
  Tembine]{BasarDjehicheTembine2026}
Tamer Ba{\c{s}}ar, Boualem Djehiche, and Hamidou Tembine.
\newblock \emph{Mean-Field-Type Game Theory {I}: Foundations and New
  Directions}.
\newblock Static \& Dynamic Game Theory: Foundations \& Applications.
  Birkh{\"a}user, Cham, 2026.
\newblock \doi{10.1007/978-3-032-07027-2}.

\bibitem[Bauso et~al.(2016{\natexlab{a}})Bauso, Tembine, and
  Ba{\c{s}}ar]{bauso2016opinion}
Dario Bauso, Hamidou Tembine, and Tamer Ba{\c{s}}ar.
\newblock Opinion dynamics in social networks through mean-field games.
\newblock \emph{SIAM Journal on Control and Optimization}, 54\penalty0 (6),
  2016{\natexlab{a}}.

\bibitem[Bauso et~al.(2016{\natexlab{b}})Bauso, Zhang, and
  Papachristodoulou]{bauso2016densitynetwork}
Dario Bauso, Xuan Zhang, and Antonis Papachristodoulou.
\newblock Density flow in dynamical networks via mean-field games.
\newblock \emph{IEEE Transactions on Automatic Control}, 62\penalty0
  (3):\penalty0 1342--1355, 2016{\natexlab{b}}.

\bibitem[Bayraktar and Kara(2025)]{bayraktar2025learning}
Erhan Bayraktar and Ali~Devran Kara.
\newblock Learning with linear function approximations in mean-field control.
\newblock \emph{Journal of Machine Learning Research}, 26\penalty0
  (192):\penalty0 1--53, 2025.

\bibitem[Bayraktar et~al.(2026)Bayraktar, Hernandez, Yan, and
  Zhu]{bayraktar2026mean}
Erhan Bayraktar, Martin Hernandez, Qinxin Yan, and Yuhua Zhu.
\newblock Mean-field {PhiBE}: Continuous-time mean-field reinforcement learning
  from discrete-time data.
\newblock \emph{arXiv preprint arXiv:2606.26498}, 2026.
\newblock \doi{10.48550/arXiv.2606.26498}.

\bibitem[Bellemare et~al.(2016)Bellemare, Ostrovski, Guez, Thomas, and
  Munos]{Bellemare}
M.~Bellemare, G.~Ostrovski, A.~Guez, P.~Thomas, and R.~Munos.
\newblock Increasing the action gap: New operators for reinforcement learning.
\newblock In \emph{Thirtieth AAAI Conference on Artificial Intelligence}, 2016.

\bibitem[Bensoussan et~al.(2013)Bensoussan, Frehse, and
  S.C.Yam]{BensoussanFrehseYam}
A.~Bensoussan, J.~Frehse, and S.C.Yam.
\newblock \emph{Mean field games and mean field type control theory}.
\newblock Springer Verlag, New York, NY, 2013.

\bibitem[Bensoussan et~al.(2018)Bensoussan, Huang, and
  Lauri{\`e}re]{bensoussan2018populations}
Alain Bensoussan, Tao Huang, and Mathieu Lauri{\`e}re.
\newblock Mean field control and mean field game models with several
  populations.
\newblock \emph{Minimax Theory and its Applications}, 3\penalty0 (2):\penalty0
  173--209, 2018.
\newblock URL
  \url{https://journalmta.com/index.php/jmta/article/download/44/42}.

\bibitem[Bernasconi et~al.(2023)Bernasconi, Vittori, Trovo, and
  Restelli]{BernasconiVittoriTrovoRestelli}
M.~Bernasconi, E.~Vittori, F.~Trovo, and M.~Restelli.
\newblock Dealer markets: A reinforcement learning mean field game approach.
\newblock \emph{North American Journal of Economics and Finance}, 68:\penalty0
  101974, 2023.

\bibitem[Bertsekas and Shreve(1996)]{BertsekasShreve}
Dimitri Bertsekas and Steven~E Shreve.
\newblock \emph{Stochastic optimal control: the discrete-time case}, volume~5.
\newblock Athena Scientific, 1996.

\bibitem[Blackwell and Dubins(1983)]{BlackwellDubins}
D.~Blackwell and L.~Dubins.
\newblock An extension of {S}korohod's almost sure representation theorem.
\newblock \emph{Proceedings of the American Mathematical Society}, 89:\penalty0
  691--692, 1983.

\bibitem[Bof et~al.(2018)Bof, Carli, and Schenato]{bof2018lyapunov}
Nicoletta Bof, Ruggero Carli, and Luca Schenato.
\newblock Lyapunov theory for discrete time systems.
\newblock \emph{arXiv preprint arXiv:1809.05289}, 2018.

\bibitem[Bongini et~al.(2017)Bongini, Fornasier, Junge, and
  Scharf]{BonginiFornasierJunge2017}
M.~Bongini, M.~Fornasier, O.~Junge, and B.~Scharf.
\newblock Sparse control of multiagent systems.
\newblock In \emph{Active Particles, Volume 1}, pages 173--228. Springer, 2017.

\bibitem[Brunnbauer et~al.(2024)Brunnbauer, Lemmel, Babaiee, Neubauer, and
  Grosu]{BrunnbauerLemmelBabaieeNeubauerGrosu}
A.~Brunnbauer, J.~Lemmel, Z.~Babaiee, S.~Neubauer, and R.~Grosu.
\newblock Scalable offline reinforcement learning for mean field games, 2024.

\bibitem[Buckdahn et~al.(2017)Buckdahn, Li, and Ma]{buckdahn2017partial}
Rainer Buckdahn, Juan Li, and Jin Ma.
\newblock A mean-field stochastic control problem with partial observations.
\newblock \emph{The Annals of Applied Probability}, 27\penalty0 (5):\penalty0
  3201--3245, 2017.
\newblock \doi{10.1214/17-AAP1280}.

\bibitem[Burger et~al.(2013)Burger, Di~Francesco, Markowich, and
  Wolfram]{burger2013mean}
Martin Burger, Marco Di~Francesco, Peter Markowich, and Marie-Therese Wolfram.
\newblock Mean field games with nonlinear mobilities in pedestrian dynamics.
\newblock \emph{Discrete and Continuous Dynamical Systems - Series B}, 19,
  2013.
\newblock \doi{10.3934/dcdsb.2014.19.1311}.

\bibitem[Bu{\c{s}}oniu et~al.(2008)Bu{\c{s}}oniu, Babu{\v{s}}ka, and
  De~Schutter]{BusoniuBabuskadeSchutter2010}
Lucian Bu{\c{s}}oniu, Robert Babu{\v{s}}ka, and Bart De~Schutter.
\newblock A comprehensive survey of multi-agent reinforcement learning.
\newblock \emph{IEEE Transactions on Systems, Man, and Cybernetics, Part C:
  Applications and Reviews}, 38\penalty0 (2):\penalty0 156–172, 2008.

\bibitem[Cabannes et~al.(2022)Cabannes, Lauri{\`e}re, Perolat, Marinier,
  Girgin, Perrin, Pietquin, Bayen, Goubault, and Elie]{cabannes2021solving}
Th{\'e}ophile Cabannes, Mathieu Lauri{\`e}re, Julien Perolat, Raphael Marinier,
  Sertan Girgin, Sarah Perrin, Olivier Pietquin, Alexandre~M. Bayen, Eric
  Goubault, and Romuald Elie.
\newblock Solving {N}-player dynamic routing games with congestion: {A}
  mean-field approach.
\newblock In \emph{Proceedings of the 21st International Conference on
  Autonomous Agents and Multiagent Systems}, pages 1557--1559, 2022.

\bibitem[Cardaliaguet(2013)]{Cardaliaguet2013}
Pierre Cardaliaguet.
\newblock Notes on mean field games, 2013.

\bibitem[Cardaliaguet and Lehalle(2018)]{cardaliaguetlehalle2018mean}
Pierre Cardaliaguet and Charles-Albert Lehalle.
\newblock Mean field game of controls and an application to trade crowding.
\newblock \emph{Mathematics and Financial Economics}, 12\penalty0 (3):\penalty0
  335--363, 2018.

\bibitem[Cardaliaguet et~al.(2023)Cardaliaguet, Daudin, Jackson, and
  Souganidis]{cardaliaguet2023rate}
Pierre Cardaliaguet, Samuel Daudin, Joe Jackson, and Panagiotis~E. Souganidis.
\newblock An algebraic convergence rate for the optimal control of
  {McKean--Vlasov} dynamics.
\newblock \emph{SIAM Journal on Control and Optimization}, 61\penalty0
  (6):\penalty0 3341--3369, 2023.
\newblock \doi{10.1137/22M1486789}.

\bibitem[Carmona et~al.(2023)Carmona, Lauri\`ere, and Tan]{CarmonaLauriere_AAP}
R.~Carmona, M.~Lauri\`ere, and Z.~Tan.
\newblock Model-free mean-field reinforcement learning: Mean-field {MDP} and
  mean-field {Q}-learning.
\newblock \emph{The Annals of Applied Probability}, 33\penalty0 (6B):\penalty0
  5334--5381, 2023.
\newblock \doi{10.1214/23-AAP1949}.

\bibitem[Carmona(2021)]{carmona2020applicationsmfgfinanceeconsurvey}
Ren{\'e} Carmona.
\newblock Applications of mean field games in financial engineering and
  economic theory.
\newblock In \emph{Mean Field Games}, pages 165--218. American Mathematical
  Society, 2021.

\bibitem[Carmona and Delarue(2018)]{CarmonaDelarue_book_I}
Ren\'e Carmona and Fran\c{c}ois Delarue.
\newblock \emph{Probabilistic theory of mean field games with applications.
  {I}}, volume~83 of \emph{Probability Theory and Stochastic Modelling}.
\newblock Springer, Cham, 2018.
\newblock ISBN 978-3-319-56437-1; 978-3-319-58920-6.
\newblock Mean field FBSDEs, control, and games.

\bibitem[Carmona and Lauriere(2025)]{CarmonaLauriere_SIREV}
Rene Carmona and Mathieu Lauriere.
\newblock Reconciling discrete-time mixed policies and continuous-time relaxed
  controls in reinforcement learning and stochastic control.
\newblock \emph{To appear in SIREV (arXiv preprint arXiv:2504.21793)}, 2025.

\bibitem[Carmona et~al.(2019)Carmona, Lauri{\`e}re, and
  Tan]{CarmonaLauriereTan2019}
Ren{\'e} Carmona, Mathieu Lauri{\`e}re, and Zongjun Tan.
\newblock Linear-quadratic mean-field reinforcement learning: Convergence of
  policy gradient methods, 2019.

\bibitem[Carmona et~al.(2020)Carmona, Hamidouche, Lauri{\`e}re, and
  Tan]{carmona2020policyCDC}
Ren{\'e} Carmona, Kenza Hamidouche, Mathieu Lauri{\`e}re, and Zongjun Tan.
\newblock Policy optimization for linear-quadratic zero-sum mean-field type
  games.
\newblock In \emph{2020 59th IEEE Conference on Decision and Control (CDC)},
  pages 1038--1043. IEEE, 2020.

\bibitem[Carmona et~al.(2021)Carmona, Hamidouche, Lauri{\`e}re, and
  Tan]{carmona2021linear}
Ren{\'e} Carmona, Kenza Hamidouche, Mathieu Lauri{\`e}re, and Zongjun Tan.
\newblock Linear-quadratic zero-sum mean-field type games: Optimality
  conditions and policy optimization.
\newblock \emph{Journal of Dynamics \& Games}, 8\penalty0 (4), 2021.

\bibitem[Cecchin(2021)]{cecchin2021finite}
Alekos Cecchin.
\newblock Finite state {$N$}-agent and mean field control problems.
\newblock \emph{ESAIM: Control, Optimisation and Calculus of Variations},
  27:\penalty0 31, 2021.
\newblock \doi{10.1051/cocv/2021032}.

\bibitem[Chan and Sircar(2015)]{ChanSircar}
Patrick Chan and Ronnie Sircar.
\newblock {B}ertrand and {C}ournot mean field games.
\newblock \emph{Applied Mathematics \& Optimization}, 71\penalty0 (3):\penalty0
  533--569, 2015.

\bibitem[Chen et~al.(2020)Chen, Qi, Zhuang, Wang, Liao, and Han]{chen2020mean}
Dezhi Chen, Qi~Qi, Zirui Zhuang, Jingyu Wang, Jianxin Liao, and Zhu Han.
\newblock Mean field deep reinforcement learning for fair and efficient {UAV}
  control.
\newblock \emph{IEEE Internet of Things Journal}, 8\penalty0 (2):\penalty0
  813--828, 2020.

\bibitem[Chen et~al.(2021{\natexlab{a}})Chen, Li, Wang, Yang, Wang, and
  Zhao]{chen2021offlinemeanfield}
Minshuo Chen, Yan Li, Ethan Wang, Zhuoran Yang, Zhaoran Wang, and Tuo Zhao.
\newblock Pessimism meets invariance: Provably efficient offline mean-field
  multi-agent {RL}.
\newblock In \emph{Advances in Neural Information Processing Systems},
  volume~34, 2021{\natexlab{a}}.
\newblock URL
  \url{https://proceedings.neurips.cc/paper/2021/hash/9559fc73b13fa721a816958488a5b449-Abstract.html}.

\bibitem[Chen et~al.(2021{\natexlab{b}})Chen, Zhang, Giannakis, and
  Ba{\c{s}}ar]{chen2021communication}
Tianyi Chen, Kaiqing Zhang, Georgios~B Giannakis, and Tamer Ba{\c{s}}ar.
\newblock Communication-efficient policy gradient methods for distributed
  reinforcement learning.
\newblock \emph{IEEE Transactions on Control of Network Systems}, 9\penalty0
  (2):\penalty0 917--929, 2021{\natexlab{b}}.

\bibitem[Choromanski et~al.(2018)Choromanski, Rowland, Sindhwani, Turner, and
  Weller]{choromanski2018structured}
Krzysztof Choromanski, Mark Rowland, Vikas Sindhwani, Richard Turner, and
  Adrian Weller.
\newblock Structured evolution with compact architectures for scalable policy
  optimization.
\newblock In \emph{International Conference on Machine Learning}, pages
  970--978. PMLR, 2018.

\bibitem[Conn et~al.(2009)Conn, Scheinberg, and Vicente]{MR2487816}
Andrew~R. Conn, Katya Scheinberg, and Luis~N. Vicente.
\newblock \emph{Introduction to derivative-free optimization}, volume~8 of
  \emph{MPS/SIAM Series on Optimization}.
\newblock Society for Industrial and Applied Mathematics (SIAM), Philadelphia,
  PA; Mathematical Programming Society (MPS), Philadelphia, PA, 2009.
\newblock ISBN 978-0-898716-68-9.
\newblock \doi{10.1137/1.9780898718768}.

\bibitem[Cui et~al.(2022)Cui, Tahir, Ekinci, et~al.]{CuiTahirEkinci2022}
Kai Cui, Ahsan Tahir, Gokhan Ekinci, et~al.
\newblock A survey on large-population systems and scalable multi-agent
  reinforcement learning, 2022.

\bibitem[Cui et~al.(2023)Cui, Li, Fabian, and Koeppl]{cui2023scalable}
Kai Cui, Mengguang Li, Christian Fabian, and Heinz Koeppl.
\newblock Scalable task-driven robotic swarm control via collision avoidance
  and learning mean-field control.
\newblock In \emph{2023 IEEE International Conference on Robotics and
  Automation (ICRA)}, pages 1192--1199. IEEE, 2023.

\bibitem[Cui et~al.(2024{\natexlab{a}})Cui, Hauck, Fabian, and
  Koeppl]{cui2023learningpomfc}
Kai Cui, Sascha~H Hauck, Christian Fabian, and Heinz Koeppl.
\newblock Learning decentralized partially observable mean field control for
  artificial collective behavior.
\newblock In \emph{The Twelfth International Conference on Learning
  Representations}, 2024{\natexlab{a}}.

\bibitem[Cui et~al.(2024{\natexlab{b}})Cui, Hauck, Fabian, and
  Koeppl]{cui2024partially}
Kai Cui, Sascha~H Hauck, Christian Fabian, and Heinz Koeppl.
\newblock Partially observable multi-agent reinforcement learning using mean
  field control.
\newblock In \emph{ICML 2024 Workshop: Aligning Reinforcement Learning
  Experimentalists and Theorists}, 2024{\natexlab{b}}.

\bibitem[Daudin et~al.(2024)Daudin, Delarue, and Jackson]{daudin2024rate}
Samuel Daudin, Fran{\c{c}}ois Delarue, and Joe Jackson.
\newblock On the optimal rate for the convergence problem in mean field
  control.
\newblock \emph{Journal of Functional Analysis}, 287\penalty0 (12):\penalty0
  110660, 2024.
\newblock \doi{10.1016/j.jfa.2024.110660}.

\bibitem[Dayan{\i}kl{\i} et~al.(2024)Dayan{\i}kl{\i}, Lauri{\`e}re, and
  Zhang]{dayanikli2024deep}
G{\"o}k{\c{c}}e Dayan{\i}kl{\i}, Mathieu Lauri{\`e}re, and Jiacheng Zhang.
\newblock Deep learning for population-dependent controls in mean field control
  problems with common noise: Extended abstract.
\newblock In \emph{Proceedings of the 23rd International Conference on
  Autonomous Agents and Multiagent Systems}, pages 2231--2233, 2024.
\newblock URL
  \url{https://www.ifaamas.org/Proceedings/aamas2024/pdfs/p2231.pdf}.

\bibitem[Delarue and Vasileiadis(2025)]{delarue2025exploration}
Fran{\c{c}}ois Delarue and Athanasios Vasileiadis.
\newblock Exploration noise for learning linear-quadratic mean field games.
\newblock \emph{Mathematics of Operations Research}, 50\penalty0 (3):\penalty0
  1762--1831, 2025.
\newblock \doi{10.1287/moor.2021.0157}.

\bibitem[Djete et~al.(2022)Djete, Possama{\"i}, and Tan]{djete2022limit}
Mao~Fabrice Djete, Dylan Possama{\"i}, and Xiaolu Tan.
\newblock {McKean--Vlasov} optimal control: Limit theory and equivalence
  between different formulations.
\newblock \emph{Mathematics of Operations Research}, 47\penalty0 (4):\penalty0
  2891--2930, 2022.
\newblock \doi{10.1287/moor.2021.1232}.

\bibitem[Doncel et~al.(2022)Doncel, Gast, and Gaujal]{doncel2022meansir}
Josu Doncel, Nicolas Gast, and Bruno Gaujal.
\newblock A mean field game analysis of {SIR} dynamics with vaccination.
\newblock \emph{Probability in the Engineering and Informational Sciences},
  36\penalty0 (2):\penalty0 482--499, 2022.

\bibitem[Elie et~al.(2019)Elie, Hubert, Mastrolia, and
  Possama{\"i}]{elie2019mean}
Romuald Elie, Emma Hubert, Thibaut Mastrolia, and Dylan Possama{\"i}.
\newblock Mean-field moral hazard for optimal energy demand response
  management.
\newblock \emph{Mathematical Finance}, 2019.

\bibitem[Elie et~al.(2020{\natexlab{a}})Elie, Hubert, and
  Turinici]{elie2020contact}
Romuald Elie, Emma Hubert, and Gabriel Turinici.
\newblock Contact rate epidemic control of {COVID-19}: an equilibrium view.
\newblock \emph{Mathematical Modelling of Natural Phenomena},
  2020{\natexlab{a}}.

\bibitem[Elie et~al.(2020{\natexlab{b}})Elie, Perolat, Lauri{\`e}re, Geist, and
  Pietquin]{elie2020convergence}
Romuald Elie, Julien Perolat, Mathieu Lauri{\`e}re, Matthieu Geist, and Olivier
  Pietquin.
\newblock On the convergence of model free learning in mean field games.
\newblock \emph{Proceedings of the AAAI Conference on Artificial Intelligence},
  34\penalty0 (05):\penalty0 7143--7150, 2020{\natexlab{b}}.

\bibitem[Elliott et~al.(2013)Elliott, Li, and Ni]{ElliottLiNi2013}
Robert~J. Elliott, Xun Li, and Yuan-Hua Ni.
\newblock Discrete time mean-field stochastic linear-quadratic optimal control
  problems.
\newblock \emph{Automatica}, 49\penalty0 (11):\penalty0 3222--3233, 2013.
\newblock \doi{10.1016/j.automatica.2013.08.012}.

\bibitem[Even-Dar and Mansour(2003)]{EvenDarMansour}
E.~Even-Dar and Y.~Mansour.
\newblock Learning rates for {Q}-learning.
\newblock \emph{Journal of Machine Learning Research}, 5:\penalty0 1--25, 2003.

\bibitem[Fabian et~al.(2025)Fabian, Cui, and Koeppl]{fabian2025sparse}
Christian Fabian, Kai Cui, and Heinz Koeppl.
\newblock Learning mean field control on sparse graphs.
\newblock In \emph{Proceedings of the 42nd International Conference on Machine
  Learning}, volume 267 of \emph{Proceedings of Machine Learning Research},
  pages 15637--15660, 2025.
\newblock URL \url{https://proceedings.mlr.press/v267/fabian25a.html}.

\bibitem[Farahmand(2011)]{Farahmand}
A.~Farahmand.
\newblock Action-gap phenomenon in reinforcement learning.
\newblock \emph{Advances in Neural Information Processing Systems}, pages
  172--180, 2011.

\bibitem[Fazel et~al.(2018)Fazel, Ge, Kakade, and Mesbahi]{fazel2018global}
Maryam Fazel, Rong Ge, Sham Kakade, and Mehran Mesbahi.
\newblock Global convergence of policy gradient methods for the linear
  quadratic regulator.
\newblock In \emph{Proceedings of the 35th International Conference on Machine
  Learning}, volume~80, 2018.

\bibitem[Flaxman et~al.(2005)Flaxman, Kalai, and McMahan]{MR2298287}
Abraham~D. Flaxman, Adam~Tauman Kalai, and H.~Brendan McMahan.
\newblock Online convex optimization in the bandit setting: gradient descent
  without a gradient.
\newblock In \emph{Proceedings of the {S}ixteenth {A}nnual {ACM}-{SIAM}
  {S}ymposium on {D}iscrete {A}lgorithms}, pages 385--394. ACM, New York, 2005.

\bibitem[Foerster et~al.(2016)Foerster, Assael, de~Freitas, and
  Whiteson]{FoersterAssaelFreitasWhiteson}
J.~Foerster, Y.~Assael, N.~de~Freitas, and S.~Whiteson.
\newblock Learning to communicate with deep multi-agent reinforcement learning.
\newblock In \emph{Advances in Neural Information Processing Systems}, pages
  2137--2145, 2016.

\bibitem[Frikha et~al.(2024)Frikha, Pham, and Song]{frikha2024full}
Noufel Frikha, Huy{\^e}n Pham, and Xuanye Song.
\newblock Full error analysis of policy gradient learning algorithms for
  exploratory linear quadratic mean-field control problem in continuous time
  with common noise.
\newblock \emph{arXiv preprint arXiv:2408.02489}, 2024.

\bibitem[Frikha et~al.(2025)Frikha, Germain, Lauri{\`e}re, Pham, and
  Song]{frikha2025actor}
Noufel Frikha, Maximilien Germain, Mathieu Lauri{\`e}re, Huy{\^e}n Pham, and
  Xuanye Song.
\newblock Actor-critic learning for mean-field control in continuous time.
\newblock \emph{Journal of Machine Learning Research}, 26\penalty0
  (127):\penalty0 1--42, 2025.

\bibitem[Fu et~al.(2019)Fu, Yang, Chen, and Wang]{fu2019actorcriticMFG}
Zuyue Fu, Zhuoran Yang, Yongxin Chen, and Zhaoran Wang.
\newblock Actor-critic provably finds {N}ash equilibria of linear-quadratic
  mean-field games.
\newblock In \emph{International Conference on Learning Representations}, 2019.

\bibitem[Gao and Pavel(2017)]{GaoPavel}
B.~Gao and L.~Pavel.
\newblock On the properties of the softmax function with application in game
  theory and reinforcement learning, 2017.

\bibitem[Gao and Caines(2020)]{gao2020graphoncontrol}
Shuang Gao and Peter~E. Caines.
\newblock Graphon control of large-scale networks of linear systems.
\newblock \emph{IEEE Transactions on Automatic Control}, 65\penalty0
  (10):\penalty0 4090--4105, 2020.
\newblock \doi{10.1109/TAC.2019.2955976}.

\bibitem[Gast and Gaujal(2011)]{gast2011mean}
N.~Gast and B.~Gaujal.
\newblock A mean field approach for optimization in discrete time.
\newblock \emph{Discrete Event Dynamic Systems}, 21:\penalty0 63--101, 2011.

\bibitem[Gast et~al.(2012)Gast, Gaujal, and Boudec]{gast2012mean}
N.~Gast, B.~Gaujal, and J.Y.~Le Boudec.
\newblock Mean field for {M}arkov decision processes: from discrete to
  continuous optimization.
\newblock \emph{IEEE Transactions on Automatic Control}, 57:\penalty0
  2266--2280, 2012.

\bibitem[Gomes and Sa{\'u}de(2020)]{gomes2020mean}
Diogo~A. Gomes and Jo{\~a}o Sa{\'u}de.
\newblock A mean-field game approach to price formation.
\newblock \emph{Dynamic Games and Applications}, 2020.

\bibitem[Gomes et~al.(2013)Gomes, Mohr, and Souza]{GomesMohrSouza2013}
Diogo~A. Gomes, Jan Mohr, and Rafael~R. Souza.
\newblock Continuous time finite state mean field games.
\newblock \emph{Applied Mathematics \& Optimization}, 68\penalty0 (1):\penalty0
  99--143, 2013.

\bibitem[Graber(2016)]{Graber2016}
P.~Jameson Graber.
\newblock Linear quadratic mean field type control and mean field games with
  common noise, with application to production of an exhaustible resource.
\newblock \emph{Appl. Math. Optim.}, 74\penalty0 (3):\penalty0 459--486, 2016.
\newblock \doi{10.1007/s00245-016-9385-x}.

\bibitem[Graber and Bensoussan(2018)]{graber2018existence}
P.~Jameson Graber and Alain Bensoussan.
\newblock Existence and uniqueness of solutions for {B}ertrand and {C}ournot
  mean field games.
\newblock \emph{Applied Mathematics \& Optimization}, 77\penalty0 (1), 2018.

\bibitem[Grover et~al.(2018)Grover, Bakshi, and
  Theodorou]{grover2018meanhomogeneousflocking}
Piyush Grover, Kaivalya Bakshi, and Evangelos~A. Theodorou.
\newblock A mean-field game model for homogeneous flocking.
\newblock \emph{Chaos: An Interdisciplinary Journal of Nonlinear Science},
  28\penalty0 (6):\penalty0 061103, 2018.

\bibitem[Gu et~al.(2023)Gu, Guo, Wei, and Xu]{gu2019dynamicmfc}
Haotian Gu, Xin Guo, Xiaoli Wei, and Renyuan Xu.
\newblock Dynamic programming principles for mean-field controls with learning.
\newblock \emph{Operations Research}, 71\penalty0 (4):\penalty0 1040--1054,
  2023.

\bibitem[Gu et~al.(2025{\natexlab{a}})Gu, Guo, Wei, and
  Xu]{gu2024meanmarldecentralized}
Haotian Gu, Xin Guo, Xiaoli Wei, and Renyuan Xu.
\newblock Mean-field multiagent reinforcement learning: {A} decentralized
  network approach.
\newblock \emph{Mathematics of Operations Research}, 50\penalty0 (1):\penalty0
  506--536, 2025{\natexlab{a}}.
\newblock \doi{10.1287/moor.2022.0055}.

\bibitem[Gu et~al.(2021)Gu, Guo, Wei, and Xu]{gu2021meanQ}
H.T. Gu, X.~Guo, X.~Wei, and R.~Xu.
\newblock Mean-field controls with {Q}-learning for cooperative {MARL}:
  Convergence and complexity analysis.
\newblock \emph{SIAM Journal on Mathematics of Data Science}, 3:\penalty0
  1168--1196, 2021.
\newblock \doi{10.1137/20M1360700}.

\bibitem[Gu et~al.(2025{\natexlab{b}})Gu, Zhao, Liang, Zheng, Wong, and
  Chae]{gu2025task}
Huixian Gu, Liqiang Zhao, Kai Liang, Gan Zheng, Kai-Kit Wong, and Chan-Byoung
  Chae.
\newblock Task offloading and position optimization for large scale unmanned
  aerial vehicle networks: A mean field learning approach.
\newblock \emph{IEEE Transactions on Cognitive Communications and Networking},
  2025{\natexlab{b}}.

\bibitem[Gu{\'e}ant et~al.(2011)Gu{\'e}ant, Lasry, and Lions]{GueantLasryLions}
O.~Gu{\'e}ant, J.M. Lasry, and P.L. Lions.
\newblock Mean field games and applications.
\newblock In \emph{Paris-Princeton lectures on mathematical finance 2010},
  pages 205--266. Springer Verlag, 2011.

\bibitem[Guo et~al.(2019)Guo, Hu, Xu, and Zhang]{guo2019learning}
Xin Guo, Anran Hu, Renyuan Xu, and Junzi Zhang.
\newblock Learning mean-field games.
\newblock \emph{Advances in Neural Information Processing Systems},
  32:\penalty0 4966--4976, 2019.

\bibitem[Hamidouche et~al.(2016)Hamidouche, Saad, Debbah, and
  Poor]{hamidouche2016mean}
Kenza Hamidouche, Walid Saad, M{\'e}rouane Debbah, and H.~Vincent Poor.
\newblock Mean-field games for distributed caching in ultra-dense small cell
  networks.
\newblock In \emph{2016 American Control Conference}. IEEE, 2016.

\bibitem[He and Liu(2025)]{HeLiu}
J.~He and A.~Liu.
\newblock A hybrid mean field framework for aggregators participating in
  wholesale electricity markets, 2025.

\bibitem[Hu et~al.(2023)Hu, Wei, Yan, and Zhang]{hu2023graphoncontrol}
Yuanquan Hu, Xiaoli Wei, Junji Yan, and Hengxi Zhang.
\newblock Graphon mean-field control for cooperative multi-agent reinforcement
  learning.
\newblock \emph{Journal of the Franklin Institute}, 360\penalty0 (18):\penalty0
  14783--14805, 2023.
\newblock \doi{10.1016/j.jfranklin.2023.09.002}.

\bibitem[Huang et~al.(2015)Huang, Li, and Yong]{HuangLiYong2015}
Jianhui Huang, Xun Li, and Jiongmin Yong.
\newblock A linear-quadratic optimal control problem for mean-field stochastic
  differential equations in infinite horizon.
\newblock \emph{Math. Control Related Fields}, 5\penalty0 (1):\penalty0
  97--127, 2015.
\newblock \doi{10.3934/mcrf.2015.5.97}.

\bibitem[Huang et~al.(2024)Huang, Yardim, and He]{huang2024statistics}
Jiawei Huang, Batuhan Yardim, and Niao He.
\newblock On the statistical efficiency of mean-field reinforcement learning
  with general function approximation.
\newblock In \emph{Proceedings of the 27th International Conference on
  Artificial Intelligence and Statistics}, volume 238 of \emph{Proceedings of
  Machine Learning Research}, pages 289--297, 2024.
\newblock URL \url{https://proceedings.mlr.press/v238/huang24a.html}.

\bibitem[Huang et~al.(2017)Huang, Di, Du, and Chen]{huang2019game}
Kuang Huang, Xuan Di, Qiang Du, and Xi~Chen.
\newblock A game-theoretic framework for autonomous vehicles velocity control:
  bridging microscopic differential games and macroscopic mean field games.
\newblock \emph{Discrete \& Continuous Dynamical Systems - B}, 22\penalty0
  (11), 2017.
\newblock \doi{10.3934/dcdsb.2020131}.

\bibitem[Huang et~al.(2006)Huang, Malham{\'e}, and
  Caines]{HuangCainesMalhame_NCE}
Minyi Huang, Roland~P. Malham{\'e}, and Peter~E. Caines.
\newblock Nash certainty equivalence in large population stochastic dynamic
  games: connections with the physics of interacting particle systems.
\newblock In \emph{Proceedings of the 45th IEEE conference on decision and
  control}, pages 4921--4926. IEEE, 2006.

\bibitem[Huang et~al.(2007)Huang, Caines, and
  Malham{\'e}]{HuangCainesMalhame2007}
Minyi Huang, Peter~E. Caines, and Roland~P. Malham{\'e}.
\newblock Large-population cost-coupled {LQG} problems with nonuniform agents.
\newblock \emph{IEEE Transactions on Automatic Control}, 52\penalty0
  (9):\penalty0 1560--1571, 2007.

\bibitem[Huang et~al.(2012)Huang, Caines, and Malham{\'e}]{huang2012social}
Minyi Huang, Peter~E. Caines, and Roland~P. Malham{\'e}.
\newblock Social optima in mean field {LQG} control: Centralized and
  decentralized strategies.
\newblock \emph{IEEE Transactions on Automatic Control}, 57\penalty0
  (7):\penalty0 1736--1751, 2012.
\newblock \doi{10.1109/TAC.2012.2183439}.

\bibitem[Hubert and Turinici(2018)]{hubert2018nash}
Emma Hubert and Gabriel Turinici.
\newblock Nash-{MFG} equilibrium in a {SIR} model with time dependent newborn
  vaccination.
\newblock \emph{Ricerche di Matematica}, 67\penalty0 (1), 2018.

\bibitem[Jackson and Lacker(2025)]{jackson2025distributed}
Joe Jackson and Daniel Lacker.
\newblock Approximately optimal distributed stochastic controls beyond the mean
  field setting.
\newblock \emph{The Annals of Applied Probability}, 35\penalty0 (1):\penalty0
  251--308, 2025.
\newblock \doi{10.1214/24-AAP2114}.

\bibitem[Jeloka et~al.(2025)Jeloka, Guan, and Tsiotras]{jeloka2025learning}
Bhavini Jeloka, Yue Guan, and Panagiotis Tsiotras.
\newblock Learning large-scale competitive team behaviors with mean-field
  interactions.
\newblock In \emph{The Seventeenth Workshop on Adaptive and Learning Agents},
  2025.

\bibitem[Jin et~al.(2021)Jin, Yang, and Wang]{jin2021pessimism}
Ying Jin, Zhuoran Yang, and Zhaoran Wang.
\newblock Is pessimism provably efficient for offline {RL}?
\newblock In \emph{Proceedings of the 38th International Conference on Machine
  Learning}, volume 139 of \emph{Proceedings of Machine Learning Research},
  pages 5084--5096, 2021.
\newblock URL \url{https://proceedings.mlr.press/v139/jin21e.html}.

\bibitem[Jusup et~al.(2024)Jusup, P{\'a}sztor, Janik, Zhang, Corman, Krause,
  and Bogunovic]{jusup2024safe}
Matej Jusup, Barna P{\'a}sztor, Tadeusz Janik, Kenan Zhang, Francesco Corman,
  Andreas Krause, and Ilija Bogunovic.
\newblock Safe model-based multi-agent mean-field reinforcement learning.
\newblock In \emph{Proceedings of the 23rd International Conference on
  Autonomous Agents and Multiagent Systems}, pages 973--982, 2024.

\bibitem[Kailath(1980)]{Kailath}
T.~Kailath.
\newblock \emph{Linear Systems}.
\newblock Prentice Hall, 1980.

\bibitem[Kallenberg(2002)]{Kallenberg2002foundations}
O.~Kallenberg.
\newblock \emph{Foundations of modern probability}.
\newblock Springer Science \& Business Media, 2002.

\bibitem[Kallenberg(2017)]{Kallenberg_RM}
O.~Kallenberg.
\newblock \emph{Random measures, theory and applications}.
\newblock Springer Verlag, 2017.

\bibitem[Karimi et~al.(2016)Karimi, Nutini, and Schmidt]{karimi2016linear}
Hamed Karimi, Julie Nutini, and Mark Schmidt.
\newblock Linear convergence of gradient and proximal-gradient methods under
  the {P}olyak-{\l}ojasiewicz condition.
\newblock In \emph{Joint European Conference on Machine Learning and Knowledge
  Discovery in Databases}, pages 795--811. Springer, 2016.

\bibitem[Kizilkale et~al.(2019)Kizilkale, Salhab, and
  Malham{\'e}]{kizilkale2019integral}
Arman~C. Kizilkale, Rabih Salhab, and Roland~P. Malham{\'e}.
\newblock An integral control formulation of mean field game based large scale
  coordination of loads in smart grids.
\newblock \emph{Automatica}, 100, 2019.

\bibitem[Kolokoltsov and Bensoussan(2016)]{kolokoltsov2016mean}
Vassili~N. Kolokoltsov and Alain Bensoussan.
\newblock Mean-field-game model for botnet defense in cyber-security.
\newblock \emph{Applied Mathematics \& Optimization}, 74\penalty0 (3):\penalty0
  669--692, 2016.
\newblock \doi{10.1007/s00245-016-9389-6}.

\bibitem[Kolokoltsov and Malafeyev(2018)]{kolokoltsov2018corruption}
Vassili~N. Kolokoltsov and Oleg~A. Malafeyev.
\newblock Corruption and botnet defense: a mean field game approach.
\newblock \emph{International Journal of Game Theory}, 47\penalty0 (3), 2018.

\bibitem[Lachapelle et~al.(2016)Lachapelle, Lasry, Lehalle, and
  Lions]{lachapelle2016efficiency}
Aim{\'e} Lachapelle, Jean-Michel Lasry, Charles-Albert Lehalle, and
  Pierre-Louis Lions.
\newblock Efficiency of the price formation process in presence of high
  frequency participants: a mean field game analysis.
\newblock \emph{Mathematics and Financial Economics}, 10\penalty0 (3), 2016.

\bibitem[Lacker(2017)]{lacker2017control}
Daniel Lacker.
\newblock Limit theory for controlled {McKean--Vlasov} dynamics.
\newblock \emph{SIAM Journal on Control and Optimization}, 55\penalty0
  (3):\penalty0 1641--1672, 2017.
\newblock \doi{10.1137/16M1095895}.

\bibitem[Laguzet and Turinici(2015)]{laguzet2015individualsir}
Laetitia Laguzet and Gabriel Turinici.
\newblock Individual vaccination as {N}ash equilibrium in a {SIR} model with
  application to the 2009--2010 influenza {A} {(H1N1)} epidemic in {F}rance.
\newblock \emph{Bulletin of Mathematical Biology}, 77\penalty0 (10):\penalty0
  1955--1984, 2015.

\bibitem[Lasry and Lions(2007)]{LasryLions2007}
Jean-Michel Lasry and Pierre-Louis Lions.
\newblock Mean field games.
\newblock \emph{Japanese Journal of Mathematics}, 2\penalty0 (1):\penalty0
  229--260, 2007.

\bibitem[Lauri{\`e}re(2025)]{lauriere2025heterogeneity}
Mathieu Lauri{\`e}re.
\newblock An overview of some extensions of mean field games beyond perfect
  homogeneity and anonymity.
\newblock \emph{Science China Information Sciences}, 68:\penalty0 210202, 2025.
\newblock \doi{10.1007/s11432-025-4658-5}.

\bibitem[Lauri{\`e}re et~al.(2022)Lauri{\`e}re, Perrin, Geist, and
  Pietquin]{lauriere2022learning}
Mathieu Lauri{\`e}re, Sarah Perrin, Matthieu Geist, and Olivier Pietquin.
\newblock Learning in mean field games: A survey.
\newblock \emph{arXiv preprint arXiv:2205.12944}, 2022.

\bibitem[Lauriere et~al.(2022)Lauriere, Perrin, Girgin, Muller, Jain, Cabannes,
  Piliouras, P{\'e}rolat, Elie, Pietquin, et~al.]{LaurierePerrinGirgin2022}
Mathieu Lauriere, Sarah Perrin, Sertan Girgin, Paul Muller, Ayush Jain,
  Theophile Cabannes, Georgios Piliouras, Julien P{\'e}rolat, Romuald Elie,
  Olivier Pietquin, et~al.
\newblock Scalable deep reinforcement learning algorithms for mean field games.
\newblock In \emph{International conference on machine learning}, pages
  12078--12095. PMLR, 2022.

\bibitem[Lee et~al.(2021)Lee, Liu, Tembine, Li, and
  Osher]{lee2021controllingepidemics}
Wonjun Lee, Siting Liu, Hamidou Tembine, Wuchen Li, and Stanley Osher.
\newblock Controlling propagation of epidemics via mean-field control.
\newblock \emph{SIAM Journal on Applied Mathematics}, 81\penalty0 (1):\penalty0
  190--207, 2021.

\bibitem[Li et~al.(2016)Li, Malham{\'e}, and Le~Ny]{li2016mean}
Feng Li, Roland~P. Malham{\'e}, and Jerome Le~Ny.
\newblock Mean field game based control of dispersed energy storage devices
  with constrained inputs.
\newblock In \emph{2016 IEEE 55th Conference on Decision and Control}. IEEE,
  2016.

\bibitem[Li et~al.(2024)Li, Li, and Xu]{li2024policy}
Na~Li, Xun Li, and Zuo~Quan Xu.
\newblock Policy iteration reinforcement learning method for continuous-time
  linear--quadratic mean-field control problems.
\newblock \emph{IEEE Transactions on Automatic Control}, 70\penalty0
  (4):\penalty0 2690--2697, 2024.
\newblock \doi{10.1109/TAC.2024.3494656}.

\bibitem[Li et~al.(2021)Li, Wang, Yang, Wang, Wang, Zhao, and
  Zha]{li2021permutationmarl}
Yan Li, Lingxiao Wang, Jiachen Yang, Ethan Wang, Zhaoran Wang, Tuo Zhao, and
  Hongyuan Zha.
\newblock Permutation invariant policy optimization for mean-field multi-agent
  reinforcement learning: A principled approach.
\newblock \emph{arXiv preprint arXiv:2105.08268}, 2021.

\bibitem[Lillicrap et~al.(2016)Lillicrap, Hunt, Pritzel, Heess, Erez, Tassa,
  Silver, and Wierstra]{lillicrap2015continuous}
Timothy~P. Lillicrap, Jonathan~J. Hunt, Alexander Pritzel, Nicolas Heess, Tom
  Erez, Yuval Tassa, David Silver, and Daan Wierstra.
\newblock Continuous control with deep reinforcement learning.
\newblock In \emph{International Conference on Learning Representations}, 2016.

\bibitem[Lojasiewicz(1963)]{lojasiewicz1963topological}
Stanislaw Lojasiewicz.
\newblock A topological property of real analytic subsets.
\newblock \emph{Coll. du CNRS, Les {\'e}quations aux d{\'e}riv{\'e}es
  partielles}, 117\penalty0 (87-89):\penalty0 2, 1963.

\bibitem[Lowe et~al.(2017)Lowe, Wu, Tamar, Harb, Abbeel, and
  Mordatch]{LoweWuTamar2017}
R.~Lowe, Y.~Wu, A.~Tamar, J.~Harb, P.~Abbeel, and I.~Mordatch.
\newblock Multi-agent actor-critic for mixed cooperative-competitive
  environments.
\newblock In \emph{Advances in Neural Information Processing Systems}, pages
  6379--6390, 2017.

\bibitem[Ma et~al.(2024)Ma, Li, Du, Dong, and Yang]{MaLiDuDongYang}
Chengdong Ma, Aming Li, Yali Du, Hao Dong, and Yaodong Yang.
\newblock Efficient and scalable reinforcement learning for large-scale network
  control.
\newblock \emph{Nature Machine Intelligence}, 6\penalty0 (9):\penalty0
  1006--1020, 2024.

\bibitem[Matignon et~al.(2007)Matignon, Laurent, and
  Fort-Piat]{MatignonLaurentLeFortPiat2007}
L.~Matignon, G.~Laurent, and N.~Le Fort-Piat.
\newblock Hysteretic {Q}-learning: An algorithm for decentralized reinforcement
  learning in cooperative multi-agent teams.
\newblock In \emph{IEEE/RSJ International Conference on Intelligent Robots and
  Systems}, pages 64--69. IEEE, 2007.

\bibitem[M{\'e}riaux et~al.(2012)M{\'e}riaux, Varma, and
  Lasaulce]{meriaux2012mean}
Fran{\c c}ois M{\'e}riaux, Vineeth Varma, and Samson Lasaulce.
\newblock Mean field energy games in wireless networks.
\newblock In \emph{2012 Conference Record of the Forty Sixth Asilomar
  Conference on Signals, Systems and Computers}. IEEE, 2012.

\bibitem[Meunier et~al.(2026)Meunier, Pham, and Reisinger]{meunier2026model}
Matthieu Meunier, Huy{\^e}n Pham, and Christoph Reisinger.
\newblock Model-free policy gradient for discrete-time mean-field control.
\newblock \emph{arXiv preprint arXiv:2601.11217}, 2026.

\bibitem[Mondal et~al.(2022{\natexlab{a}})Mondal, Agarwal, Aggarwal, and
  Ukkusuri]{mondal2022heterogeneous}
Washim~Uddin Mondal, Mridul Agarwal, Vaneet Aggarwal, and Satish~V. Ukkusuri.
\newblock On the approximation of cooperative heterogeneous multi-agent
  reinforcement learning ({MARL}) using mean field control ({MFC}).
\newblock \emph{Journal of Machine Learning Research}, 23\penalty0
  (129):\penalty0 1--46, 2022{\natexlab{a}}.
\newblock URL \url{https://jmlr.org/papers/v23/21-1312.html}.

\bibitem[Mondal et~al.(2022{\natexlab{b}})Mondal, Aggarwal, and
  Ukkusuri]{mondal2022nonuniform}
Washim~Uddin Mondal, Vaneet Aggarwal, and Satish~V. Ukkusuri.
\newblock Can mean field control ({MFC}) approximate cooperative multi agent
  reinforcement learning ({MARL}) with non-uniform interaction?
\newblock In \emph{Proceedings of the Thirty-Eighth Conference on Uncertainty
  in Artificial Intelligence}, volume 180 of \emph{Proceedings of Machine
  Learning Research}, pages 1371--1380, 2022{\natexlab{b}}.
\newblock URL \url{https://proceedings.mlr.press/v180/mondal22a.html}.

\bibitem[Mondal et~al.(2024)Mondal, Aggarwal, and
  Ukkusuri]{MondalAggarwalUkkusuri}
W.U. Mondal, V.~Aggarwal, and S.V. Ukkusuri.
\newblock Mean-field approximation of cooperative constrained multi-agent
  reinforcement learning (cmarl).
\newblock \emph{Journal of Machine Learning Research}, 25:\penalty0 1--33,
  2024.

\bibitem[Motte and Pham(2022)]{motte2022mean}
M{\'e}d{\'e}ric Motte and Huy{\^e}n Pham.
\newblock Mean-field {M}arkov decision processes with common noise and
  open-loop controls.
\newblock \emph{The Annals of Applied Probability}, 32\penalty0 (2):\penalty0
  1421--1458, 2022.
\newblock \doi{10.1214/21-AAP1713}.

\bibitem[Nesterov and Spokoiny(2017)]{MR3627456}
Yurii Nesterov and Vladimir Spokoiny.
\newblock Random gradient-free minimization of convex functions.
\newblock \emph{Found. Comput. Math.}, 17\penalty0 (2):\penalty0 527--566,
  2017.
\newblock ISSN 1615-3375.
\newblock \doi{10.1007/s10208-015-9296-2}.

\bibitem[Ni et~al.(2015)Ni, Elliott, and Li]{NiElliottLi2015}
Yuan-Hua Ni, Robert~J. Elliott, and Xun Li.
\newblock Discrete-time mean-field stochastic linear-quadratic optimal control
  problems, {II}: {I}nfinite horizon case.
\newblock \emph{Automatica}, 57:\penalty0 65--77, 2015.
\newblock \doi{10.1016/j.automatica.2015.04.002}.

\bibitem[Nourian et~al.(2011)Nourian, Caines, and
  Malham{\'e}]{NourianCainesMalhame}
M.~Nourian, P.~Caines, and R.~Malham{\'e}.
\newblock Mean field analysis of controlled {C}ucker-{S}male type flocking:
  Linear analysis and perturbation equations.
\newblock In \emph{IFAC Proceedings}, volume~44, pages 4471--4476, 2011.

\bibitem[Nourian et~al.(2010)Nourian, Caines, and
  Malham{\'e}]{nourian2010synthesisflocking}
Mojtaba Nourian, Peter~E. Caines, and Roland~P. Malham{\'e}.
\newblock Synthesis of {C}ucker-{S}male type flocking via mean field stochastic
  control theory: {N}ash equilibria.
\newblock In \emph{2010 48th Annual Allerton Conference on Communication,
  Control, and Computing}, pages 814--819. IEEE, 2010.

\bibitem[Parise and Ozdaglar(2023)]{parise2023graphon}
Francesca Parise and Asuman Ozdaglar.
\newblock Graphon games: a statistical framework for network games and
  interventions.
\newblock \emph{Econometrica}, 91\penalty0 (1):\penalty0 191--225, 2023.

\bibitem[P{\'a}sztor et~al.(2023)P{\'a}sztor, Krause, and
  Bogunovic]{pasztor2023efficient}
Barna P{\'a}sztor, Andreas Krause, and Ilija Bogunovic.
\newblock Efficient model-based multi-agent mean-field reinforcement learning.
\newblock \emph{Transactions on Machine Learning Research}, 2023.

\bibitem[Perrin et~al.(2021)Perrin, Lauri{\`e}re, P{\'e}rolat, Geist, Elie, and
  Pietquin]{PerrinLaurierePerolatGeistEliePietquin}
S.~Perrin, M.~Lauri{\`e}re, J.~P{\'e}rolat, M.~Geist, R.~Elie, and O.~Pietquin.
\newblock Mean field games flock! the reinforcement learning way.
\newblock In \emph{Proceedings of the Thirtieth International Joint Conference
  on Artificial Intelligence (IJCAI-21)}, pages 356 -- 362, 2021.

\bibitem[Perrin et~al.(2020)Perrin, Lauri{\`e}re, P{\'e}rolat,
  et~al.]{PerrinLaurierePerolat2020}
Sarah Perrin, Mathieu Lauri{\`e}re, Julien P{\'e}rolat, et~al.
\newblock Fictitious play for mean field games: Continuous time analysis and
  applications.
\newblock In \emph{Advances in Neural Information Processing Systems}, pages
  13199--13213, 2020.

\bibitem[Pham and Warin(2025)]{PhamWarin}
H.~Pham and X.~Warin.
\newblock Actor-critic learning algorithms for mean-field control with moment
  neural networks.
\newblock In \emph{Methodology and Computations in Applied Probability}, volume
  27 (13), 2025.

\bibitem[Polyak(1963)]{polyak1963gradient}
Boris~Teodorovich Polyak.
\newblock Gradient methods for minimizing functionals.
\newblock \emph{Zhurnal Vychislitel'noi Matematiki i Matematicheskoi Fiziki},
  3\penalty0 (4):\penalty0 643--653, 1963.

\bibitem[Qiu et~al.(2023)Qiu, Wang, Dong, Wang, and
  Strbac]{QiuWangDongWangStrbac}
D.~Qiu, J.~Wang, J.~Dong, X.~Wang, and G.~Strbac.
\newblock Mean-field multi-agent reinforcement learning for peer-to-peer
  multi-energy trading.
\newblock \emph{IEEE Transactions on power systems}, 38\penalty0 (5):\penalty0
  4853 -- 4866, 2023.

\bibitem[Qu et~al.(2022)Qu, Wierman, and Li]{QuWiermanLi}
Guannan Qu, Adam Wierman, and Na~Li.
\newblock Scalable reinforcement learning for multiagent networked systems.
\newblock \emph{Operations Research}, 70\penalty0 (6):\penalty0 3601--3628,
  2022.
\newblock \doi{10.1287/opre.2021.2226}.

\bibitem[Ran and Vreugdenhil(1988)]{ran1988existence}
ACM Ran and R~Vreugdenhil.
\newblock Existence and comparison theorems for algebraic riccati equations for
  continuous-and discrete-time systems.
\newblock \emph{Linear Algebra and its applications}, 99:\penalty0 63--83,
  1988.

\bibitem[Rashid et~al.(2018)Rashid, Samvelyan, Schroeder, Farquhar, Foerster,
  and Whiteson]{RashidSamvelyanSchroeder2018}
T.~Rashid, M.~Samvelyan, C.~Schroeder, G.~Farquhar, J.~Foerster, and
  S.~Whiteson.
\newblock Qmix: Monotonic value function factorisation for decentralised
  multi-agent reinforcement learning.
\newblock In \emph{International Conference on Machine Learning}, pages
  4295--4304. PMLR, 2018.

\bibitem[Recht(2018)]{recht2018tour}
Benjamin Recht.
\newblock A tour of reinforcement learning: The view from continuous control.
\newblock \emph{Annual Review of Control, Robotics, and Autonomous Systems},
  2018.

\bibitem[Salhab et~al.(2018)Salhab, Le~Ny, and
  Malham{\'e}]{salhab2018meanroute}
Rabih Salhab, Jerome Le~Ny, and Roland~P. Malham{\'e}.
\newblock A mean field route choice game model.
\newblock In \emph{2018 IEEE Conference on Decision and Control}, pages
  1005--1010. IEEE, 2018.

\bibitem[Salimans et~al.(2017)Salimans, Ho, Chen, Sidor, and
  Sutskever]{salimans2017evolution}
Tim Salimans, Jonathan Ho, Xi~Chen, Szymon Sidor, and Ilya Sutskever.
\newblock Evolution strategies as a scalable alternative to reinforcement
  learning.
\newblock \emph{arXiv preprint arXiv:1703.03864}, 2017.

\bibitem[Samarakoon et~al.(2015)Samarakoon, Bennis, Saad, Debbah, and
  Latva-Aho]{samarakoon2015energy}
Sumudu Samarakoon, Mehdi Bennis, Walid Saad, M{\'e}rouane Debbah, and Matti
  Latva-Aho.
\newblock Energy-efficient resource management in ultra dense small cell
  networks: a mean-field approach.
\newblock In \emph{IEEE Global Communications Conference}, 2015.

\bibitem[Shao et~al.(2025)Shao, Shen, and Lauriere]{shao2025reinforcement}
Kai Shao, Jiacheng Shen, and Mathieu Lauriere.
\newblock Reinforcement learning for finite space mean-field type game.
\newblock In \emph{Reinforcement Learning Conference}, 2025.

\bibitem[Singh et~al.(2019)Singh, Jain, and Sukhbaatar]{SinghJainSukhbaatar}
A.~Singh, T.~Jain, and S.~Sukhbaatar.
\newblock Learning when to communicate at scale in multiagent cooperative and
  competitive environments.
\newblock In \emph{International Conference on Learning Representations}, 2019.

\bibitem[Stella et~al.(2013)Stella, Bagagiolo, Bauso, and
  Como]{stella2013opinion}
Leonardo Stella, Fabio Bagagiolo, Dario Bauso, and Giacomo Como.
\newblock Opinion dynamics and stubbornness through mean-field games.
\newblock In \emph{52nd IEEE Conference on Decision and Control}. IEEE, 2013.

\bibitem[Subramanian et~al.(2019)Subramanian, Seraj, and
  Mahajan]{SubramanianSerajMahajan}
J.~Subramanian, R.~Seraj, and A.~Mahajan.
\newblock Reinforcement learning for mean-field teams, 2019.

\bibitem[Subramanian and Mahajan(2019)]{subramanian2019reinforcement}
Jayakumar Subramanian and Aditya Mahajan.
\newblock Reinforcement learning in stationary mean-field games.
\newblock In \emph{Proceedings of the 18th International Conference on
  Autonomous Agents and MultiAgent Systems}, pages 251--259, 2019.

\bibitem[Subramanian et~al.(2022)Subramanian, Sinha, Seraj, and
  Mahajan]{SubramanianSinhaSeraj2022}
Jayakumar Subramanian, Amit Sinha, Raihan Seraj, and Aditya Mahajan.
\newblock Approximate information state for approximate planning and
  reinforcement learning in partially observed systems.
\newblock \emph{Journal of Machine Learning Research}, 23:\penalty0 1--83,
  2022.

\bibitem[Subramanian et~al.(2020)Subramanian, Poupart, Taylor, and
  Hegde]{subramanian2020multitype}
Sriram~Ganapathi Subramanian, Pascal Poupart, Matthew~E. Taylor, and Nidhi
  Hegde.
\newblock Multi type mean field reinforcement learning.
\newblock In \emph{Proceedings of the 19th International Conference on
  Autonomous Agents and Multiagent Systems}, pages 411--419, 2020.
\newblock URL \url{https://ifaamas.org/Proceedings/aamas2020/pdfs/p411.pdf}.

\bibitem[Sukhbaatar et~al.(2016)Sukhbaatar, Szlam, and
  Fergus]{SukhbaatarSzlamFergus2016}
S.~Sukhbaatar, A.~Szlam, and R.~Fergus.
\newblock Learning multiagent communication with backpropagation.
\newblock In \emph{Advances in Neural Information Processing Systems}, pages
  2244--2252, 2016.

\bibitem[Sznitman(1991)]{Sznitman1991}
Alain-Sol Sznitman.
\newblock Topics in propagation of chaos.
\newblock In Paul-Louis Hennequin, editor, \emph{{\'E}cole d'{\'E}t{\'e} de
  Probabilit{\'e}s de Saint-Flour XIX---1989}, volume 1464 of \emph{Lecture
  Notes in Mathematics}, pages 165--251. Springer, Berlin, 1991.

\bibitem[Tanaka et~al.(2020)Tanaka, Nekouei, Pedram, and
  Johansson]{tanaka2020linearly}
Takashi Tanaka, Ehsan Nekouei, Ali~Reza Pedram, and Karl~Henrik Johansson.
\newblock Linearly solvable mean-field traffic routing games.
\newblock \emph{IEEE Transactions on Automatic Control}, 66\penalty0
  (2):\penalty0 880--887, 2020.

\bibitem[Tropp(2012)]{MR2946459}
Joel~A. Tropp.
\newblock User-friendly tail bounds for sums of random matrices.
\newblock \emph{Found. Comput. Math.}, 12\penalty0 (4):\penalty0 389--434,
  2012.
\newblock ISSN 1615-3375.
\newblock \doi{10.1007/s10208-011-9099-z}.

\bibitem[Vershynin(2018)]{vershynin2018high}
Roman Vershynin.
\newblock \emph{High-dimensional probability: An introduction with applications
  in data science}, volume~47.
\newblock Cambridge University Press, 2018.

\bibitem[Vuckovic(2026)]{vuckovic2026propagation}
James Vuckovic.
\newblock Propagation of chaos for nonlinear {M}arkov chains.
\newblock \emph{arXiv preprint arXiv:2602.07537}, 2026.

\bibitem[Wang et~al.(2014)Wang, Zhang, and Zhang]{wang2014partial}
Guangchen Wang, Chenghui Zhang, and Weihai Zhang.
\newblock Stochastic maximum principle for mean-field type optimal control
  under partial information.
\newblock \emph{IEEE Transactions on Automatic Control}, 59\penalty0
  (2):\penalty0 522--528, 2014.
\newblock \doi{10.1109/TAC.2013.2273265}.

\bibitem[Wang et~al.(2020{\natexlab{a}})Wang, Zariphopoulou, and
  Zhou]{wang2020reinforcement}
Haoran Wang, Thaleia Zariphopoulou, and Xun~Yu Zhou.
\newblock Reinforcement learning in continuous time and space: A stochastic
  control approach.
\newblock \emph{Journal of Machine Learning Research}, 21\penalty0
  (198):\penalty0 1--34, 2020{\natexlab{a}}.

\bibitem[Wang et~al.(2020{\natexlab{b}})Wang, Yang, and Wang]{wang2020breaking}
Lingxiao Wang, Zhuoran Yang, and Zhaoran Wang.
\newblock Breaking the curse of many agents: Provable mean embedding
  q-iteration for mean-field reinforcement learning.
\newblock In \emph{International Conference on Machine Learning}, pages
  10092--10103. PMLR, 2020{\natexlab{b}}.

\bibitem[Wang et~al.(2021)Wang, Han, Yang, and Wang]{wang2021globallqmfcg}
Weichen Wang, Jiequn Han, Zhuoran Yang, and Zhaoran Wang.
\newblock Global convergence of policy gradient for linear-quadratic mean-field
  control/game in continuous time.
\newblock In \emph{International Conference on Machine Learning}, pages
  10772--10782. PMLR, 2021.

\bibitem[Watkins and Dayan(1992)]{watkins1992q}
Christopher~JCH Watkins and Peter Dayan.
\newblock {Q}-learning.
\newblock \emph{Machine learning}, 8\penalty0 (3):\penalty0 279--292, 1992.

\bibitem[Wei and Yu(2025)]{wei2025continuous}
Xiaoli Wei and Xiang Yu.
\newblock Continuous time q-learning for mean-field control problems.
\newblock \emph{Applied Mathematics \& Optimization}, 91\penalty0 (1):\penalty0
  10, 2025.

\bibitem[Wu et~al.(2024)Wu, Lauri{\`e}re, Chua, Geist, Pietquin, and
  Mehta]{wu2024population}
Zida Wu, Mathieu Lauri{\`e}re, Samuel Jia~Cong Chua, Matthieu Geist, Olivier
  Pietquin, and Ankur Mehta.
\newblock Population-aware online mirror descent for mean-field games by deep
  reinforcement learning.
\newblock In \emph{Proceedings of the 23rd International Conference on
  Autonomous Agents and Multiagent Systems}, pages 2561--2563, 2024.

\bibitem[Xu et~al.(2025{\natexlab{a}})Xu, Chen, Wang, Wu, and
  Shen]{xu2025robust}
Zhenhui Xu, Jiayu Chen, Bing-Chang Wang, Yuhu Wu, and Tielong Shen.
\newblock Robust mean field social control: a unified reinforcement learning
  framework.
\newblock \emph{arXiv preprint arXiv:2502.20029}, 2025{\natexlab{a}}.

\bibitem[Xu et~al.(2025{\natexlab{b}})Xu, Wang, and Shen]{xu2025mean}
Zhenhui Xu, Bing-Chang Wang, and Tielong Shen.
\newblock Mean field {LQG} social optimization: {A} reinforcement learning
  approach.
\newblock \emph{Automatica}, 172:\penalty0 111924, 2025{\natexlab{b}}.

\bibitem[Yang et~al.(2017)Yang, Li, Sheng, Anpalagan, and Xiao]{yang2017mean}
Chungang Yang, Jiandong Li, Min Sheng, Alagan Anpalagan, and Jia Xiao.
\newblock Mean field game-theoretic framework for interference and energy-aware
  control in {5G} ultra-dense networks.
\newblock \emph{IEEE Wireless Communications}, 25\penalty0 (1), 2017.

\bibitem[Yang et~al.(2021)Yang, Ma, Li, Zheng, Zhang, Huang, Yang, and
  Zhao]{yang2021icq}
Yiqin Yang, Xiaoteng Ma, Chenghao Li, Zewu Zheng, Qiyuan Zhang, Gao Huang, Jun
  Yang, and Qianchuan Zhao.
\newblock Believe what you see: Implicit constraint approach for offline
  multi-agent reinforcement learning.
\newblock In \emph{Advances in Neural Information Processing Systems},
  volume~34, 2021.
\newblock URL
  \url{https://papers.nips.cc/paper_files/paper/2021/hash/550a141f12de6341fba65b0ad0433500-Abstract.html}.

\bibitem[Yang et~al.(2019)Yang, Chen, Hong, and Wang]{yang2019provably}
Zhuoran Yang, Yongxin Chen, Mingyi Hong, and Zhaoran Wang.
\newblock Provably global convergence of actor-critic: A case for linear
  quadratic regulator with ergodic cost.
\newblock \emph{Advances in neural information processing systems}, 32, 2019.

\bibitem[Yong(2013)]{Yong2013SICON}
Jiongmin Yong.
\newblock Linear-quadratic optimal control problems for mean-field stochastic
  differential equations.
\newblock \emph{SIAM J. Control Optim.}, 51\penalty0 (4):\penalty0 2809--2838,
  2013.
\newblock \doi{10.1137/120892477}.

\bibitem[Zajkowski(2020)]{zajkowski2020bounds}
Krzysztof Zajkowski.
\newblock Bounds on tail probabilities for quadratic forms in dependent
  sub-{G}aussian random variables.
\newblock \emph{Statistics \& Probability Letters}, 167:\penalty0 108898, 2020.
\newblock \doi{10.1016/j.spl.2020.108898}.

\bibitem[Zaman et~al.(2024{\natexlab{a}})Zaman, Koppel, Lauriere, and
  Basar]{uz2024independent}
Muhammad Aneeq~uz Zaman, Alec Koppel, Mathieu Lauriere, and Tamer Basar.
\newblock Independent {RL} for cooperative-competitive agents: A mean-field
  perspective.
\newblock \emph{arXiv preprint arXiv:2403.11345}, 2024{\natexlab{a}}.

\bibitem[Zaman et~al.(2024{\natexlab{b}})Zaman, Lauriere, Koppel, and
  Ba{\c{s}}ar]{ZamanLauriereKoppel2024}
Muhammad Aneeq~Uz Zaman, Mathieu Lauriere, Alec Koppel, and Tamer Ba{\c{s}}ar.
\newblock Robust cooperative multi-agent reinforcement learning: A mean-field
  type game perspective.
\newblock In \emph{6th Annual Learning for Dynamics \& Control Conference},
  pages 770--783. PMLR, 2024{\natexlab{b}}.

\bibitem[Zhang et~al.(2018)Zhang, Yang, Liu, Zhang, and
  Basar]{ZhangYangLiuZhangBasar}
K.~Zhang, Z.~Yang, H.~Liu, T.~Zhang, and T.~Basar.
\newblock Fully decentralizedmulti-agent reinforcement learningwith networked
  agents.
\newblock In \emph{Proceedings of the 35th Conference on Machine Learning
  PMLR}, volume~80, pages 5872 -- 5881, 2018.

\bibitem[Zhang et~al.(2021)Zhang, Yang, and Ba{\c{s}}ar]{ZhangYangBasar2021}
Kaiqing Zhang, Zhuoran Yang, and Tamer Ba{\c{s}}ar.
\newblock Multi-agent reinforcement learning: A selective overview of theories
  and algorithms.
\newblock In \emph{Handbook of Reinforcement Learning and Control}, pages
  321--384. Springer, 2021.

\end{thebibliography}

\cleardoublepage
\phantomsection
\addcontentsline{toc}{chapter}{Notation Index}
\printindex[not] 
\cleardoublepage
\phantomsection
\addcontentsline{toc}{chapter}{Subject Index}
\printindex[sub] 

\end{document}